\pdfoutput=1
 \documentclass[10pt,reqno]{amsart}

\usepackage{amsfonts}
\usepackage{graphicx}

\usepackage{cancel}

\usepackage{graphicx}
\usepackage[colorlinks=true, pdfstartview=FitV, linkcolor=blue, citecolor=blue,urlcolor=blue]{hyperref}

\usepackage{framed,color}
\definecolor{shadecolor}{rgb}{0.95, 0.95, 0.85}

\newcommand{\ba}{\begin{array}}
\newcommand{\ea}{\end{array}}

\newtheorem{corol}{Corollary}[section]
\newcommand{\bcr}{\begin{corol}}
\newcommand{\ecr}{\end{corol}}

\newtheorem{hypo}{Assumption}[section]
\newcommand{\bhypo}{\begin{hypo}}
\newcommand{\ehypo}{\end{hypo}}

\newtheorem{lemma}{Lemma}[section]
\newcommand{\ble}{\begin{lemma}}
\newcommand{\ele}{\end{lemma}}

\newtheorem{definition}{Definition}[section]
\newcommand{\bde}{\begin{definition}}
\newcommand{\ede}{\end{definition}}

\newtheorem{example}{Example}[section]
\newcommand{\bex}{\begin{example}}
\newcommand{\eex}{\end{example}}

\newtheorem{prop}{Proposition}[section]
\newcommand{\epr}{\end{prop}}
\newcommand{\bpr}{\begin{prop}}
\newtheorem{teo}{Theorem}[section]
\newcommand{\bth}{\begin{teo}}
\newcommand{\eth}{\end{teo}}
\newtheorem{rema}{Remark}[section]
\newcommand{\bre}{\begin{rema} \rm}
\newcommand{\ere}{ \end{rema}}

\newcommand{\bea}{\begin{eqnarray}}
\newcommand{\eea}{\end{eqnarray}}
\newcommand{\beas}{\begin{eqnarray*}}
\newcommand{\eeas}{\end{eqnarray*}}
\newcommand{\ee}{\end{equation}}
\newcommand{\be}{\begin{equation}}

\numberwithin{equation}{section} 

\begin{document}

\centerline{\Large\bf  Isomonodromy Deformations 
 } 
\vskip 0.3 cm 
\centerline{\Large\bf   at an Irregular Singularity with Coalescing Eigenvalues}

\vskip 0.5 cm
\centerline{\large Giordano Cotti, Boris Dubrovin, Davide Guzzetti}
\vskip 0.2 cm

\vskip 0.5 cm
\noindent
Giordano Cotti's ORCID ID: 0000-0002-9171-2569

\vskip 0.2 cm
\noindent
Boris Dubrovin's   ORCID ID: 0000-0001-9856-5883

\vskip 0.2 cm
\noindent
Davide Guzzetti's   ORCID ID: 0000-0002-6103-6563

\tableofcontents



\vskip 0.2 cm 
\noindent
\centerline{\bf Abstract} 
\vskip 0.2 cm 
\noindent 
We consider an $n\times n$ linear system of ODEs with an irregular singularity of
 Poincar\'e rank 1 at $z=\infty$, holomorphically depending on parameter $t$  within a polydisc in $\mathbb{C}^n$ centred at  $t=0$.  The eigenvalues of the leading matrix at $z=\infty$ coalesce  along a locus $\Delta$ contained in the polydisc,  passing through  $t=0$.   Namely,   $z=\infty$ is  a {\it resonant irregular singularity  for  $t\in \Delta$}.   
 We analyse the case when the leading matrix remains diagonalisable at  $\Delta$. We discuss the existence of fundamental matrix 
  solutions, their asymptotics, Stokes phenomenon and monodromy data as $t$ varies in the polydisc, and their limits for $t$ tending to 
  points of $\Delta$.  When the deformation is isomonodromic away from $\Delta$,  it is well known that a fundamental   matrix solution  has  singularities at $\Delta$.  When the system also has a Fuchsian singularity at $z=0$,  we show under  minimal vanishing conditions on the residue matrix at
   $z=0$ that isomonodromic deformations can be extended to the whole polydisc, including   $\Delta$,  in such
    a way that the fundamental matrix solutions and the constant monodromy data are well defined in the whole polydisc.  These data 
    can be computed just by considering  the system at  fixed $t=0$.    Conversely, if the $t$-dependent  system  is  isomonodromic   in a small domain 
     contained in the polydisc not intersecting $\Delta$,  if the  entries of the Stokes matrices with indices 
     corresponding to coalescing  eigenvalues vanish, then we  show that $\Delta$ is not a branching locus for the 
     fundamental matrix solutions.  The importance of these results  for the analytic theory of  
 Frobenius Manifolds is explained. An application to Painlev\'e equations is discussed.

\vskip 0.3 cm 
\hrule


\section{Introduction}
\label{17luglio2016-5}

 We study  deformations of   linear differential systems,  playing an important role in applications,  with a resonant irregular singularity at $z=\infty$.  The $n\times n$ linear (deformed) system  depends on  parameters $t=(t_1,...,t_m)\in\mathbb{C}^m$, (here $n, m\in\mathbb{N}\backslash\{0\}$) and has the following form:
 \be
\label{22novembre2016-3}
\frac{dY}{dz}=A(z,t)Y, \quad\quad\quad A(z,t):=A_0(t)+\sum_{k=1}^\infty A_k(t)z^{-k},
\ee
with  singularity of Poincar\'e rank 1 at $z=\infty$. The series $A(z,t)$ is uniformly convergent in a neighborhood of $z=\infty$  for $|z|\geq N_0>0$ sufficiently large, and the coefficients $A_0(t)$ and  $A_k(t)$, $k\geq 1$,  are  holomorphic matrix valued functions    on an  open connected  domain of $\mathbb{C}^m$. 
We  take  the Poincar\'e rank equal to 1 in view of  the  important  applications which motivate our work, as it is explained below in this Introduction.   

\vskip 0.2 cm 

The deformation theory  is well understood when $A_0(t)$ has distinct eigenvalues  $u_1(t)$, $u_2(t)$, ..., $u_n(t)$ for $t$ in the domain. On the other hand,   there are  important cases for applications (see below) when $A_0(t)$ is {\it diagonalisable}, but   two or more eigenvalues may coalesce when $t$ reaches a certain locus $\Delta$ in the $t$-domain,   called  the {\it coalescence locus}.  This means that  $u_a(t) =u_b(t)$  for some indices $a\neq b\in\{1,...,n\}$  whenever  $t$ belongs to $\Delta$, while 
  $u_1(t)$, $u_2(t)$, ..., $u_n(t)$ are pairwise distinct otherwise\footnote{$\Delta$ is a discrete set for $m=1$, otherwise it is a continuous locus for $m\geq 2$. For example, 
for the matrix $\hbox{diag}(t_1,t_2,...,t_n)$, 
the coalescence locus is 
 the union of the {\it diagonals}  $t_i=t_j$,  $i\neq j\in\{1,2,...,n\}$.}.  Points of $\Delta$ will be called {\it coalescence points}. The point $z=\infty$ for $t\in \Delta$ is usually called a {\it resonant irregular singularity}, but we will not use this nomenclature throughout the paper.  
 To the best of our knowledge, the analysis of fundamental matrix solutions and their monodromy  when  $A_0(t)$ is {\it diagonalisable}  with  coalescing eigenvalues for $t\in \Delta$,  is missing from the existing literature, as we will shortly review later. This  is the main problem  which we address in the present paper, both in the non-isomonodromic and isomonodromic cases. The main results of the paper are contained in: 
 
 \begin{itemize}
 \item[--] Theorem \ref{pincopallino}, Corollaries \ref{6marzo2016-2} and \ref{19feb2016-1}, and in Theorem \ref{8gen2017-4}, for the non-isomonodromic case;
 
 \item[--] Theorem \ref{16dicembre2016-1} (Th. \ref{3nov2015-5}), Corollary \ref{17dicembre2016-2} (Corol. \ref{17dicembre2016-5}) and Theorem \ref{9gen2017-1}, for the isomonodromic case. 
 \end{itemize}

\vskip 0.2 cm 
 
    For the sake of the  local analysis at coalescence points, we  can restrict  to the case when the domain is a polydisk
  $$
   \mathcal{U}_{\epsilon_0}(0):=\bigl\{t\in \mathbb{C}^m
   \quad\hbox{ such that }\quad 
   |t|\leq \epsilon_0 \bigr\},\quad\quad|t|:=\max_{1\leq i\leq m}|t_i|
  $$
   for suitable $\epsilon_0>0$,  being  $t=0$  a point of the coalescence locus.  We will again denote by $\Delta$ the coalescence locus in  $\mathcal{U}_{\epsilon_0}(0)$. 
    It is well known that  the eigenvalues  $u_1(t),...,u_n(t)$  are branches of one or more functions of  $t$, with algebraic branching at $\Delta$ (see \cite{Kato}). A matrix $G_0(t)$ which  diagonalises   $A_0(t)$ for $t\not \in \Delta$, namely such that  $G_0^{-1}(t)A_0(t)G_0(t)=\Lambda(t)$, where 
       \be
\label{27dic2015-1}
\Lambda(t):=
\hbox{diag}(u_1(t),...,u_n(t)),
\ee 
 is  generally  {\it   singular} when $t$ approaches  $\Delta$.  Example will be given in Section \ref{17dicembre2016-1}.  
Consider a fundamental  solution\footnote{A fundamental matrix solution will be simply called a {\it fundamental solution}.}  that for $t$ belonging to a sufficiently small domain $\mathcal{V}
\subset \mathcal{U}_{\epsilon_0}(0)$, with $\mathcal{V}\cap \Delta=\emptyset$,  has a canonical asymptotic representation (see \cite{Sh3})
 $$ 
 Y(z,t)\sim G_0(t)\Bigl(I+\sum_{k=1}^\infty F_k(t) z^{-k}\Bigr) z^{B_1}e^{\Lambda(t)z},~~~\hbox{  $z\to \infty$},
 $$
  in a suitable sector $\mathcal{S}(\overline{\mathcal{V}})$ depending on $\mathcal{V}$, explained after  formula (\ref{4gen2016-6}) below. 
  Here $I$ stands for the identity matrix  and $B_1$ is a diagonal matrix, given in   formula (\ref{4gen2016-6}).  
  Then, the $t$-analytic continuation of $Y(z,t)$  inherits  the singularities of $G_0(t)$ as $t$ tends to  $\Delta$. Thus,  in order to extend the deformation theory when $t$ approaches $\Delta$, we need the following:   

\vskip 0.2 cm 
{\bf Assumption 1:} The holomorphic matrix $A_0(t)$   is {\it holomorphically similar in  $\mathcal{U}_{\epsilon_0}(0)$} to a diagonal matrix $\Lambda(t)$ as in (\ref{27dic2015-1}), namely  there exists an invertible  matrix  $G_0(t)$ {\it holomorphic on $\mathcal{U}_{\epsilon_0}(0)$}  such that 
\be
\label{3gen2017-1}
G_0^{-1}(t)A_0(t)G_0(t)=\Lambda(t).
\ee

\vskip 0.2 cm 

Assumption 1, which is basically {\it the assumption} of the paper,   holds for example  for   Frobenius manifolds  remaining semisimple at the locus of coalescent canonical coordinates, and in applications to the sixth Painlev\'e  transcendents  holomorphic at a fixed singularity of the Painlev\'e equation (see Sections \ref{16dicembre2016-5}  and \ref{31luglio2016-2} below).  
\vskip 0.2 cm 

 Given Assumption 1, the  transformation $Y\mapsto G_0(t) Y$ changes $A(z,t)$ to a matrix valued function 
\be
\label{4gen2017-1}
 \widehat{A}(z,t):=G_0(t)^{-1}A(z,t) G_0(t),
 \ee
holomorphic on $\{|z| \geq N_0\}\times \mathcal{U}_{\epsilon_0}(0)$ for sufficiently large $N_0>0$,  so that system  (\ref{22novembre2016-3})  becomes  
\be
\label{17luglio2016-1}
\frac{dY}{dz}=\widehat{A}(z,t)Y,
\quad 
\quad
\widehat{A}(z,t)=\Lambda(t)+\sum_{k=1}^\infty \widehat{A}_k(t)z^{-k}.
\ee
where $\widehat{A}_k(t)$, $k\geq 1$,  and $\Lambda(t)$  are holomorphic on $\mathcal{U}_{\epsilon_0}(0)$.  

\vskip 0.2 cm 
When $\Delta$ is not empty, the  dependence on $t$ of fundamental solutions  of  (\ref{17luglio2016-1})   near $z= \infty$   is quite delicate. If $t\not\in \Delta$, then the system     
(\ref{17luglio2016-1}) has a  unique  formal solution  (see \cite{Sh3}),
\be
\label{4gen2016-6}
Y_F(z,t):=\Bigl(I+\sum_{k=1}^\infty F_k(t) z^{-k}\Bigr) z^{B_1(t)}e^{\Lambda(t)z},
\quad\quad
B_1(t):=\hbox{diag}(\widehat{A}_1(t)),
\ee
where  the matrices $F_k(t)$ are uniquely determined by the equation and are holomorphic on $\mathcal{U}_{\epsilon_0}(0)\backslash \Delta$.

In order to find actual solutions, and their domain of definition in the space of parameters $t$,  one can refer to the local existence results of Sibuya \cite{Sh4} \cite{Sh3} (see Theorems \ref{ShybuyaOLD} and \ref{ShybuyaOLDbis} below), which guarantees that, given $t_0\in \mathcal{U}_{\epsilon_0}(0)\backslash \Delta$, there exists a sector and a fundamental solution $Y(z,t)$ holomorphic for $|z|$ large and  $|t-t_0|<\rho$ , where $\rho$ is sufficiently  small, such that $Y(z,t)\sim Y_F(z,t)$ for $z\to \infty$ in the  sector. The condition $|t-t_0|$ is restrictive, since $\rho$ is expected to be very small.  In the present paper, we prove this result for $t$ in a wider domain   $\mathcal{V}\subset \mathcal{U}_{\epsilon_0}(0)$, extending  $|t-t_0|<\rho$. $\mathcal{V}$ is constructed as follows. Let $t=0$ and consider the {\it Stokes rays}   associated with the matrix $\Lambda(0)$, namely  rays in  the universal covering  $\mathcal R$ of the $z$-punctured plane  $\mathbb{C}\backslash \{0\}$,    defined by the condition that  $\Re e[(u_a(0)-u_b(0))z]=0$, with $u_a(0)\neq u_b(0)$ ($1\leq a\neq b \leq n$).  Then, consider  an {\it admissible ray}, namely a ray  in $\mathcal R$,  with a certain direction  $\widetilde{\tau}$, that  does not contain  any of the Stokes rays above, namely $\Re e[(u_a(0)-u_b(0))z]\neq 0$ for  any $u_a(0)\neq u_b(0)$ and $\arg z=\widetilde{\tau}$. Define the locus $X(\widetilde{\tau})$ to be the set of points $t\in \mathcal{U}_{\epsilon_0}(0)$ such that some Stokes rays $\{z\in\mathcal{R}~|~ \Re e[(u_a(t)-u_b(t))z]=0\}$
   associated with $\Lambda(t)$, $t\not \in \Delta$,  coincide with the admissible ray $\arg z=\widetilde{\tau}$.  
   Finally, define a {\it $\widetilde{\tau}$-cell } to be any connected component of 
  $\mathcal{U}_{\epsilon_0}(0)\backslash \left( \Delta\cup X(\widetilde{\tau})\right)$ (see Section \ref{onCELLS} for a thorough study of the cells).  Then,  we take an open  connected open domain  
  $\mathcal{V}$ such that its closure $\overline{\mathcal{V}}$ is contained in a $\widetilde{\tau}$-cell. 
  
    \bde
  \label{12maggio2017-2}
  The deformation of the linear system  (\ref{ourcase12}), such that  $t$ varies in  an open connected  domain $ \mathcal{V}\subset \mathcal{U}_{\epsilon_0}(0)$, such $\overline{\mathcal{V}}$ is contained in a $\widetilde{\tau}$-cell, is  called  an {\bf admissible deformation}\footnote{The definition of  {\it admissible deformation} of a linear system is in accordance with the definition given in \cite{Its}.}. For simplicity, we will  just say that  $t$ is an admissible deformation. 
  \ede
  
 By definition, an admissible deformation means that as long as $t$ varies within $\overline{\mathcal{V}}$, then {\it no Stokes rays of $\Lambda(t)$ cross the admissible ray} of  direction $\widetilde{\tau}$.

 \vskip 0.2 cm 
 If $t$ belongs to a domain  $\mathcal{V}$ as above, then we prove in Section \ref{14feb2016-10}    that there is  a family of  actual fundamental solutions  $Y_r(z,t)$, labelled by $r\in\mathbb{Z}$, uniquely determined by the canonical asymptotic representation 
$$Y_r(z,t)\sim Y_F(z,t),
$$
 for  $z\to\infty$ in  suitable sectors $\mathcal{S}_r(\overline{\mathcal{V}})$ of the universal covering $ \mathcal{R}$ of $\mathbb{C}\backslash\{0\}$.  Each $Y_r(z,t)$  is  holomorphic in $\{z\in \mathcal{R}~|~|z|\geq N\}\times \mathcal{V}$, for a suitably large $N$. The asymptotic series  $I+\sum_{k=1}^\infty F_k(t) z^{-k}$ is uniform in $\overline{\mathcal{V}}$. 
 
\begin{figure}
\centerline{\includegraphics[width=0.8 \textwidth]{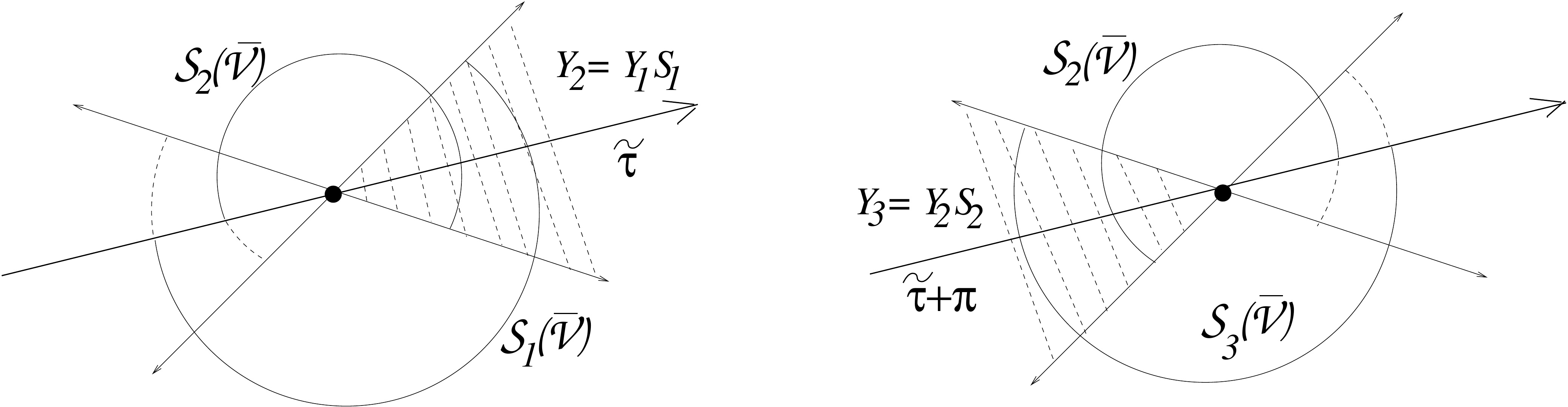}}
\caption{Stokes phenomenon of formula (\ref{2maggio2016-4}). In the left figure is represented 
the sheet of the universal covering $\widetilde{\tau}-\pi<\arg z<\widetilde{\tau}+\pi$ containing  
$\mathcal{S}_1(\overline{\mathcal{V}})\cap \mathcal{S}_2(\overline{\mathcal{V}})$,  and in the right figure
  the sheet $\widetilde{\tau}<\arg z<\widetilde{\tau}+2\pi$ containing  
  $\mathcal{S}_2(\overline{\mathcal{V}})\cap \mathcal{S}_3(\overline{\mathcal{V}})$. 
  The rays $\arg z= \widetilde{\tau}$  and $\widetilde{\tau} +\pi$ (and then $\widetilde{\tau}+k\pi$ for 
  any $k\in\mathbb{Z}$) are  {\it admissible rays},   such that   
  $\Re e\Bigl[(u_a(0)-u_b(0))z\Bigr]\neq 0$ along these rays,  for any $u_a(0)\neq u_b(0)$. 
  Moreover,  $\Re e\Bigl[(u_a(t)-u_b(t))z\Bigr]\neq 0$ for 
  any $t\in \mathcal{V}$ and any $1\leq a\neq b\leq n$.} 
\label{17aprile-1}
\end{figure}

 The sectors  $\mathcal{S}_r(\overline{\mathcal{V}})$ are constructed as follows: take for example
  the `` half plane'' $\Pi_1:=\{z\in\mathcal{R}~|~\widetilde{\tau}-\pi<\arg z<\widetilde{\tau}\}$. 
  The open sector containing $\Pi_1$ and  extending up to  the closest Stokes rays 
  of $\Lambda(t)$ outside $\Pi_1$ will be called $\mathcal{S}_1(t)$. 
  Then, we define $\mathcal{S}_1(\overline{\mathcal{V}}):=\bigcap_{t\in \overline{\mathcal{V}}}\mathcal{S}_1(t)$.  
  Analogously, we consider the ``half-planes''  $\Pi_r:=\{z\in\mathcal{R}~|~\widetilde{\tau}+(r-3)\pi<\arg z<\widetilde{\tau}+(r-1)\pi\}$ and repeat the same construction for $\mathcal{S}_r(\overline{\mathcal{V}})$. 
The sectors $\mathcal{S}_r(\overline{\mathcal{V}})$  have central opening angle
  greater than $\pi$ and their  successive intersections do not contain Stokes 
  rays $\Re e[(u_a(t)-u_b(t))z]=0$ associated with the eigenvalues of $\Lambda(t)$, $t\in \overline{\mathcal{V}}$.  The sectors  $\mathcal{S}_r(\overline{\mathcal{V}})$  for $r=1,2,3$ are represented in 
  Figure \ref{17aprile-1}. 
   An admissible ray $\arg z= \widetilde{\tau}$  in $\mathcal{S}_1(\overline{\mathcal{V}})\cap \mathcal{S}_2(\overline{\mathcal{V}})$  is also represented.

\vskip 0.2 cm 

 If the $t$-analytic continuation  of $Y_r(z,t)$  exists outside $\mathcal{V}$, then the delicate points emerge, as follows. 
 
 \begin{itemize}
 \item
  The expression $\Re e\left[(u_a(t)-u_b(t))z\right]$, $1\leq a\neq b\leq n$, has constant sign  in the  $\widetilde{\tau}$-cell containing $\mathcal{V}$, but it vanishes when a Stokes ray   $\Re e\left[(u_a(t)-u_b(t))z\right]=0$ crosses the admissible direction $\widetilde{\tau}$. This corresponds to the fact that $t$ crosses the boundary of the cell. Then, it  changes sign for $t$ outside of the cell. Hence,  the asymptotic representation $Y_r(z,t)\sim Y_F(z,t)$ for $z\to\infty$  in  $\mathcal{S}_r(\overline{\mathcal{V}})$  does no longer hold for $t$ outside the $\widetilde{\tau}$-cell containing $\mathcal{V}$. 
   
   \item
   The coefficients $F_k(t)$  are in general 
 divergent  at $\Delta$.
 
 \item
  The locus  $\Delta$ is expected to be  a locus of singularities for the  $Y_r(z,t)$'s (see Example  \ref{5dic2016-1} below). 
 
\item 
 The  {\it Stokes matrices} $\mathbb{S}_r(t)$, defined for     $t\in \mathcal{V}$ by the relations  (see Figure \ref{17aprile-1}) 
  \be
  \label{2maggio2016-4}
  Y_{r+1}(z,t)=Y_r(z,t)~\mathbb{S}_r(t),
  \ee
    are expected to be singular as $t$ approaches $\Delta$. 
  \end{itemize}
  
  \bre
  It is well known that, in order to completely describe the Stokes phenomenon,  it suffices to consider only three fundamental solutions, for example $Y_r(z,t)$ for $r=1,2,3$, and   $\mathbb{S}_1(t)$, $\mathbb{S}_2(t)$.
  \ere
 
 The matrix $A(z,t)$ may have other singularities at finite values of $z$.  In the  isomonodromic case, we will 
 consider $A(z,t)$ with a simple pole at $z=0$, namely
   \be
  \label{1novembre2016-1} 
  A(z,t)=A_0(t)+\frac{A_1(t)}{z}.
  \ee  
 An isomonodromic  system of type (\ref{1novembre2016-1}), with  antisymmetric $A_1$, is at the core of the analytic approach  to   semisimple Frobenius manifolds  \cite{Dub1} \cite{Dub3} \cite{Dub2}  (see also \cite{Sai1} \cite{Sai2} \cite{Sai3}  \cite{Man99} \cite{Sab08}). Its monodromy data  play the role of local moduli. Coalescence of eigenvalues of $A_0(t)$ occurs in important cases,  such as quantum cohomology (see \cite{Cotti1} \cite{OURPAPER} and Section \ref{16dicembre2016-5} below). 
  For $n=3$, a special case of system (\ref{1novembre2016-1})  gives an isomonodromic  description of  the general  sixth  Painlev\'e equation, according to  \cite{MazzoIrr} (see also \cite{harnad}). This description was given also in \cite{Dub1} \cite{Dub2} for a sixth  Painlev\'e equation associated with Frobenius manifolds. 
Coalescence occurs at the critical points of the   Painlev\'e equation (see Section \ref{31luglio2016-2}). 

\vskip 0.2 cm 
The transformation $Y\mapsto G_0(t) Y$ changes   (\ref{1novembre2016-1})  into 
\be
\label{ourcase12}
\frac{dY}{dz}=\widehat{A}(z,t)Y,
\quad\quad \widehat{A}(z,t)=\Lambda(t)+\frac{\widehat{A}_1(t)}{z}.
\ee 
For given $t$, a matrix $G^{(0)}(t)$ (not to be confused with $G_0(t)$ in (\ref{3gen2017-1})),  puts  $\widehat{A}_1(t)$ in  Jordan form 
$$
J^{(0)}(t):=(G^{(0)}(t))^{-1}\widehat{A}_1(t)~G^{(0)}(t).
$$
  Close to the Fuchsian singularity $z=0$, and for a given $t$, the system  (\ref{ourcase12}) has a    fundamental   solution  
   \be
\label{24gen2016-1}
Y^{(0)}(z,t)=G^{(0)}(t) \left(I+\sum_{l=1}^\infty\Psi_l(t) z^l\right) z^{D^{(0)}(t)}z^{S^{(0)}(t)+R^{(0)}(t)}, 
\ee   
in standard Birkhoff--Levelt  normal  form, whose  behaviour in $z$ and $t$ is not affected by the coalescence phenomenon.
 The matrix coefficients  $\Psi_l(t)$ of the convergent expansion are constructed by a recursive procedure.   $D^{(0)}(t)=\hbox{diag}(d_1(t),...,d_n(t))$ is a diagonal matrix of integers, piecewise constant in $t$,  $S^{(0)}(t)$ is a Jordan matrix whose eigenvalues  $\rho_1(t),...,\rho_n(t)$ have real part in $[0,1[$, and the nilpotent matrix $R^{(0)}(t)$ has non-vanishing entries only if  some eigenvalues of $\widehat{A}_1(t)$  differ by  non-zero integers. If  some eigenvalues   differ by  non-zero integers, we say that $\widehat{A}_1(t)$ is {\it resonant}. The sum 
$$J^{(0)}(t)=D^{(0)}(t)+S^{(0)}(t)$$ 
is the Jordan form of  $\widehat{A}_1(t)$ above.   
 Under the assumptions of our Theorem \ref{16dicembre2016-1} below, the solution (\ref{24gen2016-1}) turns out to be  holomorphic in $t\in \mathcal{U}_{\epsilon_0}(0)$. \footnote{If $Y^{(0)}(z,t)$ is chosen, with given $G^{(0)}(t)$,  $\Psi_l(t) $'s and $R^{(0)}(t)$,  then there is a class of suitable  matrices $\mathfrak{D}(t)$  such that   $Y^{(0)}(z,t)\mathfrak{D}(t)$ also has the standard form  (\ref{24gen2016-1}) with new $G^{(0)}(t)$,  $\Psi_l(t)$'s and $R^{(0)}(t)$.  More details are in  Section \ref{5gen2016-2}.}

  In order to completely describe the monodromy of the system (\ref{ourcase12}), we need its {\it essential} monodromy data (the adjective ``essential''  is inspired by a similar definition in  \cite{JMU}). We recall that it suffices to consider  three fundamental solutions, for example $Y_r(z,t)$ for $r=1,2,3$, and consequently  the Stokes matrices $\mathbb{S}_1(t)$ and $\mathbb{S}_2(t)$.  Moreover, chosen a solution $Y^{(0)}(z,t)$ with normal form (\ref{24gen2016-1}),  a  {\it central connection matrix} $C^{(0)}$ is defined by the relation 
\be
\label{2maggio2016-5}
Y_1(z,t)=Y^{(0)}(z,t)~C^{(0)}(t),
\quad
\quad
z\in \mathcal{S}_1(\overline{\mathcal{V}}). 
\ee    
Then, the  {\it  essential monodromy data} of the system  (\ref{ourcase12}) are defined to be
\be
\label{17april2016-1-tristezza}
\mathbb{S}_1(t) ,\quad  \mathbb{S}_2(t), \quad B_1(t)=\hbox{diag}(\widehat{A}_1(t)), \quad C^{(0)}(t), \quad J^{(0)}(t), \quad R^{(0)}(t).
\ee
   Now, when $t$ tends to a point $t_\Delta\in\Delta$, the limits of the above data may not exist. If the limits exist, they   {\it do not in general give the
   monodromy data  of the  system  $\widehat{A}(z,t_\Delta)$}. The latter  have in general  different nature, as it is clear from the results of \cite{BJL2}, and from Section \ref{31marzo2015-11} below.\footnote{See 
   for example the solution (\ref{7dic2015-1}), where it is evident that the monodromy datum $L$, defined at $t=0$,  
   is not the limit for $t\to 0$ of $B_1(t)$ as in (\ref{4gen2016-6}).}

           \bde  
      \label{19maggio2017-7}
        If the deformation is admissible in a  domain $\mathcal{V}$, as in Definition \ref{12maggio2017-2}, we say that it is {\bf isomonodromic in $\mathcal{V}$} if the essential monodromy data (\ref{17april2016-1-tristezza}) do not depend on $t\in \mathcal{V}$. 
   \ede     
        
        When this definition holds, the classical theory of Jimbo-Miwa-Ueno \cite{JMU} applies.\footnote{ Notice that in \cite{JMU} it is also assumed that $A_1(t)$ is {diagonalisable with  eigenvalues not differing by  integers}. We do not make this assumption here.} 
  We are interested in extending the deformation theory to the whole $\mathcal{U}_{\epsilon_0}(0)$, including the coalescence locus  $\Delta$. 
        
        \subsection{Main Results}

{\bf \large  a] The case of systems (\ref{22novembre2016-3})  and (\ref{17luglio2016-1}).}  Up to Section \ref{3marzo2017-2}, we  study system (\ref{22novembre2016-3})  without requiring that the deformation is isomonodromic.  First, we give the general formal and actual solutions for $z\to \infty$ of system (\ref{22novembre2016-3}) when  $t=0$ (or $t_\Delta\in\Delta$), without Assumption 1.\footnote{We give an explicit construction  of the formal solutions; their structure  can also be derived from  \cite{BJL2}.} Then,  under Assumption 1, in Proposition \ref{8dic2015-3}  we give necessary and sufficient   conditions such   that  the coefficients $F_k(t)$ of a formal solution of (\ref{22novembre2016-3})
\be
\label{3marzo2017-1}
Y_F(z,t)=G_0(t)\Bigl(I+\sum_{k=1}^\infty F_k(t) z^{-k}\Bigr) z^{B_1(t)}e^{\Lambda(t)z},\quad\quad t\not \in  \Delta,
\ee
 are actually  holomorphic also at $t\in \Delta$.  Notice that our result cannot be derived from \cite{Babbitt1} and \cite{Schafke1},  where holomorphic confluence for $t\to 0$ of formal solutions is studied, since  $\Lambda(t)z$ is in general not ``well-behaved'' (condition  (4.2) of \cite{Schafke1} is violated). In Section \ref{14feb2016-10}, we prove that fundamental solutions  $Y_r(z,t)$, $r\in\mathbb{Z}$,   of (\ref{22novembre2016-3}) can be $t$-analytically continued to a whole $\widetilde{\tau}$-cell containing the domain $\mathcal{V}$ of Definition \ref{12maggio2017-2}, preserving the asymptotic representation (\ref{3marzo2017-1}).  In Theorem \ref{pincopallino}  we give sufficient  conditions  such that fundamental solutions $Y_r(z,t)$,  together with  their Stokes matrices $ \mathbb{S}_r(t)$,   
  are actually   holomorphic also at  $\Delta$ and in the whole $\mathcal{U}_{\epsilon_0}(0)$,  in such a way that  {\it the  asymptotic representation $Y_r(z,t)\sim Y_F(z,t)$ continues to hold}, for $z\to \infty$  in wider  sectors $\widehat{\mathcal{S}}_r$ containing $\mathcal{S}_r(\overline{\mathcal{V}})$, to be introduced below (see (\ref{24maggio2017-1})).  We show in this case that the limits 
  \begin{shaded}\be
\label{3marzo2017-4}
\lim_{t\to t_\Delta} \mathbb{S}_r(t),  \quad \quad t_\Delta\in \Delta,
\ee
\end{shaded}
\noindent
exist and are finite. They  give the  Stokes matrices for the system (\ref{22novembre2016-3}) with matrix coefficient ${A}(z,t_\Delta)$    (see Corollary \ref{6marzo2016-2} and  \ref{19feb2016-1}). 
   In the analysis of the above issues,  wall crossing phenomena  and cell decompositions of $\mathcal{U}_{\epsilon_0}(0)$  will be  studied. Another result on the analytic ocntinuation of fundamental solutions, with vanishing conditions on the Stokes matrices, is given in Theorem \ref{8gen2017-4}.

  \vskip 0.2 cm 
  
We compare our results with  the existing  literature,  where sometimes  the irregular singular point is taken at $z=0$ (equivalent to $z=\infty$ by a change $z\mapsto 1/z$). One considers a ``folded'' system  $A(z,0)=z^{-k-1}\sum_{j=0}^\infty A_j(0)z^{j}$, with an irregular singularity of Poincar\'e rank $k$ at $z=0$ and studies its holomorphic  unfolding  $A(z,t)=p(z,t)^{-1}\sum_{j=0}^\infty A_j(t)z^{j}$, where $p(z,t)=(z-a_1(t))\cdots(z-a_{k+1}(t))$ is a polynomial. Early studies  on the relation between monodromy data of the ``folded'' and the ``unfolded'' systems  were started by Garnier \cite{Garnier}, and the problem was again raised by  V.I. Arnold in 1984  and studied by many authors in the '80's and '90's of the XX century, for example  see \cite{Ramis}, \cite{Duval},  \cite{Bolibruch2}.  Under suitable conditions, some results have been recently established  regarding  the convergence  for $t\to 0$ ($t$ in   sectors 
   or suitable ramified domains) of  fundamental solutions  and  monodromy data (transition or connection  matrices) of the ``unfolded'' system to the Stokes 
   matrices of the ``folded'' one  \cite{Ramis}, \cite{Duval},  \cite{Bolibruch2}, \cite{Babbitt1}, \cite{Schafke1}, \cite{Glut1}, \cite{Glut2},  \cite{Hurt-Lam-Rouss},  \cite{Lam-Rouss},  \cite{Klimes1}.  Nevertheless, to our knowledge, the   case when
    $A_0(0)$ is {\it diagonalisable with coalescing eigenvalues} has not yet been studied.  For example,    in \cite{Glut1} (see also references therein) and 
    \cite{Hurt-Lam-Rouss} \cite{Lam-Rouss}, it is assumed that  the leading matrix $A_0(0)$  has  
    distinct eigenvalues.  In \cite{Glut2},  $A_0(0)$ is a single
     $n\times n$  Jordan block (only one eigenvalue), with a generic condition on $A(z,t)$. Moreover,  the irregular singular point is required to split 
     into non-resonant Fuchsian singularities $a_1(t),...,a_{k+1}(t)$. The case when $A_0(0)$ is a $2\times 2$   Jordan block and $k=1$ 
      is thoroughly described in    \cite{Klimes1}, again under a generic condition on $A(z,t)$, with  no conditions on the polynomial $p(z,t)$. Explicit normal forms for the unfolded systems are given (including an explanation of the change of order of Borel summability when $z=0$ becomes a resonant irregular singularity  as $t\to 0$).  Nevertheless,  both in \cite{Glut2} and    \cite{Klimes1} 
       the system at $t=0$ is ramified and  the fundamental matrices  $Y_r(z,t)$ 
       diverge  when $t\to 0$, together with  the corresponding 
        Stokes matrices. Therefore,  our results  on the extension of the asymptotic representation at $\Delta$ and the  existence of the limit  (\ref{3marzo2017-4}),    for a system with diagonalisable $A_0(t_\Delta)$,  seem to be missing in the literature.   
   
      \vskip 0.3 cm 
        
 \noindent 
 {\bf \large  b] Isomonodromic  case of system (\ref{ourcase12}).} Let the deformation be isomonodromic in $\mathcal{V}$, as in Definition \ref{19maggio2017-7}, so that the classical theory of Jimbo-Miwa-Ueno applies. As a result of  \cite{JMU}, the eigenvalues  can be chosen as the independent deformation parameters.  This means that we can assume\footnote{ This assumption will be used in the paper  starting from  Section \ref{2dic2016-1}.} linearity   in $t\in \mathcal{U}_{\epsilon_0}(0)$, as follows:
    \be 
    \label{4gen2016-2}
    u_a(t)=u_a(0)+t_a, \quad 1\leq  a \leq n  \quad \Longrightarrow~m=n.
    \ee 
    Therefore, 
    $$ 
    \Lambda(t)=\Lambda(0)+\hbox{diag}(t_1,~...,~t_n)$$
        with 
  \be
  \label{17maggio2017-1}
     \Lambda(0)=\Lambda_1\oplus \cdots \oplus \Lambda_s, 
     \quad s< n,
  \quad \quad 
     \Lambda_i=\lambda_i I_{p_i},
\ee
    where $\lambda_1$, ..., $\lambda_s$ are the $s<n$ {\it distinct} eigenvalues of $\Lambda(0)$, of respectively multiplicities $p_1,...,p_s$ ($p_1+\cdots + p_s=n$). Here, $I_{p_i}$ is the $p_i\times p_i$ identity matrix. 
    Now, the size $\epsilon_0$ of $\mathcal{U}_{\epsilon_0}(0)$ is  taken  sufficiently small so  that we can write 
     \be
      \label{19maggio2017-6}
     \Lambda(t)=\Lambda_1(t)\oplus \cdots \oplus \Lambda_s(t),
     \ee
     with the properties that 
     $ 
\lim_{t\to 0}     \Lambda_j(t)=\lambda_j I_{p_j}
$
, and that 
 $\Lambda_i(t)$ and $\Lambda_j(t)$ have no common eigenvalues for $i\neq j$. 
 Thus,  $\Delta$ is represented as 
    $$ 
    \Delta=\bigcup_{i=1}^s \Delta_i,
    $$
    where 
    $$
    \Delta_i:=\{t\in\mathcal{U}_{\epsilon_0}(0)~|~t_a=t_b~\hbox{ with } u_a(0)=u_b(0)=\lambda_i\}.
    $$
    Our problem is  to extend the isomonodromy deformation theory from $\mathcal{V}$ to the whole $\mathcal{U}_{\epsilon_0}(0)$ in this case.

As it will be reviewed below after Theorem \ref{9gen2017-1},    the existing literature on isomonodromy deformations does not seem to solve our problem. 
  We give a solution in  the following Theorem \ref{16dicembre2016-1} and Corollary  \ref {17dicembre2016-2}  (equivalently,  see Theorem \ref{3nov2015-5} and Corollary \ref{17dicembre2016-5}  in the main body of the paper). 
  
  In order  state Theorem \ref{16dicembre2016-1} in a precise way, we need a last technical remark on the radius $\epsilon_0$ of the polydisc.  As explained above, $\epsilon_0$ is sufficiently small to ensure that   $\Lambda_i(t)$ has no eigenvalues in common with $\Lambda_j(t)$, for $i\neq j$ (see  (\ref{19maggio2017-6})). Moreover, we require that it satisfies the following  constraint
 \be
 \label{21feb2016-1-BBB}
 \epsilon_0<\min_{1\leq j\neq k\leq s} \delta_{jk},
 \ee
 where 
 $$ 
 \delta_{jk}:=\frac{1}{2}\min_{\rho\in\mathbb{R}}\Bigl\{   
 \left|
 \lambda_k-\lambda_j+i\rho \exp\{-i\tilde{\tau}\}
 \right| 
 \Bigr\}
 $$
(here $i$ is the imaginary unit). This condition has a geometrical  reason. If we represent $\lambda_1,...,\lambda_s$ in the same $\lambda$-plane, we can easily verify that the distance between the two parallel lines through $\lambda_j$ and $\lambda_k$ of angular direction $3\pi/2-\widetilde{\tau}$ is exactly $2\delta_{jk}$. Let us consider   Stokes rays  $\{z\in\mathcal{R}~|~\Re(z(u_a(t)-u_b(t))=0\}$ associated with couples $u_a(t)$, $u_b(t)$, $a,b\in\{1,2,...,n\}$,  such that $u_a(0)=\lambda_j$ and $u_b(0)=  \lambda_k$, with  $1\leq j\neq k \leq s$.  None  of these rays crosses the admissible directions $\widetilde{\tau}+k\pi$, $k\in\mathbb{Z}$, when $t$ varies in $\mathcal{U}_{\epsilon_0}(0)$ with $\epsilon_0$ as in (\ref{21feb2016-1-BBB}).  For a given $t$, let  $\mathfrak{R}(t)$ be the set of all the above rays for all $j\neq k$.   
We construct a sector $\widehat{\mathcal{S}}_r(t)$ containing  the ``half-plane'' $\Pi_r$ (defined above), and  extending up to the closest Stokes rays of $\mathfrak{R}(t)$ lying outside  $\Pi_r$. Clearly, $\widehat{\mathcal{S}}_r(t) \supset \mathcal{S}_r(t)$. Then, define 
\be
\label{24maggio2017-1}
\widehat{\mathcal{S}}_r:=\bigcap_{t\in \mathcal{U}_{\epsilon_0}(0)} \widehat{\mathcal{S}}_r(t).
\ee
By construction,  if $\epsilon_0$  is as in (\ref{21feb2016-1-BBB}), then  this sector has central opening angle greater than $\pi$.   Note that  $\mathcal{S}_r(\overline{\mathcal{V}}) \subseteq \widehat{\mathcal{S}}_r$. 
   
  \begin{shaded}
\bth
\label{16dicembre2016-1}
   Consider the system (\ref{ourcase12}), with  eigenvalues  of  $\Lambda(t)$ linear in $t$ as in   (\ref{4gen2016-2}), and with ${A}_1(t)$ holomorphic  on a closed polydisc $\mathcal{U}_{\epsilon_0}(0)$ centred at  $t=0$, with sufficiently small  radius $\epsilon_0$ as in (\ref{21feb2016-1-BBB}). Let  $\Delta$ be the coalescence locus  in   $\mathcal{U}_{\epsilon_0}(0)$,  passing through $t=0$.  
   Let the dependence on  $t$ be isomonodromic in a domain  $\mathcal{V}$ as in Definition \ref{19maggio2017-7}.  
    
    \vskip 0.2 cm 
    \noindent
    If the matrix entries of $\widehat{A}_1(t)$ satisfy in $\mathcal{U}_{\epsilon_0}(0)$ the vanishing conditions 
\be
\label{31luglio2016-1}
\Bigl( \widehat{A}_1(t)\Bigr)_{ab}=\mathcal{O}(u_a(t)-u_b(t)), ~~~~1\leq a\neq b \leq n,
\ee
whenever $u_a(t)$ and $u_b(t)$ coalesce as $t$ tends to a point of $\Delta$, then the following results hold: 
 
 \begin{itemize}
 \item
  The
 formal solution $Y_F(z,t)$ of (\ref{ourcase12}) as given in (\ref{4gen2016-6}) is  holomorphic on the whole $\mathcal{U}_{\epsilon_0}(0)$.  
 
 \item  The three  fundamental matrix solutions $Y_r(z,t)$,  $r=1,2,3$, of the system of (\ref{ourcase12}), which are  defined on $\mathcal{V}$,  with asymptotic representation $Y_F(z,t)$ for $z\to \infty$ in sectors $\mathcal{S}_r(\overline{\mathcal{V}})$ introduced above, can be $t$-analytically continued as single-valued holomorphic functions on $\mathcal{U}_{\epsilon_0}(0)$, with asymptotic representation 
 $$Y_r(z,t)\sim Y_F(z,t),
 \quad
 \quad
 z\to \infty \hbox{ in } \widehat{\mathcal{S}}_r,
 $$ 
 for any $t\in \mathcal{U}_{\epsilon_1}(0)$, and any $0<\epsilon_1<\epsilon_0$. In particular, they are defined at any $t_\Delta\in \Delta$ with asymptotic representation $Y_F(z,t_\Delta)$. The fundamental matrix solution $Y^{(0)}(z,t)$ is also $t$-analytically continued as a single-valued  holomorphic function on $\mathcal{U}_{\epsilon_0}(0)$
 
  \item The  constant Stokes matrices   $\mathbb{S}_1$, $\mathbb{S}_2$, and a central connection matrix $C^{(0)}$,  initially defined for  $t\in \mathcal{V}$,   are  actually {\rm globally defined  on $  \mathcal{U}_{\epsilon_0}(0)$}. They  coincide with  the Stokes and connection matrices   of the fundamental solutions  $Y_r(z,0)$ and $Y^{(0)}(z,0)$  of the system
\be
\label{4gen2016-8}
  \frac{dY}{dz}=\widehat{A}(z,0)Y,
  \quad 
  \quad 
  \widehat{A}(z,0)=\Lambda(0)+\frac{\widehat{A}_1(0)}{z}.
  \ee 
  Also the remaining $t$-independent  monodromy data in (\ref{17april2016-1-tristezza}) coincide with those of  (\ref{4gen2016-8}). 

  \item  The entries $(a,b)$ of the Stokes matrices are characterised by the following vanishing property: 
\be
\label{16dicembre2016-2}
(\mathbb{S}_1)_{ab}=(\mathbb{S}_1)_{ba}=(\mathbb{S}_2)_{ab}=(\mathbb{S}_2)_{ba}=0 
\quad \hbox{ whenever } u_a(0)=u_b(0),
\quad  1\leq a\neq b \leq n.
 \ee
 \end{itemize}
    \eth
      \end{shaded}

Theorem \ref{16dicembre2016-1}   allows to  holomorphically define the fundamental solutions and the monodromy  data on the whole $\mathcal{U}_{\epsilon_0}(0)$,  under the only condition (\ref{31luglio2016-1}).  This fact is remarkable. 
  Indeed, according to \cite{Miwa}, in general the solutions $Y^{(0)}(z,t)$, $Y_r(z,t)$ 
    and $\widehat{A}(z,t)$, $t\in \mathcal{V}$,  of monodromy preserving deformation equations  can be analytically continued as 
    meromorphic matrix valued functions  on the universal covering of 
    $\mathbb{C}^n\backslash \Delta_{\mathbb{C}^n}$, where $\Delta_{\mathbb{C}^n}=\bigcup_{a\neq b}^n \{u_a(t)=u_b(t)\}$ is  the coalescence 
    locus  in $\mathbb{C}^n$. They have  fixed singularities at the branching locus $\Delta_{\mathbb{C}^n}$, and so   at $\Delta\subset\Delta_{\mathbb{C}^n}$.  Moreover,  the $t$-analytic continuation  on $\mathcal{U}_{\epsilon_0}(0)$ of a the solutions $Y_r(z,t)$ 
 are expected to lose  their asymptotic representation $Y_r(z,t)\sim Y_F(z,t)$  in $\mathcal{S}_r(\overline{\mathcal{V}})$,   when $t$ moves sufficiently far away from $\mathcal{V}$, namely when Stokes rays cross and admissible ray of direction $\widetilde{\tau}$.   Under the assumptions of Theorem \ref{16dicembre2016-1} these singular behaviours do not occur. 

\vskip 0.2 cm

  Let  the assumptions of Theorem \ref{16dicembre2016-1} hold.  Then,  the system (\ref{4gen2016-8}) has  a  formal solutions (here we denote objects $Y$, $\mathbb{S}$ and $C$ referring to the system    (\ref{4gen2016-8})  with the symbols  $\mathring{Y}$, $\mathring{\mathbb{S}}$ and $\mathring{C}$) with behaviour\footnote{If the vanishing condition  (\ref{31luglio2016-1}) fails, formal solutions are more complicated (see Theorem \ref{31marzo2015-12}).}
\be
\label{5gen2016-1} 
\mathring{Y}_F(z)=\Bigl( 
I+\sum_{k=1}^\infty \mathring{F}_k z^{-k}
\Bigr) z^{B_1(0)} e^{\Lambda(0)z},  \quad \quad B_1(0)=\hbox{diag}(\widehat{A}_1(0)).
\ee 
The matrix-coefficients   $\mathring{F}_k$ are recursively constructed from the equation (\ref{4gen2016-8}), but not uniquely determined. Actually, there is a family of formal solutions  as above, depending on a finite number of complex parameters.  
To each element of the family,  there correspond unique actual solutions 
$\mathring{Y}_1(z)$, $\mathring{Y}_2(z)$, $\mathring{Y}_3(z)$ such that $\mathring{Y}_r(z)\sim \mathring{Y}_F(z)$ for $z\to \infty$ in a sector $\mathcal{S}_r\supset \mathcal{S}_r(\overline{\mathcal{V}})$, $r=1,2,3$, with  Stokes matrices defined by 
$$
\mathring{Y}_{r+1}(z)=\mathring{Y}(z)~\mathring{\mathbb{S}}_r,  \quad \quad r=1,2.
$$ 
Only one element of the family of formal solutions  (\ref{5gen2016-1})   satisfies the condition  $\mathring{F}_k=F_k(0)$ for any $k\geq 1$, and by  Theorem \ref{16dicembre2016-1} the relations  $\mathbb{S}_r=\mathring{\mathbb{S}}_r$ hold. 
Let us choose a solution $\mathring{Y}^{(0)}(z)$ close to $z=0$ in the Birkhoff-Levelt normal form, and define the corresponding  central connection matrix $\mathring{C}^{(0)}$ such that 
 $$
\mathring{Y}_{1}(z)=\mathring{Y}^{(0)}(z)~\mathring{C}^{(0)}.
$$  
We will  prove that the class of formal solutions (\ref{5gen2016-1}) reduces to only one element (thus the formal solution is unique) if and only if the diagonal entries of $\widehat{A}_1(0)$ do not differ by non-zero integers. This fact implies the following
  
  \begin{shaded}
  \bcr
  \label{17dicembre2016-2}
Let the assumptions of Theorem \ref{16dicembre2016-1} hold. If the diagonal entries of $\widehat{A}_1(0)$  do not differ by non-zero integers, then there is a unique formal solution (\ref{5gen2016-1}) of the system (\ref{4gen2016-8}), whose coefficients  necessarily  satisfy the condition 
$$
\mathring{F}_k\equiv F_k(0).
$$
 Hence,   (\ref{4gen2016-8})  only has at $z=\infty$  canonical fundamental solutions $\mathring{Y}_1(z)$,  $\mathring{Y}_2(z)$,  $\mathring{Y}_3(z)$,  which coincide with the canonical solutions $Y_1(z,t)$,    $Y_2(z,t)$,  $Y_3 (z,t)$ of (\ref{ourcase12})  evaluated at $t=0$, namely:
 $$
 Y_1(z,0)=\mathring{Y}_1(z),    
 \quad 
 Y_2(z,0)=\mathring{Y}_2(z),  
 \quad 
 Y_3 (z,0)=\mathring{Y}_3(z).$$
 Moreover, for any $\mathring{Y}^{(0)}(z)$ there exists $Y^{(0)}(z,t)$ such that  $Y^{(0)}(z,0)=\mathring{Y}^{(0)}(z)$. The following equalities hold:
 $$ 
 \mathbb{S}_1=\mathring{ \mathbb{S}}_1,  \quad  \mathbb{S}_2=\mathring{ \mathbb{S}}_2, \quad C^{(0)}=\mathring{C}^{(0)}.
 $$
   \ecr
  \end{shaded}

   Corollary  \ref{17dicembre2016-2}  has  a  practical computational  importance: the {\it constant monodromy data} (\ref{17april2016-1-tristezza})  of the system (\ref{ourcase12}) on    the whole  $\mathcal{U}_{\epsilon_0}(0)$ are computable just  by  considering the system (\ref{4gen2016-8}) at the coalescence point $t=0$.  This  is useful for applications in   the following two cases. 
     
  a)  When $\widehat{A}_1(t)$ is known in a whole neighbourhood of a coalescence point, but the computation of monodromy 
  data, which is  highly transcendental,  can be explicitly done (only)  at a coalescence point,  where    (\ref{ourcase12}) simplifies  due to (\ref{31luglio2016-1}). 
    An example is given in Section \ref{31luglio2016-2} for the sixth Painlev\'e equation.  Another example will be given in \cite{OURPAPER} for the $A_3$-Frobenius manifold. 
 
 b) When $\widehat{A}_1(t)$ is {\it explicitly known  only} at a coalescence point. This may happen in the case of Frobenius manifolds. So far, the theory of semisimple Frobenius manifolds has never been extended to semisimple coalescence points, which  appear frequently in important cases, such as for example the quantum cohomology of Grassmannians  \cite{Cotti1}, \cite{OURPAPER}. Our result is at the basis of the extension of the theory, as it will be thoroughly exposed in \cite{OURPAPER}.    Theorem \ref{16dicembre2016-1} and Corollary \ref{17dicembre2016-2}  allows  the computation of local moduli (monodromy data) of a semisimple Frobenius manifold just by considering a coalescence point.  The link between the present paper and  \cite{OURPAPER} will be established in Section \ref{16dicembre2016-5}.

\vskip 0.2 cm 
 In the present paper, we also  prove Theorem \ref{9gen2017-1} below, which is the converse of  Theorem \ref{16dicembre2016-1}.  Assume that the system is isomonodromic on a simply connected domain  $\mathcal{V}\subset \mathcal{U}_{\epsilon_0}(0)$ as in Definition \ref{12maggio2017-2}.  Note that now {\it we are not assuming} that  $\widehat{A}_1(t)$ is holomorphic in the whole $\mathcal{U}_{\epsilon_0}(0)$, contrary to what has been done so far.   As a result of  \cite{Miwa}, the fundamental solutions $Y_r(z,t)$, $r=1,2,3$, and $A_1(t)$  can be   analytically continued as 
    meromorphic matrix valued functions  on the universal covering of 
    $\mathcal{U}_{\epsilon_0}(0)\backslash \Delta$,  with  movable poles at the Malgrange divisor \cite{Palmer} \cite{Malg1} \cite{Malg2} \cite{Malg3}. The coalescence locus  $\Delta$ is in general  a fixed branching locus.  Moreover,  although for $t\in \mathcal{V}$ the fundamental  solutions $Y_r(z,t)$ have in $\mathcal{S}_r(\overline{\mathcal{V}})$  the canonical asymptotic behavior $Y_F(z,t)$ as in (\ref{4gen2016-6}),  in general this is no longer true when $t$ moves sufficiently far away from $\mathcal{V}$. 
  
  Nevertheless, if   the vanishing condition  (\ref{16dicembre2016-2}) on Stokes matrices holds, then  we can prove that the fundamental solutions $Y_r(z,t)$ and  $\widehat{A}_1(t)$ have {\it single-valued}  meromorphic continuation on  $\mathcal{U}_{\epsilon_0}(0)\backslash\Delta$,  so that {\it $\Delta$ is not a branching locus}. Moreover, the asymptotic behaviour is preserved, according to the following

 \begin{shaded}
 \bth
 \label{9gen2017-1}
 Let $\epsilon_0$ be as in (\ref{21feb2016-1-BBB}). Consider the  system (\ref{ourcase12}). Let the matrix   $\widehat{A}_1(t)$ be  holomorphic on an open simply connected domain $\mathcal{V}\subset\mathcal{U}_{\epsilon_0}(0)$ such that    the deformation is admissible and  isomonodromic  as in Definitions \ref{12maggio2017-2}  and \ref{19maggio2017-7}.  Assume that the entries of the   constant Stokes matrices  satisfy the vanishing condition
$$
(\mathbb{S}_1)_{ab}=(\mathbb{S}_1)_{ba}=(\mathbb{S}_{2})_{ab}=(\mathbb{S}_{2})_{ba}=0 ~~\hbox{ whenever } u_a(0)=u_b(0),~~1\leq a\neq b\leq n.
$$
Then, as functions of $t$,   the fundamental solutions $Y_{r}(z,t)$ and $\widehat{A}_1(t)$ admit single-valued  meromorphic continuation on $\mathcal{U}_{\epsilon_0}(0)\backslash\Delta$.   Moreover,     for any  $t\in \mathcal{U}_{\epsilon_0}(0)\backslash\Delta$ which is not a pole of $Y_{r}(z,\tilde{t})$ (i.e. which is not a point of the Malgrange divisor), we have 
 $$
 Y_{r}(z,t) ~\sim Y_F(z,t)~\hbox{ for $z\to \infty$ in $\widehat{\mathcal{S}}_{r}(t)$, }~~~r=1,2,3,
$$
 and 
$$ 
  Y_{r+1}(z,t)=Y_r(z,t)~\mathbb{S}_r,  \quad r=1,2.
  $$
   The $\widehat{\mathcal{S}}_{r}(t)$'s are the wide sectors described after the inequality (\ref{21feb2016-1-BBB}) above.   

\eth
  \end{shaded}

\vskip 0.2 cm 
We compare our results with the existing literature on isomonodromic deformations. 
The case when  $\Delta$ is empty  and  $\widehat{A}_1(t)$ is any matrix does not add additional  
difficulties 
 to the theory developed in  \cite{JMU}. Indeed,  in the definition  of isomonoromic deformations given above,  not only we require that the monodromy matrix at $z=0$ is independent of $t$,  but also the {\it monodromy exponents} $J^{(0)}$, $R^{(0)}$  and the connection matrix $C^{(0)}$ in (\ref{17april2016-1-tristezza}) are constant (this is an {\it isoprincipal deformation},  in the language of \cite{KV}). Given these conditions on the exponents, and assuming that $\Delta=\emptyset$, one can  essentially repeat the proofs given in \cite{JMU}. For example,  the case when $\Delta$ is empty and  $\widehat{A}_1(t)$   is skew-symmetric and diagonalisable has been studied in \cite{Dub1}, \cite{Dub2}.  We also recall that in case of  Fuchsian singularities only, isomonodromic deformations were completely studied\footnote{In    \cite{Bolibruch1} it is only assumed that the monodromy matrices are constant. This generates non-Schlesinger deformations. On the other hand, an isopricipal deformation always leads to Schlesinger deformations \cite{KV}.} in  \cite{Bolibruch1} and \cite{KV}. 

Isomonodromy deformations  at irregular singular points with    leading matrix admitting a Jordan form {\it independent of $t$}  were studied in \cite{BertMo}  (with some minor {\it Lidskii } generic conditions). For example, if the singularity is at $z=\infty$ as in (\ref{17luglio2016-1}), the results of \cite{BertMo} apply to  $\widehat{A}(z,t)=z^{k-1}(J+\sum_{j=1}^\infty \widehat{A}_j(t)z^{-j})$, with Jordan form  $J $ and Poincar\'e rank $k\geq 1$.  Although the eigenvalues of $J$  have in general algebraic multiplicity greater than 1, $J$ is ``rigid", namely  $u_1,...,u_n$ do not depend on $t$. 

Other investigations of isomonodromy deformations at irregular singularities can be found in \cite{fedo} and \cite{bibi}. Nevertheless, these  results do not apply to our coalescence problem. For example, the third  admissibility conditions of definition 10 of \cite{bibi} is not satisfied in our case. In  \cite{fedo}  the system with $A(z,t)=z^{r-1}B(z,t)$, $r\in \mathbb{Q}$, is considered, such that $B(\infty,t)$ has distinct eigenvalues;  $z=\infty$ satisfying this condition is called a simple irregular singular point. This simplicity condition does not apply in our case. 

The results of \cite{Klimes1}, cited above, are applied in \cite{Klimes2} to the $3\times 3$ isomonodromic
       description of the  Painlev\'e 6 equation and its coalescence to Painlev\'e 5. In this case, the limiting system for $t\to 0$ has leading matrix with a $2\times 2$ Jordan block, so that the fundamental matrices $Y_r(z,t)$ diverge. 

Isomonodromic deformations of a system such as our (\ref{ourcase12})  (with  $z\mapsto 1/z$, $\widehat{A}_0\mapsto Z$, $\widehat{A}_1\mapsto f$) appears also in \cite{BridgeTole}. Nevertheless, the deformations in Section  3  of  \cite{BridgeTole} are of a very particular kind. Indeed, the  eigenvalues  $u_1,...,u_n$  of the matrix $Z$ in \cite{BridgeTole}, which is the analogue of  our  $\widehat{A}_0$, are deformation parameters, but always satisfy the condition
\begin{align}
\label{17maggio2017-2}
u_1&=\dots=u_{p_1},\\
u_{p_1+1}&=\dots=u_{p_1+p_2},\\
\dots
\\
\label{17maggio2017-3}
u_{p_1+\cdots+p_{s-1}+1}&=\dots=u_{p_1+\cdots + p_s},
\end{align}
with $p_1+\cdots +p_s=n$. 
Thus, no splitting of coalescences occurs, so that the deformations are always inside the same "stratum" of the coalescence locus.  Moreover, the matrix $f=f(Z)$ in \cite{BridgeTole}, which is   the analogue of our  $\widehat{A}_1$,  satisfies quite restrictively requirements that the diagonal is zero and   $(\widehat{A}_1)_{ab}=0$ whenever $u_a=u_b$, $1\leq a\neq b \leq n$. 
 These conditions are always satisfied along the deformation ``stratum''  of  \cite{BridgeTole}; they are a particular case or the more general conditions of   Proposition \ref{6ott2015-5} in our paper below.   For these reasons, an adaptation of the classical Jimbo-Miwa-Ueno results \cite{JMU} (and those of \cite{Boalch1} for a connection on a G-bundle, with G a complex and reductive group) can be done \emph{verbatim}, in order to describe the isomonodromicity condition for such a very particular kind of deformations.
In the present paper, we studied \emph{general} isomonodromic deformations of the system (\ref{ourcase12}), not necessarily the simple decomposition of the spectrum as in (\ref{17maggio2017-2})-(\ref{17maggio2017-3}).

\subsection{Plan of the Paper}

\begin{itemize}

\item In Part I,  we study formal and  fundamental solutions of the system  (\ref{22novembre2016-3}) as $z\to \infty$, both at coalescence points and away from them. We give necessary and sufficient conditions for a formal solution, computed away from coalescence points, to admit  holomorphic continuation to the coalescence locus (see Proposition \ref{8dic2015-3}). 

\item In Part II, we study the Stokes phenomenon at $z=\infty$  for the system  (\ref{22novembre2016-3}), both at coalescence and non-coalescence points. We show existence and uniqueness results at coalescence points.  

\item In Part III,  under Assumption 1 we discuss the analytic continuation of fundamental solutions of 
(\ref{22novembre2016-3}). We show that $\mathcal{U}_{\epsilon_0}(0)$  splits into topological cells, determined by the fact 
that Stokes rays associated with $\Lambda(t)$  cross a fixed  admissible ray.  In Theorem 
\ref{pincopallino} and Corollary \ref{6marzo2016-2} we give sufficient conditions such that fundamental solutions can be 
analytically continued to  the whole $\mathcal{U}_{\epsilon_0}(0)$, preserving their asymptotic representation, so that 
the Stokes matrices admit the limits (\ref{3marzo2017-4}).  Notice that for the results  in Parts I--III no 
isomonodromicity is required.

\item  In Part IV, we formulate the monodromy preserving deformation theory  for system  (\ref{ourcase12}). We prove  Theorem \ref{16dicembre2016-1} ,  Corollary \ref{17dicembre2016-2} and  Theorem \ref{9gen2017-1}.  

\item In Part V, we show how  Theorem \ref{16dicembre2016-1} and  Corollary \ref{17dicembre2016-2} can be applied to Frobenius Manifolds and to the sixth Painlev\'e equation. 

\end{itemize}

\bre
In the main body of the paper, the matrices $Y_r$, sectors $\mathcal{S}_r$  and Stokes matrices $\mathbb{S}_r$  will be labelled differently as $Y_{\nu+(r-1)\mu}$,  $\mathcal{S}_{\nu+(r-1)\mu}$ and  $\mathbb{S}_{\nu+(r-1)\mu}$, $\nu, \mu \in\mathbb{Z}$. This labelling will be explained. 
\ere

\vskip 0.3 cm 
\noindent 
{\it Acknowledgements:} We thank Marco Bertola for helpful discussions concerning the proof of Theorem  \ref{9gen2017-1}. D. Guzzetti remembers with gratitude   Andrei Kapaev for insightful discussions at the time when  this work was initiated.

\vskip 0.5 cm
\hrule
\vskip 0.5 cm 

\centerline{\bf \Large PART I: Structure of Fundamental Solutions}

\vskip 0.3 cm 
{\bf Notational Remarks:} {\small 
If  $\alpha<\beta$ are real numbers, an open sector  and a closed sector with central opening angle $\beta-\alpha>0$ are respectively denoted by 
$$ 
 S(\alpha,\beta):=\bigr\{z\in \mathcal{R}~\bigl|  ~~\alpha< \arg z < \beta ~\bigr\},
 \quad \quad 
 \overline{S}(\alpha,\beta):=\bigr\{z\in \mathcal{R}~\bigl|  ~~\alpha\leq \arg z \leq \beta ~\bigr\}. 
$$
The rays with directions $\alpha$ and $\beta$   will be called {\it the  right and left boundary rays} respectively. If  $\overline{S}(\theta_1,\theta_2)\subset S(\alpha,\beta)$, then  $\overline{S}(\theta_1,\theta_2)$ is called a {\it proper (closed) subsector}. 

Given a  function $f(z)$ holomorphic on a sector containing  $\overline{S}(\alpha,\beta)$, we say that it admits an asymptotic expansion 
$
f(z)\sim \sum_{k=0}^\infty a_k z^{-k}$ for $z\to\infty \hbox{ in }\overline{S}(\alpha,\beta),
$
 if for any $m\geq 0$,  $\lim_{z\to \infty} {z^m} \bigl(f(z)-\sum_{k=0}^m a_k z^{-k}\bigr)= 0$,  $z\in \overline{S}(\alpha,\beta)$.  If $f$ depends on parameters $t$,     the asymptotic representation $f(z,t)\sim\sum_{k=0}^\infty a_k(t)z^{-k}$ is said to be uniform in $t$   belonging to a compact subset $K\subset \mathbb{C}^m$, if the limits above are uniform in $ K$.  
In  case  the sector is open, we write $f(z)\sim \sum_{k=0}^\infty a_k z^{-k}$ as $z\to \infty$ in  $S(\alpha,\beta)$ if the limits above are zero  {\it in every proper closed  subsector}  of $S(\alpha,\beta)$.  
 When we take the limits above for matrix valued functions  $A=(A_{ij}(z,t))_{i,j=1}^n$,   we use the norm $|A|:=\max_{ij}|A_{ij}|$. $\Box$ 

 }

 \section{Deformation of a Differential System with Singularity of the Second  Kind}
\label{preRAY}

\vskip 0.3 cm 
  We consider system  (\ref{22novembre2016-3}) of the Introduction, namely 
\be 
\label{9giugno-2}
{dY\over dz}=A(z,t) Y,
\quad \quad
t=(t_1,t_2,...,t_m)\in \mathbb{C}^m,
\ee
depending on $m$  complex parameters\footnote{Later, we will take $n=m$, as in (\ref{4gen2016-2}).}  $t$. The $n\times n$  matrix $A(z,t)$ is holomorphic in $(z,t)$   for $|z|\geq N_0>0$  and $|t|\leq  \epsilon_0$, for some positive constants $N_0$ and $\epsilon_0$, with uniformly convergent Taylor expansion
\be
\label{5gen2016-3}
A(z,t)=\sum_{j=0}^\infty A_j(t)z^{-j}.
\ee
The coefficients $A_j(t)$ are holomorphic for $|t|\leq \epsilon_0$.  We assume that  $A_0(0)$ is diagonalisable,   with {\it  distinct eigenvalues} $\lambda_1$, ..., $\lambda_s$, $s\leq n$. We are interested in the case when  $s$ is strictly less than $n$.  Up to a constant gauge transformation, there is no loss of generality in assuming that 
\be
\label{31marzo2015-8}
A_0(0)=\Lambda:=\Lambda_1\oplus\cdots\oplus \Lambda_s, \quad \quad \Lambda_i:=\lambda_i I_{p_i},  \quad i=1,2,...,s\leq n,
\ee
  being $I_{p_i}$ the $p_i\times p_i$ identity matrix.  If  $A_0(t)$ is holomorphically similar to $\Lambda(t)$, as in (\ref{3gen2017-1}), then $\Lambda=\Lambda(0)$. However, at this stage of the discussion  we do not  assume holomorphic similarity, so we keep the notation $\Lambda$ instead of  $\Lambda(0)$. 

\bre

A result due to Kostov \cite{Kostov} states that, if system (\ref{9giugno-2}) is such that $A(z,0)= A_0(0)+A_1(0)/z$, and if the matrix $A_1(0)$ has no eigenvalues differing by a non-zero integers, than there exists a gauge transformation $Y=W(z,t)\widetilde{Y}$, with  $W(z,t)$ holomorphic  at $z=\infty$ and $t=0$, such that   (\ref{9giugno-2}) becomes a system like (\ref{1novembre2016-1}):
\be
\label{13giugno2017-1}
 \frac{d\widetilde{Y}}{dz}=\left(\widetilde{A}_0(t)+\frac{\widetilde{A}_1(t)}{z} \right) \widetilde{Y}.
\ee
Nevertheless,  since $A_0(0)$ has non-distinct eigenvalues, we cannot find in general  a gauge transformation holomorphic at $z=\infty$ which transforms $A(z,0)$ of  the system (\ref{9giugno-2})  into  $ A_0(0)+A_1(0)/z$ (see also \cite{Bolibruch2}  and references therein).  Therefore the system (\ref{9giugno-2}) --  namely the system (\ref{22novembre2016-3}) --   is more general than system (\ref{13giugno2017-1}), namely than (\ref{1novembre2016-1}). 

\ere

\subsection{Sibuya's Theorem}
\label{17dicembre2016-1}
General facts  about eigenvalues and eigenvectors of a matrix $M(t)$, depending holomorphically on $t$ in  a domain $\mathcal{D}\subset \mathbb{C}^m$,  such that $M(0)$ has eigenvalues  $\lambda_1$, ..., $\lambda_s$, $s\leq n$, can be found in  \cite{Lax} and at page 63-87 of \cite{Kato}.  If $s$ is strictly smaller than $n$, then $t=0$ is a coalescence point. 
For $\mathcal{D}\subset   \mathbb{C}^m$ and $m=1$ the coalescence points are isolated, while for  $m\geq 2$ they form the {\it coalescence locus}. 
  Except for the special case when $M(t)$ is   holomorphically  similar to a Jordan form  $J(t)$, which means that there exists an invertible holomorphic matrix $G_0(t)$ on $\mathcal{D}$ such that $(G_0(t))^{-1} M(t) G_0(t)=J(t)$, in general the eigenvectors of $M(t)$ are holomorphic in the neighborhood of a non-coalescence point, but  their analytic continuation  is singular at the coalescence locus.  For example,
$$
 M(t)=
 \left(
 \begin{array}{cc}
  0 & 1 
  \\ t & 0
  \end{array}
  \right),\quad\quad t\in\mathbb{C},
  $$
  has eigenvalues $\lambda_\pm = \pm\sqrt{t}$, which are branches of $f(t)=t^{1/2}$, with ramification  at $\Delta=\{t=0\}$. The  eigenvectors can be chosen to be either
  $$ 
   \vec{\xi}_\pm=(\pm1/\sqrt{t}, 1),~~~\hbox{ or }~ \vec{\xi}_\pm=(\pm 1, \sqrt{t}).
  $$
  The matrix $G_0(t):=[\vec{\xi}_+(t),\vec{\xi}_-(t)]$ puts $M(t)$ in diagonal form $G_0(t)^{-1}A_0(t) G_0(t)={\rm diag}(\sqrt{t},-\sqrt{t})$, for $t\neq 0$, while $M(0)$ is in Jordan non-diagonal form. Either $G_0(t)$ or  $G_0(t)^{-1}$ is singular at $t=0$. The branching could be eliminated by changing  deformation parameter to $s=t^{1/2}$. Nevertheless, this would not cure the singularity of $G_0$ or $G_0^{-1}$ at $s=0$.
   Another example is  
$$ 
M(t)=
\left(
\begin{array}{cc}
1 & t 
\\ 
0 & 1
\end{array}
\right),
\quad \quad t\in\mathbb{C}. 
$$
  The eigenvalues $u_1=u_2=1$ are always coalescing. The Jordan types at $t\neq 0$ and $t=0$ are different. Indeed, $ 
M(0)={\rm diag}(1,1)$, while  for $t\neq 0$, 
$$ 
G_0(t)^{-1}M(t) ~G_0(t)=  \left(
\begin{array}{cc}
1 & 1 
\\ 
0 & 1
\end{array}
\right), \quad  G_0(t):= \left(
\begin{array}{cc}
t & 0
\\ 
0 & 1
\end{array}
\right). 
$$ 
Now,  $G_0(t)$ is not invertible and $G_0(t)^{-1}$ diverges  at $t=0$.  

In the above examples, the Jordan type of $M(t)$ changes. In the next example, the Jordan form remains diagonal, and nevertheless $G_0(t)$ is singular. Consider 
$$
M(t)=
\left( \begin {array}{cc} 1+t_1& t_2
\\ 0&1
-t_2\end {array} \right), \quad \quad t=(t_1,t_2)\in\mathbb{C}^2.
$$
The eigenvalues coalesce at $t=0$, where $M(0)=I$. Moreover, there exists a diagonalizing matrix  $G_0(t)$ such that  
$$
G_0(t)^{-1} M(t) G_0(t)=\left( \begin {array}{cc} 1+t_1& 0
\\ 0&1
-t_2\end {array} \right) \hbox{ is diagonal},
\quad
G_0(t)=\left( \begin {array}{cc} a(t) &- t_2~b(t)
\\ 0&(t_1+t_2)~b(t)
\end {array} \right),
$$
for arbitrary non-vanishing holomorphic functions $a(t), b(t)$. At $t=0$ the matrix $G_0(t)$ has zero determinant and $G_0(t)^{-1}$ diverges. 

 Although  $M(t)$ is not in general holomorphically similar to a Jordan form, holomorphic similarity can always be realised  between $M(t)$   and a block-diagonal matrix $\widehat{M}(t)$ having the same block structure of a Jordan form of $M(0)$, as follows.

\noindent
\ble
\label{SibLEM1}
{\bf [LEMMA 1 of \cite{Sh4}]:}  Let $M(t)$ be a $n\times n$ matrix holomorphically depending on  $t\in \mathbb{C}^m$,  with $|t|\leq \epsilon_0,
$  
where $\epsilon_0$ is a positive constant.  Let $\lambda_1,\lambda_2,...,\lambda_s$ be the distinct eigenvalues of $M(0)$, with multiplicities $p_1, p_2,...,p_s$, so that $p_1+p_2+\cdots +p_s=n$. Assume that $M(0)$ is in Jordan form 
$$ 
 M(0)=M_1(0)\oplus\cdots \oplus M_s(0) 
 $$
 where 
 $$ 
 M_j(0)=\lambda_j I_{p_j}+\mathcal{H}_j, 
 \quad \quad
 \mathcal{H}_j= \left[
 \begin{array}{cccccc}
 0& \mathfrak{h}_{j1} &  & && \\
                   &0 &\mathfrak{h}_{j2}     &        & &\\
                    & & \ddots & \ddots & &\\
                    &          &         & 0 &\mathfrak{h}_{jp_j-1}\\
                    &           &         &  & 0 
 \end{array}
\right], \quad 1\leq j \leq s,
$$
$\mathfrak{h}_{jk}$ being equal to 1 or 0. Then, for sufficiently small $0<\epsilon\leq \epsilon_0$ there exists a matrix $G_0(t)$, holomorphic in $t$ for $|t|\leq \epsilon$, such that 
$$ 
G_0(0)=I,
$$ 
and $
\widehat{M}(t)=(G_0(t))^{-1}M(t)G_0(t)$ has block diagonal  form 
\be 
\label{SH4continuous}
\widehat{M}(t)=\widehat{M}_{1}(t)\oplus\cdots\oplus \widehat{M}_{s}(t),
\ee 
where $\widehat{M}_{j}(t)$ are $p_j \times p_j$ matrices. For $|t|\leq \epsilon$,   $\widehat{M}_{i}(t)$ and $\widehat{M}_{j}(t)$ have no common eigenvalues for any  $i\neq j$.
\ele

\bre
The lemma also  holds when $t\in\mathbb R^m$  in the continuous (not necessarily holomorphic) setting. 
\ere

  Lemma \ref{SibLEM1} can be applied to  $M(t)\equiv A_0(t)$ in (\ref{5gen2016-3}), with $A_0(0)=\Lambda$. Therefore\footnote{Given a $n\times n$ matrix $A_0$, partitioned into $s^2$ blocks ($s\leq n$),  we use the notation $A_{ij}^{(0)}$, $1\leq i,j\leq s$,   to denote the block  in position $(i,j)$.  Such a block has dimension $p_i\times p_j$, with $p_1+...+p_n=n$.}  
\be
\label{25set2015-3}
\widehat{A}_0(t):=G_0(t)^{-1}A_0(t)G_0(t)=\widehat{A}_{11}^{(0)}(t)\oplus\cdots\oplus \widehat{A}_{ss}^{(0)}(t),
\ee
$$
G_0(0)=I, \quad \quad
\widehat{A}_0(0)=A_0(0)=\Lambda.
$$

\bre 
\label{5dic2015-1}
  $G_0(t)$ is determined up to $G_0\mapsto G_0(t)\Delta_0(t)$, where $\Delta_0(t)$ is any  block-diagonal matrix solution of  $[\Delta_0(t),\widehat{A}_0(t)]=0$.   Sibuya's normalization condition $G_0(0)=I$ can be softened to $ 
G_0(0)=\Delta_0$. 
\ere

We define a family of sectors $\mathcal{S}_\nu$ in $\mathcal{R}$ and state Sibuya's theorem.  Let $\arg_p(\lambda_j-\lambda_k)$ be the principal determination. 
  Let $\eta\in \mathbb{R}$ be  an {\bf admissible direction} for $\Lambda$ in the $\lambda$-plane  (we borrow this name and the following definition of the $\eta_\nu$'s and $\tau_\nu$'s  from  \cite{BJL1} and \cite{BJL4}). By definition, this means that, 
$$
 \eta\neq \arg_p (\lambda_j-\lambda_k) \hbox{ mod}(2\pi),  \quad \forall ~1\leq j\neq k \leq s.
$$
Introduce another determination $\widehat{\arg}$  as  follows:
\be 
\label{erni}
 \eta-2\pi < \widehat{\arg} (\lambda_j-\lambda_k)<\eta, \quad  \quad  1\leq j\neq k \leq s.
 \ee

 Let  $2\mu$, $\mu\in \mathbb{N}$, be the  number of values $\widehat{\arg} (\lambda_j-\lambda_k)$,  when $(j,k)$ spans all the indices $1\leq j\neq k \leq s$.\footnote{ 
 $
 2\mu \leq s(s-1) 
 $,
  with ``$=$" occurring when   $\arg (\lambda_j-\lambda_k)\neq  \arg (\lambda_r-\lambda_s)$ mod $2\pi$  for any $(j,k)\neq (r,s)$.}   
 Denote the $2\mu$ values of $\widehat{\arg}(\lambda_j-\lambda_k)$ with   $\eta_0,\eta_1,...,\eta_{2\mu-1}$,  according to the following ordering:
\be 
\label{deteeta}
 \eta> \eta_0>~\cdots~>\eta_{\mu-1}>\eta_\mu>~\cdots~> \eta_{2\mu-1}>\eta-2\pi.
 \ee
  Clearly
 \be 
 \label{deteeta1}
  \eta_{\nu+\mu}=\eta_\nu-\pi, 
  \quad \quad
  \nu=0,1,...,\mu-1.
  \ee
Consider the following  directional angles in the $z$-plane
\be
\label{TAU} 
\tau:=   {3\pi\over 2} -\eta,
\quad \quad
\tau_\nu:= {3\pi \over 2}-\eta_\nu, \quad \quad 0\leq \nu \leq 2\mu -1.
\ee 
From (\ref{deteeta}) if follows that,
\be 
\label{TAUn}
\tau< \tau_0<~\cdots~ < \tau_{\mu-1}<\tau_\mu<~\cdots~<\tau_{2\mu-1}<\tau+2\pi.
\ee
From (\ref{deteeta1}) if follows that,
$$ 
\tau_{\nu+\mu}= \tau_\nu+\pi,  \quad \quad \nu=0,1,...,\mu-1. 
$$ 
The extension of the above to directions in $\mathcal{R}$ is obtained by the following definition:
$$
\tau_{\nu+k\mu}:=\tau_\nu+k\pi, \quad \quad k\in\mathbb{Z}.
$$
This allows to speak of directions $\tau_\nu$ for any $\nu\in\mathbb{Z}$. 

\bde[Sector $\mathcal{S}_\nu$] 
We define the following   sectors of central opening angle greater than $\pi$:
\be
\label{SECT00bis}
\mathcal{S}_\nu:=S\bigl(\tau_\nu-\pi,\tau_{\nu+1}\bigr)\equiv S\bigl(\tau_{\nu-\mu},\tau_{\nu+1}\bigr), \quad \quad 
\nu\in\mathbb{Z}.
\ee
\ede
%
%

\begin{shaded} 
 \bth[Sibuya \cite{Sh4} \cite{Sh3}]
 \label{ShybuyaOLD}
 Let $A(z,t)$ be  holomorphic in $(z, t)$ for $|z|\geq  N_0 > 0$ and $|t| \leq  \epsilon_0$ as in (\ref{5gen2016-3}), such that $ 
A_0(0)=\Lambda=\Lambda_1\oplus\cdots\oplus \Lambda_s
$, as in (\ref{31marzo2015-8}).  
 Pick up a sector  $\mathcal{S}_\nu$ as in (\ref{SECT00bis}). Then, for any proper closed subsector $
\overline{S}(\alpha,\beta)=\{z~|~\tau_\nu-\pi<\alpha\leq \arg z\leq \beta <\tau_{\nu+1}\}\subset \mathcal{S}_\nu 
$, 
there exist a sufficiently large positive number $N\geq N_0$, a sufficiently small positive number $\epsilon \leq \epsilon_0$,  and  matrices $G_0(t)$ and ${ G}(z,t)$ with the following properties:

i) $G_0(t)$ is holomorphic for $|t|\leq \epsilon$ and
$$
G_0(0)=I,
\quad\quad
\widehat{A}_0(t):=G_0(t)^{-1}A_0(t)G_0(t)  \hbox{ is  block-diagonal as in (\ref{25set2015-3})}.
$$

 ii) $G(z,t)$ is holomorphic in $(z,t)$ for $|z|\geq N$, $z\in \overline{S}(\alpha,\beta)$, $|t|\leq \epsilon$;

 iii) ${ G}(z,t)$  has a uniform asymptotic expansion   for $|t|\leq \epsilon$, with holomorphic coefficients  $G_k(t)$:
 $$ 
 {G}(z,t)\sim I+\sum_{k=1}^\infty  {G}_k(t)z^{-k},
 \quad z\to \infty, 
 \quad
 z\in \overline{S}(\alpha,\beta),
 $$

 iv) The gauge transformation 
 $$ 
  Y(z,t)=G_0(t){G}(z,t)\widetilde{Y}(z,t),
  $$ 

  reduces the initial system to a block diagonal form
\be 
\label{9giugno-1}
  {d\widetilde{Y}\over dz}=
{B}(z,t)\widetilde{Y},
\quad
\quad
B(z,t)=B_1(z,t)\oplus\cdots\oplus B_s(z,t),
\ee

where ${B}(z,t)$ is holomorphic in $(z,t)$ in the domain $|z|\geq N$, $ z\in \overline{S}(\alpha,\beta)$, $|t|\leq \epsilon$, 

and has a uniform asymptotic expansion   for $|t|\leq \epsilon$, with holomorphic coefficients  $B_k(t)$, 
\be
\label{1gen2016-1}
{B}(z,t)\sim  \widehat{A}_0(t)+\sum_{k=1}^\infty {B}_k(t) z^{-k}, 
\quad
z\to \infty, 
\quad
z\in \overline{S}(\alpha,\beta).
\ee
In particular, setting  $\widehat{A}_1(t):=G_0^{-1}(t)A_1(t)G_0(t)$, then  $
B_1(t)=\widehat{A}_{11}^{(1)}(t)\oplus\cdots \oplus \widehat{A}_{ss}^{(1)}(t)$. 
     \eth
  \end{shaded}

\bre In the theorem above, $\epsilon$ is such that   $\widehat{A}_{ii}^{(0)}(t)$ and $\widehat{A}_{jj}^{(0)}(t)$ have no common eigenvalues for any  $i\neq j$ and  $|t|\leq \epsilon$.  Observe that one can always choose $
\beta-\alpha>\pi$. 
\ere
\bre
 $\mathcal{S}_\nu$ coincides with a sector 
  $\left\{
~z\in \mathcal{R}~\left|-3\pi/2-\omega_{-}<r \arg z<3\pi/2-\omega_+\right\}\right.
$, 
 introduced by Sibuya in \cite{Sh3}. A closed subsector $\overline{S}(\alpha,\beta)$  is  a sector $\mathcal{D}(N,\gamma)$ introduced by Sibuya in  \cite{Sh4}.
\ere

\bre
If $\Lambda=\lambda_1 I$, Theorem \ref{ShybuyaOLD} gives no new information, being  $G_0(t)=G(z,t)\equiv I$ and $\mathcal{S}_\nu=\mathcal{R}$.  
\ere


\noindent
{\it {\bf --} A Short Review of the Proof:} 
 The $z$-constant gauge transformation  
$ 
Y(z,t)=G_0(t) \widehat{Y}(z,t) 
$ 
transforms  (\ref{9giugno-2}) into
\be
\label{SeK2}
{d\widehat{Y}\over dz}=\widehat{A}(z,t)~\widehat{Y},
\quad\quad
\widehat{A}(z,t)=\sum_{i=0}^\infty \widehat{A}_i(t) z^{-i},\quad\quad
\widehat{A}_i(t):=G_0^{-1}(t)A_i(t)G_0(t).
\ee
Another gauge transformation   $ 
\widehat{Y}(z,t)=G(z,t)\widetilde{Y}(z,t)
$
yields  (\ref{9giugno-1}). Substitution into (\ref{SeK2}) gives the differential equation
\be 
\label{computblocked}
G^\prime+GB=~\widehat{A}(z,t)G,
\ee
with unknowns $G(z,t),B(z,t)$. If  formal series $G(z,t)= I +\sum_{j=1}^\infty G_j(t) z^{-j}$ and $B(z,t)=\widehat{A}_0(t)+\sum_{j=1}^\infty B_j(t)z^{-j}$ are inserted into  (\ref{computblocked}), the following recursive equations ($t$ is understood) are found:
\vskip 0.2 cm 
\noindent
\underline{For $l=0$}: ~ $B_0(t)=\widehat{A}_0(t).
$

\vskip 0.2 cm 
\noindent
\underline{For $l=1$}: 
\be 
\label{rec27marzo1}
\widehat{A}_0G_1-G_1\widehat{A}_0=-\widehat{A}_1+B_1.
\ee
\noindent
\underline{For $l\geq 2$}: 
\be 
\label{rec27marzo2}
\widehat{A}_0G_l-G_l\widehat{A}_0=\Bigr[\sum_{j=1}^{l-1}\Bigl(G_jB_{l-j}-\widehat{A}_{l-j}G_j\Bigr)-\widehat{A}_l\Bigr]~-(l-1)G_{l-1}~+B_l.
\ee
Once $G_0(t)$ has been fixed,  the recursion equations can  be solved. A solution $\{G_l(t)\}_{l=1}^\infty$, $\{B_l(t)\}_{l=1}^\infty$  is not unique in general. 
The following choice is possible: 
\be
\label{cho1}
 G^{(l)}_{jj}(t)=0, \quad
 \quad
 1\leq j \leq s,
 \quad
 \quad
 \hbox{[diagonal blocks are zero]},
\ee
and 
\be
\label{cho2}
B_l(t)=B^{(l)}_{1}(t)\oplus\cdots\oplus B^{(l)}_{s}(t),
\quad \quad
\hbox{[off-diagonal blocks are zero]}.
\ee
Then, the  $G_l(t)$'s and $B_l(t)$'s are determined by the recursion relations, because for a diagonal block $[j,j]$  the l.h.s of (\ref{rec27marzo1}) and (\ref{rec27marzo2}) is equal to 0 and  the r.h.s determines the  only unknown  variable  $B_{jj}^{(l)}$. For off-diagonal blocks $[i,j]$ there is no unknown in the r.h.s while in the l.h.s  the following expression appears
$$ 
\widehat{A}_{ii}^{(0)}(t)G_{ij}^{(l)}-G_{ij}^{(l)}\widehat{A}_{jj}^{(0)}(t), \quad \quad 1\leq i\neq j \leq s.
$$
{\it For $|t|\leq \epsilon $ small enough, $\widehat{A}_{ii}^{(0)}(t)$ and $\widehat{A}_{jj}^{(0)}(t)$ have no common eigenvalues}, so the equation is solvable for $G_{ij}^{(l)}$. 
 With the above choice, Sibuya \cite{Sh4}  proves  that  there exist actual solutions $G(z,t)$ and $B(z,t)$ of  (\ref{computblocked}) with  asymptotic expansions $I+\sum_j G_j(t)z^{-j}$ and $\widehat{A}_0+\sum_jB_j(t)z^{-j}$ respectively.  We remark that the proof relies on the above choice. It is evident that this choice also ensures that all the coefficients $G_j(t)$'s and $B_j(t)$'s are holomorphic where the $\widehat{A}_j(t)$'s are. 
 Note that (\ref{rec27marzo1}) yields $
B_1(t)=\widehat{A}_{11}^{(1)}\oplus\cdots \oplus \widehat{A}_{ss}^{(1)}(t)
$. $\Box$

   
   
   \section{Fundamental Solutions of (\ref{9giugno-1})} 
   \label{31marzo2015-11}

The system (\ref{9giugno-1}) admits  block-diagonal fundamental solutions $ 
\widetilde{Y}(z,t)=\widetilde{Y}_1(z,t)\oplus\cdots \oplus \widetilde{Y}_s(z,t)$. Here, 
 $\widetilde{Y}_i(z,t)$ is a $p_i\times p_i$ fundamental matrix of the $i$-th diagonal block of (\ref{9giugno-1}). 
The problem  is reduced to solving a system whose leading matrix has only one eigenvalue. The  case when $A_0(t)$  has distinct eigenvalues for  $|t|$ small is well known (see  \cite{Sh3}, and also \cite{BJL1} for the $t$-independent case).  
 The case when $A_0(0)=\Lambda$ is diagonalisable,  with $s\leq n$  distinct eigenvalues, will be studied here and in the subsequent sections.  \\
We do another  gauge transformation
  \be
\label{redDU}
  \widetilde{Y}(z,t)=e^{\Lambda z}~ Y_{red}(z,t),
  \ee
  where the subscript {\it red} stand for ``rank reduced". We substitute into (\ref{9giugno-1}) and find 
  $$
 \cancel{ e^{\Lambda z}}(\Lambda Y_{red}+Y_{red}^\prime)
 =B(z,t)\cancel{ e^{\Lambda z}}Y_{red}.
$$  
  The exponentials cancel because  $B(z,t)$ is block diagonal with the same structure as $\Lambda$. Thus, we obtain
  \be 
  \label{31marzo2015-2}
  {dY_{red}\over dz}= \frac{1}{z}B_{red}(z,t) ~Y_{red},
  \ee
with 
  \bea
&&  B_{red}(z,t):=z(B(z,t)-\Lambda)~=B_1^{(red)}(z,t)\oplus\cdots\oplus B_s^{(red)}(z,t),
\\
&&
B_{red}(z,t)\sim  z(\widehat{A}_0(t)-\Lambda)+\sum_{k=1}^\infty B_k(t) z^{-k+1}.
\eea
 Fundamental solutions can be taken with block diagonal structure, 
  $$
  Y_{red}(z,t)=Y^{(red)}_1(z,t)\oplus \cdots \oplus Y^{(red)}_s(z,t).
  $$
   where  $Y^{(red)}_i(z,t)$ solves 
  $$
  {dY^{(red)}_i \over dz}
  = \frac{1}{z}B^{(red)}_{i}(z,t) ~Y^{(red)}_i .
  $$
 The exponential $e^{\Lambda z}$  commutes with the above matrices, hence  a fundamental solution of  (\ref{9giugno-2}) exists in the form
$$ 
Y(z,t)=G_0(t)G(z,t)Y_{red}(z,t)~e^{\Lambda z}.
$$
 
 We proceed as follows. In Section \ref{hosforato} we describe the structure of fundamental solutions of (\ref{9giugno-2}) for  $t=0$ fixed. In Section \ref{11marzo2016-1} we describe the structure of fundamental solutions at other points $t\in\mathcal{U}_{\epsilon_0}(0)$.

  \section{A Fundamental Solution of (\ref{9giugno-2}) at $t=0$}
  \label{hosforato}

At $t=0$, {\it the rank is reduced}, since the system  (\ref{31marzo2015-2}) becomes a Fuchsian system in $\overline{S}(\alpha,\beta)$,
    \be 
 \label{PRE31marzo2015-6}
  {dY_{red}\over dz}= {1\over z} B_{red}(z,0) ~Y_{red},
  \ee
with $
B_{red}(z,0)\sim \sum_{k=1}^\infty B_k(0) z^{-k+1}$ for $z\to \infty$ in $\overline{S}(\alpha,\beta)$. 
 Let    $J_i$ be a  Jordan form of the $i$-th block $B_i^{(1)}(0)=\widehat{A}^{(1)}_{ii}(0)\equiv A^{(1)}_{ii}(0)$,    $1\leq i \leq s$. Following  \cite{Wasow}, we choose  $J_i$ arranged  into  $h_i\leq p_i$  Jordan blocks $J_1^{(i)}$, ..., $
J_{h_i}^{(i)}$
\be
\label{24settembre2015-1}
J_i= J_1^{(i)}\oplus \cdots \oplus
J_{h_i}^{(i)}.
\ee    
Each block $J_j^{(i)}$, $1\leq j\leq h_i$,  has dimension $r_j\times r_j$, with $r_j\geq 1$, $r_1+\cdots+r_{h_i}=p_i$.  Each $J_j^{(i)}$ has only one eigenvalue $\mu_j^{(i)}$, with structure,
$$
 J_j^{(i)}=\mu_j^{(i)} I_{r_j} +H_{r_j},
 \quad
 \quad
 I_{r_j}=\hbox{ $r_j\times r_j$ identity matrix},
 $$
  $$ \hbox{\rm $H_{r_j}=0$ ~if $r_j=1$}, \quad \quad 
 H_{r_j}=\left[
 \begin{array}{cccccc}
 0& 1 &  &  \cr
                   &0 &1      &        & &\cr 
                    & & \ddots & \ddots & &\cr
                    &          &         & 0 &1\cr
                    &           &         &  & 0 
 \end{array}
\right]~\hbox{  if $r_j\geq 2$}.
$$
Note that $\mu_1^{(i)}$, ..., $\mu_{h_i}^{(i)}$ are not necessarily distinct.  
One can choose a $t$-independent matrix $\Delta_0=\Delta_1^{(0)}\oplus\cdots\oplus \Delta_s^{(0)}$, in the block-diagonal of Remark \ref{5dic2015-1}, such that $
(\Delta_i^{(0)})^{-1}\widehat{A}^{(1)}_{ii}(0)\Delta_i^{(0)}=J_i$. 
Hence, 
$$ 
\Delta_0^{-1} \widehat{A}_1(0) \Delta_0\equiv \Delta_0^{-1} A_1(0) \Delta_0=\left[
\begin{array}{cccc}
J_1 &*   &  *    &  *
\cr
   *   & J_2 &         & *
\cr
    * &     * & \ddots & *
\cr
    * &     * &         * & J_s
\end{array}
\right].
$$
The  transformation $ 
Y_{red}=\Delta_0 X_{red}
$  of the system (\ref{PRE31marzo2015-6})  yields\footnote{
The gauge transformation $
\widetilde{Y}(z,0)=\Delta_0 X(z),
$ of the system  (\ref{9giugno-1}) at $t=0$ yields,
$$
\frac{dX}{dz}=\mathcal{B}(z)\widetilde{X},
\quad\quad
\mathcal{B}(z):=\Delta_0^{-1}B(z,0)\Delta_0,
\quad\quad
\mathcal{B}(z)\sim \Lambda+\frac{J}{z}+\sum_{k=2}^\infty \frac{\mathcal{B}_k}{z^k},
\quad\quad
\mathcal{B}_k:=\Delta_0^{-1}B_k(0)\Delta_0.
$$} 
\be 
\label{senzatred}
{dX_{red}\over dz}= {1\over z} \mathcal{B}_{red}(z) ~X_{red},
\quad\quad
\mathcal{B}_{red}(z):=\Delta_0^{-1}B_{red}(z,0)\Delta_0,
\ee
$$
\mathcal{B}_{red}(z) \sim J+\sum_{k=1}^\infty \frac{\mathcal{B}_{k+1}}{z^k},
\quad\quad
\mathcal{B}_k=\Delta_0^{-1}B_k(0)\Delta_0. 
$$
The system  (\ref{senzatred}) has block-diagonal  fundamental solutions $X_{red}=X^{(red)}_1\oplus\cdots\oplus X^{(red)}_1$, each block satisfying 
\be
\label{3maggio2016-1}
{dX_i^{(red)}\over dz}= {1\over z} \mathcal{B}_i^{(red)}(z) ~X_i^{(red)},
\quad\quad
1\leq i \leq s.
\ee
Now, $J_i$ has the unique decomposition 
  \bea 
  \label{31marzo2015-9}
 &&  J_i=D_i+S_i,
 \quad\quad
 D_i=\hbox{ diagonal matrix of integers},
   \\
   \label{31marzo2015-10}
&&   S_i=\hbox{ Jordan form with diagonal elements of real part $\in[0,1)$}.
\eea    
For $i=1,2,...,s$, let $m_i \geq 0$ be the  maximum integer difference between couples of eigenvalues of   $J_i$ ($m_i=0$ if eigenvalues do not differ by integers). Let $\overline{m}:=\hbox{max}_{i=1,..,s}m_i$.   The general theory of Fuchsian systems assures that (\ref{3maggio2016-1}) has a fundamental matrix solution 
$$
  X_i^{(red)}(z)=  K_i(z)~z^{D_i}z^{L_i} 
,
\quad\quad
  K_i(z)\sim I+\sum_{j=1}^\infty K_j^{(i)} z^{-j},
  \quad\quad
  z\to\infty \hbox{ in } \overline{S}(\alpha,\beta).
 $$
Here $L_i:=S_i+R_i$,  where  the  matrix $R_i$ is a sum $ 
R_i=R_{(1),i}+\cdots R_{(m_i),i}$, 
whose terms satisfy 
 \be 
 \label{RbLOk} 
 [R_{(l),i}]_{block~a,b}\neq 0  \quad 
 \hbox{ only if } 
 \quad
 \mu_b^{(i)} - \mu_a^{(i)}=l >0 \hbox{ integer}
 .
 \ee
  Let
 \be 
 \label{DSR}
 D:=D_1\oplus\cdots\oplus D_s,
 \quad\quad
 S:=S_1\oplus\cdots\oplus S_s,
 \quad\quad
 R:=R_1\oplus\cdots\oplus R_s,
 \quad\quad
 L:=R+S.
 \ee
 Observe now that $R$ has a sum decomposition 
 \be
 \label{18april2016-1} 
R=R_{(1)}+R_{(2)}+\cdots +R_{(\overline{m})},
\ee 
where $ 
R_{(l)}=R_{(l),1}\oplus\cdots\oplus R_{(l),s}.
$ 
Here it is understood that $R_{(l),i}=0$ if $m_i<l\leq \overline{m}$.  We conclude that 
\beas
&&
 X_{red}(z)=K(z)~z^Dz^L,
 \quad\quad
 K(z)\sim I+\sum_{j=1}^\infty K_jz^{-j},
 \quad\quad
 z\to\infty \hbox{ in } \overline{S}(\alpha,\beta),
\\
&&
K(z):=K_1(z)\oplus\cdots \oplus K_s(z),
\quad\quad
 K_j=K_1^{(j)}\oplus\cdots\oplus K_s^{(j)}.
\eeas
 Hence, there is a fundamental  solution of (\ref{9giugno-2}) at $t=0$, of the form  
$$
 \mathring{Y}(z):=G(z,0)~\Delta_0K(z)~z^Dz^L~e^{\Lambda z}
 .
 $$ 
This is rewritten as, 
$$ 
 \mathring{Y}(z)=\Delta_0 \mathcal{G}(z) ~z^Dz^L~e^{\Lambda z}, 
$$
where 
$
\mathcal{G}(z):=  \Delta_0^{-1}G(z,0)\Delta_0K(z) 
$. 
Clearly, 
 \be
  \label{31marzo2015-7} 
 \mathcal{G}(z) \sim 
I+\sum_{k=1}^\infty \mathring{F}_k z^{-k}:=\Bigl(I+\sum_{k=1}^\infty   
 \Delta_0^{-1}G_k (0) \Delta_0
 \Bigr)~\Bigl(I+\sum_{k=1}^\infty K_kz^{-k}\Bigr),
 \quad
 \quad
 z\to \infty \hbox{ in }\overline{S}(\alpha,\beta).
 \ee
 
 The results above can be summarized in the following theorem: 
 \begin{shaded}
  \bth
  \label{31marzo2015-12}
 Consider  the system (\ref{9giugno-2})  satisfying the assumptions of Theorem \ref{ShybuyaOLD}. There exist  an invertible block-diagonal matrix $\Delta_0$  and a matrix  $\mathcal{G}(z)$,  holomorphic for $|z|>N$, $z\in\bar{S}(\alpha,\beta)$,  with asymptotic expansion 
\be
\label{7dic2015-2}
\mathcal{G}(z)\sim I+\sum_{k=1}^\infty \mathring{F}_k z^{-k}, \quad z\to\infty,
\quad z\in\bar{S}(\alpha,\beta),
\ee
 such that the gauge transformation $Y(z,0)=\Delta_0\mathcal{G}(z) \mathcal{Y}(z)$  transforms (\ref{9giugno-2}) at $t=0$ into a blocked-diagonal system
\be
\label{6dic2015-1}
\frac{d \mathcal{Y}}{dz}= 
\left[\Lambda+\frac{1}{z}\left( J+\frac{R_{(1)}}{z}+\cdots + 
\frac{R_{(\overline{m})}}{z^{\overline{m}}}
 \right)
\right]
\mathcal{Y},\quad \quad
J=J_1\oplus\cdots\oplus J_s,
\ee
where $J_i$ is a Jordan form of $A_{ii}^{(1)}(0)=\widehat{A}_{ii}^{(1)}(0)$, $1\leq i \leq s$, and the $R_{(l)}$, $1\leq l \leq \overline{m}$ are defined in (\ref{RbLOk})-(\ref{18april2016-1}). The system (\ref{6dic2015-1})  has a fundamental solution $\mathcal{Y}(z)=z^Dz^L~e^{\Lambda z}$, hence  (\ref{9giugno-2})  restricted at $t=0$ has a fundamental solution,
\be
\label{7dic2015-1}
\mathring{Y}(z)=\Delta_0 \mathcal{G}(z) ~z^Dz^L~e^{\Lambda z}.
\ee
The matrices $D$, $L$ are defined in (\ref{31marzo2015-9}), (\ref{31marzo2015-10}) and (\ref{DSR}).  
 The matrix $\Delta_0$ satisfies
$$
\Delta_0^{-1} A_1(0)~ \Delta_0 =\left[
\begin{array}{cccc}
J_1 &*   &  *    &  *
\cr
   *   & J_2 &         & *
\cr
    * &     * & \ddots & *
\cr
    * &     * &         * & J_s
\end{array}
\right].
$$
\eth

\end{shaded}

\bre
Observe that (\ref{7dic2015-1}) does not  solve (\ref{9giugno-2}) for $t\neq 0$. 
  \ere

\bde 
\label{26maggio2017-1}
The  matrix 
\be
\label{16gen2016-3}
\mathring{Y}_F(z):=\Delta_0F(z)~z^Dz^L~e^{\Lambda z},
\quad
\quad
F(z):=I+\sum_{k=1}^\infty \mathring{F}_k z^{-k}
~,
\ee
is called a {\bf formal solution} of (\ref{9giugno-2})  for $t=0$ and $A_0(0)=\Lambda$. 
\ede
 Notice that we use the notation $\mathring{Y}$ for solutions of the system with $t=0$. For fixed $\Delta_0$, $D$, $L$ and  $\Lambda$  the formal solution  is in general not unique. See Corollary \ref{8gen2016-10}. 

\vskip 0.2 cm 
We note that (\ref{16gen2016-3}) can be transformed into a formal solution with the structure  described in \cite{BJL2}, but the specific form (\ref{16gen2016-3}) is more  refined  and is obtainable by  an explicit  construction from the differential system (see also Section \ref{18settembre2016-1} below).     
 

\subsection{Explicit computation of the $\mathring{F}_k$'s and $ R$ of (\ref{7dic2015-2}) and (\ref{6dic2015-1}). Uniqueness of Formal Solutions}
\label{18settembre2016-1}

We present the computation of the $\mathring{F}_k$'s in (\ref{7dic2015-2}) and $R$ in (\ref{18april2016-1}). This serves for two reasons. First, the details of the computation in itself will be used later, starting from section \ref{11marzo2016-3}.  Second,  it  yields the  Corollary  \ref{8gen2016-10}  below concerning the (non-)uniqueness of formal solutions. Consider the gauge transformation  $Y=\Delta_0 \widehat{X}$ at $t=0$, which transforms (\ref{9giugno-2}) into 
\beas
&&
\frac{d\widehat{X}}{dz}= \Bigl(\Delta_0^{-1}A(z,0)\Delta_0\Bigr) \widehat{X}(z),
\\
&&
\Delta_0^{-1}A(z,0)\Delta_0= \Lambda+\sum_{j=1}^\infty \mathcal{A}_j z^{-j},
\quad\quad
\mathcal{A}_j:=\Delta_0^{-1}A_j(0)\Delta_0.
\eeas
The recurrence equations (\ref{rec27marzo1}), (\ref{rec27marzo2}) become (using $F_l$ instead  of $G_l$),
\be
\label{rec27marzo1-QUAT}
\Lambda F_1-F_1\Lambda=-\mathcal{A}_1+B_1,
\quad\quad
\hbox{ with }~\hbox{diag}(\mathcal{A}_1)=J,
\ee
\be
\label{rec27marzo2-QUAT}
\Lambda F_l-F_l\Lambda=\Bigr[\sum_{j=1}^{l-1}\Bigl(F_jB_{l-j}-\mathcal{A}_{l-j}F_j\Bigr)-\mathcal{A}_l\Bigr]~-(l-1)F_{l-1}~+B_l.
\ee

\begin{shaded}
 \bpr 
 \label{31marzo2015-4-BiS}
(\ref{rec27marzo1-QUAT})-(\ref{rec27marzo2-QUAT}) admit a solution $\{F_k\}_{k\geq 1}$, $\{B_k\}_{k\geq 1}$ which satisfies, 
\beas
&& 
 B_1=J,
\\
&&
B_2=R_{(1)}, \quad...~,\quad B_{\overline{m}+1}=R_{(\overline{m})},
\\&&
  B_k=0~\hbox{ for any  }k\geq \overline{m}+2,
\eeas 
where $ 
R_{(l)}=R_{(l),1}\oplus\cdots\oplus R_{(l),s}
$,  and each $R_{(l),i}$ is as in (\ref{RbLOk}). The $F_k$'s so obtained are exactly the coefficients $\mathring{F}_k$ of  the asymptotic expansion of the gauge transformation  (\ref{7dic2015-2}), which yields (\ref{6dic2015-1}). 
  \epr 
\end{shaded}

\vskip 0.2 cm 
\noindent
{\it Proof:} Let $
  \mathcal{K}_l:=\Bigr[\sum_{j=1}^{l-1}\Bigl(F_jB_{l-j}-\mathcal{A}_{l-j}F_j\Bigr)-\mathcal{A}_l\Bigr]
  $, 
  and rewrite (\ref{rec27marzo1-QUAT}) and (\ref{rec27marzo2-QUAT}) in blocks $i,j$:
  
\vskip 0.2 cm 
\noindent
$\bullet$ \underline{For $l=1$} ($[i,j]$ is the block index, $1\leq i,j\leq s$): 
$$
\Lambda F_1-F_1\Lambda=-\mathcal{A}_1+B_1
\quad
\Longrightarrow
\quad
(\lambda_i-\lambda_j)F^{(1)}_{ij}=
-\mathcal{A}^{(1)}_{ij}~+B^{(1)}_{ij}.
$$
$\bullet$ \underline{For $l\geq 2$}: 
$$
\Lambda F_l-F_l\Lambda=\mathcal{K}_j-(l-1)F_{l-1}+B_l 
\quad
\Longrightarrow
\quad
(\lambda_i-\lambda_j)F^{(l)}_{ij}=\mathcal{K}^{(l)}_{ij}-(l-1)F^{(l-1)}_{ij}~+B^{(l)}_{ij}.
$$

\vskip 0.2 cm 
\noindent
$\bullet$ For $l=1 $ we find: 

\vskip 0.2 cm 
\noindent
-- If $i=j$:
$$
B^{(1)}_{ii}=\mathcal{A}^{(1)}_{ii}\equiv J_i,
\quad
\quad
F^{(1)}_{ii}\hbox{ not determined}.
  $$
-- If $i\neq j$:
$$
F^{(1)}_{ij}=-{\mathcal{A}^{(1)}_{ij}\over \lambda_i-\lambda_j} ,
\quad
\quad
B^{(1)}_{ij}=0.
$$

  \vskip 0.2 cm 
\noindent
$\bullet$ For $l\geq  2$ we find:

 \vskip 0.2 cm 
\noindent
-- If $i \neq j$: 
  $$
F^{(l)}_{ij}=(\lambda_i-\lambda_j)^{-1}\Bigl( \mathcal{K}^{(l)}_{ij}-(l-1)F^{(l-1)}_{ij}\Bigr),
\quad
\quad
B^{(l)}_{ij}=0.
$$  
In the r.h.s. matrix entries of $F_1$, ..., $F_{l-1}$ 
appear, therefore the equation determines $F_{ij}^{(l)}$. 
  
  \vskip 0.2 cm 
\noindent
-- If $i=j$: 
\be
\label{icona} 
0= 
\mathcal{K}^{(l)}_{ii}-(l-1)F^{(l-1)}_{ii}~+B^{(l)}_{ii}.
\ee
We observe that in $\mathcal{K}^{(l)}_{ii}$ the matrix entries of $F_1$, ..., $F_{l-1}$ appear, including the entry $F^{(l-1)}_{ii}$. 
Keeping into account that $B_1=\mathcal{A}_{11}^{(1)}\oplus \cdots \oplus \mathcal{A}_{ss}^{(1)}$, we explicitly write (\ref{icona}): 
$$
(l-1)F_{ii}^{(l-1)}=\sum_{k=1}^s \Bigl(
 F_{ik}^{(l-1)}B_{ki}^{(1)}-\mathcal{A}_{ik}^{(1)}F_{ki}^{(l-1)}
 \Bigr)
+
\sum_{j=1}^{l-2}\Bigl(F_jB_{l-j}-\mathcal{A}_{l-j}F_j \Bigr)_{[i,i]}~-\mathcal{A}_{ii}^{(l)}+B_{ii}^{(l)}  =
$$
$$
=  F_{ii}^{(l-1)}\mathcal{A}_{ii}^{(1)}-\mathcal{A}_{ii}^{(1)}F_{ii}^{(l-1)}-\sum_{k\neq i}\mathcal{A}_{ik}^{(1)}F_{ki}^{(l-1)}
+
\sum_{j=1}^{l-2}\Bigl(F_jB_{l-j}-\mathcal{A}_{l-j}F_j \Bigr)_{[i,i]}~-\mathcal{A}_{ii}^{(l)}+B_{ii}^{(l)} . 
$$
  Thus, keeping into account that $\mathcal{A}_{ii}^{(1)}=J_i$, the above is rewritten as follows:
\be 
\label{2015-03-06-1}
  \Bigl(J_i+l-1\Bigr) F_{ii}^{(l-1)}-F_{ii}^{(l-1)}J_i
  = -\sum_{k\neq i}\mathcal{A}_{ik}^{(1)}F_{ki}^{(l-1)}
+
\sum_{j=1}^{l-2}\Bigl(F_jB_{l-j}-\mathcal{A}_{l-j}F_j \Bigr)_{[i,i]}~-\mathcal{A}_{ii}^{(l)}+B_{ii}^{(l)}. 
\ee
In the r.h.s. every term is determined by previous steps (diagonal elements $F_{jj}^{(k)}$ appear up to $k \leq  l-2$), except for $B_{ii}^{(l)}$, which is still undetermined. (\ref{2015-03-06-1}) splits into the blocks inherited from $J_i=J_1^{(i)}\oplus\cdots\oplus J_{h_i}^{(i)}$. 
Let the eigenvalues of $J_i$ be $\mu_1^{(i)}$, ..., $\mu_{h_i}^{(i)}$, $h_i\leq p_i$.
Then (for $l\geq 2$), 
$$
\Bigl(\mu_a^{(i)}+l-1+H_{r_a}\Bigr) [F_{ii}^{(l-1)}]_{ab}-[F_{ii}^{(l-1)}]_{ab}(\mu_b^{(i)}+H_{r_b})
  =
$$
\be
\label{6dic2015-9}
\quad\quad
= \left[-\sum_{k\neq i}\mathcal{A}_{ik}^{(1)}F_{ki}^{(l-1)}
+
\sum_{j=1}^{l-2}\Bigl(F_jB_{l-j}-\mathcal{A}_{l-j}F_j \Bigr)_{[i,i]}~-\mathcal{A}_{ii}^{(l)}+B_{ii}^{(l)}\right]_{ab}.
\ee
Here $[\cdots]_{ab}$ denotes a block, with $1\leq a,b\leq h_i$. 

$\bullet$  
If \underline{$\mu_b^{(i)}-\mu_a^{(i)}=l-1$},  the l.h.s. of (\ref{6dic2015-9}) is 
$
H_{r_a} [F_{ii}^{(l-1)}]_{ab}-[F_{ii}^{(l-1)}]_{ab}H_{r_b}$. The homogeneous equation 
$
H_{r_a} [F_{ii}^{(l-1)}]_{ab}-[F_{ii}^{(l-1)}]_{ab}H_{r_b}=0$ has non trivial solutions, depending on parameters, since the matrices $H_{r_a}$ and $H_{r_b}$ have common eigenvalue. One can then choose $F_{ii}$ to be a solution of the homogeneous equation, and determine $[B_{ii}^{(l)}]_{ab}\neq 0$ by imposing that the r.h.s. of (\ref{6dic2015-9})  is equal to  $0$. 

$\bullet$  If \underline{$\mu_b^{(i)}-\mu_a^{(i)}\neq l-1$}, the choice $[B_{ii}^{(l)}]_{ab}= 0$ is possible and $[F_{ii}^{(l-1)}]_{ab}$ is determined. 

\vskip 0.2 cm 
\noindent
We conclude that 
$$ 
[B_{ii}^{(l+1)}]_{ab}\neq 0 \quad
\hbox{ only if }
\quad
 \mu_b^{(i)}-\mu_a^{(i)}=l>0 \quad 
 \hbox{integer}.
$$
This means that $[B_{ii}^{(l+1)}]_{ab}=[R_{(l),i}]_{ab}$. 
   $\Box$
   
\begin{shaded}
\bcr [Uniqueness of Formal Solution at $t=0$] 
\label{8gen2016-10}
A formal solution (\ref{16gen2016-3}) with given $\Delta_0$, $D$, $L$, $\Lambda$ is unique if and only if  for any $1\leq i \leq s$ the eigenvalues of $\widehat{A}_{ii}^{(1)}(0)$ do not differ by a non-zero integer. 
\ecr
\end{shaded}
\noindent
{\it Proof:} Computations above show  that  $\{F_k\}_{k=1}^\infty$ is not  uniquely determined if and only if some  $\mu_b^{(i)}-\mu_a^{(i)}=l-1$, for some $ l\geq 2$, some $i\in\{1,2,...,s\}$, and some $a,b$.  $\Box$

\subsection{Special sub-case with  $R=0$, $J$ diagonal, $\Delta_0=I$}
\label{11marzo2016-3}

A sub-case is very important for the discussion to come, occurring when   $\Delta_0=I$ and $A_{ii}^{(1)}(0)$ is diagonal. Clearly,   if $\Delta_0=I$, then  $J_i=A_{ii}^{(1)}(0)$.  Hence,  if $\Delta_0=I$, then  $J$ is diagonal if and only if $\Bigl(\widehat{A}_{ii}^{(1)}(0)\Bigr)_{pq}=0$ for any $1\leq p\neq q\leq p_i$. 

\begin{shaded}
\bpr
\label{6ott2015-5} 

There exists a fundamental solution (\ref{7dic2015-1}) at $t=0$  in a simpler form 
\be
\label{6dic2015-8}
\mathring{Y}(z)=\mathcal{G}(z) z^{B_1(0)} e^{\Lambda z},
\ee
with $\Delta_0=I$,  $J=B_1(0)=\hbox{\rm diag}(A_1(0))$ diagonal, and 
\be
\label{6dic2015-7}
 \mathcal{G}(z) 
\sim I+\sum_{k=1}^\infty \mathring{F}_k z^{-k}, \quad \quad 
z\to \infty ~\hbox{ in } \overline{S}(\alpha,\beta),
\ee 
if and only if the following conditions hold: 

\vskip 0.2 cm 
\noindent
$\bullet$ For every $i\in\{1,2,...,s\}$, and every  $p,q$, with $1\leq p\neq q\leq p_i$,  then 
\be
\label{14feb2016-1}
\Bigl(\widehat{A}_{ii}^{(1)}(0)\Bigr)_{pq}=0.
\ee

\vskip 0.2 cm 
\noindent
$\bullet$ {If   $\Bigl(\widehat{A}^{(1)}_{ii}(0)\Bigr)_{pp}
    -
   \Bigl(\widehat{A}^{(1)}_{ii}(0)\Bigr)_{qq}+l-1= 0$, for some $ l\geq 2$, some $i\in\{1,2,...,s\}$}, and some diagonal entries $\Bigl(\widehat{A}^{(1)}_{ii}(0)\Bigr)_{pp}$, $\Bigl(\widehat{A}^{(1)}_{ii}(0)\Bigr)_{qq}$,  then 
   \be
\sum_{k\neq i}^s \Bigl(\widehat{A}^{(1)}_{ik}(0)~\mathring{F}^{(l-1)}_{ki}\Bigr)_{pq}
   +
   \sum_{j=1}^{l-2}\sum_{k=1}^s
   \Bigl(\widehat{A}^{(l-j)}_{ik}(0)~\mathring{F}^{(j)}_{ki}\Bigr)_{pq}
   +
   \Bigl(\widehat{A}^{(l)}_{ii}(0)\Bigr)_{pq}= 
0,
\label{Good123}
\ee
for those values of $l$, $i$, $p$ and $q$.

\epr
\end{shaded}

\vskip 0.2 cm 
\noindent
{\it Proof:}  We only need to clarify  (\ref{Good123}), while (\ref{14feb2016-1}) has already been motivated.  We solve  (\ref{rec27marzo1-QUAT}),  (\ref{rec27marzo2-QUAT}) when $\Delta_0=I$, namely (recall that $\widehat{A}_j(0)\equiv A_j(0)$) (we write $F_l$, as in (\ref{rec27marzo1-QUAT}),  (\ref{rec27marzo2-QUAT}), but it is clear that the result of the computation will be the $\mathring{F}_l$ appearing in  (\ref{6dic2015-7})):
\beas
&&
\Lambda F_1-F_1\Lambda=-\widehat{A}_1(0)+B_1,
\\
&&
\Lambda F_l-F_l\Lambda=\Bigr[\sum_{j=1}^{l-1}\Bigl(F_jB_{l-j}-\widehat{A}_{l-j}(0)F_j\Bigr)-\widehat{A}_l(0)\Bigr]~-(l-1)F_{l-1}~+B_l.
\eeas
 At level $l=1$:
$$ 
B_1=\hbox{diag}\widehat{A}_1(0),\quad\quad
F_{ij}^{(1)}=-\frac{\widehat{A}_{ij}(0)}{\lambda_i-\lambda_j}.
$$ 
At level $l\geq 2$,   
$$ 
F_{ij}^{(l)}=\frac{ \mathcal{K}_{ij}^{(l)}-(l-1)F_{ij}^{(l-1)}}{\lambda_i-\lambda_j},
\quad\quad
B_{ij}^{(l)}=0,
$$ 
where $
  \mathcal{K}_l=\Bigr[\sum_{j=1}^{l-1}\Bigl(F_jB_{l-j}-\widehat{A}_{l-j}(0)F_j\Bigr)-\widehat{A}_l(0)\Bigr]$. 
Formula (\ref{6dic2015-9}) reads 
$$
\Bigl(\mu_a^{(i)}-\mu_b^{(i)}+l-1\Bigr) [F_{ii}^{(l-1)}]_{ab}
=
 \left[-\sum_{k\neq i}\widehat{A}_{ik}^{(1)}(0)F_{ki}^{(l-1)}
+
\sum_{j=1}^{l-2}\Bigl(F_jB_{l-j}-\widehat{A}_{l-j}(0)F_j \Bigr)_{[i,i]}~-\widehat{A}_{ii}^{(l)}(0)+B_{ii}^{(l)}\right]_{ab}.
$$
Indices above are block indices. The above can  be re-written in terms of the matrix entries,
\beas
&&
\Bigl((\widehat{A}_{ii}^{(1)}(0))_{pp}-(\widehat{A}_{ii}^{(1)}(0))_{qq}+l-1\Bigr) (F_{ii}^{(l-1)})_{pq}
=
\\
&&
\quad\quad
=
 \left[-\sum_{k\neq i}\widehat{A}_{ik}^{(1)}(0)F_{ki}^{(l-1)}
+
\sum_{j=1}^{l-2}\Bigl(F_jB_{l-j}-\widehat{A}_{l-j}(0)F_j \Bigr)_{[i,i]}~-\widehat{A}_{ii}^{(l)}(0)+B_{ii}^{(l)}\right]_{entry ~pq}.
\eeas

$\bullet$ If $(\widehat{A}_{ii}^{(1)}(0))_{pp}-(\widehat{A}_{ii}^{(1)}(0))_{qq}+l-1\neq 0$, choose $B_{ii}^{(l)}=0$ and determine $(F_{ii}^{(l-1)})_{pq}$. 

$\bullet$ If $(\widehat{A}_{ii}^{(1)}(0))_{pp}-(\widehat{A}_{ii}^{(1)}(0))_{qq}+l-1 =0$,  by induction assume that the $B_{l-j}=0$. Then   the equation is satisfied for any $(F_{ii}^{(l-1)})_{pq}$ and for  
$$
(B_{ii}^{(l)})_{pq}=
\left[\sum_{k\neq i}\widehat{A}_{ik}^{(1)}(0)F_{ki}^{(l-1)}
+
\sum_{j=1}^{l-2}\Bigl(\widehat{A}_{l-j}(0)F_j \Bigr)_{block ~[i,i]}~+\widehat{A}_{ii}^{(l)}(0)\right]_{entry~pq} .
$$
Then, if we  impose that $(B_{ii}^{(l)})_{pq}=0$ we obtain   the necessary and sufficient condition (\ref{Good123}). 
The proof by induction is justified because at the first step, namely $l=2$, we need to solve
\be 
\label{equAT0}
\left((\widehat{A}^{(1)}_{ii}(0))_{pp}
    -
   (\widehat{A}^{(1)}_{ii}(0))_{qq}+1\right){(F^{(1)}_{ii})_{pq}}
  = 
   -
   \sum_{k\neq i}^n (\widehat{A}^{(1)}_{ik}(0)~F^{(1)}_{ki})_{pq}
   -
  (\widehat{A}^{(2)}_{ii}(0))_{pq}+(B^{(2)}_{ii})_{pq}.
\ee
If $(\widehat{A}^{(1)}_{ii}(0))_{pp}
    -
   (\widehat{A}^{(1)}_{ii}(0))_{qq}+1\neq 0$, the above  has a unique solution for any choice of $(B^{(2)}_{ii})_{pq}$. We choose $({B}^{(2)}_{ii})_{pq}=0$.  
   If $\widehat{A}^{(1)}_{ii}(0))_{pp}
    -
   (\widehat{A}^{(1)}_{ii}(0))_{qq}+1= 0$, the equation leaves the choice  of  ${(F^{(1)}_{ii})_{pq}}$ free, and determines 
$$
  ({B}^{(2)}_{ii})_{pq}  =
   \sum_{k\neq i}^n (\widehat{A}^{(1)}_{ik}(0)~F^{(1)}_{ki})_{pq}
  +
  (\widehat{A}^{(2)}_{ii}(0))_{pq}
=
  - \sum_{k\neq i}^n {(\widehat{A}^{(1)}_{ik}(0)\widehat{A}^{(1)}_{ki}(0))_{pq}\over \lambda_k -\lambda_i}
  +
  (\widehat{A}^{(2)}_{ii}(0))_{pq}.
$$
   We can choose $(B^{(2)}_{ii})_{pq}=0$ if and only if
\be 
\label{thatsGood}
  (\widehat{A}^{(2)}_{ii}(0))_{pq}=\sum_{k\neq i}^n {(\widehat{A}^{(1)}_{ik}(0)\widehat{A}^{(1)}_{ki}(0))_{pq}\over \lambda_k -\lambda_i},
  \ee
which is precisely (\ref{Good123}) for $l=2$. 
$\Box$



\section{Solutions for $t\in\mathcal{U}_{\epsilon_0}(0)$ with $A_0(t)$ Holomorphically Diagonalisable.}
\label{11marzo2016-1}

In the previous section, we have constructed fundamental solutions at the  coalescence point $t=0$. Now, we let $t$ vary in $\mathcal{U}_{\epsilon_0}(0)$. 
In Sibuya Theorem,  $\widehat{A}_0(t)=\widehat{A}_{11}^{(0)}(t)\oplus\cdots\oplus \widehat{A}_{ss}^{(0)}(t)$ is neither diagonal nor  in Jordan form, except for $t=0$. $A_0(t)$ admits a Jordan form at each point of $\mathcal{U}_{\epsilon_0}(t)$, but in general this similarity is not realizable by a holomorphic transformation.   In order to procede, we need the following  fundamental assumption, already stated in the Introduction.
 \vskip 0.2 cm 

{\bf Assumption 1:}  For $|t|\leq \epsilon_0$ sufficiently small and such that  Lemma \ref{SibLEM1} and  Theorem \ref{ShybuyaOLD} apply, we assume that  $A_0(t)$ is holomorphically similar to a diagonal form $\Lambda(t)$, namely there exists a holomorphic invertible  $G_0(t)$  for $|t|\leq \epsilon_0$ such that 
$$ 
G_0(t)^{-1}A_0(t)~G_0(t)=\Lambda(t)\equiv \hbox{diag} \bigl(u_1(t),u_2(t),...,u_n(t)\bigr),
$$
with $ A_0(0)=\Lambda$, $G_0(0)=I$.  

\bre

Assumption 1 is equivalent to  the assumption that $A_0(t)$ is holomorphically similar to its Jordan form. The requirement implies  by continuity that the Jordan form is diagonal, being  equal to $\Lambda=\Lambda(0)$ at $t=0$.

\ere

With Assumption 1, we can represent the eigenvalues as well defined holomorphic functions $u_1(t)$, $u_2(t)$, ..., $u_n(t)$ such that 
 \bea
\label{29gen2016-1}
u_1(0)=\cdots=u_{p_1}(0)&&=\lambda_1,
\\
\label{29gen2016-2}
u_{p_1+1}(0)=\cdots=u_{p_1+p_2}(0)&&=\lambda_2,
\\
 \vdots&&
\\
\label{29gen2016-3}
u_{p_1+\cdots+p_{s-1}+1}(0)=\cdots=u_{p_1+\cdots+p_{s-1}+p_s}(0)&&= \lambda_s.
\eea
 Moreover, 
$$
\Lambda(t)=
\Lambda_1(t)\oplus\Lambda_2(t)\oplus~\cdots~\oplus\Lambda_s(t),
$$
 where $\Lambda_1(t)$, ..., $\Lambda_s(t)$ are diagonal matrices of dimensions respectively $p_1$, ..., $p_s$, such that $\Lambda_j(t)\to \lambda_j  I_{p_j}$ for $t\to 0$, $j=1,...,s$. For example, $\Lambda_1(t)=$ diag$(u_1(t),...,u_{p_1}(t))$, and so on.  
 Any two matrices   
   $\Lambda_i(t)$ and $\Lambda_j(t)$   have no common eigenvalues  for   $i\neq j$ and small $\epsilon_0$.   
 \vskip 0.2 cm 
The { coalescence locus in  $\mathcal{U}_{\epsilon_0}(0)$ }  is explicitly written as follows
$$
\Delta:=\bigcup_{
\begin{array}{c}
a\neq b
\\
a,b=1,...,m
\end{array}}\{
t \in \mathbb{C}^m \hbox{ such that: }
 |t| \leq \epsilon_0
\hbox{ and } 
u_a(t)=u_b(t)
\Bigr\}
.
$$
We can also write 
$$
\Delta= 
\bigcup_{i=1}^s ~\Delta_i,
$$
where 
$
\Delta_i$ is the coalescence locus of $\Lambda_i(t)$. For $m=1$, $\Delta$ is a finite set of isolated points.

\begin{shaded}
\noindent
 {\bf Improvement of Theorem \ref{ShybuyaOLD}:}
{\it
With the same assumptions and notations as of Theorem \ref{ShybuyaOLD}, if Assumption $1$ holds, then
$$B(z,t)\sim \Lambda(t)+\sum_{k\geq 1} B_k(t) z^{-k},
\quad\quad
z\to \infty \hbox{ in }\overline{S}(\alpha,\beta).
$$ }
\end{shaded}

\vskip 0.2 cm 
With Assumption $1$, we can replace the gauge trasfromation (\ref{redDU}) with 
$$ 
\widetilde{Y}(z,t)=e^{\Lambda(t) z} ~Y_{red}(z,t).
$$
Since $\widehat{A}_0(t)=\Lambda(t)$,  then  $ 
B_{red}(z,t)\sim \sum_{k=1}^\infty B_k(t) z^{-k+1}
$. Hence  the reduced  system (\ref{31marzo2015-2}) is Fuchsian also for $t\neq 0$. 
The recursive relations (\ref{rec27marzo1}) and (\ref{rec27marzo2})  become   $B_0(t)=\Lambda(t)$ for  $l=0$, and: 

\vskip 0.2 cm 
\noindent
\underline{For $l=1$}: 
\be 
\label{rec27marzo1BIS}
\Lambda(t) G_1-G_1\Lambda(t)=-\widehat{A}_1(t)+B_1.
\ee
\underline{For $l\geq 2 $}: 
\be 
\label{rec27marzo2BIS}
\Lambda(t) G_l-G_l\Lambda(t) =\Bigr[\sum_{j=1}^{l-1}\Bigl(G_jB_{l-j}-\widehat{A}_{l-j}(t)G_j\Bigr)-\widehat{A}_l(t)\Bigr]~-(l-r)G_{l-r}~+B_l.
\ee
As for Theorem \ref{ShybuyaOLD}, the choice which yields holomorphic  $G_l(t)$'s and $B_l(t)$'s is (\ref{cho1}) and (\ref{cho2}).  
Generally speaking, it is not possible to choose the $B_l(t)$'s diagonal for $l\geq 2$, because such a choice would give $G_k(t)$'s   {\it diverging} at the locus $\Delta$.

\subsection{Fundamental Solution in a neighbourhood of  $t_0\not\in\Delta$, with Assumption 1}
\label{11marzo2016-2}
 Let Assumption $1$ hold. 
 Theorem \ref{ShybuyaOLD} has been formulated  in a neighbourhood of  $t=0$, with block partition of   $A_0(0)=\Lambda_1\oplus\cdots \oplus \Lambda_s$.
 Theorem \ref{ShybuyaOLD} can also be  formulated in a neighbourhood (polydisc) of a point $t_0\in \mathcal{U}_{\epsilon_0}(0)\backslash \Delta$, of the form
 \beas 
&&\mathcal{U}_{\rho_0}(t_0):=\{t\in\mathbb{C}~|~|t-t_0|\leq \rho_0\}\subset  \mathcal{U}_{\epsilon_0}(0),
\\
\\
&&\mathcal{U}_{\rho_0}(t_0)\cap \Delta=\emptyset,
\eeas 
 where  $\Lambda(t)$ has distinct eigenvalues, provided that  $\rho_0>0$ is small enough. 
In order to do this,  we need to introduce sectors. To this end, consider a fixed point $t_*$ in $\mathcal{U}_{\epsilon_0}(0)$, and the    eigenvalues $u_1(t_*)$, ..., $u_n(t_*)$ of $\Lambda(t_*)$. We introduce an admissible direction $\eta^{(t_*)}$ such that 
\be
\label{19marzo2016-1}
\eta^{(t_*)}\neq \arg_p\Bigl(u_a(t_*)-u_b(t_*)\Bigr) ~\hbox{ mod} (2\pi) ,
\quad
\forall ~   1\leq a\neq b\leq n.
\ee
There are  $2\mu_{t_*}$ determinations satisfying $\eta^{(t_*)}-2\pi<\widehat{\arg}(u_a(t_*)-u_b(t_*))<\eta^{(t_*)}$. They  will be numbered as 
$$
\eta^{(t_*)}>\eta^{(t_*)}_0>~\cdots~>\eta_{2\mu^{(t_*)}-1}>\eta^{(t_*)}-2\pi.
$$
Correspondingly, we introduce the directions
$$
\tau^{(t_*)}:=\frac{3\pi}{2}-\eta^{(t_*)},\quad\quad
\tau_\nu^{(t_*)}=\frac{3\pi}{2}-\eta_\nu^{(t_*)},
\quad\quad
0\leq \nu \leq 2\mu_{t_*}-1,
$$
satisfying 
$$
\tau^{(t_*)}<\tau^{(t_*)}_0<\tau^{(t_*)}_1<\cdots<\tau^{(t_*)}_{2\mu_{t_*}-1}<\tau^{(t_*)}+2\pi.
$$
The following  relation defines $\tau_\sigma^{(t_*)}$ for any $\sigma\in\mathbb{Z}$, represented as $\sigma=\nu+k\mu_{t_*}$:
$$
\tau_{\nu+k\mu_{t_*}}:=\tau_{\nu}^{(t_*)}+k\pi,
\quad\quad
\nu\in \{0,1,...,\mu_{t_*}-1\},
\quad\quad
k\in\mathbb{Z}.
$$
 Finally, we introduce the sectors 
$$
\mathcal{S}_\sigma^{(t_*)}:=S(\tau_\sigma^{(t_*)}-\pi,\tau_{\sigma+1}^{(t_*)}),\quad\quad
\sigma\in\mathbb{Z}. 
$$
Theorem \ref{ShybuyaOLD}  in a neighbourhood of $t_0$ becomes: 

\begin{shaded}
\bth
\label{ShybuyaOLDbis}
Let Assumption 1 hold and let $t_0\in \mathcal{U}_{\epsilon_0}(0)\backslash \Delta$.  Pick up a sector  $
\mathcal{S}_\sigma^{(t_0)}=S(\tau_\sigma^{(t_0)}-\pi,\tau_{\sigma+1}^{(t_0)})$, $\sigma\in\mathbb{Z}$, 
 as above. For any closed sub-sector 
$$
\overline{S}^{(t_0)}(\alpha,\beta):=\Bigl\{z\in\mathcal{R}~|~\tau_\sigma^{(t_0)}-\pi<\alpha\leq \arg z \leq \beta <\tau_{\sigma+1}^{(t_0)}\Bigr\}\subset \mathcal{S}_\sigma^{(t_0)},
$$
there exist a sufficiently large positive number $N$, a sufficiently small positive number $\rho$ and  an invertible matrix valued function ${G}(z,t)$ with the following properties:

 i) $G(z,t)$ is holomorphic in $(z,t)$ for $|z|\geq N$, $z\in \overline{S}^{(t_0)}(\alpha,\beta)$, $|t-t_0|\leq \rho$;

 ii) $G(z,t)$ has  uniform asymptotic expansion for $|t-t_0|\leq\rho$, with holomorphic coefficients $G_k(t)$:
 $$ 
 G(z,t)\sim I+\sum_{k=1}^\infty  G_k(t)z^{-k},
 \quad
 z\to \infty, 
 \quad
 z\in \overline{S}^{(t_0)}(\alpha,\beta),
 $$

 iii) The gauge transformation  
 $$
  Y(z,t)=G_0(t)G(z,t)\widetilde{Y}(z,t),  
  $$ 
  
reduces the initial system (\ref{9giugno-2}) to 
$$
  {d\widetilde{Y}\over dz}=
B(z,t)\widetilde{Y},
$$

where $B(z,t)$ is a {\rm diagonal} holomorphic matrix function of $(z,t)$ in the domain

 $
|z|\geq N$, $  z\in \overline{S}(\alpha,\beta)$, $|t-t_0|\leq\rho$, with uniform asymptotic expansion  and holomorphic coefficients: 
  $$
B(z,t)\sim \Lambda(t)+\sum_{k=1}^\infty {B}_k(t) z^{-k}, 
\quad
z\to \infty, 
\quad
z\in \overline{S}^{(t_0)}(\alpha,\beta).
  $$
  In particular, $ 
B_1(t)=\hbox{\rm diag}~\widehat{A}_1(t)$. 
  \eth
  \end{shaded}

\bre 
$\overline{S}^{(t_0)}(\alpha,\beta)$ is not the same   $\overline{S}(\alpha,\beta)$ of Theorem \ref{ShybuyaOLD} (the latter should be denoted  $\overline{S}^{(0)}(\alpha,\beta)$ for consistency of notations).   
The matrices $G(z,t)$ and $B(z,t)$ {\it are not the same} of Theorem \ref{ShybuyaOLD}. On the other hand, $G_0(t)$ is the same,  by  Assumption $1$.  
\ere

As before, we let $B_{red}(z,t)=z(B(z,t)-\Lambda(t))$.  Then the system (\ref{9giugno-2}) has a fundamental matrix solution 
$$
Y(z,t)=G_0(t) \mathcal{G}(z,t) z^{B_1(t)}e^{\Lambda(t)z},
$$
where $\mathcal{G}(z,t)=G(z,t)K(z,t)$, and 
$$ 
K(z,t)=\exp\left\{\int_\infty^z \frac{B_{red}(\zeta,t)-B_1(t)}{\zeta}~d\zeta\right\}
\sim \exp\left\{  \sum_{k=2}^\infty B_k(t)\frac{z^{-k+1}}{-k+1}\right\}=I+\sum_{j=1}^\infty K_j(t)z^j,
$$
$z\to\infty$  in $\overline{S}(\alpha,\beta)$. 
This result is well known, see \cite{Sh3}.  This proves the first part of the following 

\begin{shaded}
\bcr
\label{9dic2015-1}
  The analogue of Theorem \ref{ShybuyaOLDbis} holds with a new gauge transfromation  
 $\mathcal{G}(z,t)$, enjoying the same asymptotic and analytic properties, such that $Y(z,t)=G_0(t)\mathcal{G}(z,t)\widetilde{Y}(z)$  transforms the  system (\ref{9giugno-2}) into 
\be
\label{verySIMPLE}
\frac{d\widetilde{Y}}{dz}=\left(\Lambda(t)+\frac{B_1(t)}{z}\right)\widetilde{Y},
\quad\quad
B_1(t)=\hbox{\rm diag}\widehat{A}_1(t).
\ee
With the above choice, the system (\ref{9giugno-2}) has a fundamental solution,
\be 
\label{FUNDAMSOL}
Y(z,t)=G_0(t) \mathcal{G}(z,t)z^{B_1(t)}e^{\Lambda(t)z}.
\ee
and $
\mathcal{G}(z,t)$ is holomorphic for $z\in \overline{S}^{(t_0)}(\alpha,\beta)$, $|z|\geq N$ and  $|t-t_0|\leq \rho$, with expansion 
\be
\label{8dic2015-1}
\mathcal{G}(z,t)\sim I +\sum_{k=1}^\infty F_k(t) z^{-k},
\ee
for $z\to \infty$ in $\overline{S}^{(t_0)}(\alpha,\beta)$,  uniformly in $|t-t_0|\leq \rho$. The coefficients $F_k(t)$ are uniquely determined and  holomorphic on $\mathcal{U}_{\epsilon_0}(0)\backslash \Delta$.  
\ecr 
\end{shaded}

\vskip 0.2 cm 
\noindent
{\it Proof:} The statement is clear from the previous construction. It is only to be justified that the 
 $F_k(t)$'s, $k\geq 1$, are holomorphic functions of $t\not\in \Delta$  and uniquely determined. We solve   (\ref{rec27marzo1BIS}) and (\ref{rec27marzo2BIS}) for the $F_k(t)$'s, namely
\beas
&&
\Lambda(t) F_1-F_1\Lambda(t)=-\widehat{A}_1(t)+B_1,
\\
&&
\Lambda(t) F_l-F_l\Lambda(t)=\Bigr[\sum_{j=1}^{l-1}\Bigl(F_jB_{l-j}-\widehat{A}_{l-j}(t)F_j\Bigr)-\widehat{A}_l(t)\Bigr]~-(l-1)F_{l-1}~+B_l.
\eeas
 It is convenient to use the notation $u_1(t)$, ..., $u_n(t)$ for the distinct eigenvalues. Matrix entries are here denoted $a$, $b$ $\in\{1,2,...,n\}$. 
For $l=1$, 
\beas
&&
(F_1)_{a b}(t)=-\frac{(\widehat{A}_1)_{a b}(t)}{u_a(t)-u_b(t)}, 
\quad\quad
(B_1(t))_{a b}=0,\quad
a\neq b.
\\
&&
(B_1)_{a a}(t)=(\widehat{A}_1)_{a a}(t),
\quad
\Longrightarrow
\quad
 B_1(t)=\hbox{diag}(\widehat{A}_1(t)).
\eeas
Now, {\it impose} that $ 
B_l(t)=0$ for any $l\geq 2$. 
Hence, at level $l=2$ we get:
$$ 
(F_1)_{aa}(t)=-\sum_{b\neq a}(\widehat{A}_1)_{a b}(t)(F_1)_{ba}(t)-(\widehat{A}_2)_{a a}(t).
$$
For any $l\geq 2$, we find:
\beas
&(F_l)_{a b}(t)&=-\frac{1}{u_a(t)-u_b(t)}\left\{
\Bigl[
(\widehat{A}_1)_{a a}(t)-(\widehat{A}_1)_{bb}(t)+l-1
\Bigr](F_{l-1})_{a b}(t)+
\right.
\\
&&
\left.
+\sum_{\gamma\neq a}(\widehat{A}_1)_{a \gamma}(t)(F_{l-1})_{\gamma b}(t)+\sum_{j=1}^{l-2}\Bigl(\widehat{A}_{l-j}(t)F_j(t)\Bigr)_{a b}+(\widehat{A}_l)_{a b}(t)
\right\},
\quad
a\neq b.
\\
\\
&(l-1)(F_{l-1})_{a a}(t)&=-\sum_{b\neq a }(\widehat{A}_1)_{a b}(t)(F_{l-1})_{ba}(t) -\sum_{j=1}^{l-2}\Bigl(\widehat{A}_{l-j}(t)F_j(t)\Bigr)_{a a}-(\widehat{A}_l)_{a a}(t).
\eeas
The above formulae show that the $F_l(t)$ are uniquely determined, and holomorphic away from $\Delta$. $\Box$

\vskip 0.3 cm 
The above result has two corollaries:

\begin{shaded}
\bpr
\label{8dic2015-3}
The coefficients $F_k(t)$ in the expansion (\ref{8dic2015-1}) are holomorphic at a point $t_\Delta\in\Delta$ if and only if there exists a neighbourhood of $t_\Delta$ where 
\be
\label{14feb2016-3}
(\widehat{A}_1)_{a b}(t)
\ee
and
{\small
\be
\label{8dic2015-2}
\Bigl[
(\widehat{A}_1)_{a a}(t)-(\widehat{A}_1)_{bb}(t)+l-1
\Bigr](F_{l-1})_{a b}(t)
+\sum_{\gamma\neq a}(\widehat{A}_1)_{a \gamma}(t)(F_{l-1})_{\gamma b}(t)+\sum_{j=1}^{l-2}\Bigl(\widehat{A}_{l-j}(t)F_j(t)\Bigr)_{a b}+(\widehat{A}_l)_{a b}(t)
\ee}

\noindent
vanish  as fast as  $\mathcal{O}(u_a (t)-u_b(t))$ in the neighbourhood, for those indexes $a,b\in\{1,2,...,n\} $ such that $u_a(t)$ and $u_b(t)$ coalesce   when $t$ approaches a point of $\Delta$ in the neighbourhood. In particular,  the $F_k(t)$'s are holomorphic in the whole $\mathcal{U}_{\epsilon_0}(0)$  if and only if (\ref{14feb2016-3}) and (\ref{8dic2015-2}) are zero along $\Delta$. 
\epr
\end{shaded}

Remarkably, in the isomonodromic case, we will prove that if  we just require vanishing of $(A_1)_{a b}(t)$ then  all the complicated expressions (\ref{8dic2015-2}) also vanish consequently. 

\begin{shaded}
\bpr
\label{14feb2016-2}
If the holomorphic conditions of  Proposition \ref{8dic2015-3} hold at $t=0$, then  (\ref{14feb2016-1})   and (\ref{Good123}) are satisfied, with the choice 
$$ 
\mathring{F}_k=F_k(0), 
\quad
k\geq 1.
$$
 If moreover $\Bigl(\widehat{A}_1(0)\Bigr)_{aa}-\Bigl(\widehat{A}_1(0\Bigr)_{bb}+l-1\neq 0$ for every $l\geq 2$,  then the above is the unique choice of the $\mathring{F}_k$'s, according to Corollary \ref{8gen2016-10}. 
\epr
\end{shaded}

Expression (\ref{8dic2015-2}) is a rational function of the matrix entries of 
$\widehat{A}_1(t)$, ..., $\widehat{A}_l(t)$, since $F_1(t)$,...,$F_{l-1}(t)$ are expressed in terms of $\widehat{A}_1(t)$, ..., $\widehat{A}_{l}(t)$.  For example, for $l=2$, (\ref{8dic2015-2}) becomes 
\be 
\label{zzzz}
 \left((\widehat{A}_1)_{b b}(t)
    -
  (\widehat{A}_1)_{a a }(t)-1\right)\frac{(\widehat{A}_1)_{a b}(t)}{u_a (t)-u_b(t)}
  +
  (\widehat{A}_2)_{a b}(t)
-
   \sum_{\gamma\neq a}\frac{(\widehat{A}_1)_{a \gamma}(t)(\widehat{A}_1)_{\gamma b}(t)}{u_\gamma(t)-u_b(t)}.
  \ee

\bex
\label{5dic2016-1}
 {\rm The following system does not satisfy the vanishing conditions of Proposition \ref{8dic2015-3}
 \be
\label{1novembre2016-2}      
\widehat{A}(z,t)=\left(
  \begin{array}{cc}
  0 & 0 \\ 0 & t
  \end{array}
  \right)
  +\frac{1}{z}\left(
  \begin{array}{cc}
  1 & 0 \\ t & 2
  \end{array}
  \right),
  \quad
  \Delta=\{t\in\mathbb{C}~|~t=0\}\equiv\{0\}
  \ee  
 It has a fundamental solution  
  $$Y(z,t)=
  \left[
  \begin{array}{cc}
  1 & 0 \\ w(z,t) & 1
  \end{array}
  \right]\left(
  \begin{array}{cc}
  z & 0 \\ 0 & z^2e^{tz}
  \end{array}
  \right),
  $$ with 
   $$
   w(z,t):=t^2z e^{tz}\hbox{Ei}(tz)-t\sim \sum_{k=1}^\infty \frac{(-1)^k k!}{t^{k-1}} z^{-k}, 
   \quad
   z\to \infty, \quad
     -3\pi/2<\arg (tz) <
  3\pi/2.
  $$
   The above solution  has   asymptotic representation  (\ref{8dic2015-1}), namely (\ref{4gen2016-6}). Now, $t=0$ is a branch point of logarithmic type, since ${\rm Ei}(zt)=-\ln(zt) +$ holomorphic function of $zt$. Moreover, the coefficients $F_k(t)$ diverge  when $t\to 0$. 
   The reader can check that the system has also fundamental solutions which are holomorphic at $t=0$, but without the standard asymptotic representation $Y_F(z,t)$. 
We also notice a peculiarity of this particular example,  namely that $Y(z,t)$ and $Y(ze^{-2\pi i},t)$ are connected by a  Stokes matrix $\mathbb{S}=\left[\begin{array}{cc}
  1 & 0 \\ 2\pi i t^2 & 1
  \end{array}
  \right]$, which is holomorphic also at $t=0$ and coincides with the trivial Stokes matrix  $I$ of the system $\widehat{A}(z,t=0)$.  } $\Box$ 
   
   \eex

\subsection{Fundamental Solution in a neighbourhood of  $t_\Delta\in\Delta$, with Assumption 1}

Let Assumption $1$ hold. Let   $t_\Delta\in\Delta$. Since the case $t_\Delta=0$ has  already been discussed in detail, suppose that $t_\Delta\neq 0$. Then $t_\Delta\in \Delta_i$, for some  $i\in\{1,2,...,s\}$. 

Directions $\tau_\sigma^{(t_\Delta)}$, $\sigma\in \mathbb{Z}$,	and sectors $\mathcal{S}_\sigma^{(t_\Delta)}$ have been defined in section \ref{11marzo2016-2} (just put $t_*=t_\Delta$).  We leave to the reader the task to adjust the statement of  Theorem \ref{ShybuyaOLD} reformulated in a neighbourhood of $t_\Delta$, with the block partition of $\Lambda(t_\Delta)$, which is finer than that of $\Lambda(0)$. The closed sector in the theorem will be denoted  $\overline{S}^{(t_\Delta)}(\alpha ,\beta)\subset \mathcal{S}^{(t_\Delta)}_\sigma$. A solution analogous to (\ref{7dic2015-1}) is constructed at $t=t_\Delta$, with finer block partition than (\ref{7dic2015-1}).  Special cases as in Section \ref{11marzo2016-3} are very important for us, hence we state the following.  

\begin{shaded}
\noindent
{\bf Proposition \ref{6ott2015-5} generalized at $t_\Delta$:}   For $t=t_\Delta$, the fundamental solution analogous to  (\ref{7dic2015-1})  reduces to an analogous to (\ref{6dic2015-8}), namely
\beas
&&
Y_{(t_\Delta)}(z)=G_0(t_\Delta)\mathcal{G}_{(t_\Delta)}(z) z^{B_1(t_\Delta)}e^{\Lambda(t_\Delta)z}, 
\quad
\hbox{ with }\quad
B_1(t_\Delta)=\hbox{diag}(\widehat{A}_1(t_\Delta)),
\\
 &&
\mathcal{G}_{(t_\Delta)}(z)\sim I+\sum_{k=1}^\infty F_{(t_\Delta);k} z^{-k},
\quad\quad
z\to\infty 
\quad
\hbox{ in } \overline{S}^{(t_\Delta)}(\alpha,\beta),
\eeas
if and only if  the following   conditions generalising  (\ref{Good123}) hold. 
 For  those $a\neq b\in\{1,...,n\}$ such that $u_a(t_\Delta)=u_b(t_\Delta)$, 
\be
\label{8dic2015-4}
\Bigl(\widehat{A}_1(t_\Delta)\Bigr)_{ab}=0,
\ee
 and if also $\Bigl(\widehat{A}_1(t_\Delta)\Bigr)_{aa}-\Bigl(\widehat{A}_1(t_\Delta)\Bigr)_{bb}+l-1=0$ for some $l\geq 2$, the following further  conditions must hold:
\be
\label{8dic2015-5}
\sum_{\begin{array}{c}
  \gamma\in \{1,...,n\},
\cr
u_\gamma(t_\Delta)\neq (u_a(t_\Delta)=u_b(t_\Delta))
\end{array}
} \Bigl(\widehat{A}_1(t_\Delta)\Bigr)_{a\gamma}\Bigl(F_{(t_\Delta);l-1}\Bigr)_{\gamma b}+\sum_{j=1}^{l-2} \Bigl( 
\widehat{A}_{l-j}(t_\Delta)F_{(t_\Delta);j}
\Bigr)_{ab}+\Bigl(\widehat{A}_l(t_\Delta)\Bigr)_{ab}=0.
\ee
\end{shaded}
 In the notation used here,  then  $\mathring{Y}(z)$ in (\ref{6dic2015-8}) is $Y_{(0)}(z)$, while    $\mathcal{G}(z)$ in (\ref{6dic2015-7}) is $ \mathcal{G}_{(0)}(z)$. Finally,   $\mathring{F}_k$ in (\ref{16gen2016-3}) is $F_{(0);k}$. 
Keeping into account that $(\widehat{A}_1)_{a\gamma}$ vanishes in (\ref{8dic2015-2}) for $t\to t_\Delta$ and $u_\gamma(t_\Delta)= u_a(t_\Delta)=u_b(t_\Delta)$,  it is immediate to prove the following,

\begin{shaded}
\noindent
{\bf Proposition \ref{14feb2016-2} generalised:}
If the vanishing conditions for (\ref{14feb2016-3}) and (\ref{8dic2015-2}) of  Proposition \ref{8dic2015-3} hold for  $t\to t_\Delta\in\Delta$, then  (\ref{8dic2015-4})   and (\ref{8dic2015-5}) at $t=t_\Delta$ are satisfied with the choice 
\be
\label{14dicembre2016-1}
F_{(t_\Delta);k}=F_k(t_\Delta), 
\quad
k\geq 1.
\ee
 If moreover $\Bigl(\widehat{A}_1(t_\Delta)\Bigr)_{aa}-\Bigl(\widehat{A}_1(t_\Delta)\Bigr)_{bb}+l-1\neq 0$ for every $l\geq 2$, the above (\ref{14dicembre2016-1}) is the unique choice. Namely, for the system with $t=t_\Delta$  there is only the unique formal solution   $$\Bigl(I+\sum_{k=1}^\infty 
F_k(t_\Delta)\Bigr)z^{B_1(t_\Delta)}z^{\Lambda(t_\Delta)},
\quad\quad
B_1(t_\Delta)=\hbox{diag}(\widehat{A}_1(t_\Delta)).$$
\end{shaded}
   %
   %
   %
   %
   

\noindent
\hrule

\vskip 0.5 cm
\centerline{\bf \Large PART II: Stokes Phenomenon}

\vskip 0.5 cm 
\noindent 
 When Assumption 1 holds, the system (\ref{9giugno-2}) is gauge equivalent to  (\ref{SeK2}) (i.e. system (\ref{17luglio2016-1}) in the Introduction) with $G_0(t)$ diagonalizing $A_0(t)$, namely
 \be
\label{16marzo2016-1}
{d\widehat{Y}\over dz}=\widehat{A}(z,t)~\widehat{Y},
\quad\quad
\widehat{A}(z,t):=G_0^{-1}(t)A(z,t)G_0(t)=\Lambda(t)+\sum_{k=1}^\infty\widehat{A}_k(t)z^{-k}.
\ee
 At $t_0\not \in \Delta$,  $\Lambda(t_0)$ has distinct eigenvalues, the Stokes phenomenon is studied as in   \cite{BJL1}.  We describe below the analogous  results at $t=0$ and $t_\Delta\in\Delta$, namely the existence and   uniqueness of  fundamental solutions with given asymptotics (\ref{16gen2016-3}) in  wide sectors.  
 The results could be derived from the general construction of \cite{BJL3}, especially from 
 Theorem V and VI therein\footnote{Note that notations here and in \cite{BJL3} are similar, but they indicate objects that are slightly different (for example Stokes rays $\tau_\nu$ and sectors $\mathcal{S}_\nu$ are not defined in  the same way).}. Nevertheless, it seems to be more natural to us to derive them in straightforward way, which we present below.  
  First, we concentrate on the most degenerate case $\Lambda=\Lambda(0)$, for $t=0$, so that $A(z,0)=\widehat{A}(z,0)$ and  the systems (\ref{9giugno-2}) and (\ref{16marzo2016-1}) coincide.   In Section \ref{16marzo2016-2} we consider the case of any other  $t_\Delta\in\Delta$.

\section{Stokes Phenomenon at $t=0$}   
\label{16marzo2016-4}

 \subsection{Stokes Rays of $\Lambda=\Lambda(0)$}
 \label{10marzo2015-2}

\bde
\label{31MARZO2015-1}
  The {\bf Stokes rays} associated with the pair of eigenvalues $(\lambda_j,\lambda_k)$, $1\leq j \neq k \leq n$, of $\Lambda$ are the infinitely many rays contained in the universal covering $\mathcal{R}$ of $\mathbb{C}\backslash\{0\}$, oriented outwards from $0$ to $\infty$, defined by 
$$ 
\Re\Bigl((\lambda_j-\lambda_k)z\Bigr)=0,
\quad
\Im \Bigl((\lambda_j-\lambda_k)z\Bigr)<0,
\quad\quad
z\in \mathcal{R}.
$$ 
\ede

The definition above  implies that for  a couple of eigenvalues $(\lambda_j,\lambda_k)$ the associated  rays are 
\be
\label{prethetajk}
R(\theta_{jk}+2\pi N):=\Bigl\{
~z\in \mathcal{R}~\Bigl|~z=\rho e^{i(\theta_{jk}+2\pi N)},\quad
\rho>0~
\Bigr\},
\quad\quad
N\in\mathbb{Z}.
\ee
where
\be 
\label{thetajk}
\theta_{jk}:= {3\pi\over 2} - \arg_p(\lambda_j-\lambda_k) .
\ee 

$\bullet $ {\bf Labelling:} We enumerate Stokes rays with $\nu\in\mathbb{Z}$, using directions $\tau_\nu$   introduced in Section \ref{preRAY}. Indeed, by Definition \ref{31MARZO2015-1}, Stokes rays have directions $\arg z=\tau_\nu$,  ordered in counter-clockwise sense as $\nu$ increases. For any sector of central angle $\pi$  in $\mathcal{R}$,  whose boundaries are not Stokes rays,  there exists a  $\nu_0\in \mathbb{Z}$ such that the   $\mu$ Stokes rays $\tau_{\nu_0-\mu+1}< ~\cdots~<\tau_{\nu_0-1}<\tau_{\nu_0}$ are contained in the sector. All other Stokes rays have directions
\be 
\label{genrayMARCH31}
\arg z= \tau_{\nu+k\mu}:=\tau_{\nu}+{k\pi },
\quad\quad
k\in \mathbb{Z},
\quad\quad
\nu\in\{\nu_0-\mu+1,...,~\nu_0-1,~\nu_0\}.
\ee
Rays $\tau_{\nu_0-\mu+1}< ~\cdots~<\tau_{\nu_0-1}<\tau_{\nu_0}$   are  called  a set of  {\bf basic Stokes rays}, because they generate the others \footnote{Although notations are similar to \cite{BJL3},  definitions are slightly different here.}. 

\vskip 0.2 cm 
$\bullet $ {\bf Sectors $ \mathcal{S}_\nu$:} Consider  a sector $S$ of central opening less than $\pi$, with boundary rays which are not Stokes rays. 
 The first  rays encountered outside $S$ upon moving clockwise and anti-clockwise,  will be called {\bf the two nearest Stokes rays outside $S$}.  
 If $S$ contains in its interior  a set of basic rays, say  $\tau_{\nu+1-\mu}$, $\tau_{\nu+2-\mu}$, ..., $\tau_{\nu}$,  then the two nearest Stokes rays 
 outside $S$ are $\tau_{\nu-\mu}$ and   $\tau_{\nu+1}$, namely  the  boundaries rays  of  $\mathcal{S}_\nu$ in (\ref{SECT00bis}), and obviously $S\subset \mathcal{S}_\nu$.  
\vskip 0.2 cm

$\bullet $ {\bf Projections onto $\mathbb{C}$:} If $R$ is any of the  rays in $\mathcal{R}$, its projection onto $\mathbb{C}$ will be denoted $PR$. For example,  let $\overline{\lambda}_j$ be the complex conjugate of $\lambda_j$, then for any $N$ the projection of (\ref{prethetajk}) is 
$$
PR(\theta_{jk}+2\pi N)=\bigl\{ z\in\mathbb{C}~\bigr|~z=-i\rho (\overline{\lambda}_j-\overline{\lambda}_k), ~\rho>0
\bigr\}.
$$

\bde 
\label{14feb2016-5}
An {\bf admissible ray for $\Lambda(0)$} is a   ray $R(\widetilde{\tau}):=\bigl\{
z\in \mathcal{R}~\bigl|~z=\rho e^{i\widetilde{\tau}},\quad
\rho>0
\bigr\}$ in $\mathcal{R}$,  of direction $\widetilde{\tau}\in\mathbb{R}$, 
which does not  coincide with any of  the Stokes rays of $\Lambda(0)$. Let 
\beas
&&  l_+(\widetilde{\tau}):= PR(\widetilde{\tau}+2k\pi),
\quad
 l_-(\widetilde{\tau}):= PR(\widetilde{\tau}+(2k+1)\pi),\quad\quad
 k\in\mathbb{Z},
\\
\\
&& l(\widetilde{\tau}):=l_-(\widetilde{\tau})\cup \{0\}\cup l_{+}(\widetilde{\tau}).
\eeas
 We call  the oriented line $l(\widetilde{\tau})$ an  {\bf admissible line for $\Lambda(0)$}. Its positive part is $l_{+}(\widetilde{\tau})$.
\ede
Observe that there exists a suitable $\nu$ such that $
\tau_\nu<\widetilde{\tau}<\tau_{\nu+1}
$, which implies 
$$ 
R(\widetilde{\tau})\subset \mathcal{S}_{\nu}\cap \mathcal{S}_{\nu+\mu},
\quad\quad
R(\widetilde{\tau}+\pi)\subset \mathcal{S}_{\nu+\mu}\cap \mathcal{S}_{\nu+2\mu}.
$$
In particular, if $\tau$ is as in (\ref{TAU}), then $\tau_{-1}<\tau<\tau_0$, and  $l(\tau)$ is an admissible line.


\subsection{Uniqueness of the Fundamental Solution with given Asymptotics }

In case of distinct eigenvalues, it is well known that  there exists a unique fundamental solution,  determined by the asymptotic behaviour  given by  the formal solution, on a sufficiently large sector. This fact must now be proved also at coalescence points. 

Let  the diagonal form  $\Lambda=\Lambda_1\oplus\cdots\oplus\Lambda_s$  of $A_0$ be fixed. 
  Let   a formal solution $
\mathring{Y}_F(z)=\Delta_0 F(z) z^Dz^Le^{\Lambda z} 
$
  be chosen in the class of formal solutions with given $\Delta_0$,  $D$, $L$, $\Lambda$, as in Definition \ref{26maggio2017-1}.  As a consequence of Theorem \ref{ShybuyaOLD} and Theorem \ref{31marzo2015-12},  there exists at least one actual solution as in (\ref{7dic2015-1}), namely 
\be
\label{3maggio2016-3}
\mathring{Y}(z)=\Delta_0 \mathcal{G}(z) 
z^Dz^L~e^{\Lambda z},\quad\quad
\mathcal{G}(z)\sim F(z),~ \quad
z\to\infty,\quad
z\in \overline{S}(\alpha,\beta).
\ee
Observe that $\overline{S}(\alpha,\beta)$ can be chosen in Theorem \ref{ShybuyaOLD} so that {\it it contains the  set of basic Stokes rays of $\mathcal{S}_\nu$}, namely $\tau_{\nu+1-\mu}$, ..., $\tau_{\nu-1}$,  $\tau_{\nu}$. 
  The asymptotic relation in (\ref{3maggio2016-3})   is conventionally  written as follows,
$$ 
\mathring{Y}(z)\sim \mathring{Y}_F(z),
\quad\quad
z\to \infty,
\quad
z\in \overline{S}(\alpha,\beta).
$$
Now,  $\mathcal{G}(z)$ is holomorphic for $|z|$ sufficiently big in   $\overline{S}(\alpha,\beta)$. 
 Since  $A(z)$ has no singularities for $|z|\geq N_0$ large, except the point at infinity,  then {\it $\mathring{Y}(z)$ and $\mathcal{G}(z)$ have analytic continuation on $\mathcal{R}\cap \{|z|\geq N_0\}$}.

\ble
\label{5aprile2015-2}
Let $C\in GL(n,\mathbb{C})$, and $S$ an arbitrary sector. Then
$$ 
z^Dz^LCz^{-L}z^{-D}\sim I,
\quad
z\to\infty\hbox{ in }S
\quad
\Longleftrightarrow
\quad
z^Dz^LCz^{-L}z^{-D} =I
\quad
\Longleftrightarrow
\quad
C=I.
$$
\ele
\noindent
The simple proof is left as an exercise.

\ble[Extension Lemma]
\label{3giugno2016-1}
Let $\mathring{Y}(z)$ be a fundamental matrix solution with asymptotic behaviour, 
$$
\mathring{Y}(z)\sim \mathring{Y}_F(z),
\quad\quad
z\to \infty,
\quad
z\in S,
$$
in a sector $S$ of a non specified central opening angle.
Suppose that there is a sector $\widetilde{S}$ not containing Stokes rays,  such that $S\cap \widetilde{S}\neq \emptyset$. Then, 
$$
\mathring{Y}(z)\sim \mathring{Y}_F(z),
\quad\quad
z\to \infty,
\quad
\hbox{ for }~z\in S\cup\tilde{S}.
$$
\ele

\vskip 0.2 cm 
\noindent
{\it Proof:} $\widetilde{S}$ has central opening angle less than $\pi$, because it does not contain Stokes rays. 
Therefore, by Theorem \ref{ShybuyaOLD}, there exists a fundamental matrix solution $\widetilde{Y}(z)=\Delta_0\tilde{\mathcal{G}}(z)z^D z^L e^{\Lambda z}$, with asymptotic behaviour $ 
\widetilde{Y}(z)\sim \mathring{Y}_F(z)$, for $z\to \infty$, $z\in \widetilde{S}$. 
The two fundamental matrices are connected by an invertible matrix $C$, namely  $ 
 \mathring{Y}(z)=\widetilde{Y}(z)~C$, $z\in S\cap \widetilde{S}$. 
 Therefore, 
$$
\widetilde{\mathcal{G}}^{-1}(z)~\mathcal{G}(z)=z^D z^L e^{\Lambda z}~C~e^{-\Lambda z} z^{-L} z^{-D}.
$$
Since $\mathcal{G}(z) $ and $\widetilde{\mathcal{G}}^{-1}(z)$ have the same asymptotic behaviour in $S\cap \widetilde{S}$,  the l.h.s has asymptotic series equal to the identity matrix $I$, for  $z\to\infty$ in $z\in S\cap \widetilde{S}$. Thus, so must hold for the r.h.s. The r.h.s has diagonal-block structure inherited from $\Lambda$. We write the block $[i,j]$, $1\leq i,j\leq s$, of $C$ with simple notation $C_{ij}$.  The block $[i,j]$ in r.h.s. is  then,
$
e^{(\lambda_i-\lambda_j)z}z^{D_i}z^{L_i} ~C_{ij} ~ z^{-L_j}z^{-D_j}
$.
Hence, the following  must hold,
$$
e^{(\lambda_i-\lambda_j)z}z^{D_i}z^{L_i} ~C_{ij} ~ z^{-L_j}z^{-D_j}\sim \delta_{ij} ~I_i,
\quad
\quad
z\to\infty,
\quad
z\in S\cap \widetilde{S}.
$$
Here $I_i$ is the $p_i\times p_i$ identity matrix.

\vskip 0.2 cm 
\noindent
-- For $i\neq j$:  
Since there are no Stokes rays in $\widetilde{S}$, the sign of $\Re(\lambda_i-\lambda_j)z$ does not change in $\widetilde{S}$. This implies that $
e^{(\lambda_i-\lambda_j)z}z^{D_i}z^{L_i} ~C_{ij} ~ z^{-L_j}z^{-D_j}
\sim
 0$ for $z\to\infty$ in $\widetilde{S}$.

\vskip 0.2 cm 
\noindent
-- For $i= j$: We have $
z^{D_i}z^{L_i} ~C_{ii} ~ z^{-L_i}z^{-D_i}\sim I_i$ 
for $z\to\infty$ 
in $ S\cap \widetilde{S}$. 
From Lemma \ref{5aprile2015-2} it follows that  $
z^{D_i}z^{L_i} ~C_{ii} ~ z^{-L_i}z^{-D_i} = I_i
$. This holds on the whole $\widetilde{S}$.  

\vskip 0.2 cm 
\noindent
The above considerations  imply that  $
z^D z^L e^{\Lambda z}~C~e^{-\Lambda z} z^{-L} z^{-D}
\sim I$ for  $z\to\infty$  in $\widetilde{S}$. 
From the fact that $\widetilde{\mathcal{G}}(z)\sim  I+\sum_{k\geq 1}\mathring{F}_kz^{-k} $   in $\widetilde{S}$,  we conclude that also $
\mathcal{G}(z)\sim I+\sum_{k\geq 1}\mathring{F}_kz^{-k}$ for $z\to\infty$ in $\widetilde{S}$. 
Therefore, $\mathcal{G}(z)\sim  I+\sum_{k\geq 1}\mathring{F}_kz^{-k}$ in $S\cup \widetilde{S}$.  
{$\Box$}

\vskip 0.3 cm

The extension Lemma immediately implies the following: 

\begin{shaded}
\bth[Extension Theorem]
\label{EStheoremMARCH31}
Let $\mathring{Y}(z)$ be a fundamental matrix solution such that $\mathring{Y}(z)\sim \mathring{Y}_F(z)$ in a sector $S$, containing a set of $\mu$ basic Stokes rays, and no other Stokes rays. Then, the asymptotics $\mathring{Y}(z)\sim \mathring{Y}_F(z)$ holds  on the open sector  which extends up to the two  nearest Stokes rays outside $S$. This sector has central opening angle greater than $\pi$ and is a sector $\mathcal{S}_\nu$ for a suitable $\nu$. 
\eth
\end{shaded}

\noindent
{\bf Important Remark:} The above extension theorem has the important consequence that in the statement of Theorem \ref{31marzo2015-12} and  Proposition \ref{6ott2015-5}, the matrix $\mathcal{G}(z)$, which  has analytic continuation in $\mathcal{R}$ for $|z|\geq N_0$, has the prescribed   asymptotic expansion in {\it any} proper closed subsector of $\mathcal{S}_\nu$. Hence, {\it by definition,   the asymptotics holds in the open sector $\mathcal{S}_\nu$}.

\begin{shaded}
\bth[Uniqueness Theorem]
\label{UniThT}
A fundamental matrix $\mathring{Y}(z)$   as (\ref{7dic2015-1})  such that $ 
\mathring{Y}(z)\sim \mathring{Y}_F(z)$, for  $z\to\infty$   
in a sector $S$ containing a set of basic Stokes rays, is unique. 
 In particular, this applies if $\overline{S}(\alpha,\beta)$ of  Theorem 
\ref{ShybuyaOLD} contains a set of basic Stokes rays. 
\eth
\end{shaded}

\vskip 0.2 cm 
\noindent
{\it Proof:} Suppose that there are two solutions $\mathring{Y}(z)$ and $\widetilde{Y}(z)$ with asymptotic representation $\mathring{Y}_F(z)$ in a sector  $ S$, which contains $\mu$ basic Stokes rays. Then, there exists an invertible matrix $C$ such that $ 
 \mathring{Y}(z)=\widetilde{Y}(z)~C$, namely
$$ 
\widetilde{\mathcal{G}}^{-1}(z)~\mathcal{G}(z)=z^Dz^Le^{\Lambda z}~C~e^{-\Lambda z}z^{-L}z^{-D}.
$$
The l.h.s. has asymptotic series equal to $I$ as $z\to\infty$ in $S$. 
Therefore, for the block $[i,j]$, the following must hold, 
$$
e^{(\lambda_i-\lambda_j)z}z^{D_i}z^{L_i} ~C_{ij} ~ z^{-L_j}z^{-D_j}~\sim~\delta_{ij}~I_i,
\quad
\quad
\hbox{for }~z\to \infty ~\hbox{ in }~ S.
$$
Since ${ S}$ contains a set of basic Stokes rays, $\Re(\lambda_i-\lambda_j)z$ changes sign at least once in ${S}$, for any $1\leq i \neq j \leq s$. Thus, $e^{(\lambda_i-\lambda_j)z}$  diverges in some subsector of ${S}$.  For $i\neq j$ this requires that $ 
C_{ij}=0$  for $i\neq j$. 
For $i=j$, we have $
z^{D_i}z^{L_i} C_{ii} z^{-L_i}z^{-D_i}\sim I_i
$. 
Lemma \ref{5aprile2015-2} assures that $C_{ii}=I_i$. Thus, $C=I$.  
$\Box$

\vskip 0.3 cm 

\noindent
$\bullet$ {\bf [The notation $\mathring{Y}_\nu(z)$]:} 
There exist  $\nu \in \mathbb{Z}$ such that a sector $S$ of Theorem \ref {UniThT} contains the basic rays  $\tau_{\nu+1-\mu}$, ..., $\tau_{\nu-1}$,  $\tau_{\nu}$. Hence $S\subset \mathcal{S}_\nu$. 
 The unique fundamental solution of Theorem \ref {UniThT},  with asymptotics extended to  $ \mathcal{S}_\nu$ according to Theorem \ref{EStheoremMARCH31},  will be denoted $\mathring{Y}_\nu(z)$.


   \subsection{Stokes Matrices}
\label{9dic2015-3}
The definition of Stokes matrices is  standard.   Recall that the  Stokes rays associated with $(\lambda_j,\lambda_k)$ are  (\ref{prethetajk}). Consider  also the rays
$$
R(\theta_{jk}+2\pi N+\delta)=\Bigl\{
~z\in \mathcal{R}~\Bigl|~z=\rho e^{i(\theta_{jk}+2\pi N+\delta)},\quad
\rho>0~
\Bigr\},\quad\quad
N\in\mathbb{Z}.
$$
The sign of $\Re(\lambda_j-\lambda_k)z$ for $z\in R_N(\theta_{jk}+\delta)$  is: 
$$
\left\{
\begin{array}{ccc}
\Re(\lambda_j-\lambda_k)z<0, &\hbox{ for }& -{\pi }<\delta<0 \quad
\hbox{ mod }{2\pi}
\cr
\Re(\lambda_j-\lambda_k)z>0, &\hbox{ for }& 0<\delta<{\pi } \quad
\hbox{ mod }{2\pi}
\cr
\Re(\lambda_j-\lambda_k)z=0, &\hbox{ for }& \delta=0,~{\pi},~-{\pi } \quad
\hbox{ mod }{2\pi}
\end{array}
\right. 
$$

\bde[Dominance relation] 
In a sector where $\Re(\lambda_j-\lambda_k)z>0$, $\lambda_j$ is said to be {\rm dominant} over $\lambda_k$ in that sector, and we write $\lambda_j \succ \lambda_k$. In a sector  where $\Re(\lambda_j-\lambda_k)z<0$, $\lambda_j$ is said to be {\rm sub-dominant}, or dominated by $\lambda_k$, and we write $\lambda_j\prec \lambda_k$. 
\ede
 If a sector $S$  does not contain Stokes rays in its interior, it is well defined a  dominance relation in $S$, which  {\it determines an {\bf ordering  relation} among  eigenvalues, {\bf referred to the sector $S$}}. 

\vskip 0.2 cm 
Denote  by 
$$\mathring{Y}_\nu(z)~ \hbox{ and }~\mathring{Y}_{\nu+\mu}(z)
$$
 the unique fundamental solutions (\ref{7dic2015-1}) with asymptotic  behaviours $\mathring{Y}_F(z)$ on $\mathcal{S}_\nu$ and $\mathcal{S}_{\nu+\mu}$ respectively, as in Theorem \ref{UniThT}. Observe that $\mathcal{S}_\nu\cap\mathcal{S}_{\nu+\mu}=S(\tau_\nu,\tau_{\nu+1})$ is not empty and does not contain Stokes rays. 

\bde
\label{20marzo2016-2}
For any $\nu\in\mathbb{Z}$, the {\bf Stokes matrix} $\mathring{\mathbb{S}}_\nu$ is the connection matrix such that   
\be
\label{tramonto6aprile2015}
\mathring{Y}_{\nu+\mu}(z)=\mathring{Y}_\nu(z)\mathring{\mathbb{S}}_\nu,\quad\quad
z\in\mathcal{R}. 
\ee
\ede

\bpr
\label{triastr2015}
 Let $\prec$ be the dominance relation referred to the sector $ \mathcal{S}_\nu\cap \mathcal{S}_{\nu+\mu}$. Then, the Stokes matrix $\mathring{\mathbb{S}}_\nu$ has the following  block-triangular structure:
\beas
&&
\mathring{\mathbb{S}}_{jj}^{(\nu)}=I_{p_j}, 
\\
&&
\mathring{\mathbb{S}}_{jk}^{(\nu)}=0~  ~\hbox{ for } \lambda_j\succ \lambda_k  \hbox{ in } \mathcal{S}_\nu\cap \mathcal{S}_{\nu+\mu},\quad\quad
j,k\in\{1,2,...,s\}.
\eeas
\epr

\vskip 0.2 cm 
\noindent
{\it Proof:}
We re-write  (\ref{tramonto6aprile2015}) as, 
$$
\mathcal{G}_\nu^{-1}(z)~\mathcal{G}_{\nu+\mu}(z)= z^Dz^Le^{\Lambda z}~\mathring{\mathbb{S}}_\nu~e^{-\Lambda z}z^{-L}z^{-D}.
$$
For $z\in \mathcal{S}_\nu\cap \mathcal{S}_{\nu+\mu}$, the l.h.s. has asymptotic expansion equal to $I$. 
Hence, the same must hold  for the r.h.s. 
Recalling that no Stokes rays lie in $\mathcal{S}_\nu\cap \mathcal{S}_{\nu+\mu}$, we find:

\vskip 0.2 cm 
$\bullet$ For $j\neq k$, we have $
e^{(\lambda_j-\lambda_k)z}z^{D_j}z^{L_j}\mathring{\mathbb{S}}_{jk}^{(\nu)}z^{-L_k}z^{-D_k} \sim 0$   in $ \mathcal{S}_\nu\cap \mathcal{S}_{\nu+\mu}$
 if and only if 
$\mathring{\mathbb{S}}_{jk}^{(\nu)}=0$   for $ \lambda_j\succ \lambda_k$, 
where the dominance relation is referred to the sector $ \mathcal{S}_\nu\cap \mathcal{S}_{\nu+\mu}$. 

\vskip 0.2 cm 
$\bullet$ For $j =k$, we have 
$
z^{D_j}z^{L_j}\mathring{\mathbb{S}}_{jj}^{(\nu)}z^{-L_j}z^{-D_j}\sim I_{p_j}$ if and only if 
$\mathring{\mathbb{S}}_{jj}^{(\nu)}=I_{p_j}$, by Lemma \ref{5aprile2015-2}. 
This proves the Proposition. $\Box$


\subsection{Canonical Sectors, Complete Set of Stokes Matrices, Monodromy Data}

There are no Stokes rays in the intersection of successive sectors $\mathcal{S}_{\nu+k\mu}$ and $\mathcal{S}_{\nu+(k+1)\mu}$ (recall that  $\tau_\nu+{k\pi }=\tau_{\nu+k\mu}$ for any $k\in\mathbb{Z}$). 
Therefore, we can  introduce the unique fundamental matrix solutions 
\begin{equation}
\label{8agosto2016-1}
\mathring{Y}_{\nu+k\mu}(z)
\end{equation}
 with asymptotic behaviour $\mathring{Y}_F(z)$ in $\mathcal{S}_{\nu+k\mu}$, and the Stokes matrices $\mathring{\mathbb{S}}_{\nu+k\mu}$ connecting them, 
$$
\mathring{Y}_{\nu+(k+1)\mu}(z)=\mathring{Y}_{\nu+k\mu}(z)~\mathring{\mathbb{S}}_{\nu+k\mu},\quad\quad
z\in\mathcal{R}.
$$
From Proposition \ref{triastr2015}, it follows that the blocks $[j,k]$ and $[k,j]$  satisfy
$$
\mathring{\mathbb{S}}^{(\nu)}_{jk}=0 \hbox{ for }\lambda_j\succ \lambda_k \hbox{ in } \mathcal{S}_\nu\cap \mathcal{S}_{\nu+\mu}
\quad
\Longleftrightarrow
\quad
\mathring{\mathbb{S}}^{(\nu+\mu)}_{kj}=0 \hbox{ for the same $(j,k)$}.
$$
We call   $
\mathcal{S}_\nu$, $\mathcal{S}_{\nu+\mu}$, $\mathcal{S}_{\nu+2\mu}
$ 
the {\bf canonical sectors associated with $\tau_\nu$}. 


\vskip 0.3 cm 
Given a formal solution, 
a simple computation (recall that $[L,\Lambda]=0$) yields $
\mathring{Y}_F(e^{2\pi i }z)=\mathring{Y}_F(z)~e^{2\pi i L}$. 
$L$ is called {\bf exponent of formal monodromy}.

\begin{shaded}
\bth
\label{26maggio2017-2}
We introduce the  notation $z_{(\nu)}$ if  $z\in \mathcal{S}_\nu$. Thus $
z_{(\nu+2\mu)}=e^{2\pi i }z_{(\nu)}
$.   
The following equalities hold 
\beas
&&
(i)\quad
\mathring{Y}_{\nu+2\mu}(z_{(\nu+2\mu)})=\mathring{Y}_\nu(z_{(\nu)})~e^{2\pi i L},
\\
&&
(ii)\quad
\mathring{Y}_{\nu+2\mu}(z)=\mathring{Y}_\nu(z)~\mathring{\mathbb{S}}_\nu~ \mathring{\mathbb{S}}_{\nu+\mu},\quad\quad
z\in \mathcal{R},
\\
&&
(iii)\quad
\mathring{Y}_\nu(e^{2\pi i } z)= \mathring{Y}_\nu(z)~e^{2\pi i L}~\Bigl(\mathring{\mathbb{S}}_\nu~ \mathring{\mathbb{S}}_{\nu+\mu}
\Bigr)^{-1},\quad\quad
z\in \mathcal{R}.
\eeas
where $|z|\geq N_0$ is sufficiently large, in such a way that any other singularity of $A(z)$ is contained in the ball $|z|<N_0$. 
\eth
\end{shaded}
 \noindent
{\it Proof:} As in the case of distinct eigenvalues.  Alternatively, one can adapt  Proposition 4 of \cite{BJL3} to the present case.\footnote{With the warning that notations are similar but objects are slightly different here and in \cite{BJL3}.}. $\Box$.
\vskip 0.2 cm 

The  equality  (iii) provides  the {\bf  monodromy matrix $M_\infty^{(\nu)}$ of $\mathring{Y}_\nu(z)$  at $z=\infty$}:
\be
\label{inftymon2015} 
{M_\infty}^{(\nu)}:= \Bigl(\mathring{\mathbb{S}}_\nu~ \mathring{\mathbb{S}}_{\nu+\mu}
\Bigr)~ e^{-2\pi i L}.
\ee 
  corresponding to a clockwise loop   with $|z|\geq N_0$ large, in such a way that all other singularities of 
$A(z)$ are inside the loop. 


\vskip 0.3 cm 

The two  Stokes matrices $\mathring{\mathbb{S}}_\nu$, $\mathring{\mathbb{S}}_{\nu+\mu}$, and the matrix $L$   generate all the other Stokes matrices $\mathring{\mathbb{S}}_{\nu+k\mu}$, according to the following 
proposition

\begin{shaded}
\bpr 
\label{26maggio2017-3}
For any $\nu\in\mathbb{Z}$, the following holds:  $ 
\mathring{\mathbb{S}}_{\nu+2\mu}=e^{-2\pi i L} ~\mathring{\mathbb{S}}_\nu ~e^{2\pi i L}.
$
\epr
\end{shaded}

\vskip 0.2 cm 
\noindent
{\it Proof:} For simplicity, take $\nu=0$. 
A point in  $z\in \mathcal{S}_{2\mu}\cap \mathcal{S}_{3\mu}$ can represented both as 
$z_{(2\mu)}$ and $z_{(3\mu)}$, and a point in $\mathcal{S}_0\cap \mathcal{S}_\mu$ is represented both as  $z_{(0)}$ and $z_{(\mu)}$.  Therefore, the l.h.s. of the equality $
\mathring{Y}_{3\mu}(z)=\mathring{Y}_{2\mu}(z)~\mathring{\mathbb{S}}_{2\mu}$ is $
\mathring{Y}_{3\mu}(z_{(3\mu)})=\mathring{Y}_\mu(z_{(\mu)})~e^{2\pi i L}=\mathring{Y}_0(z_{(0)}) \mathring{\mathbb{S}}_0~e^{2\pi i L}
$.  
The r.h.s. is 
$
\mathring{Y}_{2\mu}(z_{(2\mu)})~\mathring{\mathbb{S}}_{2\mu}= \mathring{Y}_0(z_{(0)})e^{2\pi i L}\quad
\mathring{\mathbb{S}}_{2\mu}$. 
Thus 
$
\mathring{Y}_0(z_{(0)})~\mathring{\mathbb{S}}_0~e^{2\pi i L}~=~ \mathring{Y}_0(z_{(0)})e^{2\pi i L}\quad
\mathring{\mathbb{S}}_{2\mu}
$. 
This proves the proposition. $\Box$

\vskip 0.2 cm 
The above proposition implies that $\mathring{\mathbb{S}}_{\nu+k\mu}$ are generated by $\mathring{\mathbb{S}}_\nu, ~\mathring{\mathbb{S}}_{\nu+\mu}$,  which  therefore form a {\bf complete set of Stokes matrices}. A  complete set of Stokes matrices and the exponent of formal monodromy are necessary and sufficient to obtain the monodromy at $z=\infty$, through formula (\ref{inftymon2015}). This justifies the following definition. 

\bde
For a chosen $\nu$, $
\bigl\{\mathring{\mathbb{S}}_\nu, ~\mathring{\mathbb{S}}_{\nu+\mu},~L\bigr\}
$ 
is a set of {\bf monodromy data} at $z=\infty$ of the system (\ref{9giugno-2}) with $t=0$.
\ede

\bre
By a factorization into Stokes factors, as in the proof of Theorem \ref{8gen2017-4} below,  it can be shown that $\mathring{\mathbb{S}}_\nu, ~\mathring{\mathbb{S}}_{\nu+\mu}$ suffice to generate $\mathring{\mathbb{S}}_{\nu+1},~...,~\mathring{\mathbb{S}}_{\nu+\mu-1}$. Hence, $\mathring{\mathbb{S}}_\nu, ~\mathring{\mathbb{S}}_{\nu+\mu}$ are really sufficient to generate all Stokes matrices. This technical part will be omitted. 
\ere

\section{Stokes Phenomenon at fixed $t_\Delta\in\Delta$}
\label{16marzo2016-2}

 The results of Section \ref{16marzo2016-4} apply to any other $t_\Delta\in \Delta$.  By a permutation matrix $P$ we  arrange 
$P^{-1}\Lambda(t_\Delta)P $ in blocks, in such a way that  each block  has only one eigenvalue and  two distinct  blocks have different eigenvalues. This is achieved by the  transformation  
 $\widehat{Y}(z,t)=P~\widetilde{Y}(z,t)$ applied to the system (\ref{16marzo2016-1}). 
 Then, the procedure is exactly the same of Section \ref{16marzo2016-4}, applied to the system
\be
\label{16marzo2016-3}
 \frac{d\widetilde{Y}}{dz}=P^{-1}\widehat{A}(z,t_\Delta)P ~\widetilde{Y}.
 \ee
 The block partition of all matrices in the computations and statements is that inherited from $P^{-1}\Lambda(t_\Delta)P$.  The Stokes rays are defined  in the same way as in Definition \ref{31MARZO2015-1}, using  the eigenvalues of $\Lambda(t_\Delta)$, namely
 \beas
&& 
\Re\Bigl((u_a(t_\Delta)-u_b(t_\Delta))z\Bigr)=0,
\quad\quad
\Im \Bigl((u_a(t_\Delta)-u_a(t_\Delta))z\Bigr)<0,
\quad\quad
z\in \mathcal{R},
\\ 
&&
\hbox{for }~ 1\leq a\neq b\leq n\quad
 \hbox{ \underline{and} } \quad
 u_a(t_\Delta)\neq u_a(t_\Delta).
\eeas
Hence, the Stokes rays associated with  $u_a(t_\Delta),u_b(t_\Delta)$ are the infinitely many rays with directions
$$ 
\arg z= \frac{3\pi}{2}-\arg_p(u_a(t_\Delta)-u_b(t_\Delta))+2N\pi, \quad
N\in \mathbb{Z}.
$$
The rays associated with  $u_b(t_\Delta),u_a(t_\Delta)$ are opposite to the above, having directions 
$$ 
\arg z= \frac{3\pi}{2}-\arg_p(u_b(t_\Delta)-u_a(t_\Delta))+2N\pi.
$$
   We conclude that all Stokes rays have directions
      $$
\arg z= \tau_{\sigma}^{(t_\Delta)}, \quad\quad
\sigma\in\mathbb{Z},
$$
 analogous to (\ref{genrayMARCH31}), with directions $\tau_\sigma^{(t_\Delta)}$ defined in Section \ref{11marzo2016-2}. 
 Once the Stokes matrices for the above system are computed, in order to go back to the original arrangement corresponding to $\Lambda(t_\Delta)$ we just apply the inverse permutation. Namely, if $\mathbb{S}$ is a Stokes matrix of   (\ref{16marzo2016-3}), then $P\mathbb{S}P^{-1}$ is a Stokes matrix for (\ref{16marzo2016-1}) with $t=t_\Delta$.


\section{Stokes Phenomenon at $t_0\not\in\Delta$}
\label{29april2016-1}
 
The results of Section \ref{16marzo2016-4} (extension theorem, uniqueness theorem, Stokes matrices, etc) apply {\it a fortiori} if the eigenvalues  are distinct, namely at a point $t_0\not\in\Delta$ such that Theorem \ref{ShybuyaOLDbis} and Corollary \ref{9dic2015-1} apply.  The block partition of $\Lambda(t_0)$ is  into one-dimensional blocks, being the eigenvalues all distinct, and we are back to the well known case of \cite{BJL1}. 
The Stokes rays   are defined  in the same way as in Definition \ref{31MARZO2015-1}, using  the eigenvalues of $\Lambda(t_0)$, namely
 $$ 
\Re\Bigl((u_a(t_0)-u_b(t_0))z\Bigr)=0,
\quad\quad
\Im \Bigl((u_a(t_0)-u_b(t_0))z\Bigr)<0,
\quad\quad
z\in \mathcal{R},\quad\quad
\forall ~1\leq a\neq b\leq n .
$$ 
Since and $u_a(t_0)\neq u_b(t_0)$ for any $a\neq b$, the above definition holds for any 
$ 
1\leq a\neq b\leq n 
$. 
Hence, the Stokes rays associated with  $u_a(t_0),u_b(t_0)$ are the infinitely many rays with directions
\be
\label{18april2016-2}
\arg z= \frac{3\pi}{2}-\arg_p(u_a(t_0)-u_b(t_0))+2N\pi, \quad
\quad
N\in \mathbb{Z}.
\ee
The rays associated with  $u_b(t_0),u_a(t_0)$ are opposite to the above, having directions 
\be
\label{18april2016-3} 
\arg z= \frac{3\pi}{2}-\arg_p(u_b(t_0)-u_a(t_0))+2N\pi.
\ee
   We conclude that all Stokes rays have directions
      $$
\arg z= \tau_{\sigma}^{(t_0)}, \quad\quad
\sigma\in\mathbb{Z},
$$
 analogous to (\ref{genrayMARCH31}), being the directions $\tau_\sigma^{(t_0)}$ defined in Section \ref{11marzo2016-2}.  We stress that $t_0$ is fixed here.  The Stokes phenomenon  is studied in the standard way. The canonical sectors are the sectors  $\mathcal{S}_\sigma^{(t_0)}$  of Theorem  \ref{ShybuyaOLDbis}. The sector $\mathcal{S}_\sigma^{(t_0)}$ contains the set of {\it basic Stokes rays}
\be
\label{18marzo2016-1}
\tau^{(t_0)}_{\sigma+1-\mu_{t_0}},\quad
 \tau^{(t_0)}_{\sigma+2-\mu_{t_0}},\quad ..., 
 \quad
 \tau^{(t_0)}_{\sigma},
\ee
which serve to generate all the other rays by adding multiples of $\pi$.  The rays $\tau^{(t_0)}_{\sigma-\mu_{t_0}}$ and $\tau^{(t_0)}_{\sigma+1}$ are the nearest Stokes rays, boundaries of  $\mathcal{S}_\sigma^{(t_0)}$.  
The Stokes matrices connect  solutions of Corollary  \ref{9dic2015-1}, having the prescribed canonical asymptotics on successive sectors, for example   $\mathcal{S}_\sigma^{(t_0)}$,  $\mathcal{S}_{\sigma+\mu_{t_0}}^{(t_0)}$,   $\mathcal{S}_{\sigma+2\mu_{t_0}}^{(t_0)}$, etc. 
 
 Our purpose is now to show how the Stokes phenomenon can be described in a consistent ``holomorphic" way as $t$ varies. The definition of Stokes matrices for varying $t$ will require some steps.

\vskip 0.3 cm 

\noindent
\hrule

\vskip 0.5 cm
\centerline{\bf \Large PART III: Cell Decomposition,  $t$-analytic Stokes Matrices}

\section{Stokes Rays rotate as  $t$ varies}
\label{13giugno2017-2}

At $t=0$, Stokes rays have  directions $3\pi/2-\arg_p(\lambda_i-\lambda_j) +2N\pi$, $1\leq i\neq j\leq s$. For $t$ away from $t=0$, the  following occurs:

\vskip 0.2 cm 
1) [Splitting]  For $1\leq i\neq j \leq s$, there are rays of directions $3\pi/2-\arg_p(u_a(t)-u_b(t))\hbox{ mod}(2\pi)$, with $u_a(0)=\lambda_i$, $u_b(0)=\lambda_j$.  These  rays are the splitting of  $3\pi/2-\arg_p(\lambda_i-\lambda_j) \hbox{ mod}(2\pi)$ into more rays.  

\vskip 0.2 cm 
2) [Unfolding] For any $i=1,2,...,s$, new rays appear,  with directions  $3\pi/2-\arg_p(u_a(t)-u_b(t))$, $u_a(0)=u_b(0)=\lambda_i$. These rays are due to the {\it unfolding} of $\lambda_i$.   

\vskip 0.2 cm 
\noindent
The cardinality of a set of  basic  Stokes rays is maximal away from the coalescence locus $\Delta$, minimal at $t=0$, and  intermediate at $t_\Delta\in\Delta\backslash \{0\}$. 

If $t\not \in \Delta$, then $u_a(t)\neq u_b(t)$ for any  $a\neq b$. The direction of every Stokes ray (\ref{18april2016-2}) or (\ref{18april2016-3})   is a continuous functions of  $t\not\in\Delta$. As $t$ varies in $\mathcal{U}_{\epsilon_0}(0)\backslash\Delta$, each one of the rays (\ref{18april2016-2}) or (\ref{18april2016-3})    rotates in $\mathcal{R}$. 

\bre {\bf Problems with enumeration of moving Stokes rays.}  Apparently, we cannot  assign a coherent labelling to the rotating rays as $t$ moves in $\mathcal{U}_{\epsilon_0}(0)\backslash \Delta$.  At a given $t_0\in \mathcal{U}_{\epsilon_0}(0)\backslash \Delta$, the rays are enumerated according to the  choice of an admissible direction  $\eta^{(t_0)}$,  as in formula (\ref{19marzo2016-1}) with $t_*=t_0$. If $t$ is very close to $t_0$, we may choose  $\eta^{(t_0)}=\eta^{(t)}$,  and we can label the rays in such a way that     $\tau^{(t)}_{\sigma}$, $\sigma\in\mathbb{Z}$,  is the result of the continuous rotation of $\tau^{(t_0)}_{\sigma}$. Nevertheless,   if $t$ moves farther in $\mathcal{U}_{\epsilon_0}(0)\backslash \Delta$,  then  some rays, while rotating, may cross with each other and cross the  rays  $R(\tau^{(t_0)}+k\pi)$, $k\in \mathbb{Z}$, which are admissible for $\Lambda(t_0)$. This phenomenon  destroys the ordering. 
 Hence,  labellings are to be taken  independently at $t_0$ and at any other  $t\in \mathcal{U}_{\epsilon_0}(0)\backslash \Delta$, with respect to independent admissible directions $\eta^{(t_0)}$ and $\eta^{(t)}$. In this way,   $\tau^{(t)}_{\sigma}$ will not be the deformation of  a $\tau^{(t_0)}_{\sigma}$ with the same $\sigma$.  
 
This complication in assigning a coherent numeration to rays and  sectors as $t$ varies will be solved in Section \ref{29sett2015-1}, by introducing a new labelling,  valid for almost all $t\in\mathcal{U}_{\epsilon_0}(0)$, induced by the labelling at  $t=0$. Before that, we need some topological preparation. 
\ere
  
\section{Ray Crossing,  Wall Crossing  and Cell Decomposition} 
 \label{onCELLS}

 We consider an oriented admissible ray $R(\widetilde{\tau})$ for $\Lambda(0)$, with  direction $\widetilde{\tau}$, as in Definition \ref{14feb2016-5} and we project $\mathcal{R}$ onto $\mathbb{C}\backslash\{0\}$. For  $t\in \mathcal{U}_{\epsilon_0}(0)\backslash \Delta$, some projected rays associated with   $\Lambda(t)$  will be to the left of $l(\widetilde{\tau})$ and  some to the right. Moreover, some projected ray  may lie exactly on $l(\widetilde{\tau})$, in which case we improperly say that {\it ``the ray lies on  $l(\widetilde{\tau})$''}.  
Suppose we start at a value  $t_*\in \mathcal{U}_{\epsilon_0}(0)\backslash \Delta$  such that  no rays associated with $\Lambda(t_*)$  lie on $l(\widetilde{\tau})$. 
 If $t$ moves away from $t_*$ in $ \mathcal{U}_{\epsilon_0}(0)\backslash \Delta$, then the directions of Stokes rays  change continuously  and the  projection of two or more  rays\footnote{Crossing involves always at least two opposite projected rays, which have directions differing by $\pi$. One projection crosses the positive part   $l_{+}(\widetilde{\tau})$ of $l(\widetilde{\tau})$, and one  projection crosses  the negative part $l_{-}(\widetilde{\tau})=l_{+}(\widetilde{\tau}\pm \pi)$.}   may cross  $l(\widetilde{\tau})$ as $t$ varies, in which case  we say that {\it ``two or more rays cross  $l(\widetilde{\tau})$''}. 
 Let 
$$ 
\widetilde{\eta}:=\frac{3\pi}{2}-\widetilde{\tau}.
$$
 Two or more Stokes rays  cross $l(\widetilde{\tau})$  for $t$ belonging to the following  {\bf crossing locus}
$$ 
X(\widetilde{\tau}):=\bigcup_{1\leq a < b\leq n
}
\Bigl\{t\in \mathcal{U}_{\epsilon_0}(0)~\Bigl|~u_a(t)\neq u_b(t),\quad
\arg_p(u_a(t)-u_b(t))
=\widetilde{\eta} \hbox{ mod } \pi \Bigr\}.
$$
Let 
$$ 
W
(\widetilde{\tau}) :=\Delta\cup X(\widetilde{\tau}).
 $$

\bde  
\label{11feb2016-1}
A {\bf  $\widetilde{\tau}$-cell} is every connected component of the set $ 
\mathcal{U}_{\epsilon_0}(0)\backslash W(\widetilde{\tau}) 
 $.
\ede

 $ W
(\widetilde{\tau})$ is the {\it ``wall" of the cells}.   
For $t$ in a $\widetilde{\tau}$-cell, $\Lambda(t)$ is diagonalisable with distinct eigenvalues, and  the Stokes rays  projected onto $\mathbb{C}$ lie either to the left or to the right of  $l(\widetilde{\tau})$. If    $t$  varies and hits $W
(\widetilde{\tau})$, then either some Stokes rays disappear (when $t\in \Delta$), or some rays cross the admissible line  $l(\widetilde{\tau})$ (when $t\in X(\widetilde{\tau})$). Notice that
$$
\Delta\cap X(\widetilde{\tau})\neq\emptyset.
$$

 A cell is open, by definition. If the eigenvalues are linear in $t$, as in (\ref{4gen2016-2}),  we will show in Section \ref{7agosto2016-1} that a cell  is simply connected and convex, namely it is a topological cell, so justifying the name. Explicit examples and figures are given in the Appendix.

\subsection{Topology of $\widetilde{\tau}$-cells and hyperplane arrangements}
\label{7agosto2016-1}

In order to study the topology of the $\widetilde{\tau}$-cells, it is convenient to first extend their definition to $\mathbb{C}^n$. 
A $\widetilde{\tau}$-cells in $\mathbb{C}^n$ can be proved to be  homeomorphic to an open ball, therefore it is a {\it cell} in the topological sense.  A $\widetilde{\tau}$-cell in $\mathbb{C}^n$ is  defined to be  a connected component of $\mathbb{C}^n\backslash (\Delta_{\mathbb{C}^n}
\cup X_{\mathbb{C}^n}(\widetilde{\tau}))$, where  
\beas 
&&
 \Delta_{\mathbb{C}^n}:=\bigcup_{1\leq a<b\leq n}\Bigl\{u\in\mathbb{C}^n~\Bigl|~u_a=u_b\Bigr\} ,
 \\
 &&
 X_{\mathbb{C}^n}(\widetilde{\tau}):= \bigcup_{1\leq a<b\leq n} \Bigl\{u\in \mathbb{C}^n~\Bigl|~u_a-u_b\neq 0~\hbox{\it and }\arg_p(u_a-u_b)
=
\widetilde{\eta} \hbox{ mod } \pi \Bigr\}.
\eeas
Recall that $\widetilde{\eta}=\frac{3\pi}{2}-\widetilde{\tau}$.

We identify $\mathbb{C}^n$ with $\mathbb{R}^{2n}$. A point $u=(u_1,...,u_n)$ is identified with $({\bf x},{\bf y})=(x_1,...,x_n,~y_1,...,y_n)$, by $u_a=x_a+iy_a$, $1\leq a \leq n$.  Therefore 

a) $\Delta_{\mathbb{C}^n}$ is identified with
$$ 
A:=\bigcup_{1\leq a< b\leq n}\Bigl\{({\bf x},{\bf y})\in\mathbb{R}^{2n}~\Bigl|~x_a-x_b=y_a-y_b=0\Bigr\}.
$$

b) $X_{\mathbb{C}^n}(\widetilde{\tau})$ is identified with 
$$ 
B:=\bigcup_{1\leq a< b\leq n} \Bigl\{({\bf x},{\bf y})\in\mathbb{R}^{2n}~\Bigl|~(x_a,y_a)\neq (x_b,y_b) ~\hbox{\it and } L_{ab}({\bf x},{\bf y})=0 \Bigr\}
$$
where $L_{ab}({\bf x},{\bf y})$ is a linear function
\bea
\label{15marzo2016-1}
&& L_{ab}({\bf x},{\bf y})=(y_a-y_b)-\tan \widetilde{\eta}~ (x_a-x_b),\quad\quad
 \hbox{ for }\widetilde{\eta}\neq \frac{\pi}{2}\hbox{ mod }\pi,
\\
\label{15marzo2016-2}
&&
L_{ab}({\bf x},{\bf y})= x_a-x_b,\quad\quad
 \hbox{ for }\widetilde{\eta}= \frac{\pi}{2}\hbox{ mod }\pi.
\eea
Hence $A\cup B$  is a {\it union of  hyperplanes} $H_{ab}$:
$$ 
A\cup B= \bigcup_{1\leq a< b\leq n} H_{ab},\quad\quad
H_{ab}:=\{({\bf x},{\bf y})\in\mathbb{R}^{2n}~|~L_{ab}({\bf x},{\bf y})=0\}.
$$
Note that $L_{ab}({\bf x},{\bf y})=0$ if and only if $L_{ba}({\bf x},{\bf y})=0$, namely $H_{ab}=H_{ba}$. The  set $\mathcal{A}=\{H_{ab}\}_{a<b}$ is known as  a {\it hyperplane arrangement} in $\mathbb{R}^{2n}$.  We have proved the following lemma

\ble 
Let $u\in\mathbb{C}^n$ be represented as $u={\bf x}+i{\bf y}$,    $({\bf x},{\bf y})\in\mathbb{R}^{2n}$. 
  Then, $\Delta_{\mathbb{C}^n}
\cup X_{\mathbb{C}^n}(\widetilde{\tau})$ is the union of hyperplanes $H_{ab}\in \mathcal{A}$  defined by the linear equations  $L_{ab}({\bf x},{\bf y})=0$, $1\leq a<b\leq $,  as in (\ref{15marzo2016-1}), (\ref{15marzo2016-2}).

\ele

 Properties of finite hyperplane arrangements in $\mathbb{R}^{2n}$ are well knows. In particular, consider the set 
 $$ 
 \mathbb{R}^{2n} - \bigcup_{1\leq a < b\leq n} H_{ab}.
 $$
 A connected component of the above set is called a {\it region} of $\mathcal{A}$. It is well known that every region of $\mathcal{A}$ is open and convex, and hence homeomorphic to the interior of an $2n$-dimensional ball of $\mathbb{R}^{2n}$. It is therefore {\it a cell} in the proper sense. We have proved the following 
 
 \bpr
 A $\widetilde{\tau}$-cell in $\mathbb{C}^n$ is a cell, namely  an  open and convex subset of $\mathbb{C}^n$, homeomorphic to the open ball $\{ u\in\mathbb{C}^n~|~|u_1|^2+\cdots + |u_n|^2<1\}=\{({\bf x},{\bf y})\in \mathbb{R}^{2n}~| ~x_1^2+\cdots +y_n^2<1\}$. 
 \epr 
  
\bre
\label{24marzo2016-1}
Three hyperplanes with one index in common intersect. Indeed, let $b$ be the common index. Then, 
$$
\left\{
\begin{array}{c}
L_{ab}({\bf x},{\bf y})=0
\\
L_{bc}({\bf x},{\bf y})=0
\end{array}
\right.\quad\quad
\Longrightarrow \quad\quad
L_{ac}({\bf x},{\bf y})=0.
$$
Hence,
$$
H_{ab}\cap H_{bc}\subset H_{ac},\quad
H_{bc}\cap H_{ac}\subset H_{ab},\quad
H_{ac}\cap H_{ab}\subset H_{bc}.
$$
Equivalently
$$ 
H_{ab}\cap H_{bc}\cap H_{ac}=H_{ab}\cap H_{bc}=H_{ab}\cap H_{ac} = H_{bc}\cap H_{ac}.
$$
\ere

 \vskip 0.2 cm 
We now consider $\widetilde{\tau}$-cells in $\mathcal{U}_{\epsilon_0}(0)$ in case the eigenvalues of $\Lambda(t)$ are {\it linear in $t$} as in (\ref{4gen2016-2}).  The arguments above apply to this case, since  $u_a=u_a(0)+t_a$ is a linear translation.  
Let $u(0)=(u_1(0),...,u_n(0))$ be as in (\ref{29gen2016-1})-(\ref{29gen2016-3}), so that $u(t)=u(0)+t$. Let us split $u(t)$ into real ($\Re$) and imaginary ($\Im$) parts:
$$
u(0)={\bf x}_0+i{\bf y}_0
,
\quad
t=\Re t+i \Im t\quad\quad
\Longrightarrow
\quad\quad
u(t)=\Bigl({\bf x}_0+i{\bf y}_0\Bigr)+\Bigr(\Re t +i \Im t \Bigr).
$$ 
Here, $\Re t:=(\Re t_1,...,\Re t_n)\in \mathbb{R}^n$ and $\Im t:=(\Im t_1,...,\Im t_n)\in \mathbb{R}^n$. 
Define the hyperplanes 
\be
\label{13april2016-1}
H_{ab}^\prime:=\Bigl\{ (\Re t, \Im t)\in\mathbb{R}^n~\Bigl|~L_{ab}(\Re t , \Im t)+L_{ab}({\bf x}_0,{\bf y}_0)=0\Bigr\},\quad
\quad
1\leq a\neq b \leq n,
\ee
and
\be
\label{21marzo2016-1}
\widetilde{H}_{ab}:=H_{ab}^\prime\cap \mathcal{U}_{\epsilon_0}(0).
\ee 
Then,  
$$
\Delta\cup X(\widetilde{\tau})=\bigcup_{1\leq a<b\leq n} \widetilde{H}_{ab}.
$$
Note that $L_{ab}({\bf x}_0,{\bf y}_0)=0$ for any $a\neq b$ corresponding to a coalescence $u_a(t)-u_b(t)\to 0$ for $t\to 0$.

\bcr
\label{20marzo2016-3}
If the eigenvalues of $\Lambda(t)$ are {\it linear} in $t$ as in (\ref{4gen2016-2}), then a $\widetilde{\tau}$-cell in $\mathcal{U}_{\epsilon_0}(0)$ is  simply connected. \ecr

\vskip 0.2 cm 
\noindent
{\it Proof:}  Any of the regions of a the hyperplane arrangement with hyperplanes (\ref{13april2016-1}) is open and convex. $\mathcal{U}_{\epsilon_0}(0)$ is a polydisc, hence it is convex. The intersection of a region and  $\mathcal{U}_{\epsilon_0}(0)$  is then convex and simply connected. $\Box$

\bre The $\widetilde{H}$'s enjoy the same properties of hyperplanes $H$'s as in Remark \ref{24marzo2016-1}. In other words, if  a Stokes ray associated with the pair $u_a(t),u_b(t)$  and a Stokes ray associated with $u_b(t),u_c(t)$ cross an admissible direction $R(\widetilde{\tau} \hbox{ mod }\pi)$ at some point $t$, then also a ray associated with $u_a(t),u_c(t)$ does. 
\ere

\bre
\label{11april2016-3}
We anticipate the fact that if $\epsilon_0$ is sufficiently small as in Section \ref{16gen2016-5}, then $
\widetilde{H}_{ab}\cap \mathcal{U}_{\epsilon_0}(0)=\emptyset
$ 
for any  $a\neq b$  such that  for $t\to 0$,  $u_a(t)\to \lambda_i$ and $u_b(t)\to\lambda_j$  with $1\leq i\neq j\leq s$ (i.e. $u_a(0)\neq u_b(0)$).  See below Remark \ref{11april2016-2} for explanations. 
\ere
\section{Sectors $\mathcal{S}_\nu(t)$ and $\mathcal{S}_\nu (K)$ }
\label{29sett2015-1}

We introduce $t-$dependent sectors, which serve to define Stokes matrices of $Y(z,t)$  of Corollary \ref{9dic2015-1} in a consistent way w.r.t. matrices of $\mathring{Y}(z)$ of Theorem \ref{31marzo2015-12}. 

 \bde[Sectors $\mathcal{S}_{\nu+k\mu}(t)$]
 \label{20marzo2016-6}
  Let  $\tau_\nu<\widetilde{\tau}<\tau_{\nu+1}$, and $k\in\mathbb{Z}$. Let $t\in  \mathcal{U}_{\epsilon_0}(0)\backslash X(\widetilde{\tau})$. 
 We define   $\mathcal{S}_{\nu+k\mu}(t)$ to be the sector containing the closed sector  $\overline{S}(\widetilde{\tau}-\pi+k\pi,\widetilde{\tau}+k\pi)$, and extending up to the nearest Stokes rays of $\Lambda(t)$  outside $\overline{S}(\widetilde{\tau}-\pi+k\pi,\widetilde{\tau}+k\pi)$. 
\ede 
The definition implies that  
$$ 
\mathcal{S}_{\nu+k\mu}(t)\subset \mathcal{S}_{\nu+k\mu},\quad\quad
\mathcal{S}_{\nu+k\mu}(0)=\mathcal{S}_{\nu+k\mu}.
$$
For simplicity, put $k=0$. Note that   $\mathcal{S}_\nu(t)$ is uniquely defined and   contains   the set of basic Stokes rays  of $\Lambda(t)$ lying  in $S(\widetilde{\tau}-\pi,\widetilde{\tau})$.  We point out the following facts:

$\bullet$ Due to the continuous dependence on $t$ of the directions of Stokes rays for $t\not\in\Delta$, then $\mathcal{S}_\nu(t)$ continuously  deforms as  $t$ varies in a $\widetilde{\tau}$ cell.

$\bullet$ $\mathcal{S}_\nu(t)$ is  ``discontinuous" at $\Delta$, by which we mean that  some Stokes rays disappear at points of $\Delta$. 

$\bullet$   $\mathcal{S}_\nu(t)$ is ``discontinuous" at $X(\widetilde{\tau})$, because one or more Stokes rays cross the admissible ray $R(\widetilde{\tau})$ (this is why $\mathcal{S}_\nu(t)$  has not been defined at $X(\widetilde{\tau})$). 
More precisely, consider a continuous monotone curve $t=t(x)$, $x$ belonging to a real interval, which for one pair $(a,b)$ intersects  $\widetilde{H}_{ab}\backslash \Delta$  at $x=x_*$  (recall that $\widetilde{H}_{ab}$ is define in (\ref{21marzo2016-1})). Hence, the curve passes from one cell to another cell, which are separated by   $\widetilde{H}_{ab}$. A Stokes ray associated with $(u_a(t),u_b(t))$ crosses $R(\widetilde{\tau})$ when $t=t(x_*)$.  Then   $\mathcal{S}_\nu(t(x))$ has a discontinuous jump at $x_*$. 

 The above observations  assure that the following definition is well posed.

\bde[Sector $\mathcal{S}_\nu(K)$]
\label{12maggio2017-1}
 Let $K$ be a compact subset of a $\widetilde{\tau}$-cell. We define 
 $$ 
  \mathcal{S}_\nu(K):=\bigcap_{t\in K}\mathcal{S}_\nu(t)~\subset \mathcal{S}_\nu.
  $$ 
\ede 

By the definitions, 
   $\mathcal{S}_\nu(t)$ and $ 
  \mathcal{S}_\nu (K)$  have the angular width strictly  greater than $\pi$ and they contain  the admissible ray $R(\widetilde{\tau})$ of Definition \ref{14feb2016-5}. Moreover $ \mathcal{S}_\nu (K_1)\supset \mathcal{S}_\nu (K_2)$ for $K_1\subset K_2$, and $\mathcal{S}_\nu (K_1\cup K_2)=\mathcal{S}_\nu (K_1)\cap \mathcal{S}_\nu (K_2)$.  Below in the paper we will consider a simply connected subset $\mathcal{V}$ of a $\widetilde{\tau}$-cell, such that the closure $\overline{\mathcal{V}}$ is also contained in the cell, and take $$K=\overline{\mathcal{V}}.$$

\begin{figure}
\centerline{\includegraphics[width=0.8\textwidth]{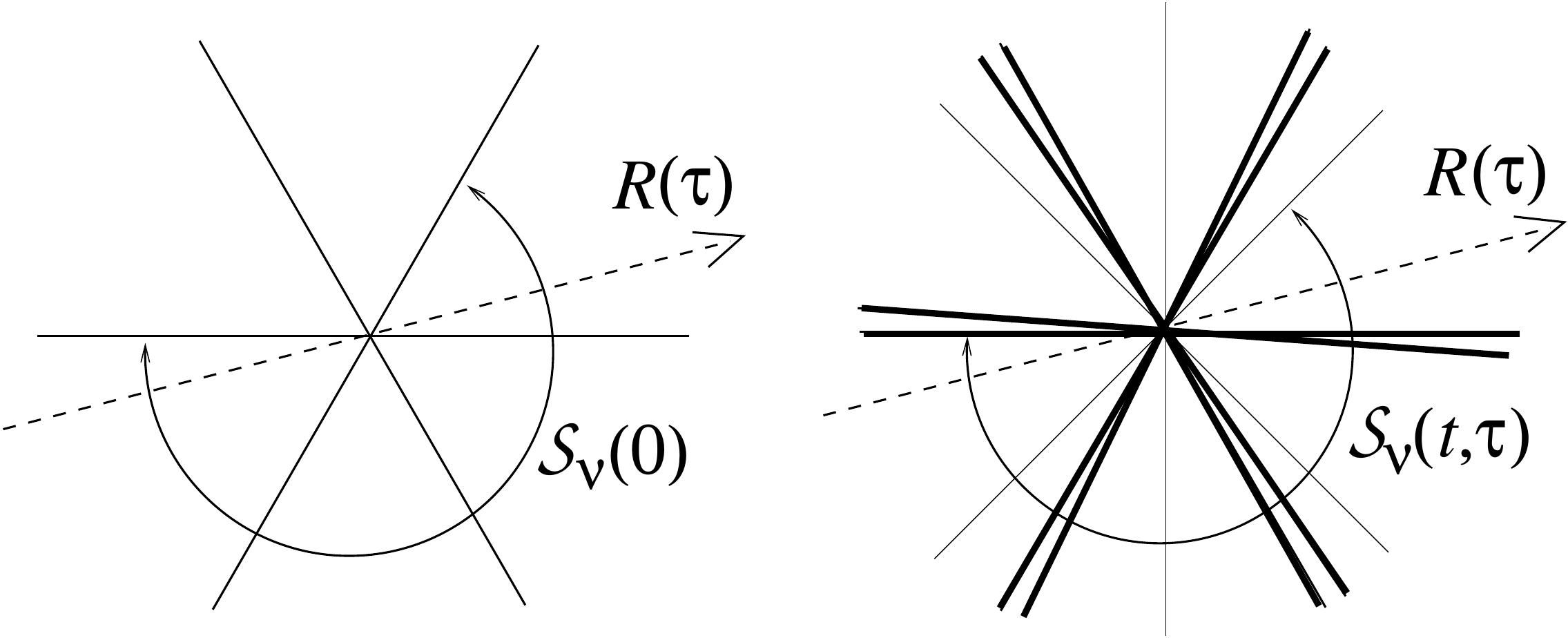}}
\caption{In the left figure $t=0$ and the sector $\mathcal{S}_\nu\equiv \mathcal{S}_\nu(0)$ is represented in a sheet of the universal covering $\mathcal{R}$. The dashed  line represents  $R(\widetilde{\tau})\cup R(\widetilde{\tau}-\pi)$ . The arrow is that of the oriented ray $R(\widetilde{\tau})$. The rays  are the Stokes rays  associated with couples $\lambda_{i},\lambda_{j}$, $1\leq i\neq j\leq s$.  
 In the right figure $t$ slightly differs from $t=0$; the rays in bold are  small deformations  of the rays  appearing in the left figure, associated with couples $u_a(t),u_b(t)$ s.t.  $u_a(0)=\lambda_{i}$, $u_b(0)=\lambda_{j}$ with $i\neq j$.   The rays in finer tone  are the  rays associated with  couples such that $u_a(0)=u_b(0)=\lambda_{i}$. The sector $\mathcal{S}_\nu(t)=\mathcal{S}_\nu(t,\widetilde{\tau})$ is represented.}
\label{sectnut}
\end{figure}

\begin{figure}
\centerline{\includegraphics[width=0.7\textwidth]{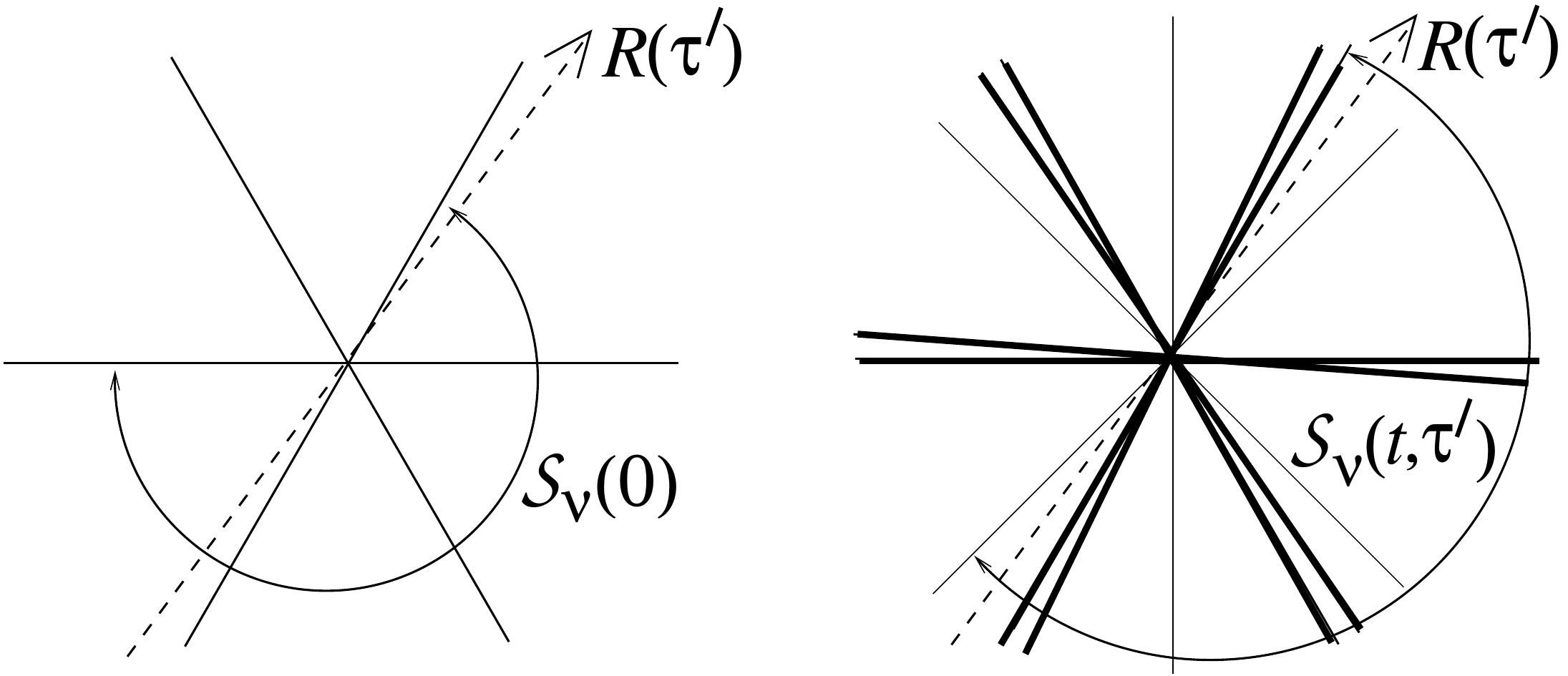}}
\caption{The explanation for this figure is the same as for Figure \ref{sectnut}, but $\widetilde{\tau}^\prime\neq \widetilde{\tau}$. $\mathcal{S}_\nu\equiv \mathcal{S}_\nu(0)$ is the same, but  $\mathcal{S}_\nu(t)=\mathcal{S}_\nu(t,\widetilde{\tau}^\prime)$  differs  from $\mathcal{S}_\nu(t,\widetilde{\tau})$ of figure \ref{sectnut}.}
\label{sectnut1}
\end{figure}

\bre
\label{16gen2016-1}
A  more precise notation could  be used as follows: 
\be
\label{21marzo2016-2}
\mathcal{S}_\nu(t)=\mathcal{S}_\nu(t;\widetilde{\tau})~,
\ee
to  
keep track of $\widetilde{\tau}$, because for given $\nu$ and two different choices of $\widetilde{\tau}\in(\tau_\nu,\tau_{\nu+1})$,  then the resulting  $\mathcal{S}_\nu(t)$'s may be different. Figures \ref{sectnut} and \ref{sectnut1} show two different  $\mathcal{S}_\nu(t)$, according to two choices of  $\widetilde{\tau}$. 
As  a consequence, while in Definition \ref{20marzo2016-6} we could well define  $\mathcal{S}_{\nu+k\mu}(t)\subset \mathcal{S}_{\nu+k\mu}$, for any $k\in\mathbb{Z}$,  we cannot define sectors  $\mathcal{S}_{\nu+1}(t)$,  $\mathcal{S}_{\nu+2}(t)$, ...,  $\mathcal{S}_{\nu+\mu-1}(t)$. 

\ere

\section{Fundamental Solutions $Y_\nu(z,t)$ and Stokes Matrices $\mathbb{S}_\nu(t)$}
\label{14feb2016-9} 
Let
 $\tau_\nu<\widetilde{\tau}<\tau_{\nu+1}$. We show that, if $t_0\not\in\Delta$ belongs to a $\widetilde{\tau}$-cell,   we can extend  the asymptotic behaviour 
 (\ref{8dic2015-1}) of  Corollary  \ref{9dic2015-1} from  $\overline{S}^{(t_0)}(\alpha,\beta)$  to $\mathcal{S}_\nu(t)$. 
 The fundamental matrix  of  Corollary  \ref{9dic2015-1} will then be denoted by $Y_\nu(z,t).$

  \begin{shaded}
\bpr[Solution $Y_\nu(z,t)$ with asymptotics on $\mathcal{S}_\nu(t)$, $t\in \mathcal{U}_{\rho}(t_0)$]
\label{14feb2016-8} 
Let Assumption 1 hold for the system (\ref{9giugno-2}).  Let $t_0$ belong to a $\widetilde{\tau}$-cell.  For any $\nu\in\mathbb{Z}$ there exists $\mathcal{U}_{\rho}(t_0)$ contained in the cell of $t_0$  and a unique fundamental solution  of  the system (\ref{9giugno-2}) as in Corollary \ref{9dic2015-1} of the form 
\be
\label{16maggio2017-4}
Y_\nu(z,t)=G_0(t)\mathcal{G}_\nu(z,t) z^{B_1(t)} e^{\Lambda(t)z},
\ee
holomorphic in $(z,t)\in \{z\in \mathcal{R}~|~|z|\geq N\}\times \mathcal{U}_{\rho}(t_0)$,  
with   asymptotic behaviour (\ref{8dic2015-1}) extended to  $\mathcal{S}_\nu(t)$, $t\in\mathcal{U}_\rho(t_0)$. Namely $\forall ~t\in \mathcal{U}_\rho(t_0)$ the following asymptotic expansion holds:
\be
\label{14feb2016-11}
\mathcal{G}_\nu(z,t)\sim I +\sum_{k=1}^\infty F_k(t) z^{-k}, \quad\quad
z\to\infty, \quad
z\in \mathcal{S}_\nu(t).
\ee
 The asymptotics (\ref{14feb2016-11})  restricted to  $z\in \mathcal{S}_\nu(\mathcal{U}_{\rho}(t_0))$ is uniform in the compact polydisc $\mathcal{U}_{\rho}(t_0)$.  
 \epr
\end{shaded}

\noindent
{\it Note}: Recall that by definition of asymptotics, the last sentence of the above Proposition means that the asymptotics (\ref{14feb2016-11}) is uniform in the compact polydisc $\mathcal{U}_{\rho}(t_0)$ when $z\to \infty$ in any {\it proper closed subsector} of $\mathcal{S}_\nu(\mathcal{U}_{\rho}(t_0))$.

\vskip 0.2 cm 
\noindent
{\it Proof:} 
 In Theorem \ref{ShybuyaOLDbis} choose $\overline{S}^{(t_0)}(\alpha,\beta)= \overline{S}(\widetilde{\tau}-\pi,\widetilde{\tau})$. This contains a set of basic Stokes rays of $\Lambda(t_0)$ and of $\Lambda(t)$ for any $t$ in the cell of $t_0$. Then, Sibuya's  Theorem \ref{ShybuyaOLDbis}  and Corollary \ref{9dic2015-1} apply, with fundamental solution $Y(z,t)$ defined for $t$ in some $\mathcal{U}_{\rho}(t_0)$. It is always possible to restrict $\rho$ so that $\mathcal{U}_{\rho}(t_0)$ is all contained in the cell. 
 
$\bullet$ [Extension to $\mathcal{S}_\nu(t)$]    For $t\in \mathcal{U}_{\rho}(t_0)$, the sector containing
 $S(\widetilde{\tau}-\pi,\widetilde{\tau})$  and extending up to the nearest Stokes rays outside is 
  $\mathcal{S}_\nu(t)$, by definition. Hence there exists a labelling as in Section \ref{11marzo2016-2}, 
  and a $\sigma\in\mathbb{Z}$, such that $\mathcal{S}_\nu(t_0)=\mathcal{S}_\sigma^{(t_0)}$. 
  The Extension Theorem and the Uniqueness Theorem can be applied to $Y(z,t)$  for any fixed $t$, because $S(\widetilde{\tau}-\pi,\widetilde{\tau})$ contains a set of basic Stokes rays. Hence,  for any $t\in\mathcal{U}_{\rho}(t_0)$ the solution 
$Y(z,t)$ is unique with the asymptotic behaviour  (\ref{8dic2015-1}) for $z\to \infty$ in $\mathcal{S}_\nu(t)$.

 $\bullet$ [Uniformity in  $\mathcal{S}_\nu(\mathcal{U}_{\rho}(t_0))$] Clearly,  $\mathcal{S}_\nu(\mathcal{U}_{\rho}(t_0))\supset 
  \overline{S}(\widetilde{\tau}-\pi,\widetilde{\tau})$.  Since $ \mathcal{S}_\nu(\mathcal{U}_{\rho}(t_0))\subset  \mathcal{S}_\nu(t) $ for any $t\in \mathcal{U}_{\rho}(t_0)$, the asymptotics  (\ref{14feb2016-11})  holds also in $ \mathcal{S}_\nu(\mathcal{U}_{\rho}(t_0))$.  Moreover,  the asymptotics is uniform in $ \mathcal{U}_{\rho}(t_0)$  if $z\to\infty$ in $ \overline{S}(\widetilde{\tau}-\pi,\widetilde{\tau})$, by Theorem \ref{ShybuyaOLDbis} and  Corollary \ref{9dic2015-1}. We apply the same proof of the
  Extension Lemma \ref{3giugno2016-1} as follows. Let $\theta_L$ and $\theta_R$ be the directions of the 
  left and right boundary rays of  $ \mathcal{S}_\nu(\mathcal{U}_{\rho}(t_0))$ 
  (i.e. $\overline{\mathcal{S}}_\nu(\mathcal{U}_{\rho}(t_0))=\overline{S}(\theta_R,\theta_L)$).
   Let $\overline{S}_1:=\overline{S}(\phi,\psi)$, for  $\theta_R+\pi <\phi<\psi<\theta_L$, and  $\overline{S}_2:=\overline{S}(\phi^\prime,\psi^\prime)$ for $\theta_R <\phi^\prime<\psi^\prime<\theta_L-\pi$. Let us consider  $\overline{S}_1$. 
  By construction, $\overline{S}_1$ does not contain Stokes rays of $\Lambda(t)$ for any $t\in \mathcal{U}_{\rho}(t_0)$, and so,  by Theorem \ref{ShybuyaOLDbis} now applied with a $\overline{S}^{(t_0)}=  \overline{S}_1$, there exists $~\widetilde{Y}(z,t)\sim Y_F(z,t)$, for $z\to\infty$ in $\overline{S}_1$, uniformly in $|t-t_0|\leq \rho_1$, for suitable $\rho_1>0$. Moreover, $Y(z,t)=\widetilde{Y}(z,t)C(t)$, where $C(t)$ is an invertible holomorphic matrix in $|t-t_0|\leq \min(\rho,\rho_1)$. The matrix entries satisfy 
  $e^{(u_a(t)-u_b(t))z}C_{ab}(t)=\widetilde{\mathcal{G}}(z,t)^{-1}\mathcal{G}(z,t)\sim \delta_{ab}$, $a,b=1,...,n$, for $|t-t_0|\leq \min(\rho,\rho_1)$ and $z\to\infty$, $z\in  \overline{S}(\widetilde{\tau}-\pi,\widetilde{\tau})\cap \overline{S}_1$.  
  Since $\Re((u_a(t)-u_b(t))z)$ does not change sign for $t$ in the cell and 
  $z\in  \overline{S}_1$, then $Y(z,t) \sim Y_F(z,t) $  also for 
  $z\in  \overline{S}(\widetilde{\tau}-\pi,\widetilde{\tau}) \cup   \overline{S}_1$, 
  uniformly in $|t-t_0|\leq \min (\rho,\rho_1)$.  The same arguments  for  $\overline{S}_2$ allow to conclude that  $Y(z,t) \sim Y_F(z,t) $   for 
  $z\in  \overline{S}(\widetilde{\tau}-\pi,\widetilde{\tau}) \cup   \overline{S}_1 \cup   \overline{S}_2$, 
  uniformly in $|t-t_0|\leq \min (\rho,\rho_1,\rho_2)$.  
  Finally, from the proof  given  by Sibuya of Theorem \ref{ShybuyaOLDbis}  (cf. \cite{Sh4}, especially from page 44 on)   it follows that $\rho_1$ and $\rho_2$ are greater or equal to $\rho$.     The proof is concluded. We denote  $Y(z,t)$ with $Y_\nu(z,t)$.  
   $\Box$

\bde[Stokes matrices $\mathbb{S}_{\nu+k\mu}(t)$] 
The {\bf Stokes matrix} $\mathbb{S}_{\nu+k\mu}(t)$, $k\in\mathbb{Z}$, is defined for $t\in  \mathcal{U}_{\rho}(t_0)$ of Proposition \ref{14feb2016-8}  by,
$$ 
Y_{\nu+(k+1)\mu}(z,t)=Y_{\nu+k\mu}(z,t)\mathbb{S}_{\nu+k\mu}(t),\quad\quad
z\in\mathcal{R},
$$
where  the $Y_{\nu+k\mu}(z,t)$ and $Y_{\nu+(k+1)\mu}(z,t)$ are as in Proposition \ref{14feb2016-8}. 
\ede

$\mathbb{S}_{\nu+k\mu}(t)$ is holomorphic in $t\in \mathcal{U}_\rho(t_0)$, because so are $Y_{\nu+(k+1)\mu}(z,t)$ and $Y_{\nu+k\mu}(z,t)$.



\section{Analytic Continuation of $Y_\nu(z,t)$ on a Cell preserving the Asymptotics}
\label{14feb2016-10} 

\begin{shaded}
\bpr[Continuation of $Y_\nu(z,t)$ preserving the asymptotics, along a curve in a cell]
\label{11dic2015-1}
 Let Assumption 1 hold for the system (\ref{9giugno-2}). 
 The fundamental solution $Y_\nu(z,t)$ of Proposition \ref{14feb2016-8}  holomorphic in  $t\in\mathcal{U}_\rho(t_0)$  admits $t$-analytic continuation along any curve contained in the $\widetilde{\tau}$-cell of $t_0$, and maintains its asymptotics (\ref{14feb2016-11}) for $z\to\infty$, $z\in \mathcal{S}_\nu(t)$, for any $t$ belonging to a neighbourhood of the curve. The asymptotics is uniform in a closed tubular neighbourhood $U$ of the curve for $z\to\infty $ in (any proper subsector of)  $\mathcal{S}_\nu(U)$. 
  \epr
\end{shaded}

\vskip 0.3 cm 
\noindent
{\it Proof:}  Let  $Y_\nu(z,t)$, $t\in \mathcal{U}_\rho(t_0)$ be as in Proposition \ref{14feb2016-8}. Join $t_0$  to a  point $t_{final}$, belonging to the  $\widetilde{\tau}$-cell of $t_0$ and not belonging to $ \mathcal{U}_\rho(t_0)$, by a curve whose support is contained in  the $\widetilde{\tau}$-cell.  
Let $t_1\in \partial   \mathcal{U}_\rho(t_0)$ be the intersection point with the curve.
 Theorem \ref{ShybuyaOLDbis} and its Corollary \ref{9dic2015-1} can be applied at $t_1$, with sector $\mathcal{S}_\sigma^{(t_1)}\equiv \mathcal{S}_\nu(t_1)$, by definition. By Proposition \ref{14feb2016-8}, there exists a unique fundamental solution, which we temporarily denote $Y_\nu^{(1)}(z,t)$, with asymptotics (\ref{14feb2016-11}) for $z\to\infty$, $z\in \mathcal{S}_\nu(t)$,  $t\in \mathcal{U}_{\rho_1}(t_1)$.  Here $\rho_1$ is possibly restricted so that  $ \mathcal{U}_{\rho_1}(t_1)$ is  contained in the cell. 
 The asymptotics is uniform in $ \mathcal{U}_{\rho_1}(t_1)$  for $z\to\infty$ in $\mathcal{S}_\nu( \mathcal{U}_{\rho_1}(t_1))$. 
Now,  when $t\in  \mathcal{U}_\rho(t_0)\cap\mathcal{U}_{\rho_1}(t_1)$, both $Y_\nu(z,t)$ and $Y_\nu^{(1)}(z,t)$  are defined, with the same asymptotic behaviour (\ref{14feb2016-11})  for $z\to\infty$, 
$ z\in\mathcal{S}_\nu\left(\mathcal{U}_\rho(t_0)\right)\cap \mathcal{S}_\nu\left(\mathcal{U}_{\rho_1}(t_1)\right)$, uniform   in $t\in \mathcal{U}_\rho(t_0)\cap \mathcal{U}_{\rho_1}(t_1)$. Moreover,  $\mathcal{S}_\nu\left(\mathcal{U}_\rho(t_0)\right)\cap \mathcal{S}_\nu\left(\mathcal{U}_{\rho_1}(t_1)\right)$  has central opening angle strictly greater than $\pi$ because both $\mathcal{U}_\rho(t_0)$ and $\mathcal{U}_{\rho_1}(t_1)$ are contained in the cell. By uniqueness it follows that 
$Y_\nu(z,t)=Y_\nu^{(1)}(z,t)$ for $t\in \mathcal{U}_\rho(t_0)\cap \mathcal{U}_{\rho_1}(t_1)$. This gives the $t$-analytic continuation of $Y_\nu(z,t)$ on   $\mathcal{U}_\rho(t_0)\cup \mathcal{U}_{\rho_1}(t_1)$.  
 The procedure can be
repeated for a sequence of neighbourhoods $\mathcal{U}_{\rho_n}(t_n)$, $n=1,2,3,...$ ($t_n$ is point of intersection of the curve with $\mathcal{U}_{\rho_{n-1}}(t_{n-1})$). Consider $U:=\bigcup_n\mathcal{U}_{\rho_n}(t_n)$. If $t_{final}$ is an internal point of $\in U$, the proof is completed and  $\mathcal{U}_{\rho_n}(t_n)$ is a finite sequence. If not, the point $t_*$ of intersection of $\partial U$ with the curve either  precedes $t_{final}$, or  $t_*=t_{final}\in\partial U$. Since $t_*$ belongs to the cell,   Proposition \ref{14feb2016-8} can be applied. The sector  $\mathcal{S}_{\sigma_*}^{(t_*)}$, $\sigma_*\in\mathbb{Z}$, prescribed by Theorem \ref{ShybuyaOLDbis} and Corollary \ref{9dic2015-1} coincides with $\mathcal{S}_\nu(t_*)$, by definition. Therefore, the analytic continuation is feasible in  a $\mathcal{U}_{\rho^*}(t_*)$, as in the construction above. We can add $\mathcal{U}_{\rho^*}(t_*)$ to $U$. In this way, $t_{final}$ is always reached by a finite sequence, and   $U$ is compact. By construction, the asymptotics is uniform in any compact subset $K\subset U$, including also $K\equiv U$,   for $z\to\infty$, $z\in \mathcal{S}_{\nu}(K)$. 
$\Box$

\begin{shaded}
\bcr {\bf (Analytic continuation of $Y_\nu(z,t)$ preserving the asymptotics on the whole cell --  case of eigenvalues (\ref{4gen2016-2}))}.
\label{20marzo2016-5} 
Let Assumption 1 hold for the system (\ref{9giugno-2}).  
If the eigenvalues of $\Lambda(t)$ are linear in $t$ as in (\ref{4gen2016-2}) then  $Y_\nu(z,t)$ of Proposition \ref{14feb2016-8} is holomorphic on the whole $\widetilde{\tau}$-cell, with  asymptotics (\ref{14feb2016-11}) for $z\to\infty$ in $\mathcal{S}_\nu(t)$, for any $t$ in the cell.  For any  compact subset $K$ of the cell, 
the asymptotics  (\ref{14feb2016-11}) for  $z\to\infty$, $z\in \mathcal{S}_\nu(K)$, is uniform in $t\in K$.
\ecr
\end{shaded}

\noindent
{\it Proof:} 
If the eigenvalues of $\Lambda(t)$ are linear in $t$ as in (\ref{4gen2016-2}), then any $\widetilde{\tau}$-cell is simply connected (see Corollary \ref{20marzo2016-3}). Hence, the continuation of $Y_\nu(z,t)$ is independent of the curve. 
$\Box$ 

\vskip 0.2 cm 
\noindent
$\bullet $ {\bf Notation:} 
If $c$ is the $\widetilde{\tau}$-cell of Corollary \ref{20marzo2016-5}, the following notation will be used 
\be
\label{22marzo2016-1}
Y_\nu(z,t)=Y_\nu(z,t;\widetilde{\tau},c),\quad\quad
t\in c.
\ee

\subsection{Analytic continuation of $Y_\nu(z,t;\widetilde{\tau},c)$ preserving the asymptotics beyond $\partial c$}

Let the eigenvalues of $\Lambda(t)$ be linear in $t$ as in (\ref{4gen2016-2}). 
The analytic continuation of Corollary \ref{20marzo2016-5} and the  asymptotics (\ref{14feb2016-11})  can be extended to values of  $t$ a little bit outside the  cell. 
 This is achieved   by a small variation $\widetilde{\tau}\mapsto\widetilde{\tau}\pm\varepsilon$, for $\varepsilon>0$ sufficiently small. 

 Recall that the Stokes rays in $\mathcal{R}$ associated with  the pair $\bigl(u_a(t),u_b(t)\bigr)$ and $\bigl(u_b(t),u_a(t)\bigr)$, $a\neq b$,  have respectively directions
$$ 
\arg z= \frac{3\pi}{2}-\arg_p(u_a(t)-u_b(t))+2N\pi\quad
\hbox{ and } 
\quad
\arg z= \frac{3\pi}{2}-\arg_p(u_b(t)-u_a(t))+2N\pi,
\quad
\quad
N\in \mathbb{Z}.
$$
Thus, their projections onto $\mathbb{C}$ are the following  opposite rays
\be
\label{15april2016-1}
PR_{ab}(t):=\bigl\{ z\in \mathbb{C}~\bigr|~z=-i\rho(\overline{u}_a(t)-\overline{u}_b(t))\bigr\},
\quad
 \quad 
 PR_{ba}(t):=\bigl\{ z\in \mathbb{C}~\bigr|~z=-i\rho(\overline{u}_b(t)-\overline{u}_a(t))\bigr\}.
 \ee
For $t\not \in W(\widetilde{\tau})$, a ray $PR_{ab}(t)$ lies either in the half plane to the left or to the right of the oriented admissible line $l(\widetilde{\tau})$. 
For $t\not \in W(\widetilde{\tau})$, the {\it finite} set of projected rays   is the union of the two disjoint subsets of (projected) rays to the left and  to the  right   of $l(\widetilde{\tau})$ respectively. Now,  for  $t$ varying  inside a cell $c$, the  projected rays never cross $l(\widetilde{\tau})$.  On the other hand, if $t$ and $t^\prime$ belong to different cells $c$ and $c^\prime$, then the two subsets  of rays to the right and the left of $l(\widetilde{\tau})$ which are associated with  $t$  do not coincide with the two subsets associated with $t^\prime$. These simple considerations imply the following: 

\begin{shaded}
\bpr
\label{19april2016-1}
 A $\widetilde{\tau}$-cell is uniquely characterised  by the subset of projected rays  
which lie to the left of $l(\widetilde{\tau})$. 
\epr
\end{shaded}

 \begin{figure}
\minipage{0.45\textwidth}
\centerline{\includegraphics[width=1\textwidth]{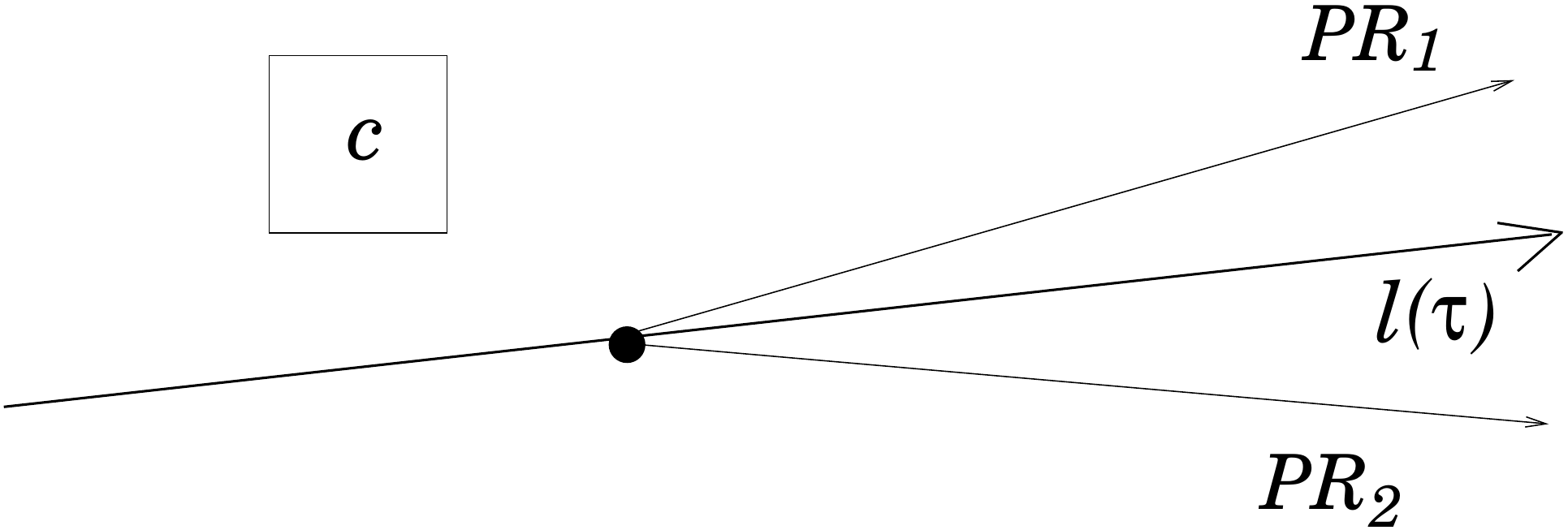}}
\caption{Configuration of rays corresponding to the cell $c$ of figures \ref{cell-april13-8} and \ref{cell-april13-9}.}
\label{cell-april13-1}
\endminipage\hfill
\minipage{0.5\textwidth}
\centerline{\includegraphics[width=0.9\textwidth]{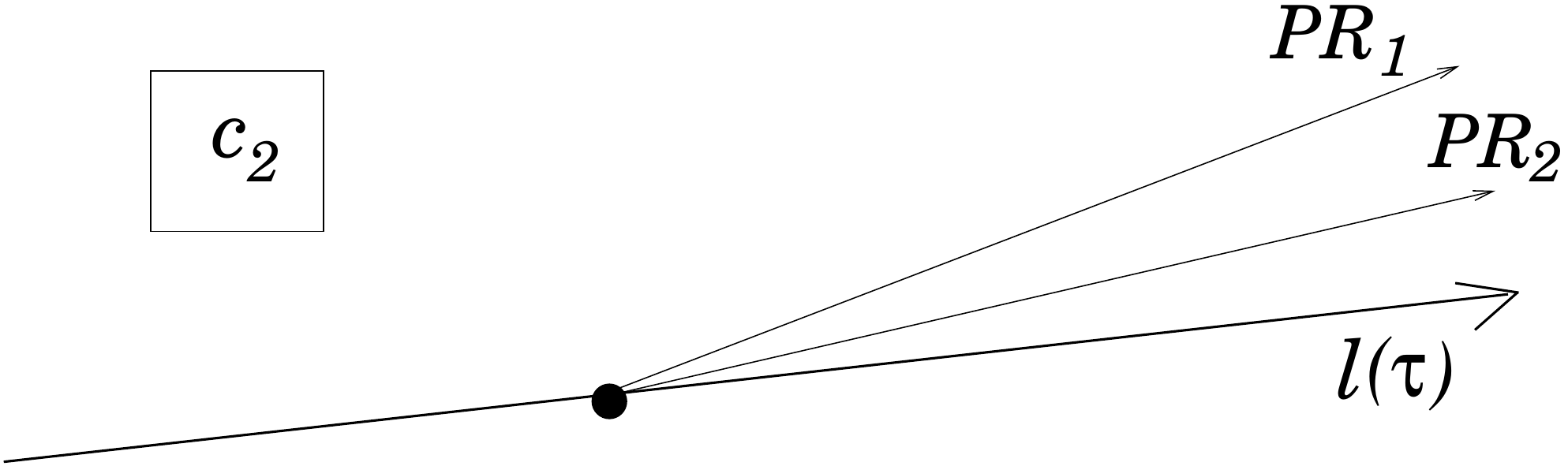}}
\caption{Configuration of rays corresponding to the cell $c_2$ of figures \ref{cell-april13-8} and \ref{cell-april13-9}.}
\label{cell-april13-2}
\endminipage
\end{figure}

 \begin{figure}
\minipage{0.45\textwidth}
\centerline{\includegraphics[width=1\textwidth]{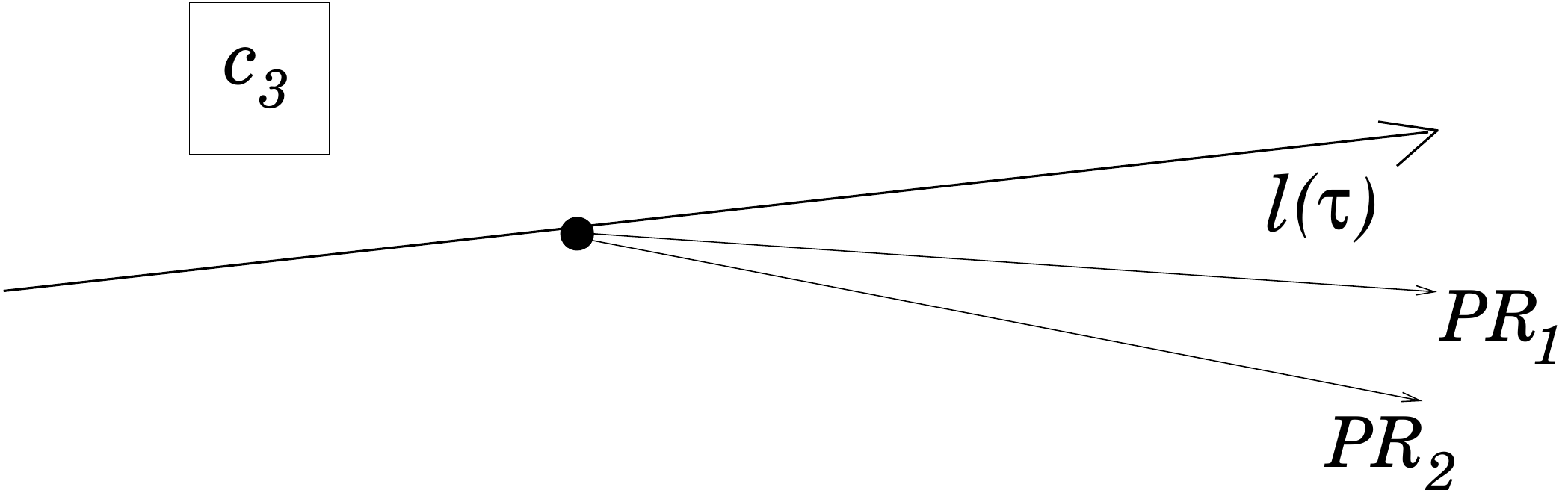}}
\caption{Configuration of rays corresponding to the cell $c_3$ of figures \ref{cell-april13-8} and \ref{cell-april13-9}.}
\label{cell-april13-3}
\endminipage\hfill
\minipage{0.45\textwidth}
\centerline{\includegraphics[width=0.9\textwidth]{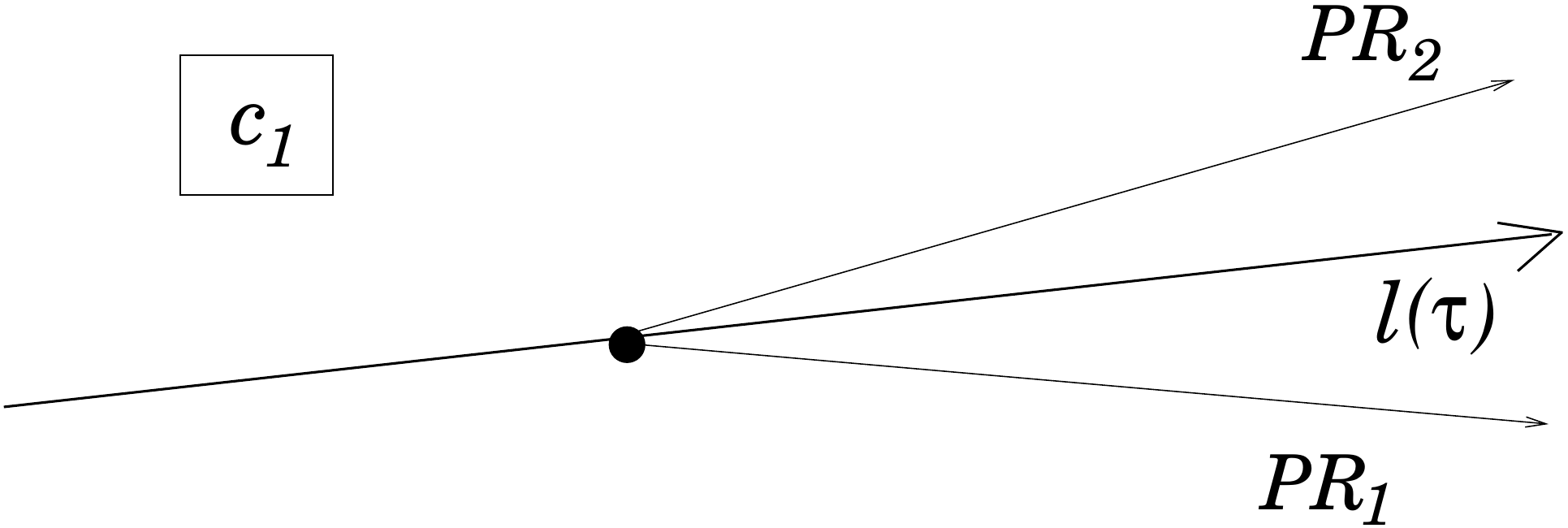}}
\caption{Configuration of rays corresponding to the cell $c_1$ of figures \ref{cell-april13-8} and \ref{cell-april13-9}.}
\label{cell-april13-4}
\endminipage
\end{figure}

\bde
\label{10april2016-1}
A point $t_*\in \widetilde{H}_{ab}\backslash \Delta$    is {\bf simple} if  $t_*\not\in \widetilde{H}_{ab}\cap  \widetilde{H}_{a^\prime b^\prime}$ for any  $(a^\prime,b^\prime)\neq(a,b)$.  
\ede

\begin{figure}
\minipage{0.40\textwidth}
\centerline{\includegraphics[width=1\textwidth]{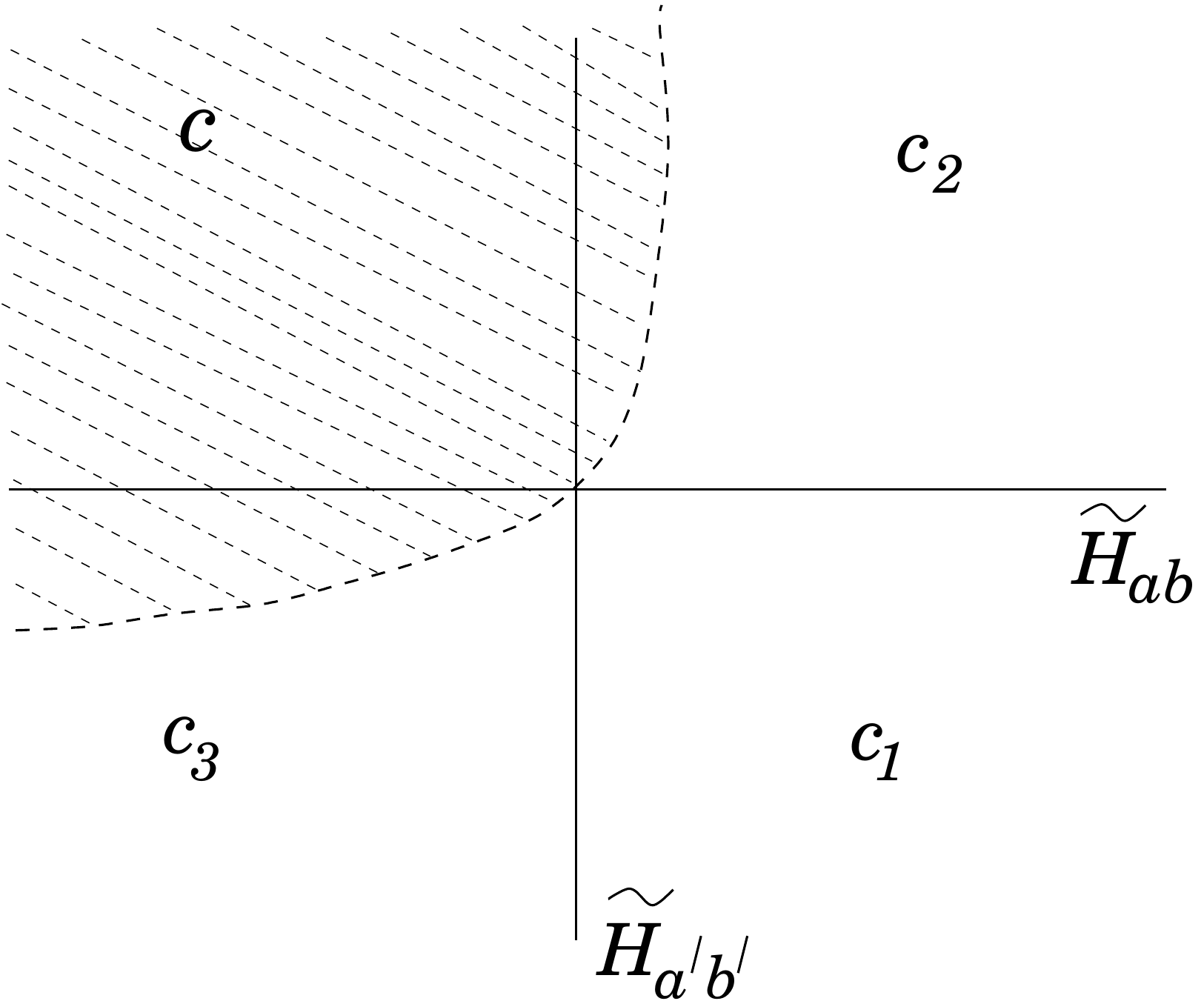}}
\caption{The cells of complex dimension $n$  (real dimension $2n$) are schematically and improperly depicted in real dimension 2. 
Boundaries $\widetilde{H}_{ab}$ and $\widetilde{H}_{a^\prime b^\prime}$ are represented as lines, 
their intersection as a point (understanding that it is not in $\Delta$). The domain of the analytic continuation of $Y_\nu(z,t;\widetilde{\tau},c)$ beyond the boundary of $c$ is the dashed region.
 The analytic continuation does not go beyond $\widetilde{H}_{ab}\cap \widetilde{H}_{a^\prime b^\prime}$, because the transition from 
 figure \ref{cell-april13-1} to figure \ref{cell-april13-4}  is obtained by a  simultaneous crossing  of  $l(\widetilde{\tau})$
   by $PR_1$ and $PR_2$ from {\it opposite sides} of  $l(\widetilde{\tau})$.}
\label{cell-april13-8}
\endminipage\hfill
\minipage{0.45\textwidth}
\centerline{\includegraphics[width=0.9\textwidth]{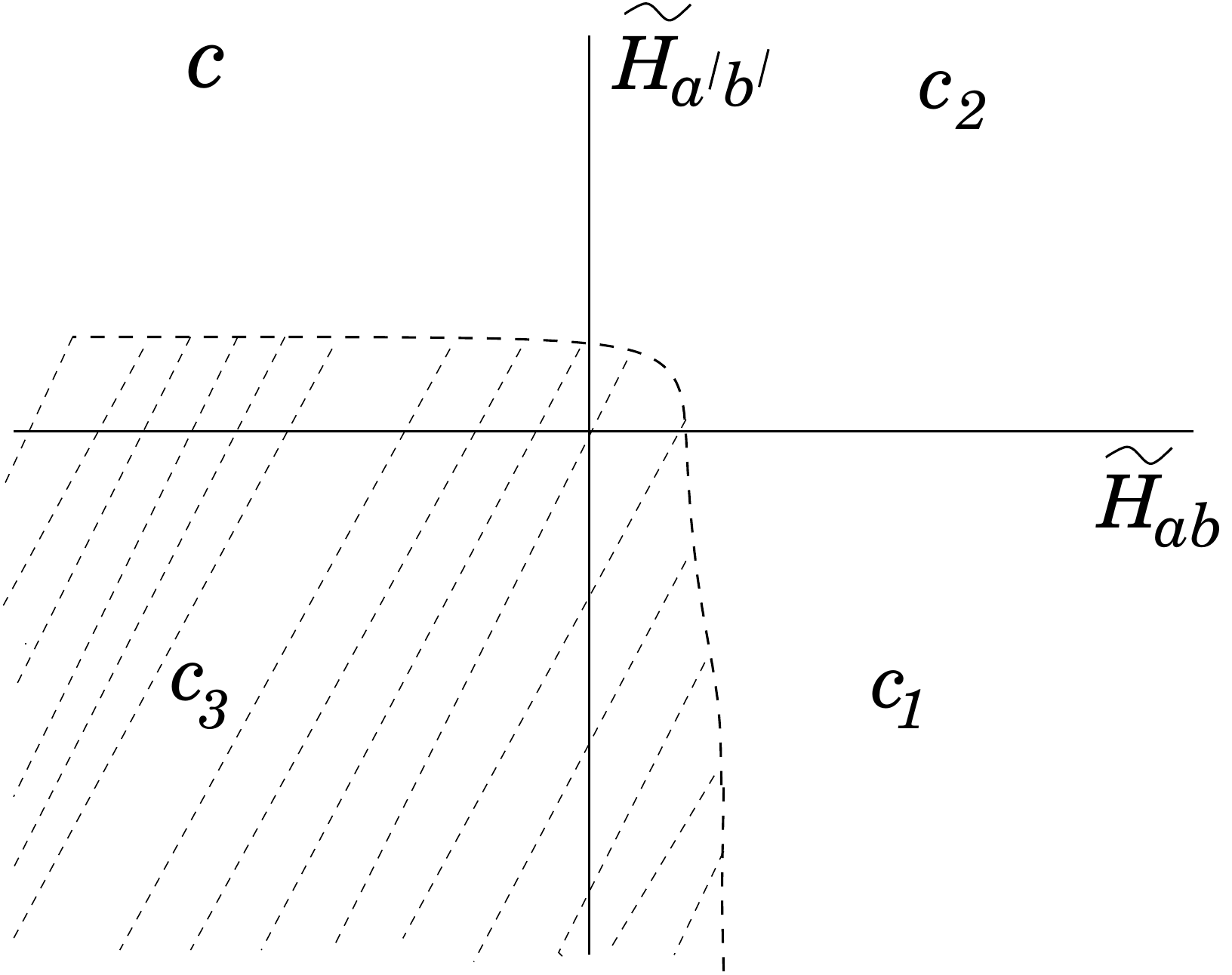}}
\caption{Analytic continuation of $Y_\nu(z,t;\widetilde{\tau},c_3)$ beyond the boundary of $c$. The continuation goes up  to the 3 neighbouring cells. This corresponds to the fact that the three transitions form figure \ref{cell-april13-3} to figures \ref{cell-april13-1} and  \ref{cell-april13-4} occur when   $PR_1$ and $PR_2$ respectively cross $l(\tau)$, while the transition from figure \ref{cell-april13-3} to figure \ref{cell-april13-2} occurs when  $PR_1$ and $PR_2$  simultaneously cross $l(\tau)$, coming from {\it the same side} of  $l(\widetilde{\tau})$ (moving  in anticlockwise sense).}
\label{cell-april13-9}
\endminipage
\end{figure}

If $t$ varies along a curve crossing the boundary  $\partial c$  of a cell $c$ at a simple point belonging to $\widetilde{H}_{ab}\backslash \Delta$, for some $a\neq b$, the ray $PR_{ab}(t)$ crosses either $l_+(\widetilde{\tau})$ or  $l_-(\widetilde{\tau})$, 
while $PR_{ba}(t)$ crosses either $l_-(\widetilde{\tau})$ or  $l_+(\widetilde{\tau})$. Since only  $PR_{ab}(t)$ and  $PR_{ba}(t)$ have crossed $l(\widetilde{\tau})$, then by Proposition \ref{19april2016-1} there is only one neighbouring cell $c^\prime$  sharing the boundary $\widetilde{H}_{ab}$ with $c$.  On the other hand, if the curve crosses $\partial c\backslash \Delta$ at a non simple point, then two or more rays simultaneously cross  $l_+(\widetilde{\tau})$ (and the opposite ones cross $l_-(\widetilde{\tau})$).  
For example, if the crossing occurs at $(\widetilde{H}_{ab}\cap  \widetilde{H}_{a^\prime b^\prime})\backslash \Delta$ then there are three cells, call them  $c_1$, $c_2$, $c_3$, sharing common boundary $(\widetilde{H}_{ab}\cap  \widetilde{H}_{a^\prime b^\prime})\backslash \Delta$ with  $c$.  Looking at the configuration of Stokes rays as in the figures \ref{cell-april13-1},   \ref{cell-april13-2},  \ref{cell-april13-3},  \ref{cell-april13-4}, we conclude that out of the three cells $c_1$, $c_2$, $c_3$,  there is one, say it is   $c_1$, such that the transition from $c$ to $c_1$  occurs with a double crossing of Stokes rays (figure \ref{cell-april13-4}), namely at a non-simple point; while for the remaining  $c_2$ and $c_3$  the transition occurs at  simple points. In figures \ref{cell-april13-1},   \ref{cell-april13-2},  \ref{cell-april13-3},  \ref{cell-april13-4},  $PR_1$ stands for $PR_{ab}(t)$ (or $PR_{ba}(t)$) and $PR_2$ stands for $PR_{a^\prime b^\prime}(t)$ 
 (or $PR_{b^\prime a^\prime}(t)$).  The transition between figure \ref{cell-april13-1} 
 and \ref{cell-april13-4}  is  between $c$ and $ c_1$ of figure \ref{cell-april13-8}, 
 through non simple  points of  $(\widetilde{H}_{ab}\cap  \widetilde{H}_{a^\prime b^\prime})\backslash \Delta$.

\bre
Recall that for any $a\neq b$, $\widetilde{H}_{ab}\cap \Delta\neq\emptyset$. Therefore, when we discuss   analytic continuation, this requires crossing of ``hyperplanes"   $\widetilde{H}_{ab}\backslash \Delta$. 
\ere

\begin{shaded}
\bpr[Continuation slightly beyond the cell, preserving asymptotics]
\label{22marzo2016-2}
 Let the assumptions of Corollary \ref{20marzo2016-5} hold.  Let $c$ and $c^\prime$ be $\widetilde{\tau}$-cells such that $\partial c\cap \partial c^\prime\neq \emptyset$. If $\partial c\cap \partial c^\prime$ does not coincide with the multiple intersection of two or more $\widetilde{H}_{ab}$'s,   then  $Y_\nu(z,t;\widetilde{\tau},c)$ has analytic continuation,  with  asymptotics (\ref{14feb2016-11}) in $\mathcal{S}_\nu(t)$, for $t$  slightly beyond $\partial c\backslash \Delta$ into $c^\prime$. The asymptotics  for $z\to\infty $ in $\mathcal{S}_\nu(K)$ is uniform  in any compact subset $K$ of the extended cell. 
 Equivalently,   $Y_\nu(z,t;\widetilde{\tau},c)$ 
  can be analytically continued along any curve crossing  $\partial c\backslash \Delta$  at a simple point and ending slightly beyond  $\partial c\backslash \Delta$  in the neighbouring cell $c^\prime$.   
  \epr
\end{shaded} 
 
 \begin{figure}
\centerline{\includegraphics[width=0.45\textwidth]{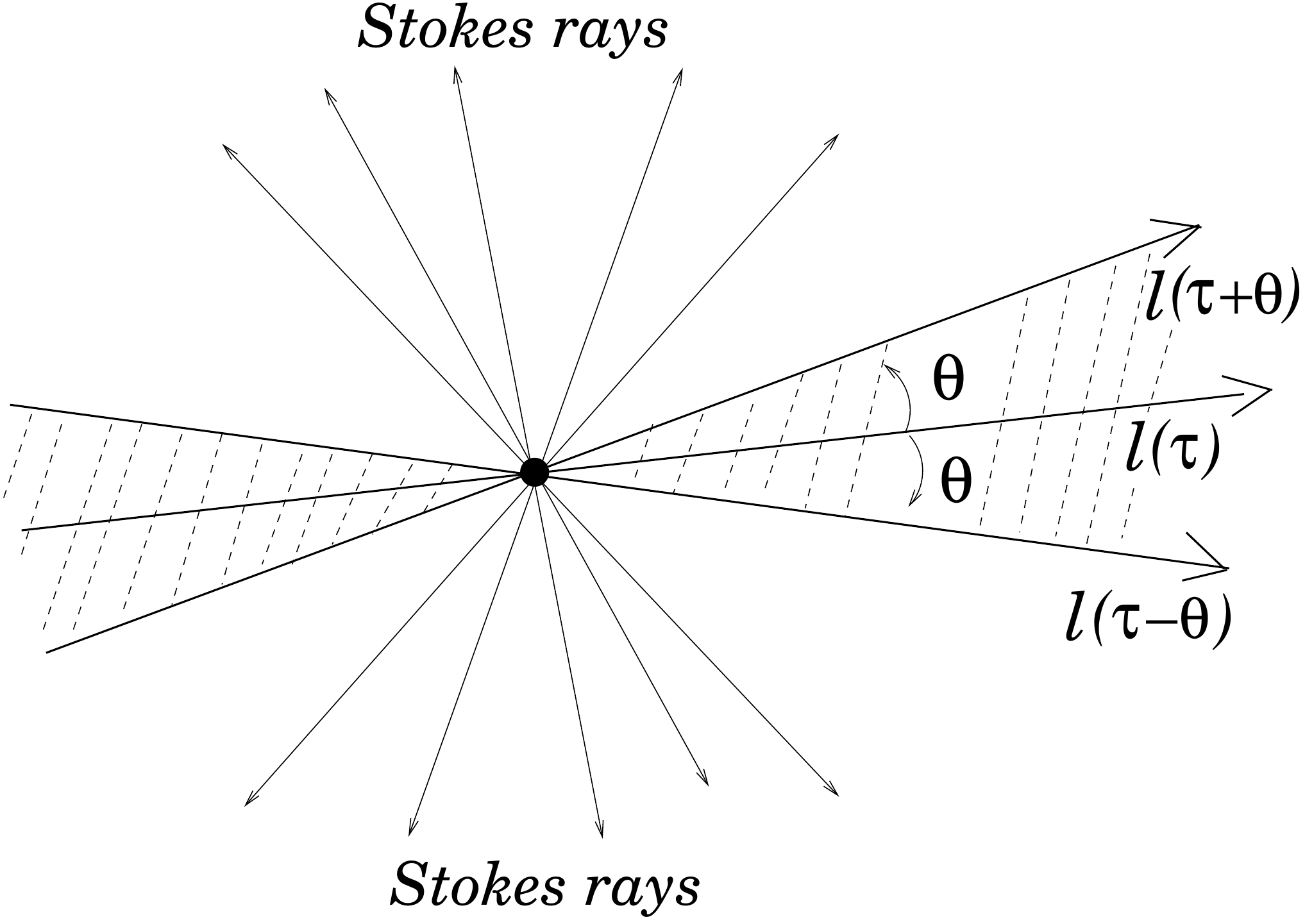}}
\caption{The two closed sectors of amplitude $2\vartheta$, not containing Stokes rays when $t\in\overline{U}$.}
\label{cell-april13-5}
\end{figure}
 
\begin{figure}
\minipage{0.45\textwidth}
\centerline{\includegraphics[width=1\textwidth]{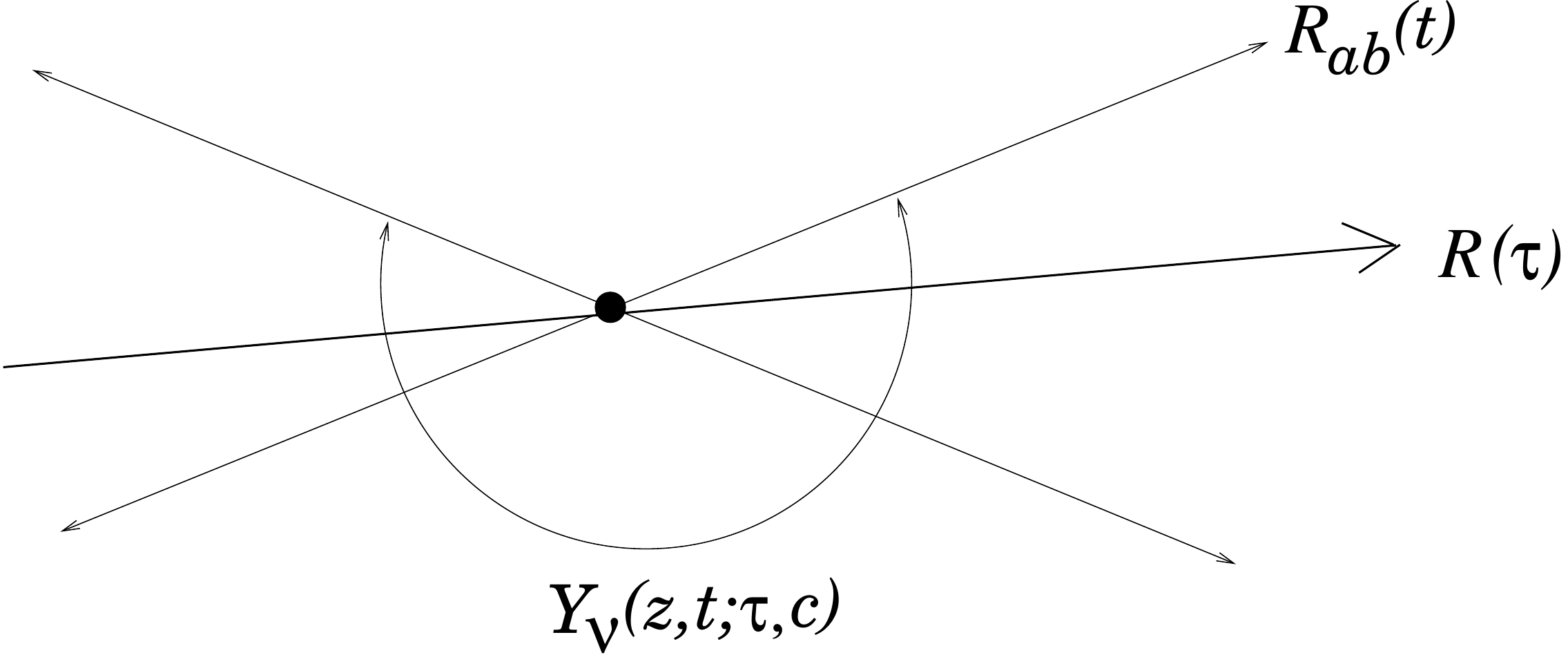}}
\caption{$Y_\nu(z,t;\widetilde{\tau},c)$ for $t\in c$. The sector  where  $Y_\nu(z,t;\widetilde{\tau},c)$ has the canonical asymptotic behaviour is represented.}
\label{cell-april13-6}
\endminipage\hfill
\minipage{0.5\textwidth}
\centerline{\includegraphics[width=0.9\textwidth]{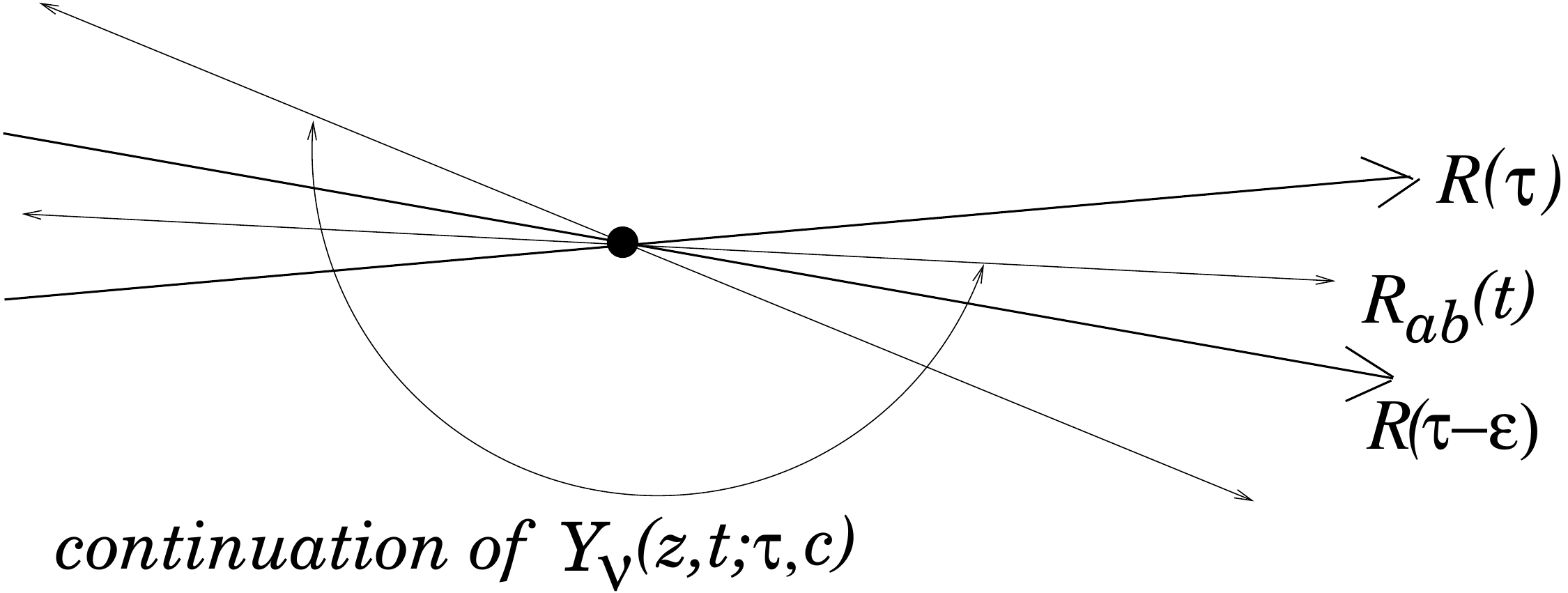}}
\caption{Analytic continuation of $Y_\nu(z,t;\widetilde{\tau},c)$ for $t$ in the neighbouring cell $c^\prime$ just after the crossing of $\partial c\backslash \Delta$, namely just after $R_{ab}(t)$ has crossed $R(\widetilde{\tau})$. The sector  where  $Y_\nu(z,t;\widetilde{\tau},c)$ has the canonical asymptotic behaviour is represented.}
\label{cell-april13-7}
\endminipage
\end{figure}

\noindent
 {\it Proof:} 
 Let $U$ be an open connected subset of the $\widetilde{\tau}$-cell $c$, such that $\overline{U}$ is contained in $c$.   There exists a small  $\vartheta=\vartheta(U)>0$ such that for any $t\in  \overline{U}$ the projected Stokes rays of $\Lambda(t)$ lie outside the two closed sectors containing  $l(\widetilde{\tau})$ and bounded by $l(\widetilde{\tau}+\theta)$ and $l(\widetilde{\tau}-\theta)$, as in figure  \ref{cell-april13-5}.  
 Let $\varepsilon\in[0,\vartheta]$.  All lines $l(\widetilde{\tau}\pm\varepsilon)$ are admissible for the Stokes rays, when $t\in \overline{U}$.  Consider the subset of projected Stokes rays to the left of $l(\widetilde{\tau})$. It uniquely identifies (cf. Proposition  \ref{19april2016-1})  the $(\widetilde{\tau}+\varepsilon)$-cell and the $(\widetilde{\tau}-\varepsilon)$-cell  obtained by deforming the boundaries of $c$  when $\widetilde{\tau}\mapsto \widetilde{\tau}+ \varepsilon$ and  $\widetilde{\tau}\mapsto \widetilde{\tau}- \varepsilon$ respectively (recall that  $L_{ab}$ in (\ref{13april2016-1}) depends on $\widetilde{\eta}=3\pi/2-\widetilde{\tau}$). Call these cells $c_{\varepsilon}$ and $c_{-\varepsilon}$.  By construction 
 \beas
 &&
   \overline{U}\subset c\cap c_{\pm \varepsilon},\quad\quad
   \varepsilon\in[0,\vartheta],
   \\
&&
Y_\nu(z,t;\widetilde{\tau},c)=Y_\nu(z,t;\widetilde{\tau}\pm \varepsilon,c_{\pm \varepsilon}) ,\quad
t\in U.   
\eeas 
   The last equality follows from the definition of $Y_\nu$, its uniqueness  and Corollary \ref{20marzo2016-5}.  Indeed, the analytic continuation explained in  the proof of Proposition \ref{11dic2015-1} can be repeated   for the function $Y_\nu(z,t;\widetilde{\tau}\pm \varepsilon,c_{\pm \varepsilon})$ initially defined in a neighbourhood of $t_0$  contained in $\overline{U}$, but {\it with  cell partition  determined by   $\widetilde{\tau} \pm \varepsilon$}. Moreover, by uniqueness of solutions with asymptotics, it follows that $Y_\nu(z,t;\widetilde{\tau},c)=Y_\nu(z,t;\widetilde{\tau}\pm \varepsilon,c_{\pm \varepsilon})$ for $t \in U$. Therefore, $Y_\nu(z,t;\widetilde{\tau},c)$ has analytic continuation to $c_{\pm \varepsilon}$.  Now, 
$$
c_{\pm \varepsilon} \cap \{\hbox{ union of  cells sharing boundary with $c$ }\} \neq \emptyset.
$$
   Then, the analytic continuation of $Y_\nu(z,t;\widetilde{\tau},c)$ obtained above is actually defined in a $t$-domain bigger than $c$.  We characterise this domain, showing that it intersect any cell  
      $c^\prime$ which is a neighbour of $c$, and  such that   $\partial c \cap \partial c^\prime$ does not coincide with the multiple
      intersection of two or more hyperplanes.      Thus, we need to show that $c_{\pm \varepsilon} \cap c^\prime\neq \emptyset$.  Notice that   $\partial c \cap \partial c^\prime = \widetilde{H}_{ab}$ for suitable $a,b$. 
      Then, suppose without loss of generality that $PR_{ab}(t)$ crosses $l_+(\widetilde{\tau})$ clockwise when $t$ crosses
       $\widetilde{H}_{ab}\backslash \Delta$ moving along a curve from $c$ to $c^\prime$.  An example of this crossing
        is the transition from  figure \ref{cell-april13-1} to figure \ref{cell-april13-3},   with the identification $c^\prime=c_3$ of Figure \ref{cell-april13-8},   
        and $PR_1=PR_{ab}$.  Then, for the small deformation $\widetilde{\tau}\mapsto \widetilde{\tau}-\varepsilon$ the
         above discussion applies.  Namely, $c_{-\epsilon}\cap c^\prime \neq \emptyset$.  See figures \ref{cell-april13-6} and \ref{cell-april13-7}.
  $\Box$

 \vskip 0.3 cm 

 If  $\partial c\cap \partial c^\prime = 
 \widetilde{H}_{ab}\cap  \widetilde{H}_{a^\prime b^\prime}$ for some $(a^\prime,b^\prime)\neq(a,b)$, there is multiple crossing 
 of $l(\widetilde{\tau})$.  The proof does not work if the crossing corresponds to a transition such as that  from figure \ref{cell-april13-1} to figure 
 \ref{cell-april13-4}, with the identification $c^\prime=c_1$. Since  $PR_1$ and $PR_2$ cross simultaneously $l_+(\widetilde{\tau})$ from opposite 
  sides, any 
 deformation $\widetilde{\tau}\mapsto \widetilde{\tau}\pm \varepsilon$ produces a cell 
 $c_{\pm \varepsilon}$ which does not intersect $c_1$.  In other words, the deformation prevents points of   $c_{\pm \varepsilon}$  from 
 getting close to  $ \widetilde{H}_{ab}\cap  \widetilde{H}_{a^\prime b^\prime}$. 
The schematic figure \ref{cell-april13-8} shows  the 4 cells corresponding to the figures  from \ref{cell-april13-1} to \ref{cell-april13-4}. 
It is shown that   $Y_\nu(z,t;\widetilde{\tau},c)$ can be continued  slightly inside $c_2$ and $c_3$, but not inside $c^\prime=c_1$. 
It is worth noticing that both  $Y_\nu(z,t;\widetilde{\tau},c_2)$ and  $Y_\nu(z,t;\widetilde{\tau},c_3)$ can be continued beyond $ \widetilde{H}_{ab}\cap  \widetilde{H}_{a^\prime b^\prime}$. See figure  \ref{cell-april13-9} for  $Y_\nu(z,t;\widetilde{\tau},c_3)$.

  \bre 
  If the eigenvalues are linear in $t$ as in (\ref{4gen2016-2}), the results of this section assures that the fundamental solutions $Y_{\nu+k\mu}(z,t;\widetilde{\tau},c)$'s are holomorphic in a $\widetilde{\tau}$-cell $c$ and a little beyond, that they maintain the  asymptotic behaviour, and then the corresponding  Stokes matrices $\mathbb{S}_{\nu+k\mu}(t)$'s  are defined and holomorphic in the whole  $\widetilde{\tau}$-cell $c$ and a little bit beyond. 
  \ere

\section{Fundamental Solutions $Y_\nu(z,t)$ and Stokes Matrices $\mathbb{S}_\nu(t)$ holomorphic at $\Delta$}
\label{2dic2016-2}

 If the fundamental solutions $Y_{\nu+k\mu}(z,t;\widetilde{\tau},c)$'s of (\ref{9giugno-2}) (with Assumption 1)  have analytic continuation to the whole $\mathcal{U}_{\epsilon_0}(0)$,  in this section   we give sufficient conditions such that the  continuations are $c$-indendent solutions  $Y_{\nu+k\mu}(z,t)$'s, which maintain the asymptotic behaviour in large sectors $\widehat{\mathcal S}_\nu$ defined below, 
so that the Stokes matrices $\mathbb{S}_{\nu+k\mu}(t)$ are well defined in the whole   $\mathcal{U}_{\epsilon_0}(0)$. Moreover, we  show that $Y_{\nu+k\mu}(z,0)\equiv \mathring{Y}_{\nu+k\mu}(z)$ and  $\mathbb{S}_{\nu+k\mu}(0)\equiv \mathring{\mathbb{S}}_{\nu+k\mu}$, where  $\mathring{Y}_{\nu+k\mu}(z)$, $
\mathring{\mathbb{S}}_{\nu+k\mu}$ have been defined in Section \ref{16marzo2016-4} for the system at fixed $t=0$.

\subsection{Restriction of $\epsilon_0$}
\label{16gen2016-5}
So far, $\epsilon_0$ has been taken so small that $\Lambda_i(t)$ and $\Lambda_j(t)$, $1\leq i\neq j\leq s$,  have no common eigenvalues for $t\in \mathcal{U}_{\epsilon_0}(0)$.  If $\Lambda=\Lambda(0)$ has at least two distinct eigenvalues, we consider a further restriction of $\epsilon_0$.    Let $\widetilde{\eta}=3\pi/2-\widetilde{\tau} $ be the admissible direction associated with the direction  $\widetilde{\tau}$ of the admissible ray $R(\widetilde{\tau})$. 
 Let $\delta_0$ be a small positive number  such that 
\be
\label{22marzo2016-11}
\delta_0<\min_{1\leq i\neq j\leq s} \delta_{ij},
\ee
where $\delta_{ij}$ is 1/2 of the distance  between two parallel lines of angular direction $\widetilde{\eta}$ in the $\lambda$-plane, one passing through $\lambda_i$ and one through $\lambda_j$; namely
\be
\label{22marzo2016-10}
\delta_{ij}:=\frac{1}{2}\min \left\{
 \left|
\lambda_i-\lambda_j +\rho e^{i\tilde{\eta}}\right|, ~\rho\in\mathbb{R}
\right\},\quad
\quad
i\neq j =1,2,...,s.
\ee
 Clearly,  $\delta_0$ depends on the choice of $\widetilde{\eta}$ (see also Remark \ref{26marzo2016-2}). 
Let  $\overline{B}(\lambda_i;\delta_0)$ be the closed ball in $\mathbb{C}$ with center $\lambda_i$ and radius $\delta_0$. Then, we choose $\epsilon_0$ so small that the eigenvalues $u_1(t)$, ..., $u_n(t)$ for $t\in \mathcal{U}_{\epsilon_0}(0)$ satisfy  
$$ 
(u_1(t),...,u_n(t))\in~  \overline{B}(\lambda_1;\delta_0)^{\times p_1}\times ~ \cdots ~\times \overline{B}(\lambda_s;\delta_0)^{\times p_s}.
$$
 As $t $ varies  in $\mathcal{U}_{\epsilon_0}(0)$ above, the Stokes rays continuously move, but the directions of the  rays   associated with   a $u_a\in\overline{B}(\lambda_i;\delta_0)$ and  a $u_b\in \overline{B}(\lambda_j;\delta_0)$, $i\neq j$,  never  cross the values $\widetilde{\eta}$ and $\widetilde{\eta}-\pi$ (mod $2\pi$), so that the projected rays $PR_{ab}(t)$ and $PR_{ba}(t)$ never cross the admissible line $l(\widetilde{\tau})$. It follows that 
 
 \begin{shaded}
 {\it 
  \noindent
  the cell decomposition only depends on the Stokes rays associated with  couples ($u_a(t)$, $u_b(t)$) such that $u_a(0)=u_b(0)=\lambda_i$, $i=1,...,s$.
  }
\end{shaded}
   
For eigenvalues {\it linear} in $t$ as in (\ref{4gen2016-2}),  we can take $\epsilon_0=\delta_0$ and
\be
\label{25marzo2016-1}
\mathcal{U}_{\epsilon_0}(0) \equiv \overline{B}(0;\delta_0)^{\times p_1}\times ~ \cdots ~\times \overline{B}(0;\delta_0)^{\times p_s},\quad\quad
\epsilon_0=\delta_0.
\ee

  \bre 
\label{11april2016-2}
If $t$ moves  from one $\widetilde{\tau}$-cell to another,  the only  Stokes rays  which may cross admissible rays  $R(\widetilde{\tau}+k\pi)$, $k\in \mathbb{Z}$, are those associated with  pairs $u_a(t),u_b(t)$ with  $u_a(0)=u_b(0)=\lambda_i$, $i=1,...,s$.    Therefore, the boundaries of the cells are only the  $\widetilde{H}_{ab}$'s  such that $u_a(0)=u_b(0)$. In this case, $L_{ab}({\bf x}_0,{\bf y}_0)=0 $,   so that 
$$
H_{ab}^\prime:=\Bigl\{ (\Re t, \Im t)\in\mathbb{R}^{2n}~\Bigl|~L_{ab}(\Re t , \Im t)=0\Bigr\}. 
$$
Remark \ref{11april2016-3} follows from the above observations.   
   \ere


\subsection{The Sectors $\widehat{\mathcal{S}}_\nu(t)$ and $\widehat{\mathcal{S}}_\nu$}
\label{2dic2016-1}

\begin{figure}
\centerline{\includegraphics[width=0.7\textwidth]{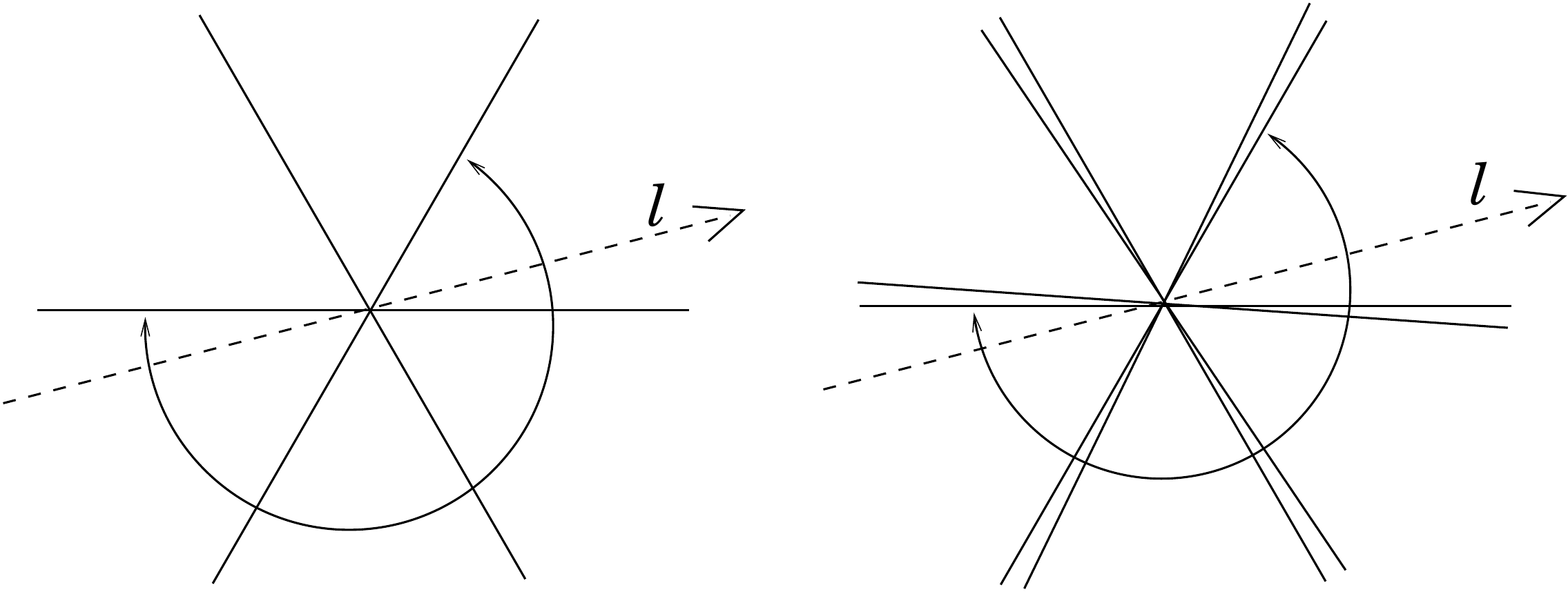}}
\caption{In the left figure  $t=0$ and the sector $\mathcal{S}_\nu$ is represented. The explanation is  as for the left part of  Figure \ref{sectnut}. In the right figure, $t\neq 0$.  Represented are only the rays associated with couples $u_a(t),u_b(t)$ with  $u_a(0)=\lambda_{i}$, $u_b(0)= \lambda_{j}$, for $i\neq j$, together with the sector $\widehat{\mathcal{S}}_\nu(t)$.}
\label{sectnut2}
\end{figure}

Let $\Lambda(t)$ be  of the form (\ref{27dic2015-1}) with eigenvalues (\ref{4gen2016-2}). Let  $\epsilon_0= \delta_0$ be as in subsection \ref{16gen2016-5}. 
We define  a  subset $\mathfrak{R}(t)$ of the set of Stokes rays  of $\Lambda(t)$ as follows:    $\mathfrak{R}(t)$ contains only those Stokes  rays $\{z\in\mathcal{R}~|~\Re(z(u_a(t)-u_b(t)))=0\}$ which are  associated with   pairs  $u_a(t)$, $u_b(t)$ satisfying the condition $u_a(0)\neq u_b(0)$ (namely,  $u_a(0)=\lambda_i$, $u_b(0)=\lambda_j$,   $i\neq j$;  see (\ref{29gen2016-1})-(\ref{29gen2016-3})). The reader may visualise the rays in $\mathfrak{R}(t)$  as being originated  by the splitting of Stokes rays of $\Lambda(0)$.
  See figure \ref{sectnut2}.

$\mathfrak{R}(t)$ has the following important property: if $t$ varies in $\mathcal{U}_{\epsilon_0}(0)$, the rays in $\mathfrak{R}(t)$ continuously move, but   since $\epsilon_0 = \delta_0$, they never cross  any admissible ray $R(\widetilde{\tau}+k\pi)$, $k\in \mathbb{Z}$.

\bde[Sectors  $\widehat{\mathcal{S}}_{\nu+k\mu}(t)$] 
\label{13giugno2017-4}
We define $\widehat{\mathcal{S}}_{\nu+k\mu}(t)$ to be the unique sector containing $S(\widetilde{\tau}-\pi+k\pi,\widetilde{\tau}+k\pi)$  and extending up to the nearest Stokes rays in $\mathfrak{R}(t)$, $t\in\mathcal{U}_{\epsilon_0}(0)$. 
\ede

Any $\widehat{\mathcal{S}}_{\nu+k\mu}(t)$ contains a set of basic Stokes rays of  $\mathfrak{R}$. Moreover, 
$$R(\widetilde{\tau})\subset \widehat{\mathcal{S}}_\nu(t)\cap
\widehat{\mathcal{S}}_{\nu+\mu}(t)\subset S(\tau_\nu,\tau_{\nu+1}),
$$
 and  
$$ 
\mathcal{S}_\nu(t)\subset \widehat{\mathcal{S}}_\nu(t),\quad\quad
\widehat{\mathcal{S}}_\nu(0)\equiv \mathcal{S}_\nu.
$$
In case $\Lambda(0)=\lambda_1 I$, then  $\widehat{\mathcal{S}}_\nu(t)$ is unbounded, namely it coincides with $\mathcal{R}$.

\bde[Sectors $\widehat{\mathcal{S}}_\nu(K)$] 
For any compact $K\subset \mathcal{U}_{\epsilon_0}(0)$ we define 
$$
\widehat{\mathcal{S}}_\nu(K):=\bigcap_{t\in K} \widehat{\mathcal{S}}_\nu(t). 
$$ 
\ede
If  $K_1\subset K_2$, then $
\widehat{\mathcal{S}}_\nu(K_2)\subset \widehat{\mathcal{S}}_\nu(K_1)$. For any $K_1$, $K_2$, we have  $\widehat{\mathcal{S}}_\nu(K_1\cup K_2)=\widehat{\mathcal{S}}_\nu(K_1)\cap \widehat{\mathcal{S}}_\nu(K_2)$.

\bde[Sectors $\widehat{\mathcal{S}}_\nu$] 
\label{26marzo2016-1}
If  $K=\mathcal{U}_{\epsilon_0}(0)$, we define 
$$ 
\widehat{\mathcal{S}}_\nu:=  \widehat{\mathcal{S}}_\nu(\mathcal{U}_{\epsilon_0}(0)).
$$
\ede
Since  $\epsilon_0=\delta_0$,   $\widehat{\mathcal{S}}_\nu$  has angular opening greater than $\pi$ and
\beas
&&
\widehat{\mathcal{S}}_\nu \subset  \widehat{\mathcal{S}}_\nu(0)\equiv \mathcal{S}_\nu,  
\\
&&R(\widetilde{\tau})\subset \widehat{\mathcal{S}}_\nu\cap
\widehat{\mathcal{S}}_{\nu+\mu}\subset S(\tau_\nu,\tau_{\nu+1}).
\eeas

\bre
\label{26marzo2016-2}
Notice that $\widetilde{\tau}\in(\tau_\nu,\tau_{\nu+1})$ determines $\delta_0$ through (\ref{22marzo2016-10})  and (\ref{22marzo2016-11}). Let $\widetilde{\tau}^\prime\in (\tau_\nu,\tau_{\nu+1})$ and let $\delta_0^\prime$ be obtained through  (\ref{22marzo2016-10})  and (\ref{22marzo2016-11}). Let $\epsilon_0=\min\{\delta_0,\delta_0^\prime\}$. We temporarily denote by $\widehat{\mathcal{S}}_\nu[\widetilde{\tau}]$ the sector $\widehat{\mathcal{S}}_\nu$ of Definition \ref{26marzo2016-1} obtained starting from $\widetilde{\tau}$.  Then  for the above $\epsilon_0$ we have 
$$ 
\widehat{\mathcal{S}}_\nu[\widetilde{\tau}]=\widehat{\mathcal{S}}_\nu[\widetilde{\tau}^\prime].
$$
\ere

\subsection{Fundamental group of $\mathcal{U}_{\epsilon_0}(0)\backslash \Delta$ and generators}
\label{8agosto2016-2}

Let the eigenvalues of $\Lambda(t)$ be linear in $t$ as in (\ref{4gen2016-2}), 
   $\tau_\nu<\widetilde{\tau}<\tau_{\nu+1}$ and $\widetilde{\eta}=3\pi/2-\widetilde{\tau}$.

 The fundamental group 
 $\pi_1(\mathcal{U}_{\epsilon_0}(0)\backslash \Delta, t_{\it base})$ is generated by    loops $\gamma_{ab}$, $1\leq a\neq b \leq n$,    which are homotopy classes of  simple curves encircling  the component   $\{t\in \mathcal{U}_{\epsilon_0}(0)~|~u_a(t)=u_b(t)\}$ of $\Delta$. The choice of the base point is  free, because   $\mathcal{U}_{\epsilon_0}(0)\backslash \Delta$ is path-wise connected, since $\Delta$ is a  braid arrangement in $\mathcal{U}_{\epsilon_0}(0)$ and the hyperplanes are {\it complex}.

 For  $\epsilon_0= \delta_0$ of Section \ref{16gen2016-5},   Stokes rays in $\mathfrak{R}(t)$ never cross the admissible rays $R(\widetilde{\tau}+k\pi)$, $k\in\mathbb{Z}$,  when $t$ goes along any  loop in $\mathcal{U}_{\epsilon_0}(0)$ (see Remark \ref{11april2016-2}).  Therefore, as far as the analytic continuation of $Y_\nu(z,t)$ is concerned, it is enough to consider $u_a(t)$ and $u_b(t)$ 
 coming from  the unfolding of an eigenvalue $\lambda_i$ of $\Lambda(0)$ (see the beginning of Section \ref{13giugno2017-2}), namely 
 \be
 \label{11april2016-1}
u_a(t)=\lambda_i+t_a,\quad\quad
u_b(t)=\lambda_i+t_b.
\ee 
    If we represent $t_a$ and $t_b$ in the same complex plane, so that $ t_a-t_b$ is a complex number,  a representative of $\gamma_{ab}$, which we also denote $\gamma_{ab}$ with abuse of notation, is represented by the following  loop around $t_a-t_b=0$, 
\be
\label{10maggio2016-1}
 t_a-t_b\longmapsto (t_a-t_b)e^{2\pi i}.
\ee
$|t_a-t_b| $ will be taken small.  
The Stokes rays associated with $u_a(t)$ and $u_b(t)$  have directions
\be
\label{3april2016-1}
\frac{3\pi}{2}-\arg(t_{a}-t_b)\quad
\hbox{\rm mod}(2\pi),\quad\quad
\frac{3\pi}{2}-\arg(t_b-t_a)\quad
\hbox{\rm mod}(2\pi).
\ee
The projection of these rays  onto $\mathbb{C}$  are the two  opposite rays $PR_{ab}(t)$ and $PR_{ba}(t)$, as in (\ref{15april2016-1}) . Along the loop (\ref{10maggio2016-1}), each of  these rays rotate clockwise and 
crosses the line $l(\widetilde{\tau})$ twice (recall Definition \ref{14feb2016-5}), once passing over  the positive half line and once over the negative half line, returning to the initial position at the end of the loop.  Hence, the support of  $\gamma_{ab}$ is contained in at least two cells, but generally in more than two, as follows.

\begin{figure}
\centerline{\includegraphics[width=0.6\textwidth]{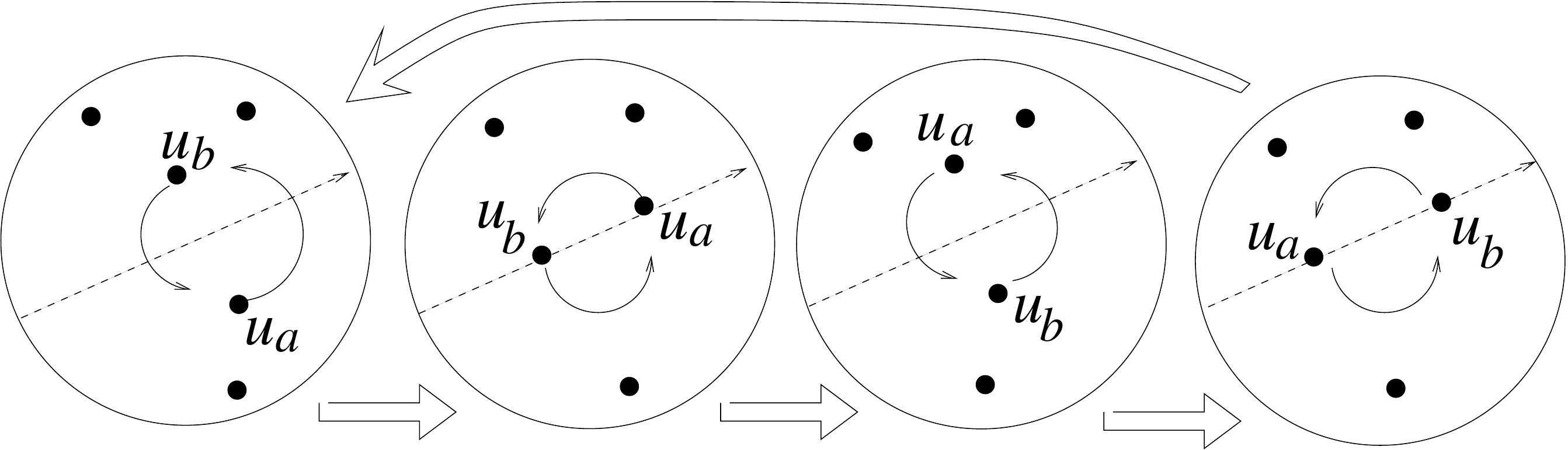}}
\caption{Loop $\gamma_{ab}$ represented in $\overline{B}(\lambda_i;\epsilon_0)$. 
The dashed oriented  line is the direction $\widetilde{\eta}$. Along the loop, $u_a$ and $u_b$ get aligned with $\widetilde{\eta}$ twice, in the second and fourth figures. The second figure corresponds to the passage from 
one initial cell $c$ to a neighbouring cell $c^\prime$ (while $PR_{ab}$ crosses clockwise a half line 
of $l(\widetilde{\tau})$)  and the fourth figure to the return to $c$ (while $PR_{ab}$ crosses clockwise  the opposite half line of 
$l(\widetilde{\tau})$). Other dots represent other eigenvalues $u_\gamma(t)$ in $\overline{B}(\lambda_i;\epsilon_0)$.}
\label{105may}
\end{figure}

$\bullet$  There exists a representative  contained in only two cells if only $PR_{ab}(t)$ and its opposite $PR_{ba}(t)$ cross  $l(\widetilde{\tau})$,  each   twice. For example,  in figure \ref{105may}  the ball $\overline{B}(\lambda_i;\epsilon_0)$ is represented with  the loop (\ref{10maggio2016-1}). The dots represent other  points $u_\gamma(t) \in \overline{B}(\lambda_i;\epsilon_0)$, $\gamma\neq a,b$.  $PR_{ab}(t)$ and $PR_{ba}(t)$ cross   $l(\widetilde{\tau})$ when $u_a(t)$ and $u_b(t)$ are aligned with the admissible direction $\widetilde{\eta}$. Along the loop, no other $u_\gamma$  aligns  with  $u_a(t)$ and $u_b(t)$.

\begin{figure}
\centerline{\includegraphics[width=0.5\textwidth]{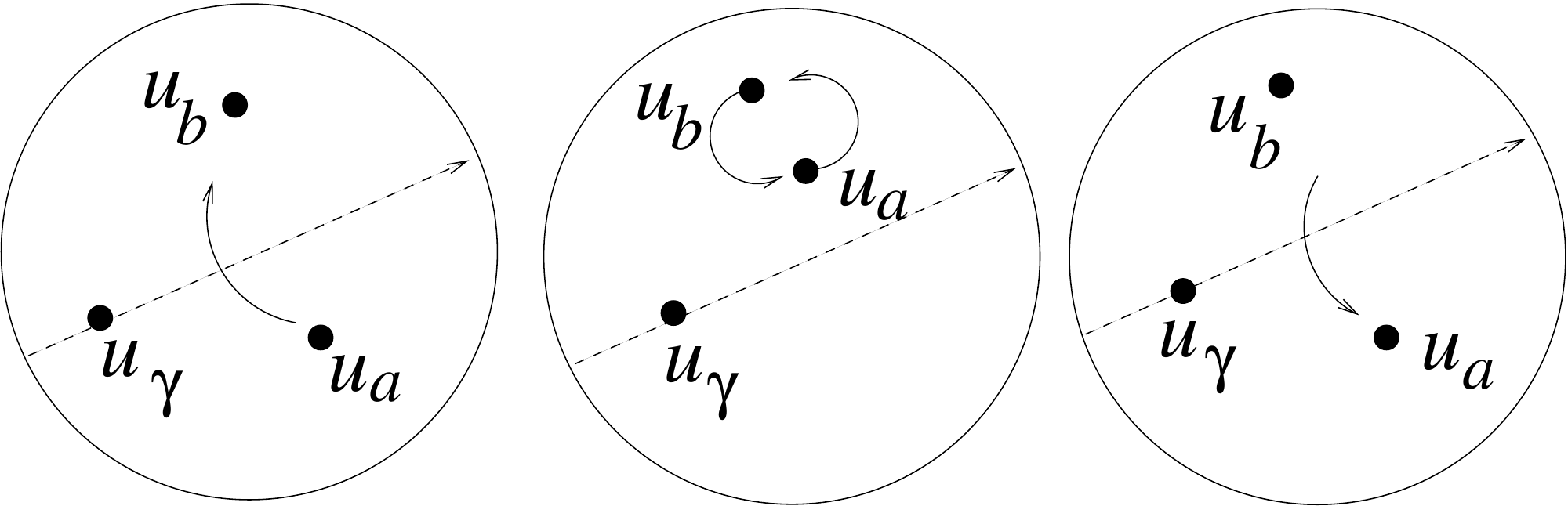}}
\caption{Loop $\gamma_{ab}$ represented in $\overline{B}(\lambda_i;\epsilon_0)$. The dashed oriented line is the direction $\widetilde{\eta}$. 
In the first figure, $u_a$ moves close to $u_b$. Along the way it gets aligned with $u_\gamma$. At this alignment, $PR_{a\gamma}$ crosses clockwise  a half line of $l(\widetilde{\tau})$ and $t$ passes from the initial cell $c$ to a cell $c^\prime$. The second figure is figure \ref{105may}. Here $t$ passes from $c^\prime $ to another cell $c^{\prime\prime}$ and then back to $c^\prime$. In the third figure, $u_a$ moves to the initial position. Along the way it gets aligned with $u_\gamma$, $PR_{a\gamma}$ crosses anti-clockwise the same  half line of $l(\widetilde{\tau})$  and  $t$ returns to the cell $c$. In this example, $\gamma_{ab}$ has support contained in three cells.}
\label{105may1}
\end{figure}

$\bullet$  In general, other (projected) rays cross $l(\widetilde{\tau})$  along any possible representative of $\gamma_{ab}$. For example,   the representative of  (\ref{10maggio2016-1}) in figure   \ref{105may1} is contained in three cells. Indeed, also $PR_{a\gamma}(t)$ and $PR_{\gamma a}(t)$ cross $l(\widetilde{\tau})$ when $u_a$ and $u_\gamma$ get aligned with $\widetilde{\eta}$. Alignment  corresponds to the passage from one cell to another.

\subsection{Holomorphic conditions such that $Y_\nu(z,t)\to\mathring{Y}_\nu(z)$ and $\mathbb{S}_\nu(t)\to \mathring{\mathbb{S}}_\nu$ for $t\to 0$, in case of linear eigenvalues  (\ref{4gen2016-2})} 
\label{1april2016-1}

The following theorem is one of the central results of the paper, and it will be used to prove Theorem \ref{16dicembre2016-1}. 

\begin{shaded}
\bth
\label{pincopallino} 
Consider the system (\ref{9giugno-2}) (i.e. system (\ref{22novembre2016-3}) of the Introduction). Let  Assumption 1 hold, so that   (\ref{9giugno-2}) is holomorphically equivalent to the system (\ref{16marzo2016-1})  (i.e. to  (\ref{17luglio2016-1}) of the Introduction).  
Let $\Lambda(t)$ be  of the form (\ref{27dic2015-1}), with eigenvalues (\ref{4gen2016-2}) and $\epsilon_0 =\delta_0$ as in  subsection  \ref{16gen2016-5}. Let $\widetilde{\tau}$ be the direction of an admissible ray $R(\widetilde{\tau})$, satisfying $\tau_\nu<\widetilde{\tau}<\tau_{\nu+1}$.  Suppose  that: 

\vskip 0.2 cm
1)   For every integer  $j\geq 1$, the  $F_j(t)$'s are holomorphic in $\mathcal{U}_{\epsilon_0}(0)$ (so necessary and sufficient conditions of Proposition \ref{8dic2015-3} hold);

\vskip 0.2 cm

 2)  For any $\widetilde{\tau}$-cell $c$ of $ \mathcal{U}_{\epsilon_0}(0)$ and any $k\in\mathbb{Z}$,  the fundamental solution $Y_{\nu+k\mu}(z,t;\widetilde{\tau},c)$ has analytic continuation as a single-valued holomorphic function on the whole $ \mathcal{U}_{\epsilon_0}(0)$. Denote the analytic continuation with the same symbol  $Y_{\nu+k\mu}(z,t;\widetilde{\tau},c)$, $t\in\mathcal{U}_{\epsilon_0}(0)$. 

\vskip 0.2 cm 
\noindent
Then:

$\bullet$ For any  $\widetilde{\tau}$-cells $c$ and $c^\prime$, 
$$
Y_{\nu+k\mu}(z,t;~\widetilde{\tau},c)=Y_{\nu+k\mu}(z,t;~\widetilde{\tau},c^\prime),
\quad
\quad
t\in\mathcal{U}_{\epsilon_0}(0).
$$ 

Therefore, we can simply write $Y_{\nu+k\mu}(z,t;\widetilde{\tau})$. 
\vskip 0.2 cm 
$\bullet$ Let 
$\mathcal{G}_{\nu+k\mu}(z,t;\widetilde{\tau}):= G_0(t)^{-1}Y_{\nu+k\mu}(z,t;\widetilde{\tau})z^{-B_1(t)} e^{-\Lambda(t)z}$. 
 For any $\epsilon_1<\epsilon_0$  the following  asymptotic expansion holds: 
\be
\label{26marzo2016-3}
\mathcal{G}_{\nu+k\mu}(z,t;\widetilde{\tau})\sim I+\sum_{k=0}^\infty F_k(t)z^{-k},
\quad
\quad
z\to\infty,
\quad
z\in\widehat{\mathcal{S}}_{\nu+k\mu},
\quad
t\in \mathcal{U}_{\epsilon_1}(0).
\ee
The asymptotic expansion is uniform  in $t$ in $ \mathcal{U}_{\epsilon_1}(0)$ and uniform in  $z$ in any closed subsector of  $\widehat{\mathcal{S}}_{\nu+k\mu}$.  

\vskip 0.2 cm 
$\bullet$ For any $t\in \mathcal{U}_{\epsilon_1}(0)$, the diagonal blocks of any Stokes matrix $\mathbb{S}_{\nu+k\mu}(t)$ are the identity matrices $I_{p_1}$, $I_{p_2}$, ...,$I_{p_s}$. Namely 
$$
(\mathbb{S}_{\nu+k\mu})_{ab}(t)=(\mathbb{S}_{\nu+k\mu})_{ba}(t)=0 
\quad
\hbox{ whenever } u_a(0)=u_b(0).
$$
\eth
\end{shaded}
\bre[Continuation of Remark \ref{26marzo2016-2}]
Since $Y_{\nu+k\mu}(z,t;~\widetilde{\tau},c)=Y_{\nu+k\mu}(z,t;~\widetilde{\tau},c^\prime)\equiv Y_{\nu+k\mu}(z,t;~\widetilde{\tau}) $, only the choice of $\widetilde{\tau}$ is relevant.
If $\widetilde{\tau}$ and $\widetilde{\tau}^\prime$ are as in Remark \ref{26marzo2016-2} , then 
$$
Y_{\nu+k\mu}(z,t;\widetilde{\tau})=Y_{\nu+k\mu}(z,t;\widetilde{\tau}^\prime),
$$
 because the rays in $\mathfrak{R}(t)$, $t\in\mathcal{U}_{\epsilon_0}(0)$, neither cross the admissible rays $R(\widetilde{\tau}+m\pi)$ nor the rays $R(\widetilde{\tau}^\prime+m\pi)$, $m\in\mathbb{Z}$. In other words, $Y_{\nu+k\mu}(z,t;\widetilde{\tau})$ depends on $\widetilde{\tau}$ only through $\epsilon_0$. Hence, we can  restore the notation  
$$Y_{\nu+k\mu}(z,t),
\quad\quad
 t\in \mathcal{U}_{\epsilon_0}(0).
$$ 

\ere

\begin{shaded}
\bcr 
\label{6marzo2016-2} 
Let the assumptions of Theorem \ref{pincopallino}  hold. 
Let $\mathring{Y}_{\nu+k\mu}(z)$, $k\in \mathbb{Z}$,  denote the unique fundamental solution (\ref{8agosto2016-1}) of the form 
(\ref{6dic2015-8}), namely
$$
\mathring{Y}_{\nu+k\mu}(z)=\mathcal{G}_{\nu+k\mu}(z)z^{B_1(0)}e^{\Lambda z},
$$  
  with the asymptotics (\ref{6dic2015-7})
  $$
  \mathcal{G}_{\nu+k\mu}(z)\sim I+\sum_{j=1}^\infty \mathring{F}_j z^{-j},\quad\quad
  z\to\infty,
  \quad
  z\in \mathcal{S}_{\nu+k\mu},
  $$
      corresponding to the  particular choice $\mathring{F}_j=F_j(0)$, $j\geq 1$.  Then, 
\beas
&&\mathcal{G}_{\nu+k\mu}(z,0)=\mathcal{G}_{\nu+k\mu}(z),
\\
&&Y_{\nu+k\mu}(z,0)=\mathring{Y}_{\nu+k\mu}(z).
\eeas
\ecr 
\end{shaded}

\noindent
{\it Proof:}  Observe that $Y_{\nu+k\mu}(z,0)$ is defined at $t=0$. Now, $\widehat{\mathcal{S}}_\nu\subset  \mathcal{S}_\nu$ and both sectors have central opening angle greater than $\pi$. Hence, the solution with given asymptotics in $\widehat{\mathcal{S}}_\nu$ is unique, namely $\mathcal{G}_\nu(z)=\mathcal{G}_{\nu}(z,0)$. $\Box$

 \begin{shaded}
 \bcr
\label{19feb2016-1}
Let the assumptions  of Theorem \ref{pincopallino} hold. 
   Let $\mathbb{S}_\nu(t)$, $\mathbb{S}_{\nu+\mu}(t)$ be a complete set of Stokes matrices  associated with fundamental solutions $
Y_{\nu}(z,t)$, $ Y_{\nu+\mu}(z,t)$,  $Y_{\nu+2\mu }(z,t)$, with canonical asymptotics, for $t$ in a $\widetilde{\tau}$-cell of  $\mathcal{U}_{\epsilon_0}(0)$, in sectors $\mathcal{S}_\nu(t)$, $\mathcal{S}_{\nu+\mu}(t)$ and $\mathcal{S}_{\nu+2\mu}(t)$ respectively, which by Theorem \ref{pincopallino} extend to $\widehat{\mathcal{S}}_\nu$, $\widehat{\mathcal{S}}_{\nu+\mu}$ and $\widehat{\mathcal{S}}_{\nu+2\mu}$ respectively for $t\in\mathcal{U}_{\epsilon_1}(0)$.    Then there exist  
 $$ 
\lim_{t\to 0} \mathbb{S}_\nu(t)= \mathring{\mathbb{S}}_\nu,
\quad\quad
\lim_{t\to 0}\mathbb{S}_{\nu+\mu}(t)= \mathring{\mathbb{S}}_{\nu+\mu},
 $$ 
 where $\mathring{\mathbb{S}}_\nu$, $\mathring{\mathbb{S}}_{\nu+\mu}$ is a complete set of Stokes matrices for the system at $t=0$, referred to three fundamental solutions $
\mathring{Y}_{\nu+k\mu}(z)$, $k=0,1,2$, of Corollary \ref{6marzo2016-2} 
  having asymptotics  in sectors $\mathcal{S}_{\nu+k\mu}$, with  $\mathring{F}_j=F_j(0)$, $j\geq 1$.
  
 \ecr
 \end{shaded}
 
 \noindent{\it Proof:} The analyticity  of $Y_{\nu+k\mu}(z,t)$ in assumption 2) of Theorem \ref{pincopallino} implies that the  Stokes matrices are holomorphic in $\mathcal{U}_{\epsilon_0}(0)$. 
Hence, for $k=1,2$, there exists  
$$
\mathbb{S}_{\nu+k\mu}(0)=\lim_{
t\to 0
}
\left(Y_{\nu+(k+1)\mu}(z,t)^{-1}Y_{\nu+k\mu}(z,t)\right)
=
 \mathring{Y}_{\nu+(k+1)\mu}(z)^{-1}\mathring{Y}_{\nu+k\mu}(z)=\mathring{\mathbb{S}}_{\nu+k\mu}.
$$
$\Box$

\subsubsection{Proof of Theorem \ref{pincopallino}}

\ble
\label{28marzo2016-2} 
Let Assumption 1 hold for the system (\ref{9giugno-2}).  Let  the eigenvalues of $\Lambda(t)$ be linear in $t$ as in (\ref{4gen2016-2}). 
Suppose that $Y_\nu(z,t;\widetilde{\tau},c)$ has $t$-analytic continuation on $\mathcal{U}_{\epsilon_0}(0)\backslash\Delta$, with $\epsilon_0 =\delta_0$ as in  subsection  \ref{16gen2016-5}. Temporarily call $Y^{cont}_\nu(z,t;\widetilde{\tau},c)$ the continuation. Also suppose that
$$  
Y^{cont}_\nu(z,t;\widetilde{\tau},c)\Bigr|_{t\in c^\prime}=Y_\nu(z,t;\widetilde{\tau},c^\prime).
$$
Then:

$a)$ Any $Y_\nu(z,t;\widetilde{\tau},c^\prime)$ has analytic continuation on $\mathcal{U}_{\epsilon_0}(0)\backslash\Delta$, coinciding with $Y^{cont}_\nu(z,t;\widetilde{\tau},c)$. 
Due to the independence of $c$, we denote this continuation by
$$
Y_\nu(z,t;\widetilde{\tau}).
$$

b)  $\mathcal{G}_{\nu}(z,t;\widetilde{\tau}):= G_0(t)^{-1}Y_{\nu}(z,t;\widetilde{\tau})z^{-B_1(t)z} e^{-\Lambda(t)z}$ has asymptotic expansion
$$
\mathcal{G}_{\nu}(z,t;\widetilde{\tau})\sim I+\sum_{k=0}^\infty F_k(t)z^{-k},
\quad
\quad
z\to\infty,
\quad
z\in\widehat{\mathcal{S}}_{\nu}(t),
\quad
t\in \mathcal{U}_{\epsilon_0}(0)\backslash \Delta.
$$
The asymptotics for $z\to\infty$ in $\widehat{\mathcal{S}}_\nu(K)$ is uniform on any compact subset $K\in \mathcal{U}_{\epsilon_0}(0)\backslash\Delta$.
\ele

\vskip 0.2 cm 
\noindent
{\it Proof of Lemma \ref{28marzo2016-2}:}  a) is obvious. We prove  b), dividing the proof into two parts. 

\vskip 0.2 cm 
{\bf Part 1} (in steps).  Chosen an arbitrary cell $c$ (all cells are equivalent, by a)) and  any $\breve{t}\in c$, we prove that the sector where  $Y_\nu(z,\breve{t};\widetilde{\tau})$ has canonical asymptotics can be extended from $\mathcal{S}_\nu(\breve{t})$  to  $\widehat{\mathcal{S}}_\nu(\breve{t})$. For clarity in the discussion below, let us still write $Y_\nu(z,\breve{t};\widetilde{\tau},c)$. 

\begin{figure}
\centerline{\includegraphics[width=0.5\textwidth]{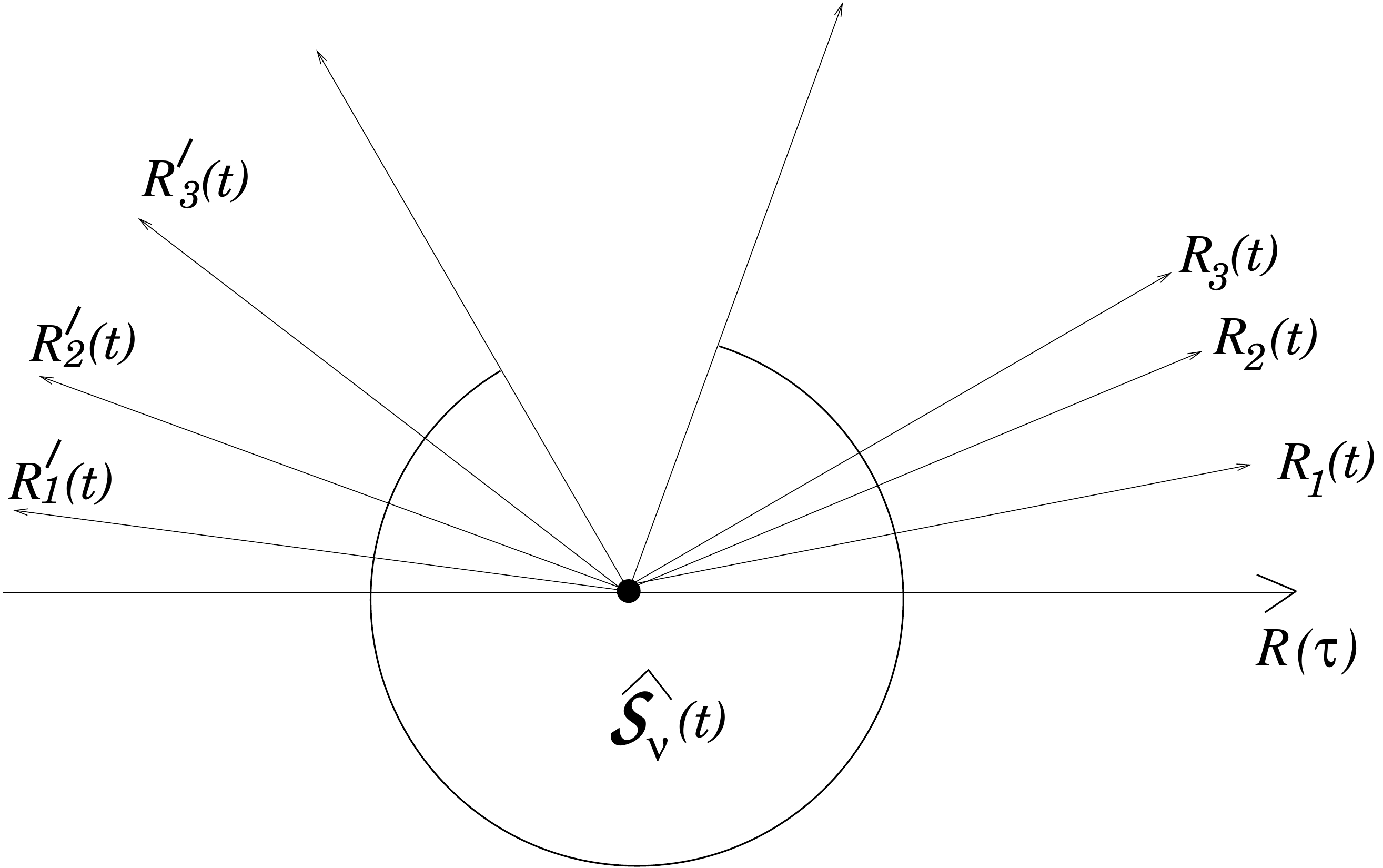}}
\caption{Rays in $\widehat{\mathcal{S}}_\nu(t)$ which are going to be erased in the proof.}
\label{pria15}
\end{figure}

\vskip 0.2 cm 
Step 1.  At $\breve{t}$, consider the Stokes rays in $\widehat{\mathcal{S}}_\nu(\breve{t})\backslash S(\widetilde{\tau}-\pi,\widetilde{\tau})$ associated with the unfolding of the
 $\lambda_i$'s.  Those with direction greater than $\widetilde{\tau}$ will be labelled in anticlockwise sense as  $R_1(\breve{t})$, $R_2(\breve{t})$, ..., etc.  Those with direction  smaller than $\widetilde{\tau}-\pi$ will be labelled in clockwise sense $R_1^\prime(\breve{t})$, $R_2^\prime(\breve{t})$, etc. Therefore,  
  $R_1(\breve{t})$ is the closest to the admissible ray $R(\widetilde{\tau})$, while  $R_1^\prime(\breve{t})$ is the closest to $R(\widetilde{\tau}-\pi)$. 
     (Warning about the notation: The dependence on $t$ is indicated in Stokes rays $R_1$, $R_2$ etc,   while for the admissible ray $R(\widetilde{\tau})$, $\widetilde{\tau}$ is the direction as in Definition \ref{14feb2016-5}).    See figure \ref{pria15}.

Let $t$ vary from $\breve{t}$ into a neighbouring cell $c_1$,  in such a way that $R_1(t)$ approaches and crosses $R(\widetilde{\tau})$ clockwise.  By Proposition \ref{22marzo2016-2},  $Y_\nu^{cont}(z,t;\widetilde{\tau},c)$ is well defined with canonical asymptotics on a sector having left boundary ray  equal to $R_1(t)$,  for values of $t\in c_1$ just after the crossing.\footnote{As long as $R_1(t)$ does not reach another Stokes ray}

 By assumption,  $Y_\nu^{cont}(z,t;\widetilde{\tau},c)=Y_\nu(z,t;\widetilde{\tau},c_1)$. For $t\in c_1$ just after the crossing, $Y_\nu(z,t;\widetilde{\tau},c_1)$ has canonical asymptotics in $\mathcal{S}_\nu(t)$, 
  which now  has left boundary ray equal to $R_2(t)$. See Figures \ref{pil1} e \ref{pil1bis}. This implies that  $Y_\nu^{cont}(z,t;\widetilde{\tau},c)$ has canonical asymptotics extended up to $R_2(t)$, $t\in c_1$ as above. See Figure \ref{pil2}. 
  
\begin{figure}
\minipage{0.45\textwidth}
\centerline{\includegraphics[width=1\textwidth]{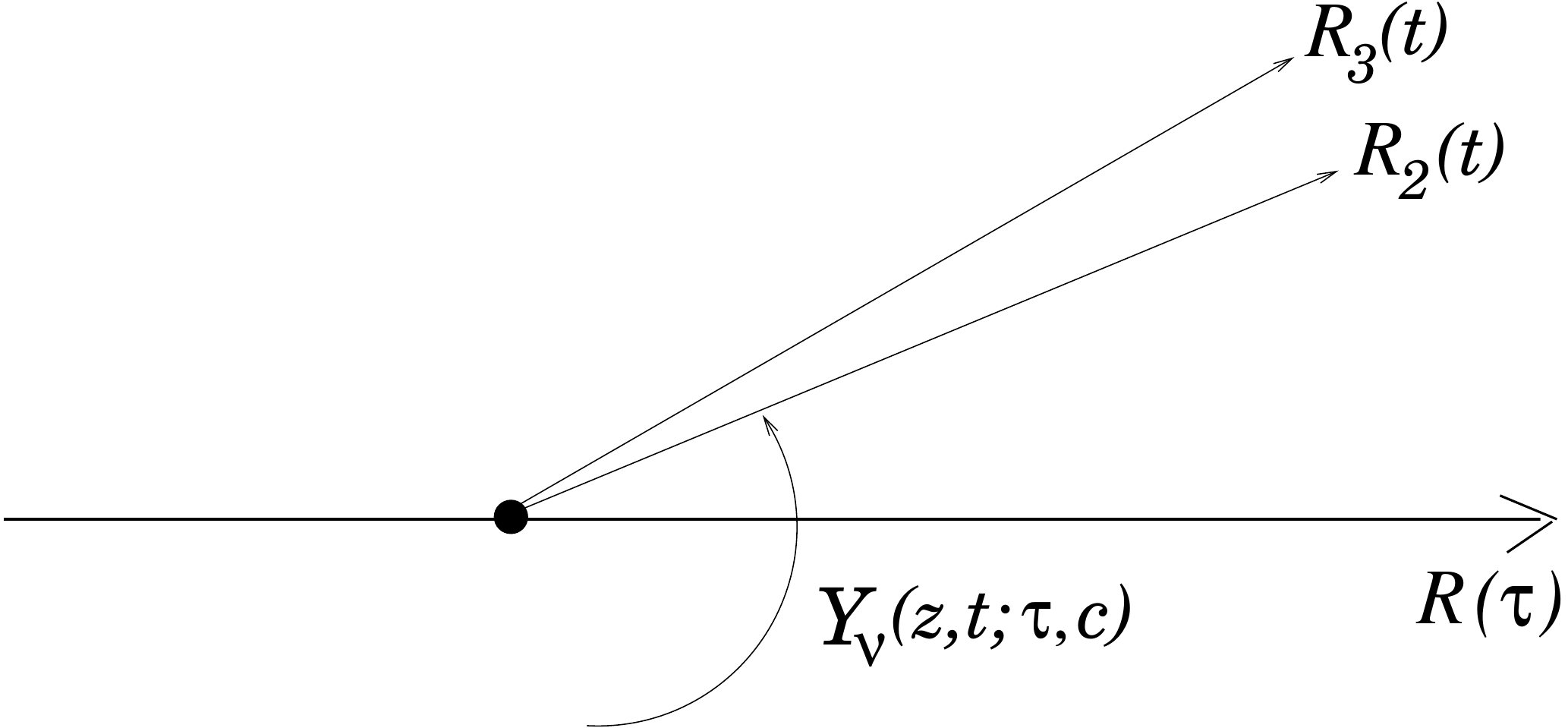}}
\caption{The sector where $Y_\nu(z,t;\widetilde{\tau},c)$ has canonical asymptotics has been extended up to $R_2(t)$, $t\in c$.}
\label{pil4}
\endminipage\hfill
\minipage{0.45\textwidth}
\centerline{\includegraphics[width=0.9\textwidth]{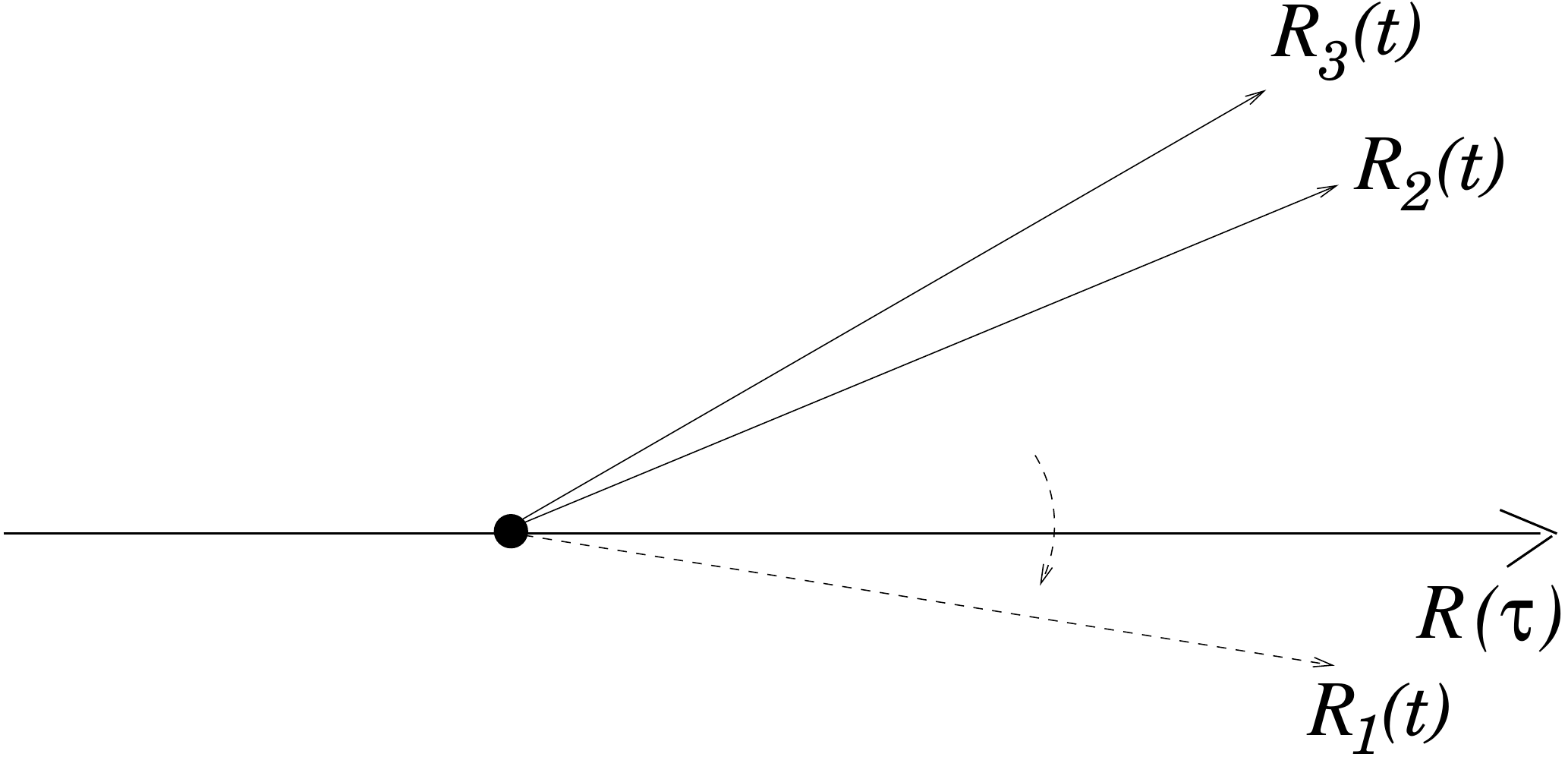}}
\caption{The dashed ``virtual ray" $R_1(t)$ crosses $R(\widetilde{\tau})$, when $t$ enters into $c_1$.}
\label{pil5}
\endminipage
\end{figure}

  Let $t$ go back along the same path, so that $R_1(t)$ crosses $R(\widetilde{\tau})$ anticlockwise. Proposition \ref{22marzo2016-2} now can be applied to $Y_\nu(z,t;\widetilde{\tau},c_1)$ for this crossing.\footnote{In the proof, deform $\widetilde{\tau}\mapsto \widetilde{\tau}+\varepsilon$.} Hence, $Y_\nu(z,t;\widetilde{\tau},c_1)$ has analytic continuation for $t$ before the crossing, certainly up to $\breve{t}$ (because $R_1(t)$ does not cross $R_2(t)$), with canonical asymptotics in a sector having $R_2(\breve{t})$ as left boundary. See Figure \ref{pil3}.  Again,  by assumption, we have that  $Y_\nu(z,\breve{t};\widetilde{\tau},c)=Y_\nu^{cont}(z,\breve{t};\widetilde{\tau},c_1)$. Hence, $Y_\nu(z,\breve{t};\widetilde{\tau},c)$ has canonical asymptotics extended up to the ray $R_2(\breve{t})$. See Figure \ref{pil4}. In conclusion, $R_1(t)$ has been erased.  
  
\begin{figure}
\minipage{0.45\textwidth}
\centerline{\includegraphics[width=1\textwidth]{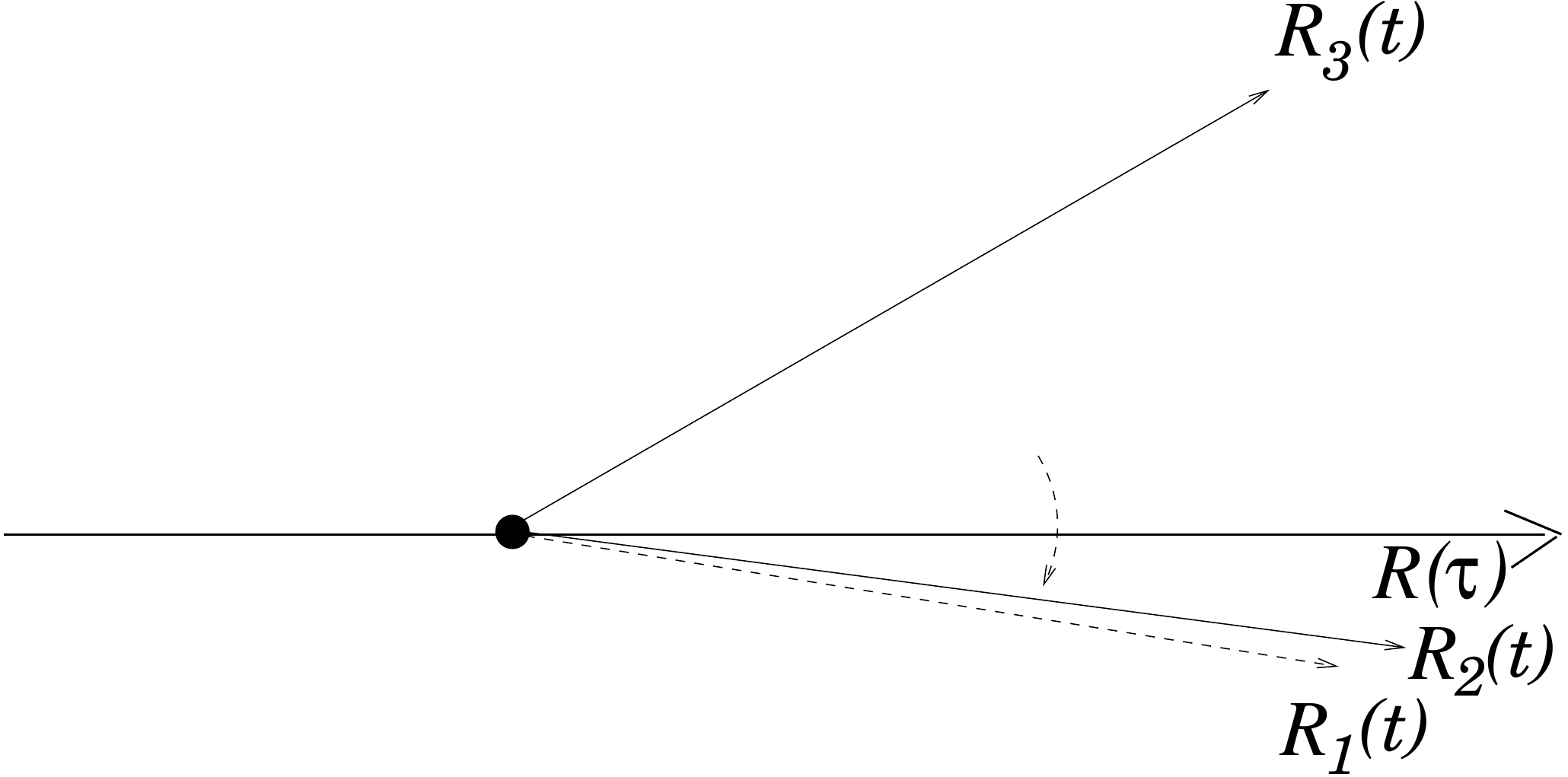}}
\caption{$R_2(t)$ crosses $R(\widetilde{\tau})$ when $t$ enters into $c_2$}
\label{pil6}
\endminipage\hfill
\minipage{0.45\textwidth}
\centerline{\includegraphics[width=0.9\textwidth]{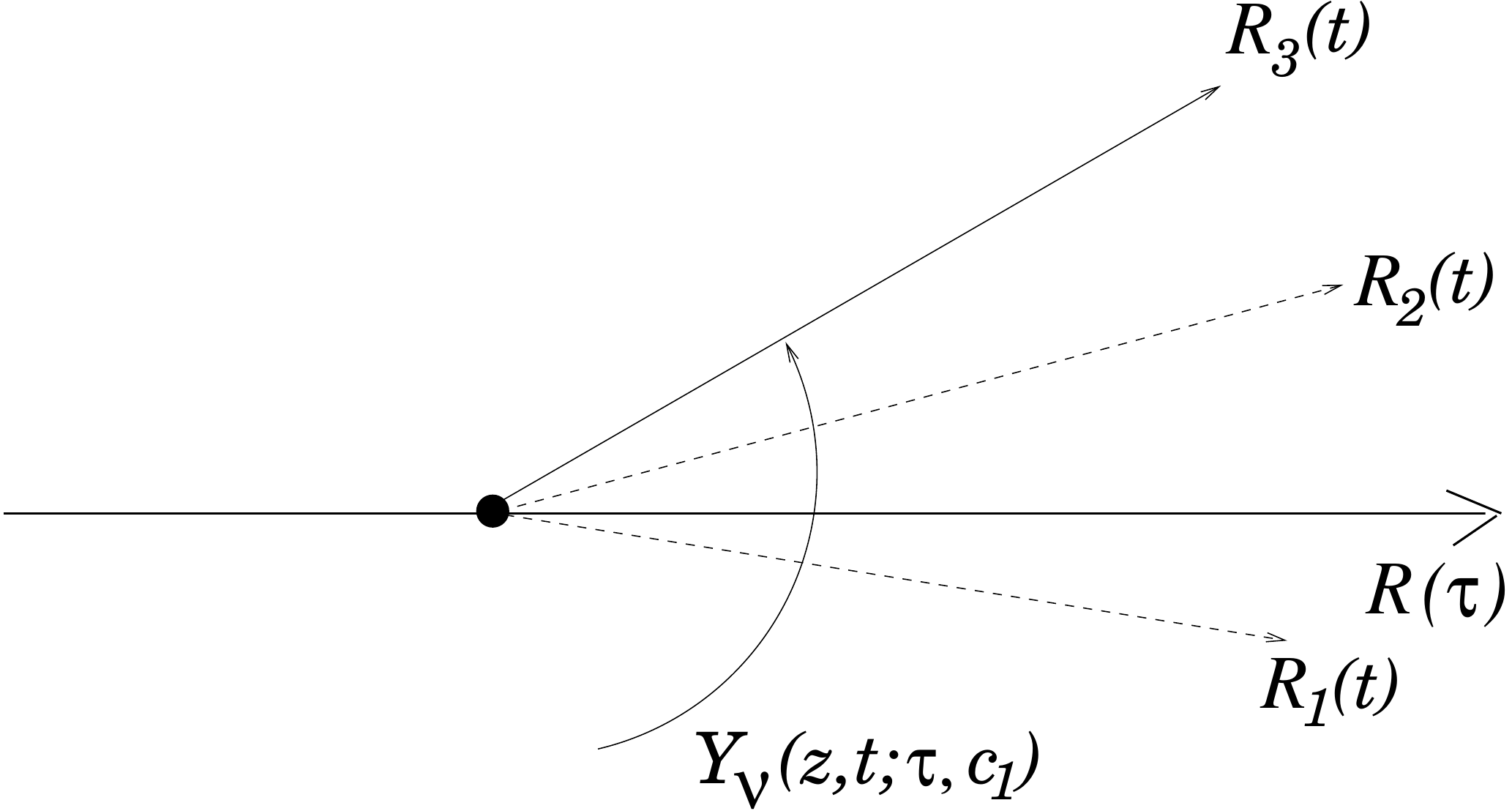}}
\caption{Extension up to $R_3(t)$ of the sector for the asymptotics of  $Y_\nu(z,t;\widetilde{\tau},c_1)$, for $t\in c_1$.}
\label{pil7}
\endminipage
\end{figure}

  \vskip 0.2 cm 
  Step 2. We repeat the arguments analogous to those of Step 1 in order to erase $R_2(t)$. Let $t$ vary in such a way that $R_1(t)$, which is now a ``virtual ray", crosses $R(\widetilde{\tau})$ clockwise, as in step 1. After the crossing, $t\in c_1$ and  $Y_\nu^{cont}(z,t;\widetilde{\tau},c)=Y_\nu(z,t;\widetilde{\tau},c_1)$. Then, let $t$ vary in such a way that also $R_2(t)$ crosses $R(\widetilde{\tau})$ clockwise. See Figures \ref{pil5}, \ref{pil6}. Just after the crossing, $t$ belongs to another cell $c_2$ (clearly, $c_2\neq c$ and $c_1$; see Proposition \ref{19april2016-1}). 
  
  The same discussion done at Step 1 for $Y_\nu(z,t;\widetilde{\tau},c)$ is repeated now for $Y_\nu(z,t;\widetilde{\tau},c_1)$. Indeed, 
  $Y_\nu^{cont}(z,t;\widetilde{\tau},c_1)=Y_\nu(z,t;\widetilde{\tau},c_2)$, for $t\in c_2$ just after $R_2(t)$ has crossed $R(\widetilde{\tau})$. The 
  conclusion, as before, is that $Y_\nu^{cont}(z,t;\widetilde{\tau},c_1)$ has canonical asymptotics extended up to $R_3(t)$ for $t\in c_1$. See Figure \ref{pil7}. 
  
\begin{figure}
\minipage{0.45\textwidth}
\centerline{\includegraphics[width=1\textwidth]{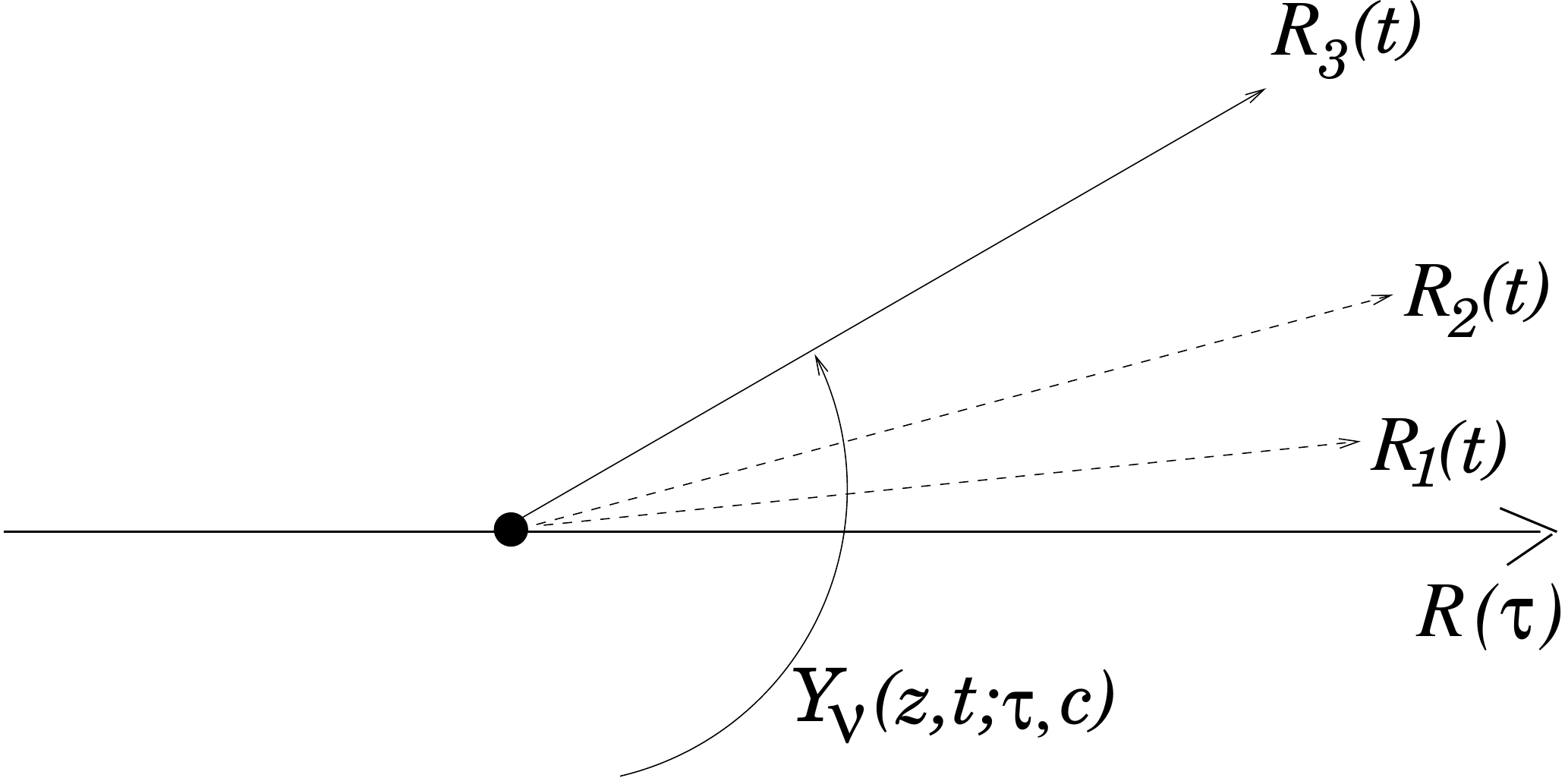}}
\caption{Extension up to $R_3(t)$ of the sector for the asymptotics of  $Y_\nu(z,t;\widetilde{\tau},c)$, for $t\in c$.}
\label{pil8}
\endminipage\hfill
\minipage{0.5\textwidth}
\centerline{\includegraphics[width=0.9\textwidth]{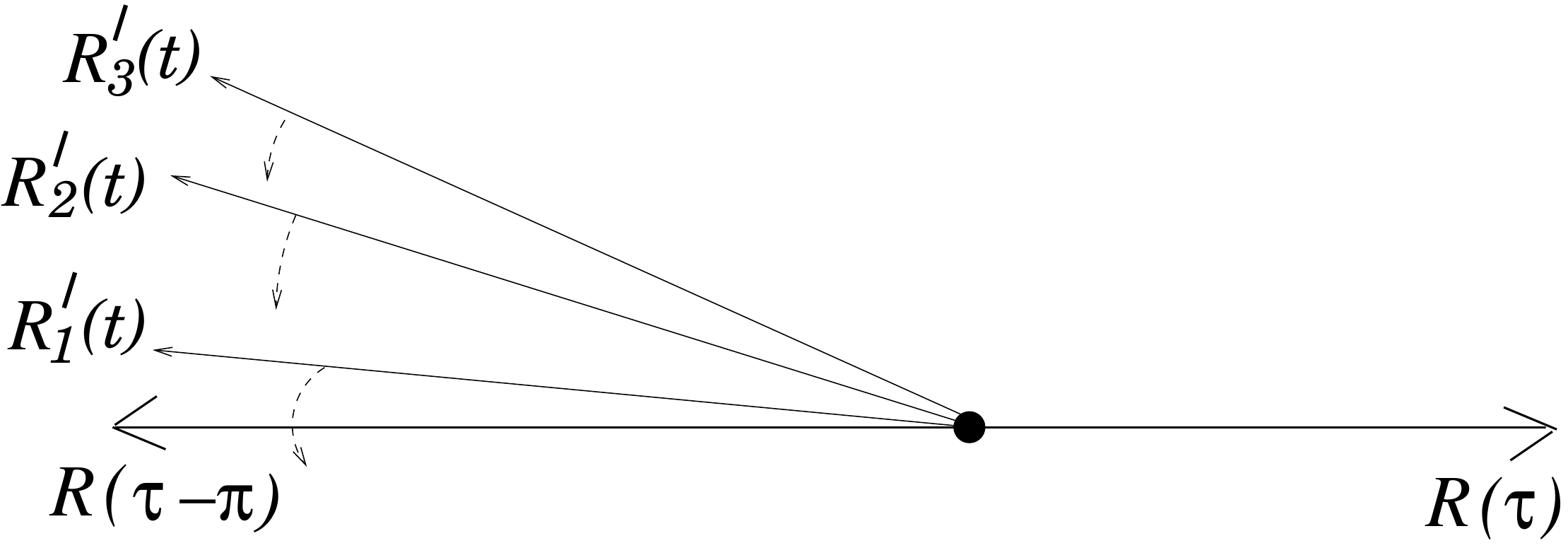}}
\caption{The extension of the sector for the asymptotics of  $Y_\nu(z,t;\widetilde{\tau},c)$ must be done as above also at $R(\widetilde{\tau}-\pi)$, considering crossings as in figure.}
\label{pil9}
\endminipage
\end{figure}

  Now, let $t$ go back along the same path up to $\breve{t}$. Also the virtual ray $R_1(t)$ comes to the initial position, and 
  $Y_\nu(z,\breve{t};\widetilde{\tau},c)=Y_\nu^{cont}(z,\breve{t};\widetilde{\tau},c_1)=Y_\nu^{cont}(z,\breve{t};\widetilde{\tau},c_2)$, with canonical asymptotics extended up to $R_3(\breve{t})$. See figure \ref{pil8}. 

\vskip 0.2 cm 
Step 3. The discussion above can be repeated for all Stokes rays $R_1$, $R_2$, $R_3$ , etc. 
\vskip 0.2 cm 
Step 4.  Observe that the right boundary ray $R_1^\prime$ of the sector where $Y_\nu(z,t;\widetilde{\tau},c)$ has asymptotics is not affected by the above construction. Once the left boundary rays $R_1$, $R_2$,... have been erased,  the same discussion must be repeated considering crossings of the admissible ray $R(\widetilde{\tau}-\pi)$ by the rays $R_1^\prime$, $R_2^\prime$, etc,  as in figure \ref{pil9}.

\vskip 0.2 cm 
In conclusion, all rays $R_1,R_2,...$, $R_1^\prime,R_2^\prime,...$ from unfolding lying in $\widehat{\mathcal{S}}_\nu(\breve{t})\backslash S(\widetilde{\tau}-\pi,\widetilde{\tau})$ are erased. Hence $Y_\nu(z,\breve{t};\widetilde{\tau},c)\equiv Y_\nu(z,\breve{t};\widetilde{\tau})$ has canonical asymptotics  extended up to the closest Stokes rays in $\mathfrak{R}(\breve{t})$ outside $S(\widetilde{\tau}-\pi,\widetilde{\tau})$, namely the asymptotics holds  in $\widehat{\mathcal{S}}_\nu(\breve{t})$. 

The above discussion can be repeated also if one of more  rays among $R_1$, $R_2$, etc. is double (i.e. it corresponds to three eigenvalues) at $\breve{t}$, because as $t$ varies the rays unfold. Thus, the above discussion  holds for any $\breve{t}\in c$ and any  $c$.  Therefore, $Y_\nu(z,t;\widetilde{\tau})$ has asymptotics in  $\widehat{\mathcal{S}}_\nu(t)$ for any $t$ belonging to the union of the 
cells.\footnote{Namely,  $t\in \mathcal{U}_{\epsilon_0}(0)\backslash \left(\Delta\cup X(\widetilde{\tau})\right)=  \mathcal{U}_{\epsilon_0}(0)\backslash \left(\bigcup\widetilde{H}_{ab}\right)$, $a,b$ from unfolding.}

\vskip 0.2 cm 
We observe that  a ray $R_1(t)$, $R_2(t)$, etc, crosses $R(\widetilde{\tau})$ for $t$ equal to a simple point $t_*$ 
 (see Definition \ref{10april2016-1}). The above proof allows to conclude that $Y_\nu(z,t_*;\widetilde{\tau})$ has asymptotics in  $\widehat{\mathcal{S}}_\nu(t_*)$ also when   $\breve t=t_*$.

\begin{figure}
\centerline{\includegraphics[width=0.45\textwidth]{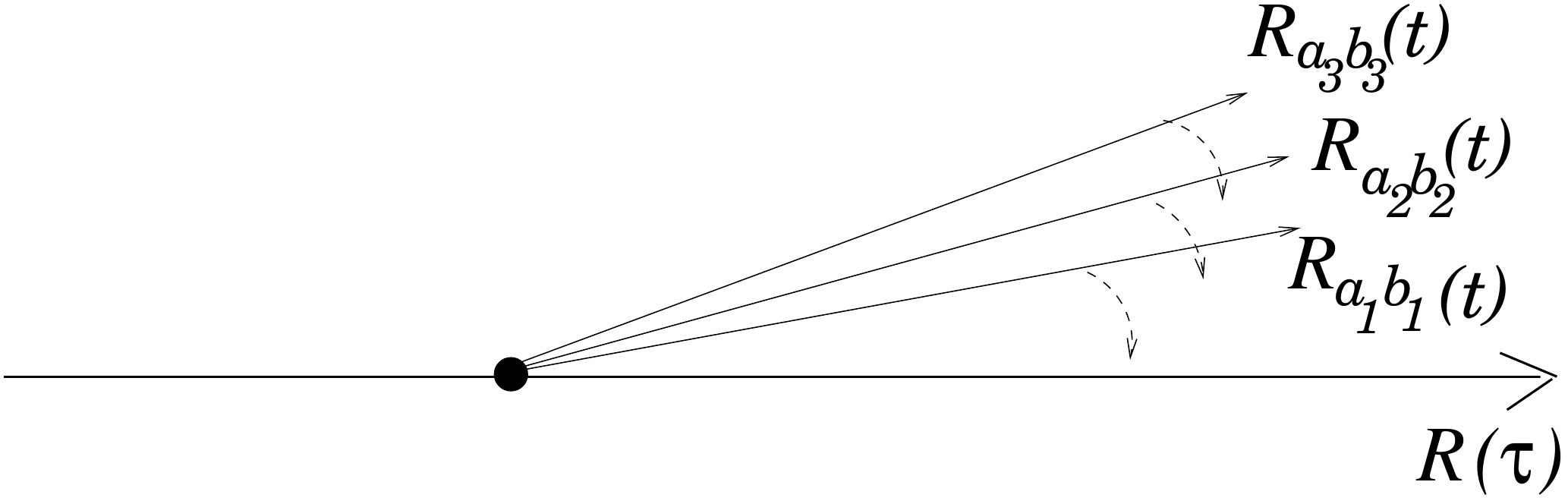}}
\caption{$t$ belongs to a cell $c$ whose boundary contains $\widetilde{H}_{a_1b_1}\cap\widetilde{H}_{a_2b_2}\cap\cdots \cap \widetilde{H}_{a_l,b_l}$, and such that the Stokes rays associated with these hyperplanes cross $R(\widetilde{\tau})$ simultaneously from the same side ($c$ can be taken so that the crossing is clockwise).}
\label{pil10-1}
\end{figure}

\begin{figure}
\minipage{0.5\textwidth}
\centerline{\includegraphics[width=1\textwidth]{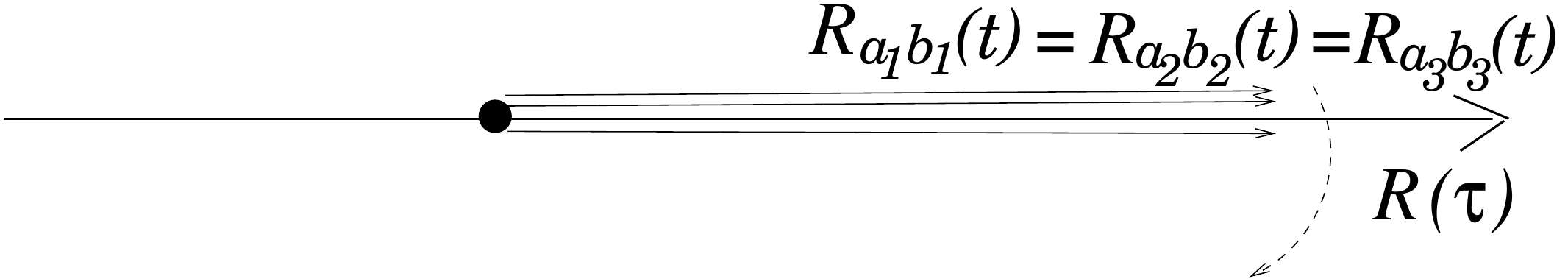}}
\caption{Simultaneous crossing for  $t\in \bigl(\widetilde{H}_{a_1b_1}\cap\widetilde{H}_{a_2b_2}\cap\cdots \cap \widetilde{H}_{a_l,b_l}\bigr)\backslash \Delta$.}
\label{pil10-2}
\endminipage\hfill
\minipage{0.5\textwidth}
\centerline{\includegraphics[width=0.9\textwidth]{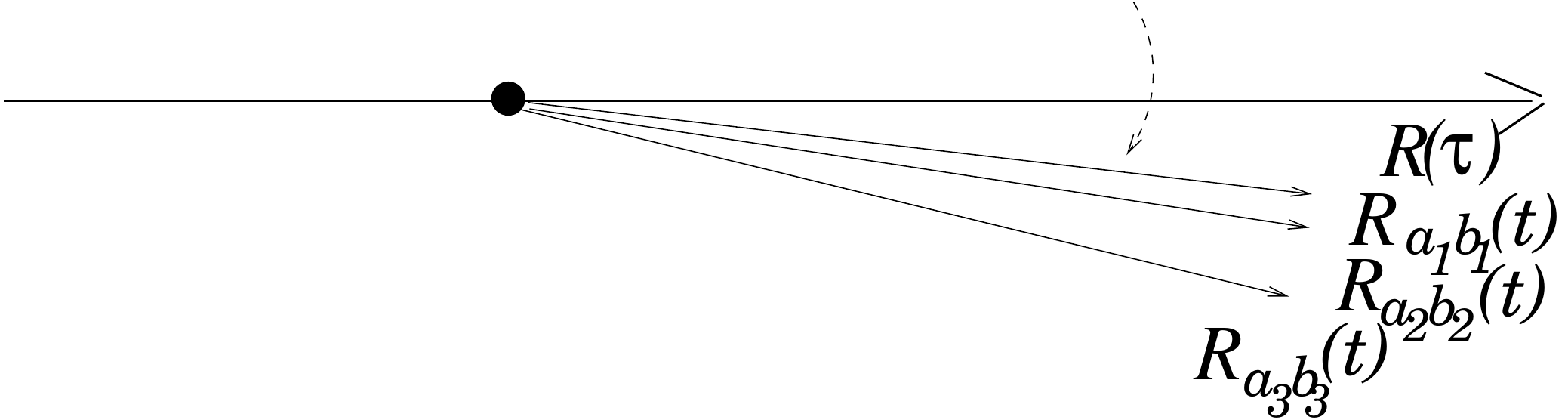}}
\caption{After the simultaneous crossing, $t\in c^\prime$.}
\label{pil10-3}
\endminipage
\end{figure}

\vskip 0.2 cm 
{\bf Part 2:} Points $\breve{t}$ internal to cells and simple points have been considered. It remains  to discuss  non simple points 
  $t_*\in\bigl( \widetilde{H}_{a_1b_1}\cap\widetilde{H}_{a_2b_2}\cap\cdots \cap \widetilde{H}_{a_l,b_l}\bigr)\backslash \Delta$, for some 
  $l\geq 2$.  Consider all the Stokes  rays associated with either 
  one of  $(u_{a_m}(t),u_{b_m}(t))$ or  $(u_{b_m}(t),u_{a_m}(t))$,  $m=1,...,l$,  and lying  in $S(\widetilde{\tau},\widetilde{\tau}+\pi)$.  There exists a cell $c$, among the cells having boundary sharing the above intersection, such that these rays   
    cross $R(\widetilde{\tau})$ clockwise and simultaneously at $t_*$, when $t$ approaches $t_*$ from $c$. Call these rays $R_{a_1b_1}(t)$, $R_{a_2b_2}(t)$, etc.  
    See figures \ref{pil10-1},   \ref{pil10-2},  \ref{pil10-3}. 

Let $t$ start from $\breve{t}\in c$ and vary, reaching $t_*$ and penetrating into a neighbouring cell $c^\prime$  through  $\bigl( \widetilde{H}_{a_1b_1}\cap\widetilde{H}_{a_2b_2}\cap\cdots \cap \widetilde{H}_{a_l,b_l}\bigr)\backslash \Delta$. At $t_*$ the above Stokes rays cross $R(\widetilde{\tau})$ clockwise and simultaneously, from the same side. Hence  $Y_\nu^{cont}(z,t;\widetilde{\tau},c)$ has analytic continuation into $c^\prime$ (here the situation is similar to the continuation from $c_3$ to $c_2$ in figure \ref{cell-april13-9}). 
After the crossing, $t\in c^\prime$ and the same discussion of Part 1 applies. Namely,
 $Y_\nu^{cont}(z,t;\widetilde{\tau},c)=Y_\nu(z,t;\widetilde{\tau},c^\prime)$. The canonical asymptotics is extended up to
  the nearest Stokes ray in $S(\widetilde{\tau},\widetilde{\tau}+\pi)$. Then,\footnote{By a small deformation $\widetilde{\tau}\mapsto \widetilde{\tau}+\varepsilon$.}  as in Proposition \ref{22marzo2016-2},  $Y_\nu(z,t;\widetilde{\tau},c^\prime)$ 
  is analytically continued for $t$ back to $c$, up to  $\breve{t}$.  Therefore,  the 
  asymptotics of $Y_\nu^{cont}(z,t;\widetilde{\tau},c)$ gets extended up to the above mentioned
   nearest Stokes ray in $S(\widetilde{\tau},\widetilde{\tau}+\pi)$. This fact holds also for $t=t_*$.  In this way, $R_{a_1b_1}(t)$, $R_{a_2b_2}(t)$, 
   etc, get erased also at $t_*$. Proceeding as in Part 1, we conclude that 
   $Y_\nu(z,t_*;\widetilde{\tau})\equiv Y_\nu^{cont}(z,t_*;\widetilde{\tau},c)$ has asymptotics in
    the sector $\widehat{\mathcal{S}}_\nu(t_*)$. 
    
    Uniformity follows from Corollary \ref{20marzo2016-5}  and Proposition \ref{22marzo2016-2} applied to any $Y_\nu(z,t;\widetilde{\tau},c^\prime)$. 
$\Box$

\bre 
If $\Lambda(0)=\lambda_1 I$, then  $\widehat{\mathcal{S}}_\nu=\mathcal{R}$, so that the asymptotics extends to $\mathcal{R}$. 
\ere

 \begin{figure}
\minipage{0.45\textwidth}
\centerline{\includegraphics[width=1\textwidth]{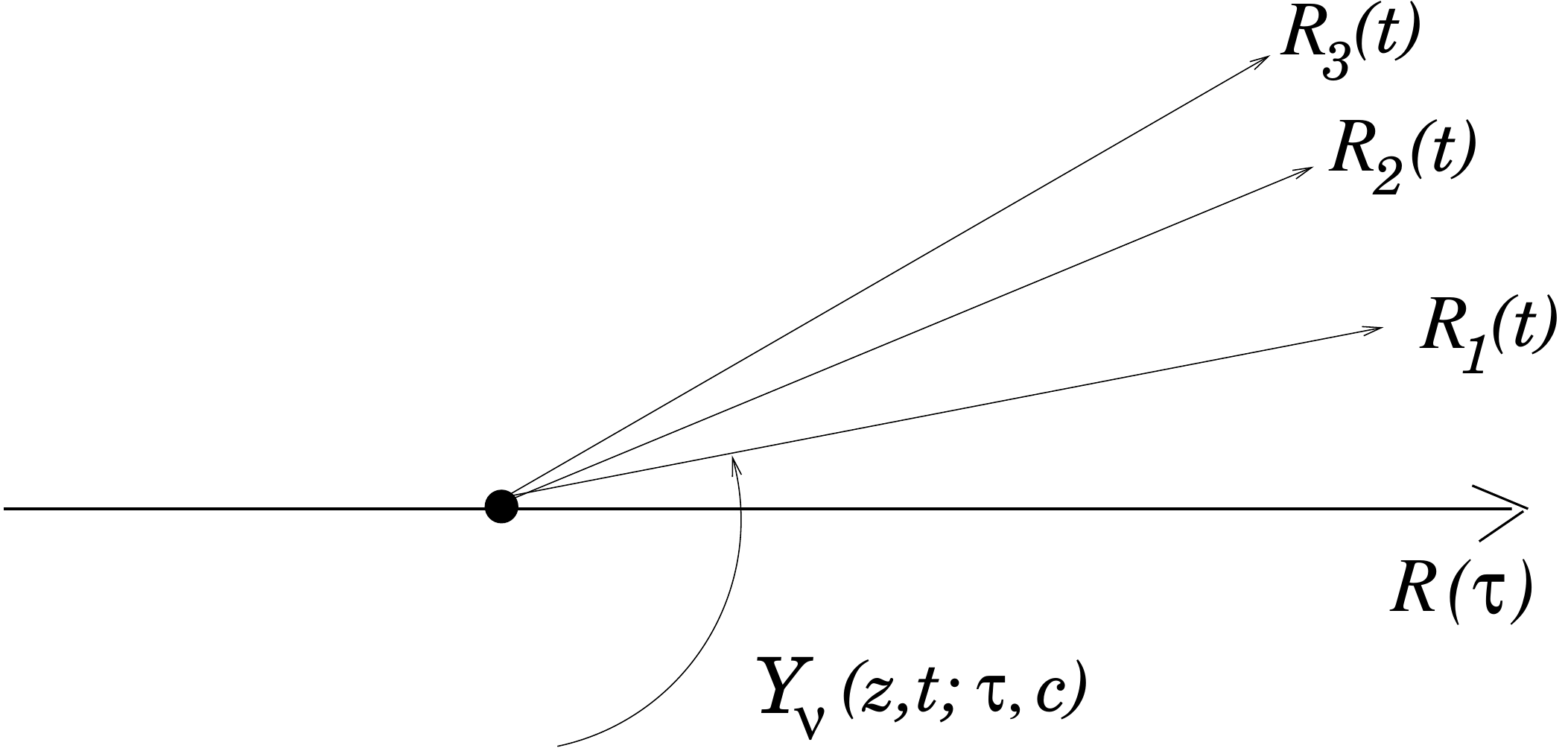}}
\caption{$Y_\nu(z,t;\widetilde{\tau},c)$ for $t\in c$, before $R_1(t)$ crosses $R(\widetilde{\tau})$. A portion of $\mathcal{S}_\nu(t)$ is represented by an arc.}
\label{pil1}
\endminipage\hfill
\minipage{0.45\textwidth}
\centerline{\includegraphics[width=0.9\textwidth]{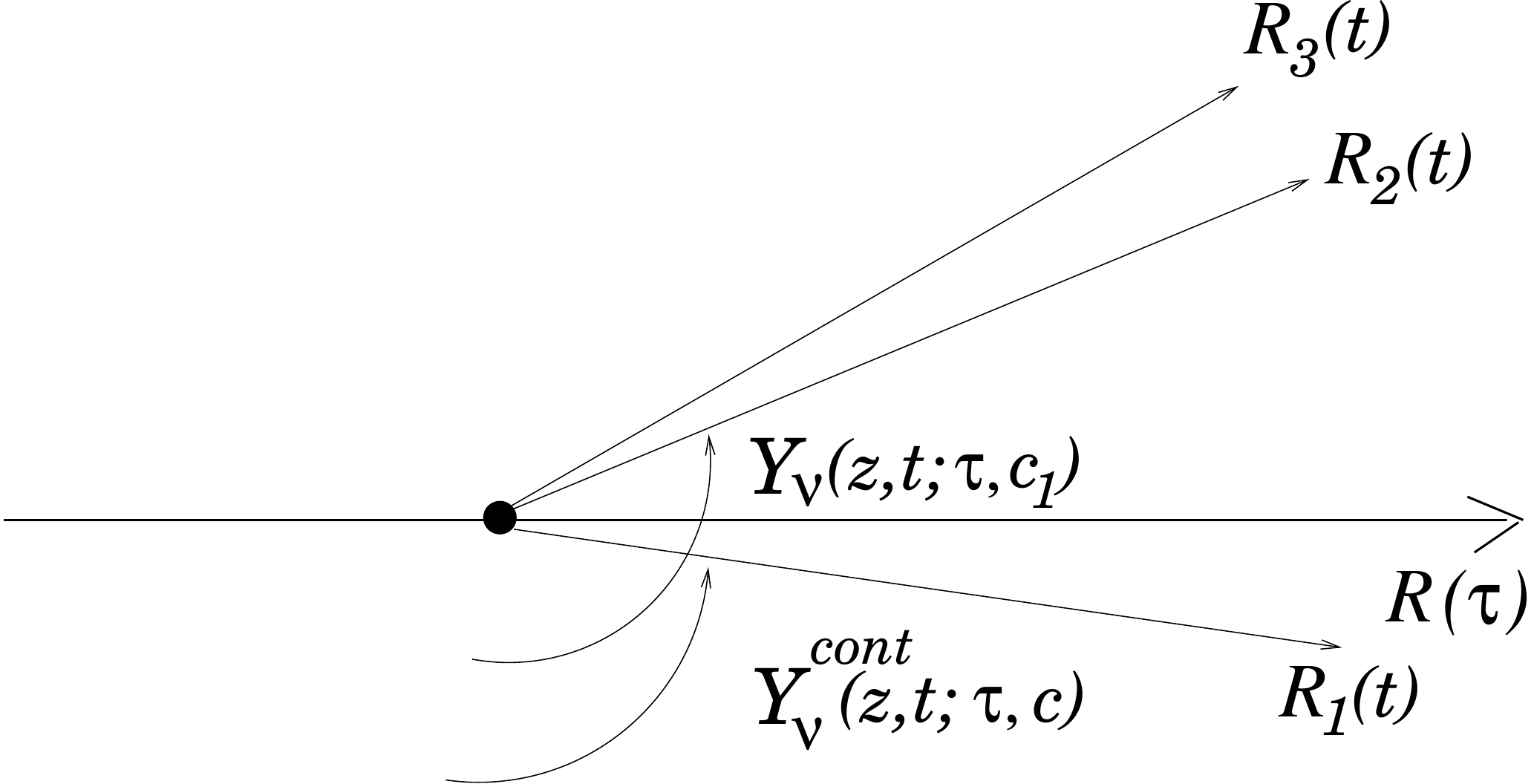}}
\caption{$Y_\nu^{cont}(z,t;\widetilde{\tau},c)$ and $Y_\nu(z,t;\widetilde{\tau},c_1)$ just after $R_1(t)$ has crossed $R(\widetilde{\tau})$. Portions of sectors where the asymptotics holds are represented.}
\label{pil1bis}
\endminipage
\end{figure}

 \begin{figure}
\minipage{0.45\textwidth}
\centerline{\includegraphics[width=1\textwidth]{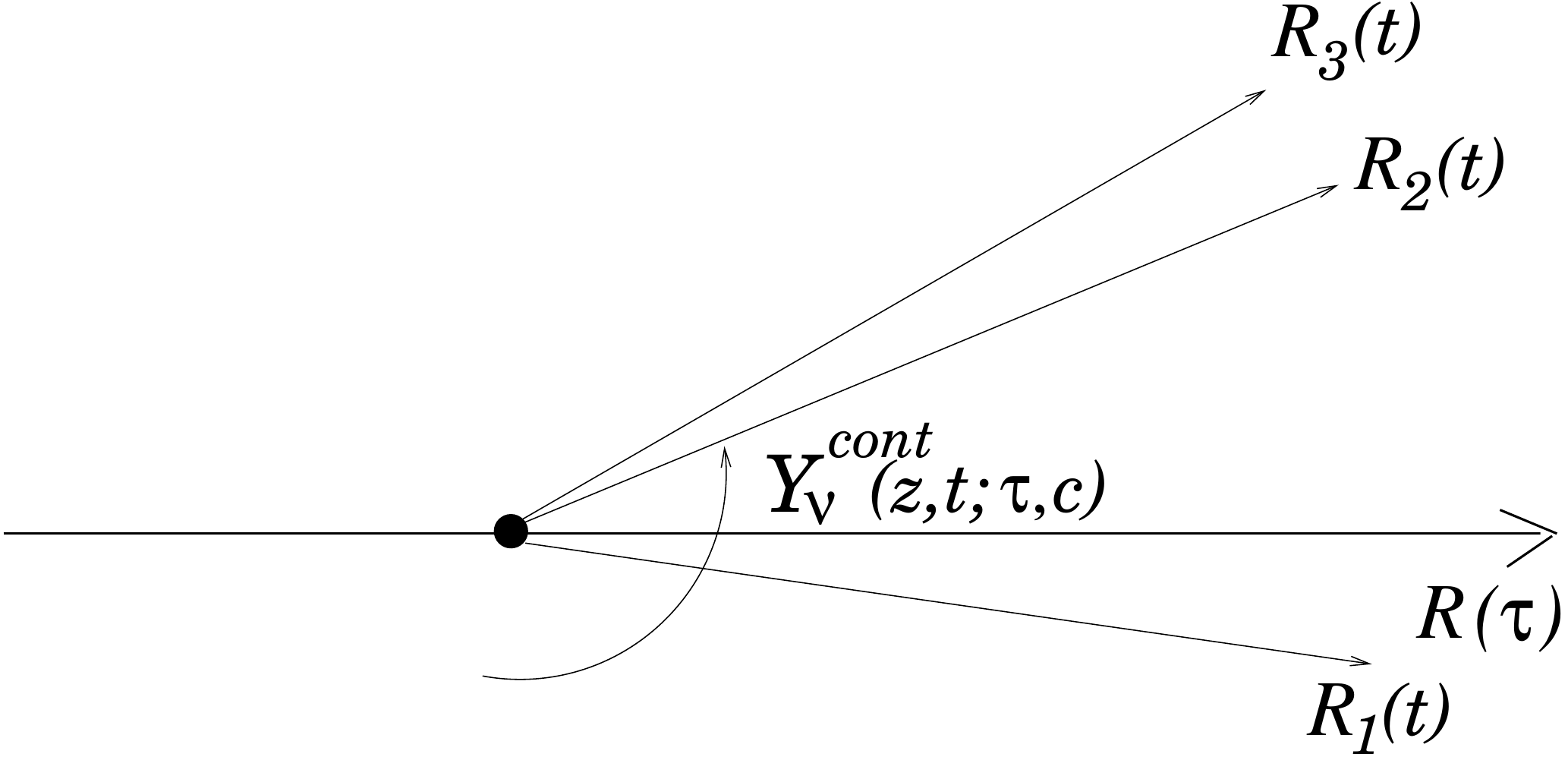}}
\caption{Extension of sector for the asymptotics of $Y_\nu^{cont}(z,t;\widetilde{\tau},c)$, $t\in c_1$.}
\label{pil2}
\endminipage\hfill
\minipage{0.45\textwidth}
\centerline{\includegraphics[width=0.9\textwidth]{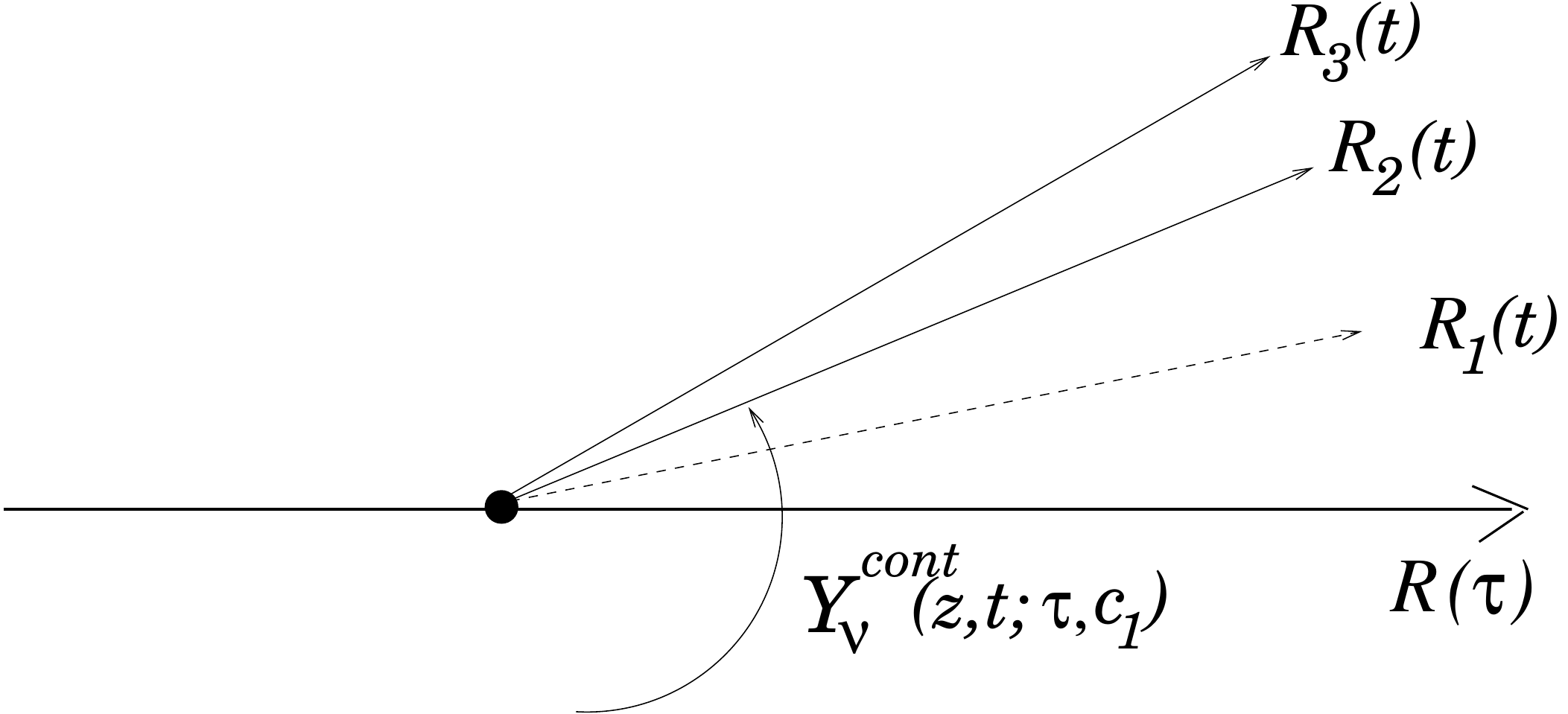}}
\caption{Continuation   $Y_\nu^{cont}(z,t;\widetilde{\tau},c_1)$, $t\in c$ before the crossing. The sector of the asymptotics is represented.}
\label{pil3}
\endminipage
\end{figure}

\begin{figure}
\centerline{\includegraphics[width=0.5\textwidth]{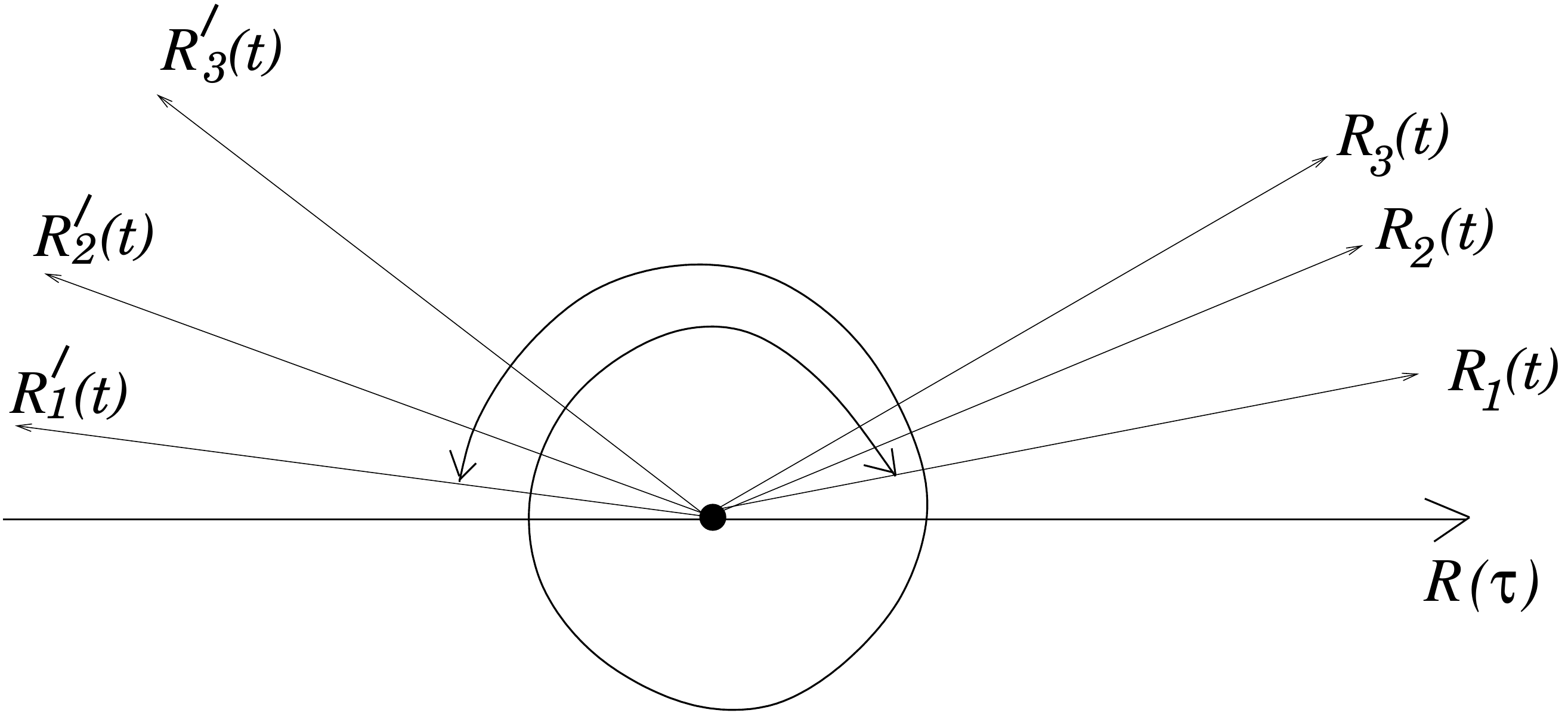}}
\caption{If $\Lambda(0)=\lambda_1 I$, the asymptotics extends  to $S(\arg(R_1(\breve{t}))-2\pi,\arg(R_1^\prime(\breve{t}))+2\pi)$.}
\label{vera}
\end{figure}

\vskip 0.3 cm 
\noindent
$\bullet$ \underline{\it Proof of Theorem \ref{pincopallino}:}   We do the proof for  $Y_{\nu+\mu}(z,t;\widetilde{\tau},c)$. For any other $Y_{\nu+k\mu}$, $k\in\mathbb{Z}$, 
 the proof is the same.    We compute the analytic continuation of $Y_{\nu+\mu}(z,t;\widetilde{\tau},c)$  
 along  loops $\gamma_{ab}$ in $\pi_1(\mathcal{U}_{\epsilon_0}(0)\backslash \Delta, t_{\it base})$,  associated with $u_a(t)$ and $u_b(t)$ in  (\ref{11april2016-1}). For these $a,b$,   only one of the infinitely many rays of directions (\ref{3april2016-1}) is contained in $S(\widetilde{\tau},\widetilde{\tau}+\pi)$  for $t\in c$. 
 We can suppose  that this is the ray 
    $$
  R_{ab}(t):=\left\{z\in\mathcal{R}~\Bigl|~\arg z= \frac{3\pi}{2}-\arg_p\bigl(  t_a-t_b\bigr)+2N_c\pi\right\},
  $$
  (recall that $\arg_p\bigl( u_a(t)-u_b(t)\bigr)=\arg_p\bigl( t_a-t_b\bigr)$) where  $N_c$ is  a suitable integer such that
  $$
\widetilde{\tau}<  \frac{3\pi}{2}-\arg_p\bigl(  t_a-t_b\bigr)+2N_c\pi<\widetilde{\tau}+\pi,\quad\quad
t\in c.
$$
If it is not the above ray,  then it is a ray with $\arg z= \frac{3\pi}{2}-\arg_p\bigl(  t_b-t_a\bigr)+2N^\prime_c$ and suitable $N^\prime_c$, so that the proof holds in the same way.  $R_{ab}(t)$ rotates clockwise as $t$ moves along the support of $\gamma_{ab}$. 

\begin{figure}
\minipage{0.45\textwidth}
\centerline{\includegraphics[width=1\textwidth]{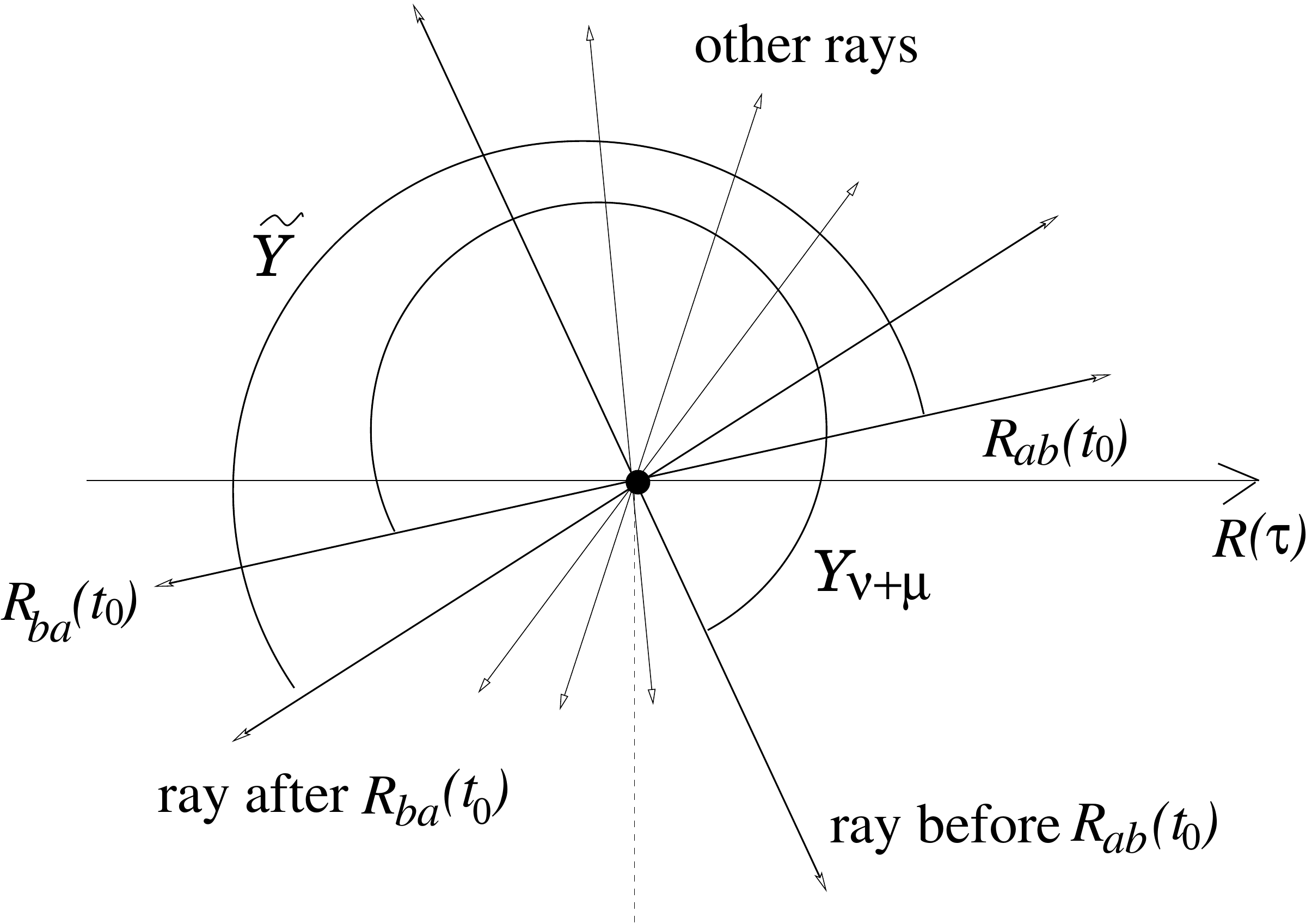}}
\caption{This and the following pictures represent the sheet $S(\widetilde{\tau}-\pi/2,\widetilde{\tau}+3\pi/2)$ (this is the meaning of the dashed vertical half-line). The Stokes rays  at the starting point $t_0$ are represented. $Y_{\nu+\mu}$ is $Y_{\nu+\mu}(z,t;\widetilde{\tau},c)$, while $\widetilde{Y}$ is  $Y_{\nu+\mu}(z,t;\widetilde{\tau},c^\prime)$}
\label{FiguraA}
\endminipage\hfill
\minipage{0.45\textwidth}
\centerline{\includegraphics[width=0.9\textwidth]{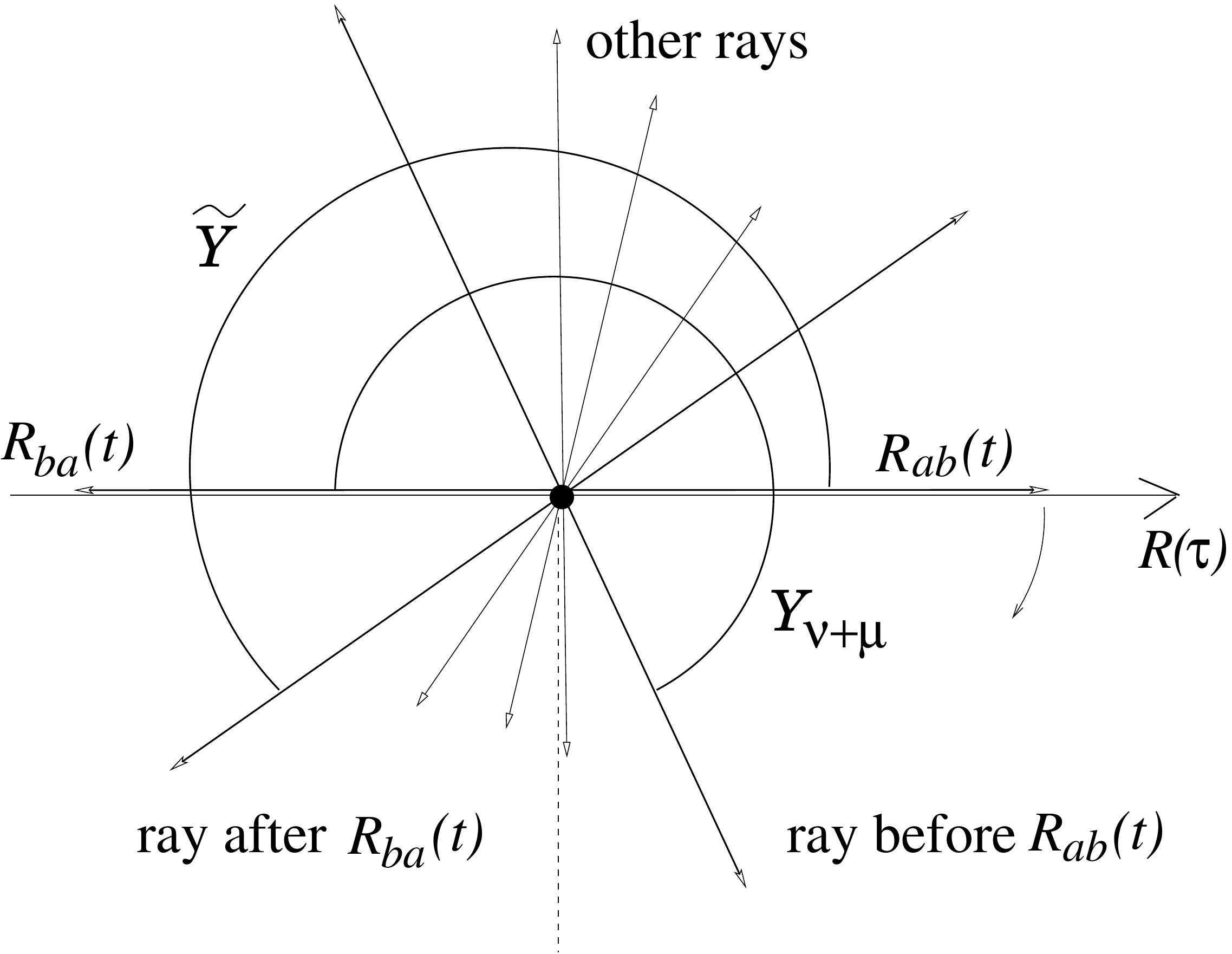}}
\caption{Crossing of $R(\widetilde{\tau})$. Note that also the other rays can move, but never cross the admissible ray $R(\widetilde{\tau})$ or $R(\widetilde{\tau}\pm\pi)$.}
\label{FiguraB}
\endminipage
\end{figure}

For the sake of this proof, if a ray $R$ has angle $\theta$ and $R^\prime$ has angle $\theta +\theta^\prime$, we agree to write $R^\prime=R+\theta^\prime$. Hence, 
   let 
 $$ 
 R_{ba}(t):=R_{ab}(t)+\pi.
$$
See Figure \ref{FiguraA}.

\begin{figure}
\minipage{0.4\textwidth}
\centerline{\includegraphics[width=1\textwidth]{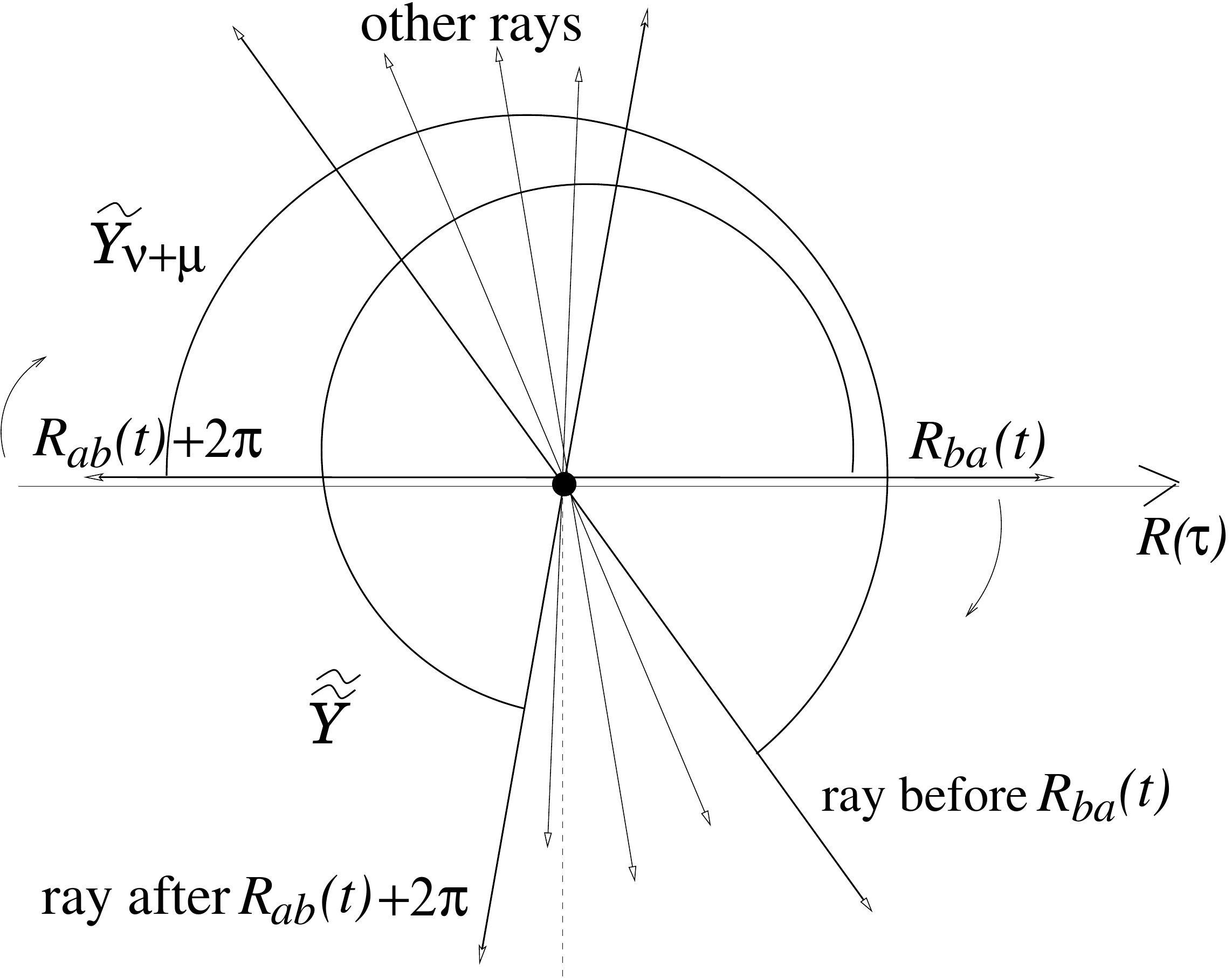}}
\caption{Second crossing. $\widetilde{Y}_{\nu+\mu}$ is   $Y_{\nu+\mu}(z,t;\widetilde{\tau},c^\prime)$ and $\widetilde{\widetilde{Y}}$ is  $Y_{\nu+\mu}(z,t;\widetilde{\tau},c)$. The other rays represented are moving, without crossing  $R(\widetilde{\tau})$ or $R(\widetilde{\tau}\pm\pi)$.}
\label{FiguraF}
\endminipage\hfill
\minipage{0.55\textwidth}
\centerline{\includegraphics[width=0.9\textwidth]{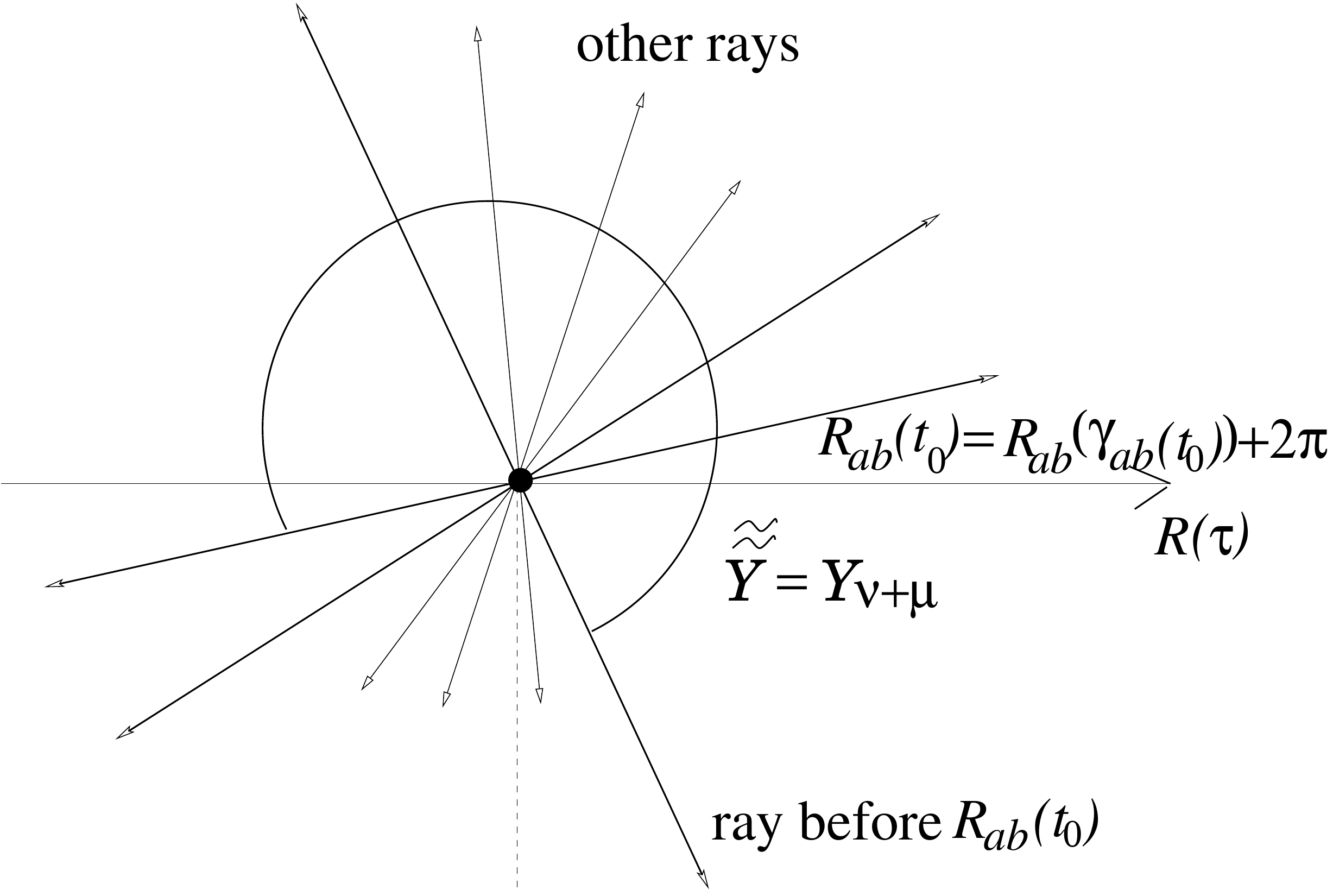}}
\caption{After the loop $\gamma_{ab}$}
\label{FiguraA-tris}
\endminipage
\end{figure}

Assume first that  $a,b$ are such that  for $t\in c$  and $|t_a-t_b|$ sufficiently small, then no projected Stokes rays other than $PR_{ab}$ and $PR_{ba}$ cross $l(\widetilde{\tau})$ when $t$ varies along $\gamma_{ab}$ (the case discussed in figure \ref{105may}). Cases when also other projected Stokes rays   cross $l(\widetilde{\tau})$, as for figure \ref{105may1},  will be discussed later.

\vskip 0.2 cm

 Step 1) As base point consider $t_0\in c$, close to $\widetilde{H}_{ab}$, in  such a way that $R_{ab}(t_0)\subset S(\widetilde{\tau},\widetilde{\tau}+\pi)$    is close to  $R(\widetilde{\tau})$,\footnote{
 $\widetilde{\tau}$ in   $R(\widetilde{\tau})$ is the direction, while $t$ in $R_{ab}(t)$ is the dependence on $t$} 
 and it is the first ray in $S(\widetilde{\tau},\widetilde{\tau}+\pi)$ encountered on moving anti-clockwise from $R(\widetilde{\tau})$.  $Y_{\nu+\mu}(z,t_0;\widetilde{\tau},c)$ has the canonical asymptotics in $\mathcal{S}_{\nu+\mu}(t_0)$, which contains $R(\widetilde{\tau})$.  
By definition,  $\mathcal{S}_{\nu+\mu}(t_0)$ contains $S(\widetilde{\tau},\widetilde{\tau}+\pi)$ and extends to the closest Stokes rays outside. These rays are:

 a) [left ray]  the ray $R_{ba}(t_0)$.
 
 b) [right ray]  the first ray encountered on moving clockwise from  $R_{ab}(t_0)$, which we call ``the ray before" $R_{ab}(t_0)$ ( see Figure \ref{FiguraA}). The name ``before" means that  this ray comes before $R_{ab}(t_0)$  in the natural anti-clockwise orientation of angles). This ray is to the right of $R(\widetilde{\tau})$.

Step 2) As $t$  moves along $\gamma_{ab}$, $R_{ab}(t)$ moves clockwise and  crosses  $R(\widetilde{\tau})$, while $R_{ba}(t)$ crosses $R(\widetilde{\tau}+\pi)$  (see Figure \ref{FiguraB}). 
The curve $\gamma_{ab}$ crosses $\widetilde{H}_{ab}\backslash \Delta$ and penetrates into another cell $c^\prime$. 
As in Proposition \ref{22marzo2016-2}, just before the intersection of the curve with $\widetilde{H}_{ab}\backslash \Delta$,  also $Y_{\nu+\mu}(z,t;\widetilde{\tau},c^\prime)$ is well defined with the same asymptotics as $Y_{\nu+\mu}(z,t;\widetilde{\tau},c)$,  but in the sector bounded by  $R_{ab}(t)$,  as right ray, and the ray coming after $R_{ba}(t)$  in anti-clockwise sense, as left ray, which we call ``the ray after" (see Figures \ref{FiguraA} and \ref{FiguraB}).  A connection matrix $\mathbb{K}^{[ab]}(t)$  (called  {\it Stokes factor}) connects  $Y_{\nu+\mu}(z,t;\widetilde{\tau},c^\prime)$ and $Y_{\nu+\mu}(z,t;\widetilde{\tau},c)$ ,
\be
\label{15dic2015-1}
  Y_{\nu+\mu}(z,t;\widetilde{\tau},c^\prime) =Y_{\nu+\mu}(z,t;\widetilde{\tau},c)~\mathbb{K}^{[ab]}(t).
\ee
$\mathbb{K}^{[ab]}(t)$ is holomorphic on $\mathcal{U}_{\epsilon_0}(0)$, because the fundamental solutions are holomorphic by assumption  2).   Again by the proof of Proposition \ref{22marzo2016-2}, just after the crossing, $Y_{\nu+\mu}(z,t;\widetilde{\tau},c)$ maintains its asymptotics between the ray before $R_{ab}(t)$, which has possibly only slightly moved, and $R_{ba}(t)$.  
Both $Y_{\nu+\mu}(z,t;\widetilde{\tau},c)$ and $Y_{\nu+\mu}(z,t;\widetilde{\tau},c^\prime)$ have the same asymptotics in successive sectors, and in particular they have the same asymptotics on  the sector having right ray $R_{ab}$ and left ray $R_{ba}$. Since $\Re [(u_a-u_b)z]>0$ on this sector, 
it follows from (\ref{15dic2015-1})  that for $t$ in a small open neighbourhood of the intersection point of the curve with $\widetilde{H}_{ab}\backslash \Delta$,  the structure of $\mathbb{K}^{[ab]}(t)$ must be as follows
$$ 
  (\mathbb{K}^{[ab]})_{ii}=1,
  \quad
  \quad
  1\leq i \leq n; 
  \quad
  \quad
  (\mathbb{K}^{[ab]})_{ij}= 0 
  \quad
  \forall
  \quad
  i\neq j \hbox{ except  for } i=b,j=a.
$$
The entry $(\mathbb{K}^{[ab]})_{ba}(t)$ may possibly be different from zero. 
Since $\mathbb{K}^{[ab]}(t)$ is holomorphic on $\mathcal{U}_{\epsilon_0}(0)$, the above structure holds  for every $t\in \mathcal{U}_{\epsilon_0}(0)$.

Step 3) As $t$ moves along $\gamma_{ab}$, $R_{ab}(t)$ continues to rotate clockwise.  It will cross other Stokes rays along the way, but  $Y_{\nu+\mu}(z,t;\widetilde{\tau},c^\prime)$  will maintain its canonical asymptotics in ${\mathcal{S}}_{\nu+\mu}(t)$, because $t\in c^\prime$,  until $R_{ab}(t)$ reaches $R(\widetilde{\tau}-\pi)$. 

Step 4) Just before $R_{ab}(t)$ crosses $R(\widetilde{\tau}-\pi)$,  ${\mathcal{S}}_{\nu+\mu}(t)$ has left ray equal to $R_{ab}(t)+2\pi$ and the
  right ray is the ray before $R_{ba}(t)$.  Again by Proposition \ref{22marzo2016-2}, $Y_{\nu+\mu}(z,t;\widetilde{\tau},c)$ is defined 
   with canonical asymptotics in the sector following $\mathcal{S}_{\nu+\mu}(t)$  anticlockwise 
   (see Figure \ref{FiguraF}). There is a Stokes factor $\widetilde{\mathbb{K}}^{[ab]}(t)$ such that, 
\be
\label{22feb2016-2}
Y_{\nu+\mu}(z,t;\widetilde{\tau},c)
=Y_{\nu+\mu}(z,t;\widetilde{\tau},c^\prime)~\widetilde{\mathbb{K}}^{[ab]}(t). 
\ee
The above relation and the common asymptotic behaviour imply that for $t$ in a neighbourhood of the crossing point  the structure must be 
$$ 
  (\widetilde{\mathbb{K}}^{[ab]})_{ii}=1,
  \quad
  1\leq i \leq n; 
  \quad
  \quad
  (\widetilde{\mathbb{K}}^{[ab]})_{ij}= 0 
  \quad
  \forall
  \quad
  i\neq j \hbox{ except  for } i=a, j=b.
$$
The entry $ (\widetilde{\mathbb{K}})^{[ab]}_{ab}(t)$ may be possibly non zero. By assumption 2), $\widetilde{\mathbb{K}}^{[ab]}(t)$ is holomoprhic on $\mathcal{U}_{\epsilon_0}(0)$, so the above structure holds for any $t\in \mathcal{U}_{\epsilon_0}(0)$.

Step 5) The rotation of $R_{ab}(t)$ continues, crossing  other Stokes rays.  Finally, $R_{ab}(t)$ reaches the position
$$
R_{ab}\Bigl(\gamma_{ab}(t_0)\Bigr)=R_{ab}(t_0)-2\pi,
$$
 after a full rotation of $-2\pi$. 
This corresponds to the full loop $t_a-t_b\mapsto (t_a-t_b)e^{2\pi i}$.

From (\ref{15dic2015-1}) and  (\ref{22feb2016-2}) we conclude that,
\be
\label{15dic2015-2}
Y_{\nu+\mu}(z,t;\widetilde{\tau},c)=Y_{\nu+\mu}(z,t;\widetilde{\tau},c)~\mathbb{K}^{[ab]}(t)\widetilde{\mathbb{K}}^{[ab]}(t),
\quad\quad
t\in \mathcal{U}_{\epsilon_0}(0).
\ee
Hence
$$
\mathbb{K}^{[ab]}(t)\widetilde{\mathbb{K}}^{[ab]}(t)=I,
\quad\quad
t\in\mathcal{U}_{\epsilon_0}(0).
$$
This implies that   
$(\mathbb{K}^{[ab]})_{ba}=(\widetilde{\mathbb{K}}^{[ab]})_{ab}=0$. Therefore, 
\be
\label{11dic2015-3}
\mathbb{K}^{[ab]}(t)=\widetilde{\mathbb{K}}^{[ab]}(t)=I, 
\quad
\quad
t\in \mathcal{U}_{\epsilon_0}(0).
\ee
We conclude from (\ref{15dic2015-1}) or (\ref{22feb2016-2}) that 
\be
\label{22marzo2016-3}
Y_{\nu+\mu}(z,t;\widetilde{\tau},c)=Y_{\nu+\mu}(z,t;\widetilde{\tau},c^\prime),
\quad\quad
t\in \mathcal{U}_{\epsilon_0}(0).
\ee
The above discussion can be repeated for all loops $\gamma_{ab}$ starting in $c$ involving  a simple crossing of $R(\widetilde{\tau})$. 

We now turn to the case when also other projected Stokes  rays,  not only  $PR_{ab}$ and $PR_{ba}$, cross $l(\widetilde{\tau})$ along $\gamma_{ab}$. In this case, the representative of $\gamma_{ab}$ can be decomposed into steps, for each of which the analytic continuation studied above and formula  (\ref{22marzo2016-3}) hold.  See for example the configuration of figure \ref{105may1}.
 In these occurrences, the analytic continuation is done first  from $c$ to $c^\prime$. The passage from $c$ to $c^\prime$   corresponds to the alignment of $u_\gamma$ and $u_a$. Hence,  $Y_{\nu+\mu}(z,t;\widetilde{\tau},c)$ is continued from $c$  to $c^\prime$ and (\ref{22marzo2016-3}) holds.  Then,  $Y_{\nu+\mu}(z,t;\widetilde{\tau},c^\prime)$ can be used in place of $Y_{\nu+\mu}(z,t;\widetilde{\tau},c)$, applying the same proof previously explained, since for  $t\in c^\prime$, if $|t_a-t_b|$ is sufficiently small, then the crossing involves only $PR_{ab}$ and $PR_{ba}$. 


Concluding,  (\ref{22marzo2016-3}) holds for any cell $c^\prime$ which has a boundary in common with $c$. 

Now,  we consider  a cell $c^\prime$ which has a boundary in common with $c$, and we do  the analytic continuation of $Y_{\nu+\mu}(z,t;\widetilde{\tau},c^\prime)$  to all cells $c^{\prime\prime}$ which have a boundary in common with $c^\prime$, in the same way it was done above. In this way, we conclude that $Y_{\nu+\mu}(z,t;\widetilde{\tau},c)=Y_{\nu+\mu}(z,t;\widetilde{\tau},c^\prime)$ and $Y_{\nu+\mu}(z,t;\widetilde{\tau},c^\prime)=Y_{\nu+\mu}(z,t;\widetilde{\tau},c^{\prime\prime})$, for $t\in \mathcal{U}_{\epsilon_0}(0)$. With this procedures, all cells can be reached, so that (\ref{22marzo2016-3}) holds for any cell $c$ and $c^\prime$ of $\mathcal{U}_{\epsilon_0}(0)$. 
For the above reasons,  we are allowed to write 
\be
\label{22marzo2016-5} 
Y_{\nu+\mu}(z,t;\widetilde{\tau}),
\quad\quad
t\in\mathcal{U}_{\epsilon_0}(0),
\ee
in place of $Y_{\nu+\mu}(z,t;\widetilde{\tau},c)$.
\vskip 0.2 cm 
The above conclusions imply that the assumptions of  Lemma \ref{28marzo2016-2} hold.   Lemma \ref{28marzo2016-2}  assures that the asymptotics extends to the closest Stokes rays in $\mathfrak{R}(t)$ outside $S(\widetilde{\tau},\widetilde{\tau}+\pi)$. Hence the asymptotics    
  \be
\label{11nov2015-1}
G_0(t)^{-1}Y_{\nu+\mu}(z,t;\widetilde{\tau})e^{-\Lambda(t)}z^{-B_1(t)}\sim I+\sum_{k=0}^\infty F_k(t)z^{-k}
\ee
holds for $z\to \infty$ in  $\widehat{\mathcal{S}}_{\nu+\mu}( t)$, and $t\in  \mathcal{U}_{\epsilon_0}(0)\backslash \Delta$. 
 {\it A fortiori}, the asymptotics holds in $\widehat{\mathcal{S}}_{\nu+\mu}=\widehat{\mathcal{S}}_{\nu+\mu}(\mathcal{U}_{\epsilon_0}(0))$. It is uniform on any compact subset $K\subset \mathcal{U}_{\epsilon_0}(0)\backslash \Delta$ for $z\to\infty$ in $\widehat{\mathcal{S}}_\nu(K) $. 

\vskip 0.2 cm 
The last  property  to be  verified is that the asymptotics in $\widehat{\mathcal{S}}_{\nu+\mu}$ holds also for $t\in \Delta$.  Let 
$$
R_k(z,t):=G_0(t)^{-1}Y_{\nu+\mu}(z,t;\widetilde{\tau})e^{-\Lambda(t)z}z^{-B_1(t)}-\left(I+\sum_{l=1}^{k-1} F_l(t)z^{-k}\right),
\quad\quad
t\in \mathcal{U}_{\epsilon_0}(0).
$$ 
Let $\left(R_k(z,t)\right)_{ls}$, $l,s=1,...,n$ be the entries of the matrix $R_k$.  Since $R_k $ is  the $k$-th {remainder} of the asymptotic expansion, it satisfies the inequality
\be
\label{30maggio2016-1}
\Bigl|R_k(z,t) \Bigr|:=\max_{l,s=1,...,n}\Bigl| \left(R_k(z,t)\right)_{ls}\Bigr|
\leq \frac{C(k;\overline{S};t)}{|z|^k} ,
\quad\quad
t\in \mathcal{U}_{\epsilon_0}(0)\backslash\Delta, 
\quad\quad
z\in \overline{S},
\ee
for $z$ belonging to  a proper closed subsector $\overline{S} \subset \widehat{\mathcal{S}}_{\nu+\mu}$. Here $C(k;\overline{S};t)$ is a constant depending on $k$, $\overline{S}$ and $t\in \mathcal{U}_{\epsilon_0}(0)\backslash\Delta$.  Our goal is to prove a similar relation for $t\in\Delta$.

We consider $n$ positive numbers $r_a\leq \epsilon_0$, $a=1,...,n$. We further require that for any $i=1,...,s$ and for any $a\neq b$, such that $u_a(0) =u_b(0)=\lambda_i$,  these numbers are distinct, i.e. $r_a\neq r_b$. We introduce the polydisc $\mathcal{U}_{r_1,....,r_n}(0):= \{t\in\mathbb{C}^n~|~|t_a|\leq r_a,~a=1,...,n\} $.
  Clearly, $\mathcal{U}_{r_1,....,r_n}(0)\subset \mathcal{U}_{\epsilon_0}(0)$. 
   Let us denote the {\it skeleton} of  $\mathcal{U}_{r_1,....,r_n}(0)$ with $\Gamma:=\{t\in\mathbb{C}^n~|~|t_a|=r_a,~a=1,...,n\}$.  The above choice of pairwise distinct $r_a$'s assures that $\Gamma\cap \Delta=\emptyset$. 

The inequality (\ref{30maggio2016-1}) holds in $\mathcal{U}_{r_1,....,r_n}(0)\backslash \Delta$ for any fixed $z\in \overline{S}$. Since $R_k(z,t) $ is holomorphic on the interior of   $\mathcal{U}_{r_1,....,r_n}(0)$ and continuous {up to} the boundary, every matrix entry of $R_k(z,t) $  attains its maximum modulus on the {\it Shilov boundary} (cf. \cite{Sa}, page 21-22)  of  $\mathcal{U}_{r_1,....,r_n}(0)$, which coincides with  $\Gamma$. Since  (\ref{30maggio2016-1}) holds on $\Gamma$, we conclude that 
\be
\label{30maggio2016-2}
\Bigl|R_k(z,t) \Bigr|\leq \frac{C(k;\overline{S};\Gamma)}{|z|^k} ,
\quad
\quad
\forall~t\in \mathcal{U}_{r_1,...,r_n}(0),
\ee
where $
C(k;\overline{S};\Gamma)=\max_{t\in \Gamma}C(k;\overline{S};t)
$. This maximum is finite, because the asymptotics is uniform on every compact subset of $\mathcal{U}_{\epsilon_0}(0)\backslash\Delta$. 
The above {estimate} (\ref{30maggio2016-2}) means that the asymptotics (\ref{11nov2015-1}) holds uniformly in $t$ on the whole $\mathcal{U}_{r_1,...,r_n}(0)$, 
including $\Delta$, for $z\to \infty $ in $\overline{S}$.  
 A fortiori, the {asymptotics} holds in $\mathcal{U}_{\epsilon_1}(0)$, with $\epsilon_1\leq \min_a r_a<\epsilon_0$. Since (\ref{30maggio2016-2}) holds for any  
 closed proper subsector $\overline{S}\subset \widehat{\mathcal{S}}_{\nu+\mu}$, by 
 definition $G_0(t)^{-1}Y_{\nu+\mu}(z,t;\widetilde{\tau})e^{-\Lambda(t)}z^{-B_1(t)}$ is asymptotic to $ I+\sum_{k=0}^\infty F_k(t)z^{-k}$ in $\widehat{\mathcal{S}}_{\nu+\mu}$. 
 
 \vskip 0.3 cm 
 It remains to comment on the structure of a Stokes matrix. 
 In the proof above, a ray $R_{ab}(t)$ associated {with} a pair  $u_a(t),u_b(t)$ with $u_a(0)=u_b(0)=\lambda_{i}$  is ``invisible" as far as the asymptotics is concerned, because  $\mathbb{K}^{[ab]}(t)=\widetilde{\mathbb{K}}^{[ab]}(t)=I$ for any $\gamma_{ab}$.  Therefore, in the factorisation of any $ 
\mathbb{S}_\nu(t)$, the Stokes  factors associated {with} rays $3\pi/2-\arg(u_{a}(t)-u_{b}(t))$ mod $2\pi$, with  $u_a(0)=u_b(0)=\lambda_{i}$, are the identity.  
 $\Box$

\section{Meromorphic Continuation}
\label{3marzo2017-2}

In Theorem \ref{pincopallino} we have assumed that  for any $\widetilde{\tau}$-cell $c$ of $ \mathcal{U}_{\epsilon_0}(0)$ and any $k\in\mathbb{Z}$,  the fundamental solution $Y_{\nu+k\mu}(z,t;\widetilde{\tau},c)$ has analytic continuation as a single-valued holomorphic function on the whole $ \mathcal{U}_{\epsilon_0}(0)$.   In this section, we assume that the above fundamental matrices have continuation on the universal covering $\mathcal{R}( \mathcal{U}_{\epsilon_0}(0)\backslash \Delta)$ of $\mathcal{U}_{\epsilon_0}(0)\backslash \Delta$ as meromorphic matrix-valued functions. We show that if the Stokes matrices satisfy a vanishing condition, then the continuation is actually holomorphic and single valued on $\mathcal{U}_{\epsilon_0}(0)\backslash \Delta$. In particular, $\Delta$ {\it is not a branching locus}.

Recall that the Stokes matrices are defined by 
$$
Y_{\nu+(k+1)\mu}(z,t;\widetilde{\tau},c)=Y_{\nu+k\mu}(z,t;\widetilde{\tau},c)~\mathbb{S}_{\nu+k\mu}(t), 
\quad
\hbox{ for } t\in c. 
$$

\begin{shaded}
\bth
\label{8gen2017-4}
Consider the system (\ref{9giugno-2}) (i.e. system (\ref{22novembre2016-3}) of the Introduction) with holomorphic coefficients and  Assumption 1. Let $\Lambda(t)$ be  of the form (\ref{27dic2015-1}), with eigenvalues (\ref{4gen2016-2}) and $\epsilon_0 =\delta_0$ as in  subsection  \ref{16gen2016-5}. Let $\widetilde{\tau}$ be the direction of an admissible ray $R(\widetilde{\tau})$, satisfying $\tau_\nu<\widetilde{\tau}<\tau_{\nu+1}$.

Assume that for any $\widetilde{\tau}$-cell $c$ of $ \mathcal{U}_{\epsilon_0}(0)$ and any $k\in\mathbb{Z}$,  the fundamental solution $Y_{\nu+k\mu}(z,t;\widetilde{\tau},c)$, defined for $t\in c$, has analytic continuation on the universal covering $\mathcal{R}( \mathcal{U}_{\epsilon_0}(0)\backslash \Delta)$ as a meromorphic matrix-valued function. Assume that the entries of the Stokes matrices  satisfy the vanishing condition
\be
\label{6gen2017-1}
(\mathbb{S}_\nu(t))_{ab}=(\mathbb{S}_\nu(t))_{ba}=(\mathbb{S}_{\nu+\mu}(t))_{ab}=(\mathbb{S}_{\nu+\mu}(t))_{ba}=0,
\quad
\quad
\forall t\in c,
\ee
for any $1\leq a\neq b \leq n$ such that $u_a(0)=u_b(0)$. 
\\
\\
Then: 
\begin{enumerate}

\item[$\bullet$] The continuation of $Y_{\nu+k\mu}(z,t;\widetilde{\tau},c)$ defines a single-valued holomorphic (matrix-valued) function on  $\mathcal{U}_{\epsilon_0}(0)\backslash \Delta$.

\item[$\bullet$]   $Y_{\nu+k\mu}(z,t;\widetilde{\tau},c)=Y_{\nu+k\mu}(z,t;\widetilde{\tau},c^\prime)$, for $t\in c$. Therefore, we write $Y_{\nu+k\mu}(z,t;\widetilde{\tau})$

\item[$\bullet$]  The asymptotics 
$$G_0^{-1}(t)Y_{\nu+k\mu}(z,t;\widetilde{\tau})e^{-\Lambda(t)z}z^{-B_1(t)}\sim I+\sum_{j\geq 1}F_j(t)z^{-j},
$$
 holds for $z\to \infty$ in $\widehat{\mathcal{S}}_{\nu+k\mu}(t)$, $t\in\mathcal{U}_{\epsilon_0}(0)\backslash\Delta$.

\end{enumerate}

\eth

\end{shaded}

\bre
\label{26maggio2017-5}
Recall that $B_1(t)=\hbox{diag}(\widehat{A}_1(t))$ is the exponent of formal monodromy,  appearing in the fundamental solutions (\ref{16maggio2017-4}). The formula $\mathbb{S}_{\nu+2\mu}=e^{-2\pi i B_1}\mathbb{S}_\nu~ e^{2\pi i B_1}$, analogous to that of Proposition \ref{26maggio2017-3}, implies that (\ref{6gen2017-1}) holds for any $\mathbb{S}_{\nu+k\mu}$. Notice that the $F_j(t)$'s are holomorphic on $\mathcal{U}_{\epsilon_0}(0)\backslash\Delta$. 
\ere

\vskip 0.2 cm 
\noindent
{\it Proof:} 
Without loss of generality, we  label the eigenvalues  as in (\ref{29gen2016-1})-(\ref{29gen2016-3}), so that $\mathbb{S}_{\nu+k\mu}(t)$ is partitioned into $p_j\times p_k$ blocks ($1\leq j,k\leq s$) such that the $p_j\times p_j$  diagonal blocks  have matrix entries  $(\mathbb{S}_{\nu+k\mu}(t))_{ab}$ corresponding to coalescing eignevalues $u_a(0)=u_b(0)$. 

We consider $Y_{\nu+\mu}(z,t;\widetilde{\tau},c)$. For any other $Y_{\nu+k\mu}(z,t;\widetilde{\tau},c)$ the discussion is analogous.  We denote the meromorphic  continuation  of $Y_{\nu+\mu}(z,t;\widetilde{\tau},c)$ on $\mathcal{R}( \mathcal{U}_{\epsilon_0}(0)\backslash \Delta)$  by $\mathbb Y_{\nu+\mu}(z,\tilde{t};\widetilde{\tau},c)$, $\tilde{t}\in \mathcal{R}( \mathcal{U}_{\epsilon_0}(0)\backslash \Delta)$. Therefore, the continuation   along a loop $\gamma_{ab}$ as in (\ref{11april2016-1}) and (\ref{10maggio2016-1}), starting in $c$,  will be denoted by $\mathbb Y_{\nu+\mu}(z,\gamma_{ab}t;\widetilde{\tau},c)$, where $\tilde{t}=\gamma_{ab}t$ is the point in $\mathcal{R}( \mathcal{U}_{\epsilon_0}(0)\backslash \Delta)$ after the loop. 

We then proceed as in  the proof of Theorem \ref{pincopallino}, up to  {eq.} (\ref{15dic2015-2}). 
Assume first that  $a,b$ are such that  for $t\in c$  and $|t_a-t_b|$ sufficiently small, then no projected Stokes rays other than $PR_{ab}$ and $PR_{ba}$ cross $l(\widetilde{\tau})$ when $t$ varies along $\gamma_{ab}$ (the case discussed in figure \ref{105may}). Cases when also other projected Stokes rays   cross $l(\widetilde{\tau})$, as for figure \ref{105may1},  can be   discussed later as we did in the proof of Theorem \ref{pincopallino}.
The intermediate steps along $\gamma_{ab}$, corresponding to the formulae (\ref{15dic2015-1}) and (\ref{22feb2016-2}), hold. Namely:
 \be
\label{15dic2015-1-BiiS}
  Y_{\nu+\mu}(z,t;\widetilde{\tau},c^\prime) =Y_{\nu+\mu}(z,t;\widetilde{\tau},c)~\mathbb{K}^{[ab]}(t)
\ee
for $t$ in a neighbourhood of the intersection of the support of $\gamma_{ab}$ with the common boundary of $c$ and $c^\prime$ (i.e. $\widetilde{H}_{ab}\backslash \Delta$) corresponding to  $R_{ab}$ crossing $R(\widetilde{\tau})$. Moreover,
\be
\label{22feb2016-2-BiiS}
Y_{\nu+\mu}(z,t;\widetilde{\tau},c)
=Y_{\nu+\mu}(z,t;\widetilde{\tau},c^\prime)~\widetilde{\mathbb{K}}^{[ab]}(t),
\ee
for $t$ in a neighbourhood of the intersection of the support of $\gamma_{ab}$ with the common boundary of $c$ and $c^\prime$  corresponding to  $R_{ab}$ crossing $R(\widetilde{\tau}-\pi)$. Note that to such $t$ there corresponds a point $\tilde{t}$ in the covering, which is reached along $\gamma_{ab}$, so that  $Y_{\nu+\mu}(z,t;\widetilde{\tau},c)$  in the right hand-side of  (\ref{15dic2015-1-BiiS}) becomes $\mathbb{Y}_{\nu+\mu}(z,\tilde{t};\widetilde{\tau},c)$.

 $\mathbb{K}^{[ab]}(t)$, $\widetilde{\mathbb{K}}^{[ab]}(t)$ have the same structure as in the proof of Theorem \ref{pincopallino}, for $t$ in a small open neighborhood of the crossing points. By assumption, $\mathbb{K}^{[ab]}(t)$, $\widetilde{\mathbb{K}}^{[ab]}(t)$ are meromorphic on $\mathcal{R}( \mathcal{U}_{\epsilon_0}(0)\backslash \Delta)$, so they preserve their structure. 

At the end of the loop, $t$ is back to the initial point, but in the universal covering the point $\tilde{t}=\gamma_{ab}t$ is reached and $Y_{\nu+\mu}(z,t;\widetilde{\tau},c)$ has been analytically continued to  $\mathbb Y_{\nu+\mu}(z,\gamma_{ab}t;\widetilde{\tau},c)$. 
 Thus,   the analogous of  formula (\ref{15dic2015-2}) now reads as follows
\be
\label{8gen2017-7}
Y_{\nu+\mu}(z,t;\widetilde{\tau},c)=\mathbb Y_{\nu+\mu}(z,\gamma_{ab}t;\widetilde{\tau},c)~\mathbb{K}^{[ab]}(t)\widetilde{\mathbb{K}}^{[ab]}(t),
\quad
\quad
t\in c.
\ee

We need to compute  the only non trivial entries $(\mathbb{K}^{[ab]}(t))_{ba}$ and $(\widetilde{\mathbb{K}}^{[ab]}(t))_{ab}$.  Let us consider $\mathbb{K}^{[ab]}(t)$. As it is well known,  $\mathbb{S}_{\nu+\mu}$ can be factorised  into Stokes factors. At the beginning of the loop $\gamma_{ab}$, just before $t$ crosses the boundary of the cell $c$ as in Figure  \ref{FiguraA}, we have
$$ 
\mathbb{S}_{\nu+\mu}=\mathbb{K}^{[ab]} \cdot \mathbb{T},
$$
where $\mathbb{K}^{[ab]}$ is a  Stokes factor and  the matrix  $\mathbb{T}$ is factorised into the remaining Stokes factors of $\mathbb{S}_{\nu+\mu}$. For simplicity, we suppose that $\mathbb{S}_{\nu+\mu}$ is upper triangular (namely $a<b$; if not, the discussion is modified in an obvious way): 
\be
\label{9gen2017-2}
\mathbb{S}_{\nu+\mu}=\left(
\begin{array}{ccccc}
I_{p_1} & *            & *    &\cdots & * 
\\
0           & I_{p_2}  &    *         & \cdots & *
\\
0          &       0       & I_{p_3} &  \cdots & * 
\\ \vdots & \vdots & \vdots &     \ddots           &          \vdots 
\\ 0      &        0        & 0      &     0           & I_{p_s}
\end{array}
\right).
\ee  
It follows that $b<a$, namely $\mathbb{K}^{[ab]} $ has entries equal to $1$'s on the diagonal, $0$ {elsewhere}, except for    a non-trivial entry $m_{ba}:=(\mathbb{K}^{[ab]})_{ba}$ above the diagonal in a block corresponding to one of the $I_{p_1}$, ..., $I_{p_s}$ in (\ref{9gen2017-2}).  Let $E_{jk}$ be the matrix with zero entries except for $(E_{jk})_{jk}=1$. Then, $\mathbb{K}^{[ab]} =I+m_{ba}E_{ba}$, {and} we factorise $\mathbb{T}$ as follows:
$$
\mathbb{S}_{\nu+\mu}=(I+m_{ba}E_{ba})\cdot \prod_{j<k
\hbox{ in $V$}
}(I+m_{jk}E_{jk})\cdot \prod_{\hbox{The others }j<k}(I+m_{jk}E_{jk}),
$$
where $V$ is the set of indices $j<k\in \{1,2,...,n\}$ such that $u_j(0)=u_k(0)$ and $(j,k)\neq (b,a)$ (the entries of the diagonal 
blocks of the matrix block partition associated {with}  $p_1,...,p_s$). 

Now, all the {numbers} $m_{ba}$ and  $m_{jk}$ are uniquely
 determined by the entries of $\mathbb{S}_{\nu+\mu}$. This fact follows from the following  result (see for example \cite{BJL1}). Let $S$ be any upper
  triangular matrix with diagonal elements equal to $1$. Label the upper triangular entries entries $(j,k)$, $j<k$, 
   in an arbitrary way,
$$
(j_1,k_1),
\quad
(j_2,k_2),
\quad
...~,
\quad
(j_{\frac{n(n-1)}{2}},k_{\frac{n(n-1)}{2}}).
$$
Then, there exists numbers $m_1$, $m_2$, ..., $m_{\frac{n(n-1)}{2}}$ which are uniquely determined by the labelling and the entries of $S$, such that 
$$ 
S=(I+m_1E_{j_1,k_1})(I+m_2 E_{j_2,k_2})\cdots (I+m_{\frac{n(n-1)}{2}}E_{j_{\frac{n(n-1)}{2}}k_{\frac{n(n-1)}{2}}}).
$$
Indeed, {a} direct computation gives
\be
S= I + \sum_{a=1}^{\frac{n(n-1)}{2}} m_aE_{j_a k_a}+ \hbox{ non linear terms in the $m_a$'s}.
\label{6gen2017-2}
\ee
The commutation relations 
$$ 
E_{ij}E_{jk}=E_{ik},
\quad
\quad
E_{ij}E_{lk}=0 \hbox{ for } j\neq l,
$$ 
imply that the non linear terms are in an upper sub-diagonal lying above the sub-diagonal where the corresponding factors appear. Hence, (\ref{6gen2017-2}) gives  uniquely solvable recursive relations, expressing the $m_a$'s in terms of the entries of $S$. 

Applying the above procedure to $S=\mathbb{S}_{\nu+\mu}$, and keeping (\ref{6gen2017-1}) into account, we obtain 
$$ 
m_{ba}=0, 
\quad
m_{jk}=0 ~\forall  j<k \hbox{ in } V. 
$$
This proves that 
$$
\mathbb{K}^{[ab]}(t) =I,
$$
 for $t$ in a small open neighborhood of the intersection point of the curve $\gamma_{ab}$ with $\widetilde{H}_{ab}\backslash \Delta$. This structure is preserved by analytic continuation. 
Analogously, we factorise into Stokes factor the (lower triangular) matrix   $\mathbb{S}_\nu=\widetilde{\mathbb{T}}\cdot \widetilde{\mathbb{K}}^{[ab]}$ and prove that 
$$\widetilde{\mathbb{K}}^{[ab]}=I.$$

We conclude that 
$$Y_{\nu+\mu}(z,t;\widetilde{\tau},c)=\mathbb Y_{\nu+\mu}(z,\gamma_{ab}t;\widetilde{\tau},c).
$$
Formulae (\ref{15dic2015-1-BiiS}) and (\ref{22feb2016-2-BiiS}) also imply that 
\be
\label{6gen2017-3}
Y_{\nu+\mu}(z,t;\widetilde{\tau},c^\prime) =Y_{\nu+\mu}(z,t;\widetilde{\tau},c)
\ee
This discussion can be repeated for any loop and any cell, as we did in the proof of Theorem \ref{pincopallino} in the paragraphs following {eq.} (\ref{22marzo2016-3}).  Since $Y_{\nu+\mu}(z,t;\widetilde{\tau},c)$ is holomorphic on $c$ by Corollary \ref{20marzo2016-5}, the above formulae imply  the analyticity of $Y_{\nu+\mu}(z,t;\widetilde{\tau},c^\prime)$ on $\mathcal{U}_{\epsilon_0}(0)\backslash \Delta$. Since (\ref{6gen2017-3}) holds, the first two statements are proved. 

{Equation} (\ref{6gen2017-3}) also implies that the rays $R_{ab}$ and $R_{ba}$ are not the boundaries of the sector where the asymptotic behaviour  of $Y_{\nu+\mu}(z,t;\widetilde{\tau})$  holds. The above discussion  repeated for all $a,b$ such that $u_a(0)=u_b(0)$ proves the third statement of the theorem. 
 $\Box$

\vskip 0.5 cm 
\noindent
\hrule

\vskip 0.5 cm
\centerline{\bf \Large PART IV: Isomonodromy Deformations of system  (\ref{ourcase12}).}
\vskip 0.2 cm 
\centerline{\bf \Large Theorem  \ref{16dicembre2016-1},  Corollary \ref{17dicembre2016-2} and Theorem  \ref{9gen2017-1}}

\vskip 0.3 cm


We have established the theory of coalescence in $\mathcal{U}_{\epsilon_0}(0)$, and the corresponding characterisation of the limiting Stokes matrices for the  system (\ref{22novembre2016-3}) -- namely system (\ref{9giugno-2}) of Section \ref{preRAY} -- under  Assumption 1, or equivalently for the system (\ref{17luglio2016-1}). 
We  now consider the system (\ref{1novembre2016-1}) {under Assumption} 1, already put in the form (\ref{ourcase12}),  namely 
$$
\frac{dY}{dz}=\widehat{A}(z,t)Y,
\quad
\quad
\widehat{A}(z,t)=\Lambda(t)+\frac{\widehat{A}_1(t)}{z},
$$  
and study its isomonodromy deformations. 
The eigenvalues are taken to be linear in $t$,  as in (\ref{4gen2016-2}):
$$
u_i(t)=u_i(0)+t_i,
\quad
1\leq i \leq n.
$$

\section{Structure of Fundamental Solutions in  Levelt form at $z=0$}
\label{5gen2016-2}

At any point $t\in \mathcal{U}_{\epsilon_0}(0)$, let $\mu_1(t),\mu_2(t),...,\mu_n(t)$ be the (non necessarily distinct)
 eigenvalues of $\widehat{A}_1(t)$, and let $J^{(0)}(t)$ be  a Jordan form   {of $\hat A_1(t)$}, with diag$(J^{(0)})=$ diag$(\mu_1,...,\mu_n)$ (see also  (\ref{9feb2016-1}) below). The eigenvalues are decomposed uniquely as, 
$$ 
\mu_i(t)=d_i^{(0)}(t)+\rho_i^{(0)}(t),
\quad
\quad
d_i^{(0)}(t)\in\mathbb{Z},
\quad
~0\leq \Re \rho_i^{(0)}(t)<1.
$$
Let $D^{(0)}(t)=\hbox{diag}(d_1^{(0)}(t),...,d_n^{(0)}(t))$, which is piecewise constant, 
so that 
$$ 
 J^{(0)}(t)=D^{(0)}(t)+S^{(0)}(t) ,
$$
where $S^{(0)}(t)$ is the Jordan matrix with $\hbox{diag}(S^{(0)})=\hbox{diag}(\rho_1^{(0)},...,\rho_n^{(0)})$. 

Let $\mathcal{V}$ be an open connected subset of $\mathcal{U}_{\epsilon_0}(0)$. 
 In order to write a solution at $z=0$ in  Levelt form which is holomorphic on $\mathcal{V}$, we need the following  assumption. 
\vskip 0.2 cm
{\bf Assumption 2:} We assume that $\widehat{A}_1(t)$ is {\it holomorphically similar} to $J^{(0)}(t)$ on $\mathcal{V}$. This means that there exists an invertible matrix $G^{(0)}(t)$ holomorphic on $\mathcal{V}$ such that 
$$
 (G^{(0)}(t))^{-1}\widehat{A}_1(t)~G^{(0)}(t)=J^{(0)}(t).
$$
\vskip 0.2 cm 
\noindent
Assumption 2  in $\mathcal{V}$ implies that the eigenvalues $\mu_i(t)$ are holomorphic on $\mathcal{V}$.  In the isomonodromic case (to be {defined} below),   Assumption 2  for $\mathcal{V}=\mathcal{U}_{\epsilon_0}(0)$ turns out to be equivalent to the vanishing condition (\ref{31luglio2016-1}). See Proposition \ref{28dic2015-4} below.

\bre
In order to realise the above assumption it is not {sufficient} to assume, for example, that the eigenvalues of $\widehat{A}_1(t)$ are independent of $t$, as the example 
$\widehat{A}_1(t)=\left(\begin{array}{cc}\mu & t \\ 0 & \mu \end{array}\right)$ shows. Sufficient conditions can be {found} in {the Wasow's book} \cite{Wasow}, Ch. VII.
\ere

\noindent
With Assumption 2, the following fundamental solutions in Levelt form are found.

\vskip 0.2 cm

A] If $\widehat{A}_1(t)$ has distinct eigenvalues at any point of $\mathcal{V}$, it is automatically holomorphically similar to 
$$
\widehat{\mu}(t):=\hbox{diag}(\mu_1(t),~...~,\mu_n(t)).
$$
A fundamental matrix exists of the form 
$$ 
Y^{(0)}(z,t)=G^{(0)}(t)\Bigl( I+\sum_{l=1}^\infty \Psi_l(t)z^l\Bigr) z^{\widehat{\mu}(t)}.
$$ 
Each matrix $\Psi_l(t)$ is holomorphic on $\mathcal{V}$,   and the series $ I+\sum_{l=1}^\infty \Psi_l(t)z^l$ is absolutely convergent for $|z|$ bounded, defining a holomorphic matrix-valued function in $(z,t)$   on $\{|z|<r\}\times \mathcal{V}$, for any $r>0$.

\vskip 0.2 cm 

B] If  $\mu_i(t)-\mu_j(t)\not \in\mathbb{Z}\backslash\{0\}$ for any $i\neq j$ and any $t\in\bar{\mathcal{V}}$,  then there exists a fundamental matrix
$$ 
Y^{(0)}(z,t)=G^{(0)}(t)\Bigl( I+\sum_{l=1}^\infty \Psi_l(t)z^l\Bigr) z^{J^{(0)}(t)},
$$ 
such that $G^{(0)}(t)$, $J^{(0)}(t)$ and each matrix $\Psi_l(t)$ are holomoprhic on  $\mathcal{V}$,   and the series $ I+\sum_{l=1}^\infty \Psi_l(t)z^l$ is absolutely convergent for $|z|$ bounded, defining a holomorphic matrix-valued function  in $(z,t)$  on $\{|z|<r\}\times \mathcal{V}$, for any $r>0$.

\vskip 0.2 cm 

The above forms of the matrix $Y^{(0)}(z,t)$ are obtained by a recursive procedure (see  \cite{Wasow}),  aimed at constructing a gauge transformation $Y=G^{(0)}(t)\Bigl( I+\sum_{l=1}^\infty \Psi_l(t)z^l\Bigr)\mathcal{Y}$  that reduces   the linear system  to a simple form $\frac{d \mathcal{Y}}{dz}=\frac{J^{(0)}}{z}\mathcal{Y}$, whose solution $z^{J^{(0)}(t)}$ can be immediately written.   
In {\it resonant cases}, namely when $\mu_i(t)-\mu_j(t)\in\mathbb{Z}\backslash\{0\}$, this procedure yields a gauge transformation $Y=G^{(0)}(t)\Bigl( I+\sum_{l=1}^\infty \Psi_l(t)z^l\Bigr)\mathcal{Y}$ that reduces the system to the  form  
 \be
 \label{26maggio2017-6}
 \frac{d \mathcal{Y}}{dz}=\frac{1}{z}
\left(
J^{(0)}(t)+R_1(t) z+\cdots +R_\kappa(t) z^\kappa \right)\mathcal{Y},
\ee
where  $1 \leq \kappa$ is the maximal integer difference of eigenvalues of $J^{(0)}$, and  the $R_j(t)$'s are certain nilpotent matrices (see (\ref{9feb2016-2}) below  for more details).  These matrix coefficients  may be discontinuous in $t$, even if Assumption 2 is made. In order to avoid this,  we need  the following

\vskip 0.2 cm
{\bf (Temporary) Assumption 3} [Resonant Case]: If for some $i\neq j$ it happens that $\mu_i(t)-\mu_j(t)\in\mathbb{Z}\backslash\{0\}$ at a point $t\in \mathcal{V}$, then we require that $\mu_i(t)-\mu_j(t)=\hbox{constant}\in\mathbb{Z}\backslash\{0\}$ all over  $\mathcal{V}$. If moreover $J^{(0)}(t)$ is not diagonal, then  we require that the $d_i$'s, $1\leq i \leq n$, are constant on $\mathcal{V}$. 

\vskip 0.2 cm 

Assumptions 3  certainly holds if the eigenvalues $\mu_1$, ..., $\mu_n$ are independent of $t$ in $\mathcal{V}$, namely in the isomonodromic case of Definition \ref{30dic2015-1} below.\footnote{In case we define a deformation to be isomonodromic when the monodromy matrices are constant, this is still true, namely $\mu_1$, ..., $\mu_n$ are independent of $t$. See Lemma 1 of \cite{Bolibruch1}.}   Hence, Assumptions  3 is only ``temporary'' here, being unnecessary in the isomonodromic case. 
 
 When Assumptions 2 and 3 hold together, fundamental matrices in Levelt form can always be constructed in such a way that they are holomorphic on $\mathcal{V}$.  Besides the cases A] and B] (which require only Assumption 2), we have the following resonant cases:

\vskip 0.2 cm 
C] If $J^{(0)}(t)\equiv \widehat{\mu}(t):=\hbox{diag}(\mu_1(t),\mu_1(t),...,\mu_n(t))$ (eigenvalues non necessarily distinct) then there exists a fundamental matrix
$$ 
Y^{(0)}(z,t)=G^{(0)}(t)\Bigl( I+\sum_{l=1}^\infty \Psi_l(t)z^l\Bigr) z^{\widehat{\mu}(t)}z^{R^{(0)}(t)},
$$ 
were the matrix $R^{(0)}(t):=R_1(t)+\cdots R_\kappa(t)$ has entries $R^{(0)}_{ij}(t)\neq 0$ only if $\mu_i(t)-\mu_j(t)\in \mathbb{N}\backslash\{0\}$. Moreover,  $G^{(0)}(t)$,  $\widehat{\mu}(t)$ $R^{(0)}(t)$ and each matrix $\Psi_l(t)$ can be chosen holomorphic on $\mathcal{V}$, and the series $ I+\sum_{l=1}^\infty \Psi_l(t)z^l$ is absolutely convergent for $|z|$ bounded, defining a holomorphic matrix-valued function in $(z,t)$  on $\{|z|<r\}\times \mathcal{V}$, for any $r>0$.

\vskip 0.2 cm

\vskip 0.2 cm 
D] If some $\mu_i(t)-\mu_j(t)\in\mathbb{Z}\backslash\{0\}$ and $J^{(0)}(t)$ is not diagonal, 
then there exists a fundamental matrix holomorphic on $\mathcal{V}$, 
\be
\label{24gen2016-1BBIS}
Y^{(0)}(z,t)=G^{(0)}(t)\Bigl( I+\sum_{l=1}^\infty \Psi_l(t)z^l\Bigr) z^{D^{(0)}} z^{L^{(0)}(t)},
\ee
where 
$$ 
L^{(0)}(t):=S^{(0)}(t)+R^{(0)}(t), 
$$ 
 $G^{(0)}$, $S^{(0)}$ are holomorphic on $\mathcal{U}_{\epsilon_0}(0)$, and $R^{(0)}$ and the $\Psi_l$'s can be chosen holomorphic on $\mathcal{V}$.   The series $ I+\sum_{l=1}^\infty \Psi_l(t)z^l$ is absolutely convergent for $|z|$ bounded, defining a holomorphic matrix-valued function in $(z,t)$  on $\{|z|<r\}\times \mathcal{V}$, for any $r>0$.
 
 The structure of $R^{(0)}$ is more conveniently described if the eigenvalues $\mu_1,\mu_2,...,\mu_n$ are re-labelled as follows. Up to a permutation $J^{(0)}\mapsto P^{-1}J^{(0)}P$, which  corresponds to $G^{(0)}\mapsto G^{(0)}P$, where $P$ is a permutation matrix,  the Jordan blocks structure can be arranged as 
\be
\label{9feb2016-1}
J^{(0)}=J^{(0)}_1\oplus\cdots\oplus J^{(0)}_{s_0},
\quad
\quad
s_0\leq n.
\ee
 For $i=1,2,...,s_0$, each $J^{(0)}_i$ has dimension $n_i$ (then $n_1+\cdots+n_{s_0}=n$) and has only one eigenvalue $\widetilde{\mu}_i$, with structure 
\be
\label{26maggio2017-7}
 J_i^{(0)}=\widetilde{\mu}_i I_{n_i} +H_{n_i},
 \quad
 \quad
 I_{n_i}=\hbox{ $n_i\times n_i$ identity matrix},
 \ee
  $$ \hbox{\rm $H_{n_i}=0$ ~if $n_i=1$},
  \quad
  \quad
 H_{n_i}=\left[
 \begin{array}{cccccc}
 0& 1 &  &  \cr
                   &0 &1      &        & &\cr 
                    & & \ddots & \ddots & &\cr
                    &          &         & 0 &1\cr
                    &           &         &  & 0 
 \end{array}
\right]~\hbox{  if $n_i\geq 2$}.
$$
$\widetilde{\mu}_1$, ..., $\widetilde{\mu}_{s_0}$ are not necessarily distinct.    Let us partition $R^{(0)}$ according to the block structure $n_1$, ..., $n_{s_0}$. Then $[R^{(0)}]_{block~i,j}\neq 0$ only if $\widetilde{\mu}_i-\widetilde{\mu}_j\in\mathbb{N}\backslash\{0\}$, for $1\leq i\neq j \leq s_0$.

\bre
\label{17april2016-3}
Also in cases A], B] and C] the fundamental solution can be written in the Levelt form  (\ref{24gen2016-1BBIS}), with $L^{(0)}=S^{(0)}$ in A] and B], and $L^{(0)}=S^{(0)}+R^{(0)}$ in C]. 
\ere

 \subsection{Freedom} Let the matrix $J^{(0)}(t)$ be fixed with the convention (\ref{9feb2016-1}). Let Assumptions 2 and 3 hold. The class of normal forms at the Fuchsian singularity $z=0$ with given $J^{(0)}$ is not unique, when some eigenvalues of  $\widehat{A}_1(t)$ differ by non-zero integers. Let $\kappa$ be the maximal integer difference. Then, if (\ref{24gen2016-1BBIS}) is a Levelt form,  there are other Levelt forms 
\beas
\widetilde{Y}^{(0)}(z,t)&&=\widetilde{G}^{(0)}(t) \Bigl(I+\sum_{l=1}^\infty\widetilde{\Psi}_l(t) z^l\Bigr) z^{D^{(0)}(t)}z^{\widetilde{L}^{(0)}(t)}
\\
&&\equiv Y^{(0)}(z,t) \mathfrak{D}(t),
\eeas
where $\mathfrak{D}(t)$ is a connection matrix. From the standard theory of equivalence of Birkhoff normal forms of a given differential system with Fuchsian singularity, it follows that $\mathfrak{D}(t)$ must have the following property 
$$
 z^{D^{(0)}(t)}z^{L^{(0)}(t)} \mathfrak{D}(t)=\mathfrak{D}_0(t)\Bigl(I+\mathfrak{D}_1(t)z+\cdots+\mathfrak{D}_\kappa(t) z^\kappa\Bigr)z^{D^{(0)}(t)}z^{\widetilde{L}^{(0)}(t)},
$$
being  $\mathfrak{D}_0$, ..., $\mathfrak{D}_\kappa$  {\it arbitrary} matrices satisfying   $[\mathfrak{D}_0,J^{(0)}]=0$, $\mathfrak{D}_{ij}^{(l)}\neq 0$ only if $\widetilde{\mu}_i-\widetilde{\mu}_j=l>0$.   The connection matrix is then 
 $$
\mathfrak{D}(t)= \mathfrak{D}_0(t)\Bigl(I+\mathfrak{D}_1(t)+\cdots+\mathfrak{D}_k(t)\Bigr).
$$
Being $\mathfrak{D}_0(t)$, ...,  $\mathfrak{D}_\kappa(t)$  arbitrary, we can choose the subclass of those connection matrices  $\mathfrak{D}(t)$ which are holomorphic in $t$. Note that $\mathfrak{D}_0$ commutes with $D^{(0)}$. 
The relation between matrices with $\widetilde{\quad}$ and without is as follows:
\beas
\widetilde{G}^{(0)}(t)\Bigl(I+\sum_{l=1}^\infty\widetilde{\Psi}_l(t) z^l\Bigr)=
\\
=G^{(0)}(t)\Bigl(I+\sum_{l=1}^\infty\Psi_l(t) z^l\Bigr) ~\Bigl[\mathfrak{D}_0(t)\Bigl(I+\mathfrak{D}_1(t)z+\cdots+\mathfrak{D}_\kappa(t) z^\kappa\Bigr)\Bigr].
\eeas
 Moreover, 
\be
\label{24gen2016-3}
\widetilde{L}^{(0)}=\mathfrak{D}^{-1} L^{(0)} \mathfrak{D},
\quad
\quad
\widetilde{R}^{(0)}=\mathfrak{D}^{-1}R^{(0)}\mathfrak{D}+\mathfrak{D}^{-1}[S^{(0)},\mathfrak{D}].
\ee
Observe  that 
$$\widetilde{G}^{(0)}(t)={G}^{(0)}(t)
\quad
~\Longleftrightarrow
\quad
~\mathfrak{D}_0(t)=I.$$

\section{Definition of Isomonodromy Deformation of the System (\ref{ourcase12}) with Eigenvalues (\ref{4gen2016-2})}
\label{5gen2016-2-bis}
\vskip 0.2 cm 
The Stokes phenomenon at $z=\infty$ has been already described. 

Let $\widetilde{\tau}$ be an admissible direction for $\Lambda(0)$.
 For the remaining part of the paper,  $\mathcal{V}$ will denote   an open simply  connected subset of a $\widetilde{\tau}$-cell, such that the closure  $\overline{\mathcal{V}}$ is also   contained in the cell.   Let  the label $\nu$   satisfy $\tau_\nu<\widetilde{\tau}<\tau_{\nu+1}$.  The holomorphic fundamental matrices of Section \ref{14feb2016-9}, namely $Y_\sigma(z,t)$, $\sigma=\nu,\nu+\mu,\nu+2\mu$, exist and satisfy Corollary \ref{20marzo2016-5} and Proposition  \ref{22marzo2016-2}. Therefore, in particular, they have canonical asymptotics on $\mathcal{S}_\sigma(\overline{\mathcal{V}})$,  with holomorphic on $\mathcal{V}$ Stokes matrices $\mathbb{S}_\nu(t)$ and $\mathbb{S}_{\nu+\mu}(t)$.

\bre[{\bf Notations}]
\label{4gen2017-2}
The notation  $Y_\nu(z,t)$ of Sections \ref{14feb2016-9}-\ref{2dic2016-2}  has been used  for the fundamental matrix solutions of the system (\ref{22novembre2016-3}),  (\ref{9giugno-2}). 
 We  consider now the system (\ref{ourcase12}) and  use the same notation  $Y_\nu(z,t)$, with the  replacement  $G_0(t) \mapsto I$ in all the formulae where $G_0(t)$ appears. 
\ere
 
\bde 
\label{13giugno2017-7}
The {\bf  central connection matrix} $C^{(0)}_\nu(t)$ is defined {by}
$$ 
Y_\nu(z,t)=Y^{(0)}(z,t)C^{(0)}_\nu(t),
\quad
\quad
z\in\mathcal{R}.
$$ 
\ede

\bde[Isomonodromic Deformation  in  $\mathcal{V}$]
\label{30dic2015-1}
 Let $\mathcal{V}$ be an open connected subset of  a $\widetilde{\tau}$-cell, such that $\overline{\mathcal{V}}$ is also contained the cell.  A  $t$-deformation  of the system (\ref{ourcase12}) satisfying Assumption 2 in $\mathcal{V}$ is said to be  {\bf isomonodromic in  $\mathcal{V}$}  if the {\it essential monodromy data}, $$ 
 \mathbb{S}_\nu
 \quad
 \mathbb{S}_{\nu+\mu},
 \quad
 B_1=\hbox{\rm diag}(\widehat{A}_1);
 \quad
 \quad
 \{\mu_1,~...,~\mu_n
\},
 $$ 
 are independent of $t\in \mathcal{V}$, and if there exists a fundamental solution (\ref{24gen2016-1BBIS}) (see Remark \ref{17april2016-3}), holomorphic in $t\in \mathcal{V}$, such that also the corresponding essential monodromy data 
$$
R^{(0)},
\quad
C^{(0)}_\nu,
$$
 are independent of $t\in \mathcal{V}$. 
\ede

\bre If $\mu_1$, ..., $\mu_n$ are independent of $t$ as in Definition \ref{30dic2015-1}, then 
{\it Assumption 2 in $B$ implies that also Assumption 3 holds in $\mathcal{V}$}. 
\ere

The existence of a fundamental solution with constant $R^{(0)}$ implies that  the system (\ref{ourcase12}) can be reduced to  a simpler form (\ref{26maggio2017-6}) which is     independent of $t\in \mathcal{V}$, namely 
\be
\label{9feb2016-2}
\frac{d \mathcal{Y}}{dz}=\frac{1}{z}
\left(
J^{(0)}+R_1 z+\cdots +R_\kappa z^\kappa \right)\mathcal{Y},
\ee
where  $1 \leq \kappa$ is the maximal integer difference of eigenvalues of $J^{(0)}$, $[R_l]_{block ~i,j}\neq 0$ only if $\widetilde{\mu}_i-\widetilde{\mu}_j=l$,  $R_1+\cdots+R_\kappa=R^{(0)}$, with all $R_l$ independent of $t\in \mathcal{V}$, and the $\widetilde{\mu}_i$'s are the eigenvalues of $\widehat{A}_1(t)$ as arranged in the  Jordan from (\ref{9feb2016-1})-(\ref{26maggio2017-7}).

\bre 
\label{24gen2016-2}
There is a freedom in the isomonodromic $R^{(0)}$ and $L^{(0)}$, as in (\ref{24gen2016-3}), for a $t$-independent  $\mathfrak{D}$ such that $\widetilde{Y}^{(0)}=Y^{(0)}\mathfrak{D}$. Hence,  there is  {a} freedom in the isomonodromic central connection matrix, according to  $$C^{(0)}=\mathfrak{D}\widetilde{C}^{(0)}.$$  We call $\mathcal{C}_0(J^{(0)},L^{(0)})$ the group of such transformations $\mathfrak{D}$ which leave $L^{(0)}$ invariant in (\ref{24gen2016-3}). This notation is a slight variation of a notation introduced in \cite{Dub2} for a particular subclass of our systems (\ref{ourcase12}), related to Frobenius manifolds.
\ere

\bre 
\label{16gen2016-2}
Definition \ref{30dic2015-1} is given with reference to some $\nu$. Nevertheless, it implies that it holds for any other $\nu^\prime$ in a suitably small $\mathcal{V}^\prime\subset \mathcal{V}$. To see this, consider another admissible  $\widetilde{\tau}^\prime\in (\tau_{\nu^\prime},\tau_{\nu^\prime+1})$, and define $\mathcal{S}_{\nu^\prime+k\mu}(t)$, $Y_{\nu^\prime+k\mu}(z,t)$ in the usual way,  for $t$ in the intersection of $\mathcal{V}$ with a $\widetilde{\tau}^\prime$-cell. 
\footnote{Note that there may be more than one choices for  $\mathcal{S}_{\nu^\prime+k\mu}$, $Y_{\nu^\prime+k\mu}(z,t)$, depending on the neighbourhood of $t$ considered. See Remark \ref{16gen2016-1}.}  Call $\mathcal{V}^\prime$ the intersection. 
 Now, there is a finite product  of {\it Stokes factors} $K_1(t)\cdots K_M(t)$ ($M\leq $ number of basic Stokes rays of
 $\Lambda(t)$) such that  $Y_{\nu}(z,t)= Y_{\nu^\prime}(z,t)K_1(t)\cdots K_M(t)$, $t\in \mathcal{V}^\prime$. The Stokes matrices $\mathbb{S}_\nu(t)$ and
 $\mathbb{S}_{\nu+\mu}(t)$ are determined uniquely by their factors, and conversely   a Stokes matrix determines uniquely 
the factors of  a factorization of {the} prescribed structure (see the proof of Theorem \ref{8gen2017-4}, or section 4 of  \cite{BJL1}, point D). Moreover,  the product $ K_1(t)\cdots K_M(t) $
 appears  in the factorization of  $\mathbb{S}_{\nu}$ or  $\mathbb{S}_{\nu+\mu}$. Hence, if $\mathbb{S}_{\nu}$ and $\mathbb{S}_{\nu+\mu}$ do not depend
 on $t\in \mathcal{V}$ for a certain $\nu$, also  $\mathbb{S}_{\nu^\prime}$ and $\mathbb{S}_{\nu^\prime+\mu}$  do not depend
 on $t\in \mathcal{V}^\prime \subset \mathcal{V}$. Thus, the same is true for $C^{(0)}_{\nu^\prime}$.
\ere

\ble
\label{7gen2016-6}
Let the deformation be isomonodromic in $\mathcal{V}$ as in Definition \ref{30dic2015-1} (here it is not necessary to suppose that $\mathcal{V}$ is in a cell, since we are considering solutions at  $z=0$). Let Assumption 2 hold in $\mathcal{U}_{\epsilon_0}(0)$, namely   let $\widehat{A}_1(t)$ be holomorphically equivalent to $J^{(0)}$ in $\mathcal{U}_{\epsilon_0}(0)$. Then: 
\begin{itemize}
\item[i)] $\mu_1,...,\mu_n$, $D^{(0)}$ , $S^{(0)}$ and $J^{(0)}$ are independent of $t$ in  $\mathcal{U}_{\epsilon_0}(0)$.

\item[ii)] Any fundamental matrix (also non-{isomonodromic} ones) in Levelt form $Y^{(0)}(z,t)=G^{(0)}(t)(I+\sum_{l}\Psi_l(t)z^l)z^Dz^{L^{(0)}(t)}$, which is  holomorphic of $t\in \mathcal{V}$,   is also holomorphic on the whole  $ \mathcal{U}_{\epsilon_0}(0)$. 

\item[iii)] If $R^{(0)}$ (i.e $L^{(0)}$) is independent of $t$ in $\mathcal{V}$, then it is independent of $t$ in $\mathcal{U}_{\epsilon_0}(0)$. 
\end{itemize}
\ele

\vskip 0.2 cm 
\noindent
{\it Proof:}  i) That $\mu_1,...,\mu_n$, $D^{(0)}$, $S^{(0)}$,  $J^{(0)}$ are constant in $\mathcal{U}_{\epsilon_0}(0)$ follows from the fact that $\mu_1,...,\mu_n$ are constant in $\mathcal{V}$, and that   $G^{(0)}(t)$, and so the $\mu_1,...,\mu_n$, are holomorphic on $\mathcal{U}_{\epsilon_0}(0)$. So $\mu_1,...,\mu_n$ are constant in $\mathcal{U}_{\epsilon_0}(0)$.  

ii) Since   $\mu_1,...,\mu_n$ are constant in $\mathcal{U}_{\epsilon_0}(0)$, and $\Lambda(t)$ and $\widehat{A}_1(t)$ are holomorphic,  the recursive standard procedure which yields the Birkhoff normal form at $z=0$ allows to choose  $\Psi_l(t)$'s and $R^{(0)}(t)$ holomorphic on $\mathcal{U}_{\epsilon_0}(0)$. 

iii) That  $R^{(0)}$  is independent of $t\in \mathcal{U}_{\epsilon_0}(0)$ follows from  the fact that $R^{(0)}(t)$ is holomoprhic on $\mathcal{U}_{\epsilon_0}(0)$ and constant on $\mathcal{V}$. 
$\Box$

\begin{shaded}
\bpr
\label{8feb2016-5}
Let the deformation of the system (\ref{ourcase12}) be isomonodromic in $\mathcal{V}$ as in Definition \ref{30dic2015-1} (here {it} is not necessary to assume that $\mathcal{V}$ is contained in a cell).  Let Assumption 2 hold in $\mathcal{U}_{\epsilon_0}(0)$, namely  let $\widehat{A}_1(t)$ be holomorphically equivalent to $J^{(0)}=D^{(0)}+S^{(0)}$ in $\mathcal{U}_{\epsilon_0}(0)$.  Consider the system 
\be
\label{7gen2016-7}
\frac{dY}{dz}=\widehat{A}(z,0) Y,
\ee
and a fundamental solution in {the} Levelt form 
\be
\label{17dicembre2016-7}
\mathring{Y}^{(0)}(z)=\mathring{G}^{(0)}\mathring{G}(z) z^{D^{(0)}} z^{\mathring{L}},
\quad
\quad
\mathring{G}(z)=I+\mathcal{O}(z),
\ee
with $\mathring{L}=S^{(0)}+\mathring{R}$. Here $\mathring{R}$ is obtained by reducing  (\ref{7gen2016-7})  to a Birkhoff normal form at $z=0$. Then, there exists  an isomonodromic fundamental solution of (\ref{ourcase12}), call it $
Y^{(0)}_{isom}(z,t)$, with the same monodromy exponent $\mathring{L}$ and Levelt form 
$$
Y_{isom}^{(0)}(z,t)=G^{(0)}(t)G_{isom}(z,t)z^{D^{(0)}}z^{\mathring{L}},
$$ 
with $G_{isom}(z,t)=I+\sum_{k=1}^\infty \Psi_l(t)z^l$, holomorphic on $\mathcal{U}_{\epsilon_0}(0)$, 
such that 
$$ 
\mathring{Y}^{(0)}(z)=Y_{isom}^{(0)}(z,0).
$$
\epr
\end{shaded}

\vskip 0.2 cm 
\noindent
{\it Proof:} We prove the proposition in two steps.

 $\bullet$ The first step is  the following 

\ble
\label{8gen2016-3}
  Let the deformation be isomonodromic in $\mathcal{V}$ as in Definition \ref{30dic2015-1}  (here {it} is not necessary to assume that $\mathcal{V}$ is contained in a cell).  Let $\widehat{A}_1(t)$ be holomorphically equivalent to $J^{(0)}$ in $\mathcal{U}_{\epsilon_0}(0)$.  For any holomorphic fundamental solution in Levelt form 
$$ 
Y(z,t)=G(t) H(z,t) z^{D^{(0)}} z^{L^{(0)}(t)},
\quad
\quad
H(z,t)=I+\sum_{l=1}^\infty h_l(t)z^l,
$$
with monodromy exponent $L^{(0)}(t)$, there exists {an} isomonodromic $Y^{(0)}(z,t)$,  with  monodromy exponent equal to  $L^{(0)}(0)$, in the Levelt form 
$$ 
Y^{(0)}(z,t)=G^{(0)}(t) G(z,t) z^{D^{(0)}} z^{L^{(0)}(0)},
\quad\quad
G(z,t)=I+\sum_{l=1}^\infty \Psi_l(t)z^l,
$$
 such that 
$ 
Y^{(0)}(z,0)=Y(z,0)$. 
\ele

To prove this Lemma, consider an isomonodromic fundamental solution, which exists by assumption, say
$$ 
\widetilde{Y}^{(0)}(z,t)=\widetilde{G}^{(0)}(t) \widetilde{G}(z,t) z^{D^{(0)}} z^{\widetilde{L}^{(0)}},
\quad\quad
\widetilde{G}(z,t)=I+\sum_{l=1}^\infty \widetilde{\Psi}_l(t)z^l, 
$$ 
with $t$-independent monodromy exponent  $\widetilde{L}^{(0)}$ and $t$-independent connection matrix  defined by 
$$Y_\nu(z,t)=\widetilde{Y}^{(0)}(z,t) ~\widetilde{C}_\nu^{(0)}.
$$  Then, there exists a holomorphic invertible connection matrix  $\mathfrak{D}(t)$ such that 
$$ 
Y(z,t)=\widetilde{Y}^{(0)}(z,t)\mathfrak{D}(t).
$$
Hence, 
\be
\label{8gen2016-1}
\mathfrak{D}_0(t)\Bigl(I+\mathfrak{D}_1(t)z+\cdots+\mathfrak{D}_\kappa(t)z^\kappa\Bigr)    
z^{D^{(0)}}z^{L^{(0)}(t)}=z^{D^{(0)}}
z^{\widetilde{L}^{(0)}}~\mathfrak{D}(t)
\ee
with 
$
\mathfrak{D}(t)=\mathfrak{D}_0(t)\Bigl(I+\mathfrak{D}_1(t)+\cdots+\mathfrak{D}_\kappa(t)\Bigr)
$. 
Observe that $z^{D^{(0)}}z^{L^{(0)}(0)}$ and $z^{D^{(0)}}
z^{\widetilde{L}^{(0)}}$ are fundamental solutions of two Birkhoff normal forms of (\ref{7gen2016-7}), related by (\ref{8gen2016-1}) with $t=0$, namely
$$
\mathfrak{D}_0(0)(I+\mathfrak{D}_1(0)z+\cdots + \mathfrak{D}_\kappa(0) z^\kappa)
z^{D^{(0)}}z^{L^{(0)}(0)}
=z^{D^{(0)}}z^{\widetilde{L}^{(0)}} ~\mathfrak{D}(0).
$$ 
Therefore, the isomonodromic  fundamental solution we are looking for is 
$$ 
Y^{(0)}(z,t):=\widetilde{Y}^{(0)}(z,t)~\mathfrak{D}(0)= Y(z,t)\mathfrak{D}(t)^{-1}\mathfrak{D}(0).
$$
 
$\bullet$ Second step.  Consider a fundamental solution of (\ref{7gen2016-7}) in  {the} Levelt form
$$ 
\mathring{Y}^{(0)}(z)=\mathring{G}^{(0)}\mathring{G}(z) z^{D^{(0)}} z^{\mathring{L}},
$$ 
where $\mathring{L}=S^{(0)}+\mathring{R}$,  $\mathring{R}=\sum_{l=1}^\kappa \mathring{R}_l$. The  $\mathring{R}_l$, $l=1,2,...,\kappa$, are coefficients of a simple gauge equivalent 
form(\ref{26maggio2017-6}), with $t=0$,  of (\ref{7gen2016-7}). It can be proved that there is a form (\ref{26maggio2017-6})  for the system  (\ref{ourcase12}), with coefficients $R_l(t)$,  such that  the $\mathring{R}_l$'s  coincide with the values $R_l(0)$'s at $t=0$. Moreover, the  $R_l(t)$'s  are holomorphic on $\mathcal{U}_{\epsilon_0}(0)$. This  fact follows from the recursive procedure which yileds the gauge transformation from  (\ref{ourcase12}) to (\ref{26maggio2017-6}).  Therefore,  there exists a holomorphic exponent $L^{(0)}(t)$ such that $L^{(0)}(0)=\mathring{L}$. Consider an isomonodromic  fundamental solution $Y^{(0)}(z,t)$ of Lemma \ref{8gen2016-3}, with exponent $L^{(0)}(0)=\mathring{L}$. Since  $Y^{(0)}(z,0)$ is a fundamental solution of (\ref{7gen2016-7}), there exists an invertible and constant connection matrix $C$ such that 
$$ 
Y^{(0)}(z,0)C=\mathring{Y}^{(0)}(z).
$$ 
Now,  $C\in\mathcal{C}_0(J^{(0)},\mathring{L})$ (cf.  Remark \ref{24gen2016-2}), because $Y^{(0)}(z,0)$ and $\mathring{Y}^{(0)}(z)$ have the same monodromy exponent. This implies that
\beas
Y^{(0)}(z,t)C&&
=G^{(0)}G(z,t)z^{D^{(0)}} z^{\mathring{L}}C=
\\
&&
=G^{(0)}G(z,t)C_0(I+C_1z+\cdots+C_\kappa z^\kappa)z^{D^{(0)}} z^{\mathring{L}},
\quad
\quad
C=C_0(I+C_1+\cdots+C_\kappa).
\eeas
Moreover, also $Y^{(0)}(z,t)C$ is  isomonodromic. Therefore, the solution we are looking for is 
 $Y^{(0)}_{isom}(z,t):=Y^{(0)}(z,t)C$. $\Box$

\section{Isomonodromy Deformation Equations}
Let 
$$
\Omega(z,t):=\sum_{k=1}^n \Omega_k(z,t) ~dt_k,
\quad
\quad
\Omega_k(z,t):=zE_k+[F_1(t),E_k].
$$
 Here $E_k$ is the matrix with all entries equal to zero, except for $(E_k)_{kk}=1$, and $(F_1)_{ab}=-(\widehat{A}_1)_{ab}/(u_a-u_b)$, so that 
\be
\label{3nov2015-2}
[F_1(t),E_k]=\left(
 \frac{(\widehat{A}_1(t))_{ab}(\delta_{ak}-\delta_{bk})}{u_a(t)-u_b(t)}
 \right)_{a,b=1..n}=\left(
\begin{array}{ccccc}
0 & 0&\frac{ -(\widehat{A}_1)_{1k}}{u_1-u_k} &  0& 0
 \\
0 &0 &\vdots & 0& 0
 \\
 \frac{ (\widehat{A}_1)_{k1}}{u_k-u_1} & \cdots&0 & \cdots&\frac{ (\widehat{A}_1)_{kn}}{u_k-u_n}
 \\
0 &0 &\vdots &0  & 0
 \\
 0 &0 &\frac{ -(\widehat{A}_1)_{nk}}{u_n-u_k} &  0& 0
\end{array}
\right)
\ee
Let $df(z,t):=\sum_{i=1}^n \partial f(z,t)/\partial t_i dt_i$. 
\begin{shaded}
\bth 
\label{8luglio2016-1}
If the deformation of the system (\ref{ourcase12}) is isomonodromic in $\mathcal{V}$ as in Definition \ref{30dic2015-1}, then an isomonodromic $Y^{(0)}(z,t)$ and  the $Y_\sigma(z,t)$'s, for $\sigma=\nu,\nu+\mu,\nu+2\mu$, satisfy the total differential system 
\be
\label{3nov2015-3}
dY=\Omega(z,t)Y.
\ee
Conversely, if the $t$-deformation satisfies assumptions 2 and  3 in $\mathcal{V}$, and if a fundamental solution $Y^{(0)}(z,t)$ in Levelt form at $z=0$, and the canonical solution $Y_\sigma(z,t)$, $\sigma=\nu,\nu+\mu,\nu+2\mu$ at $z=\infty$, satisfy the total differential system (\ref{3nov2015-3}), then the deformation is isomonodromic in $\mathcal{V}$. 
\eth
\end{shaded}

\noindent
{\it Proof:}  The proof is done in the same way as for Theorem 3.1 at page 322 in \cite{JMU}. In  \cite{JMU} the proof is given for non resonant $\widehat{A}_1(t)$, but it can be repeated in our case with no changes, except for the assumptions 2, 3. \footnote{The result was announced in  \cite{NBOKA} and not proved. It can also be proved by the methods of \cite{KV}, since the requirement that  $\mu_1,...,\mu_n$, $R^{(0)}$ and $C^{(0)}$ are constant is equivalent to having an isoprincipal deformation.}  
 The matrix valued differential form $\Omega(z,t)$ turns out to be still  as in formula (3.8) and (3.14)  of \cite{JMU}, which in our case  becomes, 
 $$ 
 \Omega(z,t)=\left[\left(I+\sum_{k=1}^\infty F_k(t)z^{-k}\right)~d\Lambda(t)z~
 \left(I+\sum_{k=1}^\infty F_k(t)z^{-k}\right)^{-1}\right]_{sing},
 $$ 
 where $[\cdots]_{sing}$ stands for the singular terms at infinity, namely the terms with powers $z^j$, $j\geq 0$, in the above formal expansion. This is 
 $$ 
 \Omega(z,t)=d\Lambda(t)z+[F_1(t),d\Lambda(t)].
 $$
 Therefore, 
 $$ 
 \Omega_k(z,t)= \frac{\partial \Lambda(t)}{\partial t_k}~z~+\left[
 F_1(t),\frac{\partial \Lambda(t)}{\partial t_k}
 \right]~= E_k+[F_1(t),E_k].
 $$
In the last step we have used the fact that  
$ 
\Lambda(t)=\hbox{ diag} (u_1(t),u_2(t),~...,~u_n(t))
$, 
 with eigenvalues  (\ref{4gen2016-2}). 
 In the domain $\mathcal{V}$ the eigenvalues are distinct,  so the off-diagonal entries of $F_1$ are: 
 $$
 (F_1)_{ab}=\frac{(\widehat{A}_1)_{ab}}{u_b-u_a},
 \quad
 \quad
 1\leq a\neq b\leq n.
 $$ 
Hence, 
$$
\Omega_k(z,t)=E_k~z~+
\left(
\frac{\widehat{A}^{(1)}_{ab}}{u_{b}(t)-u_{a}(t)}~
 \frac{\partial }{\partial t_k}\Bigl( u_{b}(t)- u_{a}(t)\Bigr)
\right)_{a,b=1}^n.
$$
Finally, observe that  $\frac{\partial }{\partial t_k}( u_{b}(t)- u_{a}(t))=\frac{\partial }{\partial t_k}( t_{b}- t_{a})=\delta_{kb}-\delta_{ka}$. The proof is concluded. 
$\Box$

\begin{shaded}
\bcr
\label{17april2016-4}
If the deformation of the system (\ref{ourcase12}) is isomonodromic in $\mathcal{V}$ as in Definition \ref{30dic2015-1}, then $G^{(0)}(t)$ satisfies 
\be
\label{28dic2015-1}
dG^{(0)}=\Theta^{(0)}(t)~G^{(0)},
\ee
where
$$ 
\Theta^{(0)}(t)=\Omega(0,t)=\sum_k[F_1(t),E_k]dt_k.
$$
More explicitly, 
$$ 
\Theta^{(0)}(t)=\left(
\frac{\widehat{A}^{(1)}_{ab}}{u_{a}(t)-u_{b}(t)}~(dt_a-dt_b)
\right)_{a,b=1}^n.
$$
\ecr
\end{shaded}

\noindent
{\it Proof:} Substitute $Y^{(0)}$ into (\ref{3nov2015-3}) an  compare coefficients of equal powers of $z$. Equation (\ref{28dic2015-1}) comes form the coefficient of $z^0$. $\Box$

\begin{shaded}
\bpr
If the deformation is isomonodromic in $\mathcal{V}$  as in Definition \ref{30dic2015-1}, then 
\be
\label{30dic2015-3}
d\widehat{A}=\frac{\partial \Omega}{\partial z} +[\Omega,\widehat{A}].
\ee 
\epr
\end{shaded}

\noindent
{\it Proof:} 
Let the deformation be isomonodromic. Then, by Theorem \ref{8luglio2016-1},  equations (\ref{ourcase12}) and (\ref{3nov2015-3}) are compatible. The compatibility condition  is (\ref{30dic2015-3}).
$\Box$

\vskip 0.2 cm 
Note that  (\ref{30dic2015-3})  is a necessary condition of isomonodromicity, but not sufficient in case of resonances (sufficiency can be proved if the eigenvalues of $\widehat{A}_1$ do not differ by integers, cf. \cite{JMU}). Explicitly, (\ref{30dic2015-3}) is 
$$
\left\{
\begin{array}{cc}
[E_k,\widehat{A}_1]=[\Lambda,[F_1,E_k]],& k=1,...,n,
\\
\\
d\widehat{A}_1=[\Theta^{(0)},\widehat{A}_1].&
\end{array}
\right.
$$
The first $n$ equations are automatically satisfied by definition of $F_1$. The last equation  in components is
\be
\label{23gen2016-4}
\frac{\partial \widehat{A}_1}{\partial t_k}=\Bigl[  [F_1,E_k],\widehat{A}_1\Bigr],
\ee
where $[F_1,E_k]$ is in (\ref{3nov2015-2}).

\section{Holomorphic Extension of Isomonodromy Deformations to $\mathcal{U}_{\epsilon_0}(0)$ and Theorem \ref{16dicembre2016-1}}
\label{17luglio2016-2}

\begin{shaded}
\ble
\label{21luglio-1}
 In case the eigenvalues of $\Lambda(t)$ are as in (\ref{4gen2016-2}) and $\widehat{A}_1(t)$ is holomorphic on $\mathcal{U}_{\epsilon_0}(0)$, then $\Omega(z,t)$  is holomoprhic ({in} $t$) on  $\mathcal{U}_{\epsilon_0}(0)$ if and only if 
\be 
\label{30giugno-3}
(\widehat{A}_1)_{ab}(t)=\mathcal{O}(u_a(t)-u_b(t))\equiv \mathcal{O}(t_a-t_b),
\ee
whenever $u_a(t)$ and $u_b(t)$ coalesce as $t$ tends to a point of $\Delta\subset \mathcal{U}_{\epsilon_0}(0)$. 

\noindent
Also  $\Theta^{(0)}(t)$ of Corollary \ref{17april2016-4} is holomorphic on  $\mathcal{U}_{\epsilon_0}(0)$ if and only if (\ref{30giugno-3}) holds. 
\ele
\end{shaded}

\vskip 0.2 cm 
\noindent
{\it Proof:}
By (\ref{3nov2015-2}),  $\Omega(z,t)$ and $\Theta^{(0)}(t)$ are  continuous at $t_\Delta\in\Delta$  if and only if (\ref{30giugno-3}) holds 
for those  $u_a(t)$, $u_b(t)$ coalescing at $t_\Delta\in \Delta$. Hence, any point of $\Delta$ is a removable singularity if and only if (\ref{30giugno-3}) holds. 
$\Box$

\bpr The system
\beas
&&
(\ref{30dic2015-3})
\quad\quad
d\widehat{A}=\frac{\partial \Omega}{\partial z} +[\Omega,\widehat{A}],
 \\
&&  (\ref{28dic2015-1}) 
\quad\quad
 dG^{(0)}=\Theta^{(0)}(t)~G^{(0)},
 \eeas 
with $\widehat{A}_1$ holomorphic satisfying condition (\ref{30giugno-3}) on $\mathcal{U}_{\epsilon_0}(0)$,  is Frobenius integrable for  $t\in \mathcal{U}_{\epsilon_0}(0)$.
 \epr 
\vskip 0.2 cm 
\noindent
The proof is as in  \cite{JMU}. It holds also in our case, because the algebraic relations are the same as in our case, no matter if  $\widehat{A}_1$ is resonant (see e.g. Example 3.2 in \cite{JMU}). 

Write $\Theta^{(0)}=\sum_k \Theta_k^{(0)} dt_k$.  Since  (\ref{28dic2015-1})  is integrable, the compatibility of equations holds: 
\be
\label{30dic2015-2}
\frac{\partial \Theta^{(0)}_j}{\partial t_i}-\frac{\partial \Theta^{(0)}_i}{\partial t_j}=
\Theta^{(0)}_i\Theta^{(0)}_j-\Theta^{(0)}_j\Theta^{(0)}_i.
\ee

\begin{shaded}
\bpr
\label{28dic2015-4}
 Let the deformation of the system (\ref{ourcase12}) be isomonodromic  in $\mathcal{V}$  as in Definition \ref{30dic2015-1}, with $\Lambda(t)$ is as in (\ref{4gen2016-2}) and $\widehat{A}_1(t)$ holomorphic on $\mathcal{U}_{\epsilon_0}(0)$. Then, $\widehat{A}_1(t)$ is holomorphically similar to $J^{(0)}$ in the whole $\mathcal{U}_{\epsilon_0}(0)$ if and only if  (\ref{30giugno-3}) holds as $t$ tends to points of   $\Delta\subset \mathcal{U}_{\epsilon_0}(0)$.  In other words, if the deformation is isomonodromic in $\mathcal{V}$ with holomorphic  $\widehat{A}_1(t)$, then Assumption 2 in the whole $\mathcal{U}_{\epsilon_0}(0)$ is equivalent to (\ref{30giugno-3}).
\epr
\end{shaded}

\noindent
{\it Proof:}  Let  $\widehat{A}_1(t)$ be holomorphic and let  (\ref{30giugno-3}) hold, so that $\Theta^{(0)}(t)$ is holomorphic on  $\mathcal{U}_{\epsilon_0}(0)$ by Lemma \ref{21luglio-1}. The {\it linear} Pfaffian systems $dG^{(0)}=\Theta^{(0)}(t)G^{(0)}$ and $d[(G^{(0)})^{-1}]=-(G^{(0)})^{-1} \Theta^{(0)}(t)$ are integrable in $\mathcal{U}_{\epsilon_0}(0)$,  with holomorphic coefficients $\Theta^{(0)}(t)$. Then, a solution $G^{(0)}(t)$ has analytic continuation {onto} $\mathcal{U}_{\epsilon_0}(0)$. We take a solution satisfying $(G^{(0)}(t))^{-1}$ $ \widehat{A}_1(t)$ $G^{(0)}(t)=J^{(0)}$ for $t\in \mathcal{V}$,  which then has analytic continuation {onto} $\mathcal{U}_{\epsilon_0}(0)$ as a holomorphic invertible  matrix.  Hence, $(G^{(0)}(t))^{-1}$ $ \widehat{A}_1(t)$ $G^{(0)}(t)=J^{(0)}$ holds in $\mathcal{U}_{\epsilon_0}(0)$ with holomorphic 
$G^{(0)}(t)$. 
 Conversely, suppose that Assumption 2 holds in $\mathcal{U}_{\epsilon_0}(0)$. Then  $G^{(0)}(t)$  and $G^{(0)}(t)^{-1}$ are holomorphic  on   $\mathcal{U}_{\epsilon_0}(0)$. Therefore, also   $\Theta^{(0)}(t)$ is holomorphic on  $\mathcal{U}_{\epsilon_0}(0)$, because $\Theta^{(0)}(t)=dG^{(0)}\cdot (G^{(0)})^{-1}$ defines  the analytic continuation of  $\Theta^{(0)}(t)$ on $\mathcal{U}_{\epsilon_0}(0)$. Then (\ref{30giugno-3}) holds, by Lemma \ref{21luglio-1}.  $\Box$

\vskip 0.2 cm 

Summarising, if $\Lambda(t)$ is as in (\ref{4gen2016-2}) and $\widehat{A}_1(t)$ is holomorphic on $\mathcal{U}_{\epsilon_0}(0)$, if the deformation is isomonodromic in a simply connected subset $\mathcal{V}$ of a cell, s.t. $\overline{\mathcal{V}}\subset $ cell, then it suffices to assume that $\widehat{A}_1(t)$ is  holomorphically similar to a Jordan form $J^{(0)}(t)$ in  $\mathcal{U}_{\epsilon_0}(0)$, or equivalently that (\ref{30giugno-3}) holds at  $\Delta\subset \mathcal{U}_{\epsilon_0}(0)$, in order to conclude that the system
\beas
(\ref{3nov2015-3}) &&
\quad
dY=\Omega(z,t)~Y,
\\
\\
(\ref{28dic2015-1}) &&
\quad
dG^{(0)}=\Theta^{(0)}(t)~ G^{(0)},
\eeas
has holomorphic coefficients on $\mathcal{R}\times \mathcal{U}_{\epsilon_0}(0)$.  
The integrability/compatibility condition of (\ref{3nov2015-3}) is 
\be
\label{8luglio2016-2}
\frac{\partial \Omega_j}{\partial t_i}-\frac{\partial \Omega_i}{\partial t_j}=
\Omega_i\Omega_j-\Omega_j\Omega_i.
\ee
 If this relation is explicitly written, it turns out to be equivalent to (\ref{30dic2015-2}). 
Hence, being (\ref{28dic2015-1}) integrable,  also the linear Pfaffian system (\ref{3nov2015-3}) is integrable, with coefficients holomorphic in $\mathcal{U}_{\epsilon_0}(0)$. 
Therefore, due to linearity, any solution $Y(z,t)$ can be $t$-analytically continued along any curve in $\mathcal{U}_{\epsilon_0}(0)$, for $z$ fixed.  

\begin{shaded}
\bcr 
\label{28dic2015-2}
 Let  the deformation be isomonodromic in a simply connected subset $\mathcal{V}$ of a cell, s.t. $\overline{\mathcal{V}}\subset $ cell. If  $\widehat{A}_1(t)$ is   holomorphically similar to a Jordan form $J^{(0)}$ in   $\mathcal{U}_{\epsilon_0}(0)$, or equivalently if (\ref{30giugno-3}) holds in   $\mathcal{U}_{\epsilon_0}(0)$, then  the 
$Y_\sigma(z,t)$'s, $\sigma=\nu,\nu+\mu,\nu+2\mu$, together with an isomonodromic  $Y^{(0)}(z,t)$, can be $t$-analytically continued as single valued  holomorphic functions on $\mathcal{U}_{\epsilon_0}(0)$. 
\ecr
\end{shaded}

\noindent
{\it Proof:} If the deformation is isomonodromic, by Theorem \ref{8luglio2016-1} the system (\ref{ourcase12}),(\ref{3nov2015-3}) is a completely integrable linear Pfaffian system (compatibility conditions  (\ref{30dic2015-3}) and (\ref{8luglio2016-2}) hold), with common solutions $Y_\sigma(z,t)$'s, $\sigma=\nu,\nu+\mu,\nu+2\mu$,  and $Y^{(0)}(z,t)$.  If  $\widehat{A}_1(t)$ is   holomorphically similar to a Jordan form $J^{(0)}$ in   $\mathcal{U}_{\epsilon_0}(0)$, or equivalently if (\ref{30giugno-3}) holds in   $\mathcal{U}_{\epsilon_0}(0)$, then  the coefficients are holomorphic in $\mathcal{U}_{\epsilon_0}(0)$, by Proposition \ref{28dic2015-4}. In particular, since $Y_\sigma(z,t)$'s, $\sigma=\nu,\nu+\mu,\nu+2\mu$,  and $Y^{(0)}(z,t)$ solve (\ref{3nov2015-3}), they can be $t$-analytically continued along any curve in $\mathcal{U}_{\epsilon_0}(0)$.  
$\Box$.

\bre
\label{8gen2017-1}
Corollary \ref{28dic2015-2} can be compared with the result of \cite{Miwa}. It is always true that   the  $Y_{\sigma}(t,z)$'s and  $Y^{(0)}(t,z)$ {\it can be $t$-analytically continued on 
$\mathcal{T}$ as a meromorphic function}, where (in our case): 
$$ 
\mathcal{T}=\hbox{ universal covering of } \mathbb{C}^n\backslash \Delta_{\mathbb{C}^n}.
$$
Here $\Delta_{\mathbb{C}^n}$ is the  locus of $\mathbb{C}^n$ where eigenvalues of $\Lambda(t)$ coalesce. It is a  locus of  ``fixed singularities" (including branch points and essential singularities) of $\Omega(z,t)$ and of any solution of $dY=\Omega Y$. The movable singularities of $\Omega(z,t)$, $Y_{\sigma}(t,z)$  and  $Y^{(0)}(t,z)$ outside the locus  are  poles and constitute  the zeros of the Jimbo-Miwa isomonodromic $\tau$-function \cite{Miwa}.
Here, we have {furthermore} assumed that $\widehat{A}_1$ is holomorphic in $\mathcal{U}_{\epsilon_0}(0)$ and that (\ref{30giugno-3}) holds. This fact has allowed us to conclude that $Y_\sigma(z,t)$'s, $\sigma=\nu,\nu+\mu,\nu+2\mu$,  and $Y^{(0)}(z,t)$ are $t$-holomorphic in  $\mathcal{U}_{\epsilon_0}(0)$. 
\ere

In order to prove Theorem \ref{16dicembre2016-1}, we need a last ingredient, namely the analyticity at $\Delta$ of the coefficients $F_k(t)$ of the formal solution computed away from $\Delta$. 
\begin{shaded}
\bpr
\label{2gen2016-1}
Let  the deformation of the system (\ref{ourcase12}) be isomonodromic in a simply connected subset $\mathcal{V}$ of a cell, s.t. $\overline{\mathcal{V}}\subset $ cell. If  $\widehat{A}_1(t)$ is   holomorphically similar to a Jordan form $J^{(0)}$ in   $\mathcal{U}_{\epsilon_0}(0)$, or equivalently if (\ref{30giugno-3}) holds in   $\mathcal{U}_{\epsilon_0}(0)$, then   the coefficients $F_k(t)$, $k\geq 1$, of a formal solution of (\ref{ourcase12})
 \be
 \label{21luglio-2}
Y_F(z,t)=\Bigl(I+\sum_{k=1}^\infty F_k(t)z^{-k}\Bigr)z^{B_1}e^{\Lambda(t)z},
\ee 
are holomorphic on $\mathcal{U}_{\epsilon_0}(0)$. 
\epr
\end{shaded}

\vskip 0.2 cm 
\noindent
{\it Proof:} 
Recall that 
\beas
&&
 (F_1)_{ab}(t)=\frac{(\widehat{A}_1)_{ab}(t)}{u_b(t)-u_a(t)},
 \quad\quad
 a\neq b,
\\
&&
(F_1)_{aa}(t)=-\sum_{b\neq a}(\widehat{A}_1)_{ab}(t)(F_1)_{ba}(t).
\eeas
If by assumption  (\ref{30giugno-3}) holds,  the above formulas  imply that $F_1(t)$ is holomorphic in $\mathcal{U}_{\epsilon_0}(0)$, because  the singularities at  $\Delta$, i.e. for  $u_a(t)-u_b(t)\to 0$, become  removable.  
 Since the asymptotics corresponding to (\ref{21luglio-2}) is uniform in a compact subset $K$ of a simply connected open subset of a cell, we substitute it into $dY=\Omega(z,t)Y$, with 
$$ 
\Omega(z,t)= z d\Lambda(t)+[F_1(t),d\Lambda(t)].
$$
By {comparing coefficients of} powers of $z^{-l}$ we obtain 
\be  
\label{21luglio-3} 
[F_{l+1}(t),d\Lambda(t)]=[F_1(t),d\Lambda(t)]F_l(t)-d F_l(t),\quad
\quad
l\geq 1.
\ee
In components of the differential $d$, this becomes a recursive relation (use $\partial \Lambda(t)/\partial t_i=E_i$):
$$
\Bigl[F_{l+1}(t),E_i\Bigr]=\Bigl[F_1(t),E_i\Bigr]F_l(t)-\frac{\partial F_l(t)}{\partial t_i},
$$
with,
$$
\Bigl[F_{l+1}(t),E_i\Bigr]
=
\left(
\begin{array}{ccccc}
0 &  & (F_{l+1})_{1i} &  & 0 
\\
  &        & \vdots &  & 
\\
-(F_{l+1})_{i1} & \cdots & 0 & \cdots &   -(F_{l+1})_{in}
\\
   &     & \vdots &  
 \\
 0 &  & (F_{l+1})_{ni} &  & 0    
\end{array}
\right),
$$
The diagonal element $(i,i)$ is zero. Therefore, (\ref{21luglio-3}) recursively determines $F_{l+1}$ as a function of $F_l,F_{l-1},...,F_1$, except for the diagonal $\hbox{ diag}(F_{l+1})$. 
On the other hand, the diagonal elements are determined by the off-diagonal elements according to the already proved  formula, 
\be
\label{3nov2015-6} 
l~(F_{l+1})_{aa}(t)=-\sum_{b\neq a}(\widehat{A}_1)_{ab}(t)(F_l)_{ba}(t).
\ee

 Let us start with $l+1=2$. Since $F_1$ is holomorphic, the above formulae (\ref{21luglio-3}), (\ref{3nov2015-6}) imply that $F_2$ is holomorphic. Then, by induction the same formulae imply  that all the $F_{l+1}(t)$ are holomorphic. $\Box$

\vskip 0.3 cm  
 Corollary \ref{28dic2015-2}    means that assumption 2) of Theorem \ref{pincopallino} applies, while  Proposition \ref{2gen2016-1} means that assumption 1) applies.  This, together with Proposition \ref{8feb2016-5}, proves the following theorem, which is indeed  our Theorem \ref{16dicembre2016-1}.

\begin{shaded}
\bth[Theorem \ref{16dicembre2016-1}.]
\label{3nov2015-5}
 Let  $\Lambda(t)$  and   $\widehat{A}_1(t)$ be  holomorphic on $\mathcal{U}_{\epsilon_0}(0)$, with eigenvalues as in (\ref{4gen2016-2}). If the deformation of the system (\ref{ourcase12}) is isomonodromic on a simply connected subset $\mathcal{V}$ of a cell, such that  $\overline{\mathcal{V}}$ is in the cell,  and if  $\widehat{A}_1(t)$ is holomorphically similar to a Jordan form $J^{(0)}$ in   $\mathcal{U}_{\epsilon_0}(0)$, or equivalently  the vanishing condition 
 $$(\widehat{A}_1)_{ab}(t)=\mathcal{O}(u_a(t)-u_b(t))\equiv \mathcal{O}(t_a-t_b),
 $$
 holds at points of $\Delta$ in $\mathcal{U}_{\epsilon_0}(0)$, then Theorem \ref{pincopallino}  and Corollary \ref{19feb2016-1} hold (with  $G_0(t)\mapsto I$, see Remark \ref{4gen2017-2}), so that  $\mathcal{G}_\sigma(z,t)=Y_\sigma(z,t)e^{\Lambda(t)}z^{-B_1(t)}$,  $\sigma=\nu,\nu+\mu,\nu+2\mu$,  maintains the  canonical asymptotics
$$ 
\mathcal{G}_\sigma(z,t)\sim I+\sum_{k=1}^\infty F_k(t) z^{-k},
\quad
\quad
z\to\infty \hbox{ in }\widehat{\mathcal{S}}_\sigma,
$$
for any $t\in \mathcal{U}_{\epsilon_1}(0)$ and any $\epsilon_1<\epsilon_0$. The   Stokes matrices,
$$ 
 \mathbb{S}_\nu,
 \quad
 ~\mathbb{S}_{\nu+\mu},
$$
are defined and constant on  the whole $\mathcal{U}_{\epsilon_0}(0)$. They  coincide with the Stokes matrices $ \mathring{\mathbb{S}}_\nu$, $ \mathring{\mathbb{S}}_{\nu+\mu}$ of the specific fundamental solutions $\mathring{Y}_\sigma(z)$ of the system (\ref{7gen2016-7}) 
$$
\frac{dY}{dz}=\widehat{A}(z,0) Y,
$$
 which satisfy  $\mathring{Y}_\sigma(z)\equiv Y_\sigma(z,0)$, according to  Corollary \ref{19feb2016-1}.   
Any  central connection matrix  $C^{(0)}_\nu$ is defined and constant on  the whole $\mathcal{U}_{\epsilon_0}(0)$, coinciding  with a  matrix  $\mathring{C}^{(0)}_\nu$ defined by the relation
$$ 
 \mathring{Y}_\nu(z)=\mathring{Y}^{(0)}(z) \mathring{C}^{(0)}_\nu,   
 $$
where $\mathring{Y}^{(0)}(z)$ is a fundamental solution of (\ref{7gen2016-7}) 
 in {the} Levelt form (\ref{17dicembre2016-7}), and $\mathring{Y}_\nu(z)=Y_\nu(z,0)$ as above. \\
The matrix entries of Stokes matrices vanish in correspondence with coalescing eigenvalues, i.e. 
 $$(\mathbb{S}_1)_{ij}=(\mathbb{S}_1)_{ji}=(\mathbb{S}_2)_{ij}=(\mathbb{S}_2)_{ji }=0 
 \quad
 \hbox{ whenever } u_i(0)=u_j(0).
 $$
\eth
\end{shaded}

\begin{shaded}
\bcr
\label{17dicembre2016-5} (Corollary  \ref{17dicembre2016-2})
If moreover the diagonal entries of $\widehat{A}_1(0)$  do not differ by non-zero integers, Corollary \ref{8gen2016-10} applies. Accordingly, there is a unique formal solution of the  system with  $t=0$, whose coefficients are necessarily  $$\mathring{F}_k\equiv F_k(0).$$ Hence,  there exists only one choice of fundamental solutions $\mathring{Y}_\sigma(z)$'s 
  with canonical asymptotics at $z=\infty$  corresponding to the unique formal solution,  which necessarily coincide with   the $Y_\sigma(z,0)$'s.
\ecr
\end{shaded}

Summarizing, the monodromy data are computable from the system with fixed $t=0$  and are: 

$\bullet$  $J^{(0)} =$ a Jordan form of  $\widehat{A}_1(0)$;    $R^{(0)}=\mathring{R}$.  See Proposition \ref{8feb2016-5}.
 
$\bullet$  $B_1=\hbox{\rm diag}(\widehat{A}_1(0))$. 

$\bullet$  $\mathbb{S}_\nu= \mathring{\mathbb{S}}_\nu$, $\mathbb{S}_{\nu+\mu}=\mathring{\mathbb{S}}_{\nu+\mu}$.

$\bullet$ $C^{(0)}_\nu=\mathring{C}^{(0)}_\nu$.

Here,  $\mathring{\mathbb{S}}_1$ and $\mathring{\mathbb{S}}_2$ are 
 the Stokes matrices  of those fundamental solutions $\mathring{Y}_1(z)$,  $\mathring{Y}_2(z)$,  $\mathring{Y}_3(z)$ of  the system (\ref{7gen2016-7})  (i.e. system (\ref{4gen2016-8})) with the specific canonical asymptotics (\ref{5gen2016-1})  satisfying   $\mathring{F}_k\equiv F_k(0)$, $k\geq 1$.  For these solutions the identity  
 $\mathring{Y}_r(z)=Y_r(z,0)$ holds. In case of Lemma \ref{17dicembre2016-5}, only these solutions exist.

\section{Isomonodromy Deformations with Vanishing Conditions on Stokes Matrices, Proof of Theorem \ref{9gen2017-1}}
\label{26maggio2017-8}

We now consider again system (\ref{ourcase12})  with eigenvalues (\ref{4gen2016-2}) coalescing at $t=0$, but we {\it give up the assumption that $\widehat{A}_1(t)$ is holomorphic in the whole $\mathcal{U}_{\epsilon_0}(0)$}. 
We assume that    $\widehat{A}_1(t)$ is  holomorphic on a simply connected open  domain $\mathcal{V}\subset \mathcal{U}_{\epsilon_0}(0)$,  as in Definition \ref{12maggio2017-2},  so that the Jimbo-Miwa-Ueno isomonodromy deformation theory\footnote{The fact that $\widehat{A}_1$ may have eigenvalues differing by integers does not constitute a problem; see   the proof of Theorem \ref{8luglio2016-1}.} is well defined $\mathcal{V}$. Therefore $Y_{\nu+k\mu}(t,z)$'s ($k\in\mathbb Z$)  and  $Y^{(0)}(t,z)$ satisfy the system
\bea
\label{ourcase12-TTRIS}
&&
\frac{dY}{dz}=\left(\Lambda(t) +\frac{\widehat{A}_1(t)}{z}\right)Y,
\\
&&dY=\Omega(z,t)Y,
\eea
and $\widehat{A}_1(t)$ solves the non-linear isomonodromy deformation equations 
\beas
&&
d\widehat{A}=\frac{\partial \Omega}{\partial z} +[\Omega,\widehat{A}],
\\
&&
dG^{(0)}=\Theta^{(0)}~ G^{(0)} .
\eeas
 Here $\Omega$ and $\Theta^{(0)}$ are the same as in the previous sections, defined for $t\in \mathcal{V}$.

Since the deformation is admissible, there exists $\widetilde{\tau}$ such that $\overline{\mathcal{V}}\subset c$, where $c$ is a $\widetilde{\tau}$-cell in $ \mathcal{U}_{\epsilon_0}(0)$. 
 The Stokes rays of $\Lambda(0)$ will be numerated so that $\tau_\nu<\widetilde{\tau}<\tau_{\nu+1}$.

As in Remark \ref{8gen2017-1}, the solutions $\widehat{A}_1(t)$,   any $Y_{\nu+k\mu}(t,z)$'s   and  $Y^{(0)}(t,z)$ of the above isomonodromy deformation equations, initially defined in $\mathcal{V}$, can be $t$-analytically continued on the  universal covering of $\mathbb{C}^n\backslash \Delta_{\mathbb{C}^n}$, 
as a meromorphic functions. 
The coalescence locus $\Delta_{\mathbb{C}^n}$ is a  locus of  fixed singularities \cite{Miwa}, so that it may be a branching 
locus for  $\widehat{A}_1(t)$ and for any of the fundamental matrices $Y(z,t)$ of  (\ref{ourcase12-TTRIS}) (i.e. of (\ref{ourcase12})).  Notice that our $\Delta$ is obviously contained in $\Delta_{\mathbb{C}^n}$.  The  movable singularities 
of $\widehat{A}_1(t)$, $Y_{\nu+k\mu}(t,z)$  and  $Y^{(0)}(t,z)$ outside $\Delta_{\mathbb{C}^n}$  are  poles and 
constitute, according to \cite{Miwa},  the locus of zeros of the Jimbo-Miwa-Ueno isomonodromic $\tau$-function. This locus can also be called  {\it Malgrange's divisor}, since it has been proved in \cite{Palmer} that it coincides with a divisor,  introduced by Malgrange (see \cite{Malg1} \cite{Malg2} \cite{Malg3}), where  a certain Riemann-Hilbert problem fails to have solution (below, we formulate  a Riemann-Hilbert problem in proving Lemma \ref{27feb2017-1}). This divisor has a complex co-dimension equal to 1, so it  does not disconnect $\mathbb{C}^n\backslash \Delta_{\mathbb{C}^n}$ and $\mathcal{U}_{\epsilon_0}(0)\backslash \Delta$.

The fundamental solutions $Y_{\nu+k\mu}(t,z)$'s above  are the unique solutions which have for $t\in \mathcal{V}$ the asymptotic behaviour 
\be
\label{8gen2017-3}
Y_{\nu+k\mu}(z,t)e^{-\Lambda(t)z}z^{-B_1}\sim I+\sum_{j\geq 1}F_j(t)z^{-j}, 
\quad
~z\to\infty \hbox{ in }\mathcal{S}_{\nu+k\mu}(t).
\ee 
The $t$-independent Stokes matrices are then defined by the relations 
$$
Y_{\nu+(k+1)\mu}(t,z)=Y_{\nu+k\mu}(t,z)\mathbb{S}_{\nu+k\mu}.
$$
Notice that also the coefficients $F_j(t)$ are analytically continued as meromorphic multivalued matrix functions.  
For the sake of the proof of the Lemma \ref{8gen2017-9} below, the analytic continuation of  $Y_{\nu+k\mu}(t,z)$ will be denoted by 
$$
\mathbb{Y}_{\nu+k\mu}(z,\tilde{t}),
$$
 where $\tilde{t}$ is a point of the universal covering $\mathcal{R}(\mathcal{U}_{\epsilon_0}(0)\backslash\Delta)$, whose projection is $t$.  The analytic continuation of  $F_j(t)$ will be simply denoted by $F_j(\tilde{t})$

By arguments similar to those in Section \ref{14feb2016-10}, it is seen that as $t$ varies in $c$ or slightly beyond the boundary $\partial c$, then $Y_{\nu+k\mu}(t,z)$ maintains its asymptotic behaviour, for $t$  away from the Malgrange's divisor. But when $t$ moves sufficiently far form  $c$, then the asymptotic representation (\ref{8gen2017-3}) is lost. 
 The following Lemma gives the sufficient condition such that the asymptotics (\ref{8gen2017-3}) is not lost by $\mathbb{Y}_{\nu+k\mu}(z,\tilde{t})$. 
\ble
\label{8gen2017-9}
Assume that  the Stokes matrices  satisfy the vanishing condition
\be
\label{6gen2017-1-BBBIS}
(\mathbb{S}_\nu)_{ab}=(\mathbb{S}_\nu)_{ba}=(\mathbb{S}_{\nu+\mu})_{ab}=(\mathbb{S}_{\nu+\mu})_{ba}=0,
\ee
for any $1\leq a\neq b\leq n$ such that $u_a(0)=u_b(0)$. 
Then the meromorphic continuation  $\mathbb{Y}_{\nu+k\mu}(z,\tilde{t})$, $k\in\mathbb{Z}$, on the universal covering $\mathcal{R}(\mathcal{U}_{\epsilon_0}(0)\backslash\Delta)$   maintains the asymptotic behaviour 
 $$
 \mathbb{Y}_{\nu+k\mu}(z,\tilde{t}) e^{-\Lambda(t)z}z^{-B_1}~\sim I+\sum_{j\geq 1}F_j(\tilde{t})~z^{-j},
$$
  for $z\to \infty$ in $\widehat{\mathcal{S}}_{\nu+k\mu}(t)$ and any $\tilde{t}\in\mathcal{R}(\mathcal{U}_{\epsilon_0}(0)\backslash\Delta)$  away from the Malgrange's divisor.  Moreover, 
  $$\mathbb{Y}_{\nu+(k+1)\mu}(z,\tilde{t}) =\mathbb{Y}_{\nu+k\mu}(z,\tilde{t}) \mathbb{S}_{\nu+k\mu}.
  $$
  Here $\widehat{\mathcal{S}}_{\nu+k\mu}(t)$ is the sector in Definition \ref{13giugno2017-4}.

\ele

\bre
 Notice that $B_1=\hbox{diag}(\widehat{A}_1(t))$ is independent of $t \in \mathcal{V}$ by assumption, and $\widehat{A}_1(t)$ is meromorphic, so $B_1$  is constant everywhere. Moreover,  the relation $\mathbb{S}_{\nu+2\mu}=e^{-2\pi i B_1}\mathbb{S}_\nu e^{2\pi i B_1}$ implies that (\ref{6gen2017-1-BBBIS}) holds for any $\mathbb{S}_{\nu+k\mu}$, $k\in\mathbb{Z}$.
\ere

\noindent
{\it Proof:}   
Since $\overline{\mathcal{V}}$ belongs to the $\widetilde{\tau}$-cell  $ c$, then $Y_{\nu+k\mu}(z,t)$ can be denoted by $Y_{\nu+k\mu}(z,t;\widetilde{\tau},c)$,  as in Theorem \ref{8gen2017-4}, for $t\in \mathcal{V}$ and for any  $t\in c$ away from the Malgrange's divisor.  Noticing that the Malgrange's divisor does not disconnect $\mathcal{U}_{\epsilon_0}(0)\backslash \Delta$, we proceed exactly as in the proof of  Theorem \ref{8gen2017-4}.  Now $\mathcal{V}$ is considered as lying on a sheet of the covering $\mathcal{R}(\mathcal{U}_{\epsilon_0}(0)\backslash \Delta)$. The relation (\ref{15dic2015-1-BiiS}) holds unchanged, and reads
\be
\label{2marzo2017-1}
Y_{\nu+\mu}(z,\tilde{t};\widetilde{\tau},c^\prime)=\mathbb{Y}_{\nu+\mu}(z,\tilde{t},\widetilde{\tau},c)\mathbb{K}^{[ab]}
.
\ee
On the other hand, the relation (\ref{22feb2016-2-BiiS}) becomes 
$$
\mathbb{X}_{\nu+\mu}(z,\tilde{t})
=Y_{\nu+\mu}(z,\tilde{t};\widetilde{\tau},c^\prime)~\widetilde{\mathbb{K}}^{[ab]}(t),
$$
where $\mathbb{X}_{\nu+\mu}(z,\tilde{t})$ is {a} solution of  the system (\ref{ourcase12-TTRIS}) with coefficient $\widehat{A}_1(\tilde{t})$, where $\tilde{t}$ is a point of the universal covering, reached   along $\gamma_{ab}$ after $R_{ab}(t)$ has crossed $R(\widetilde{\tau}-\pi)$ in Figure \ref{FiguraF}.  $\mathbb{X}_{\nu+\mu}(z,\tilde{t})$ is the unique fundamental matrix solution having asymptotic behaviour 
 $$
\mathbb{X}_{\nu+\mu}(z,\tilde{t}) e^{-\Lambda(t)z}z^{-B_1}~\sim I+\sum_{j\geq 1}F_j(\tilde{t})~z^{-j},
$$
in $\mathcal{S}_{\nu+\mu}(t)$. 
Then (\ref{8gen2017-7}) is replaced by 
$$
\mathbb{X}_{\nu+\mu}(z,\gamma_{ab}t)=\mathbb Y_{\nu+\mu}(z,\gamma_{ab}t)~\mathbb{K}^{[ab]}~\widetilde{\mathbb{K}}^{[ab]},
\quad
\quad
t\in c.
$$
Here, $\mathbb Y_{\nu+\mu}(z,\gamma_{ab}t)$ is the continuation of $ Y_{\nu+\mu}(z,t)
\equiv Y_{\nu+\mu}(z,t;\widetilde{\tau},c) $ at 
$$
\tilde{t}=\gamma_{ab}t.
$$ 
The proof that $\mathbb{K}^{[ab]}=\widetilde{\mathbb{K}}^{[ab]}=I$ holds unchanged, following from (\ref{6gen2017-1-BBBIS}). Therefore,  
$$
\mathbb{X}_{\nu+\mu}(z,\gamma_{ab}t)=\mathbb Y_{\nu+\mu}(z,\gamma_{ab}t).
$$
This proves that the analytic continuation $\mathbb Y_{\nu+\mu}(z,\tilde{t})$  along $\gamma_{ab}$ maintains the canonical  asymptotic behaviour. 
Moreover, the ray $R_{ab}$ {plays} no role in the asymptotics, as it follows from (\ref{2marzo2017-1}) with $\mathbb{K}^{[ab]}=I$.  
Repeating the construction for all possible loops $\gamma_{ab}$, as in the proof of Theorem \ref{pincopallino} and  Theorem \ref{8gen2017-4}, we conclude that $\mathbb Y_{\nu+\mu}(z,\tilde{t})$ maintains its the canonical asymptotic representation  for any $\tilde{t}$ in the universal covering ($\tilde{t}$ away from the Malgrange divisor), when $z\to\infty$  in $\widehat{S}_{\nu+\mu}(t)$.  
$\Box$

\vskip 0.3 cm 
In Lemma \ref{8gen2017-9}, we have taken into account the fact that $\Delta$ is expected to be a branching locus, so that $\mathbb{Y}(z,\tilde{t})$ is defined on $\mathcal{R}(\mathcal{U}_{\epsilon_0}(0)\backslash\Delta)$, as the result of \cite{Miwa} predicts. In fact, it turns out that (\ref{6gen2017-1-BBBIS}) implies that there is no branching at $\Delta$, as the following lemma states. 

\ble
\label{27feb2017-1}
If (\ref{6gen2017-1-BBBIS}) holds, then: 
\begin{itemize}
\item
The meromorphic continuation  on the universal covering $\mathcal{R}(\mathcal{U}_{\epsilon_0}(0)\backslash\Delta)$  of any $Y_{\nu+k\mu}(z,t)$, $k\in\mathbb Z$, and $Y^{(0)}(z,t)$ is single-valued   on $\mathcal{U}_{\epsilon_0}(0)\backslash\Delta$. 
 \item
 The meromorphic continuation of $\widehat{A}_1(t)$ is  single-valued   on $\mathcal{U}_{\epsilon_0}(0)\backslash\Delta$. 
\end{itemize}
In other words, $\Delta$ is not a branching locus.
\ele

 The single-valued continuation of $Y_{\nu+k\mu}(z,t)$ will be simply denoted by $Y_{\nu+k\mu}(z,t)$ in the remaining part of this section, so we will no longer need the notation $\mathbb Y_{\nu+k\mu}(z,\tilde{t})$.

 \vskip 0.2 cm 
 \noindent
 \underline{\it Proof of Lemma \ref{27feb2017-1}:}  Let $t\in \mathcal{V}$ be an admissible isomonodromic deformation  and  $\widehat{A}_1(t)$ be holomorphic in $\mathcal{V}$. Let   $\widetilde{\tau}$ be the direction of an admissible ray for $\Lambda(0)$ such that $\mathcal{V}$ lies in a  $\widetilde{\tau}$-cell. 
 Since  the linear relation 
(\ref{4gen2016-2})
   $$
    u_i(t)=u_i(0)+t_i,
    \quad\quad~~1\leq i \leq n,
    $$
    holds, we  will use $u$ as variable in place of $t$. Accordingly, we will write $\Lambda(u)$ instead of $\Lambda(t)$ and $Y(z,u)$ instead of $Y(z,t)$.   Now, the fundamental solutions $Y_{\nu+k\mu}(z,u)$ and $Y^{(0)}(z,u)$ are holomorphic functions of $u\in \mathcal{V}$. We construct a Riemann-Hilbert boundary value problem (abbreviated by R-H) satisfied by\footnote{Recall that $Y_{\nu+2k\mu}(ze^{2k\pi  i})=Y_\nu(z) e^{2k\pi  i B_1}$, $k\in\mathbb Z$.} $Y_{\nu-\mu}(z,u)$, $Y_\nu(z,u)$, $Y_{\nu+\mu}(z,u)$  and $Y^{(0)}(z,u)$. 
     
     The given data are the essential monodromy data (see Definition \ref{30dic2015-1})  $\mathbb{S}_{\nu-\mu}$, $\mathbb{S}_\nu$, $B_1$, $\mu_1,...,\mu_n$, $R^{(0)}$ and $C^{(0)}_\nu$. Instead of $\mu_1,...,\mu_n$, $R^{(0)}$, we can use $D^{(0)}$ and $L^{(0)}$ (see (\ref{24gen2016-1BBIS}) and Remark \ref{17april2016-3}).  They satisfy a constraint, because  the monodromy  $( C_\nu^{(0)})^{-1} ~e^{2\pi i L^{(0)}}~ C_\nu^{(0)}$ at $z=0$   can be expressed in the equivalent way $e^{2\pi i B_1} (\mathbb{S}_\nu  \mathbb{S}_{\nu+\mu})^{-1}$. Recalling  that $\mathbb{S}_{\nu+\mu}=e^{-2\pi i B_1} \mathbb{S}_{\nu-\mu} e^{2\pi i B_1}$, the constraint is
    \be
    \label{17aprile2017-1}
    \mathbb{S}_{\nu-\mu}^{-1}~e^{2\pi i B_1}~\mathbb{S}_\nu^{-1} =~( C_\nu^{(0)})^{-1} ~e^{2\pi i L^{(0)}}~ C_\nu^{(0)}.
    \ee
  The following relations hold for fundamental solutions:
    \bea
    \label{27feb2017-2}
   && Y_\nu(z,u)=Y_{\nu-\mu}(z,u) \mathbb{S}_{\nu-\mu},
    \\
 &&    Y_{\nu+\mu}(z,u)= Y_\nu(z,u) \mathbb{S}_\nu,
    \\
  &&   Y_\nu(z,u)=Y^{(0)}(z,u)C^{(0)}_\nu,
    \\
  &&   Y_{\nu+\mu}(z,u)=Y^{(0)}(z,u)C^{(0)}_\nu \mathbb{S}_\nu.
    \eea
    Since $Y_{\nu+\mu}(ze^{2\pi i})=Y_{\nu-\mu}(z)e^{2\pi i B_1}$, we can rewrite (\ref{27feb2017-2}) as 
\be
 \label{27feb2017-3}
Y_\nu(z,u)=Y_{\nu+\mu}(ze^{2\pi i},u) e^{-2\pi i B_1}\mathbb{S}_{\nu-\mu}
\ee
 We now write 
\beas 
&& Y_{\nu+k\mu}(z,u)=\mathcal{G}_{\nu+k\mu}(z,u) e^{Q(z,u)},
\quad
\quad
Q(z,u):=\Lambda(u)z+B_1\ln z,
\\
&&\mathcal{G}_{\nu+k\mu}(z,u) \sim I+\sum_{j=1}^\infty F_j(u)z^{-j}, 
\quad
\quad
z\to \infty \hbox{ in } \mathcal{S}_{\nu+k\mu}(u),
\quad
~k=0,1. 
\\
\\
&&
Y^{(0)}(z,u)=\mathcal{G}_0(z,u)~ z^{D^{(0)}} z^{L^{(0)}}
\\
&&
\mathcal{G}^{(0}(z,u)=G^{(0)}(u)+O(z)
\quad
 \hbox{  holomorphic at $z=0$.}
\eeas
Therefore, from (\ref{27feb2017-2})-(\ref{27feb2017-3}) we obtain 
\bea
\label{17marzo2017-1}
&&\mathcal{G}_\nu(z,u)=\mathcal{G}_{\nu+\mu}(ze^{2\pi i},u)~ e^{Q(z,u)}\mathbb{S}_{\nu-\mu} e^{-Q(z,u)},
\\
&&
\mathcal{G}_{\nu+\mu}(z)=\mathcal{G}_\nu(z,u)~ e^{Q(z,u)} \mathbb{S}_\nu e^{-Q(z,u)},
\\
&&
\mathcal{G}_\nu(z,u)=\mathcal{G}^{(0)}(z,u) ~z^{D^{(0)}}z^{L^{(0)}} C^{(0)}_\nu e^{-Q(z,u)},
\\
\label{17marzo2017-2}
&&
\mathcal{G}_{\nu+\mu}(z,u)=\mathcal{G}^{(0)}(z,u) ~z^{D^{(0)}}z^{L^{(0)}} C^{(0)}_\nu \mathbb{S}_\nu e^{-Q(z,u)}.
\eea

 \begin{figure}
\centerline{\includegraphics[width=0.6 \textwidth]{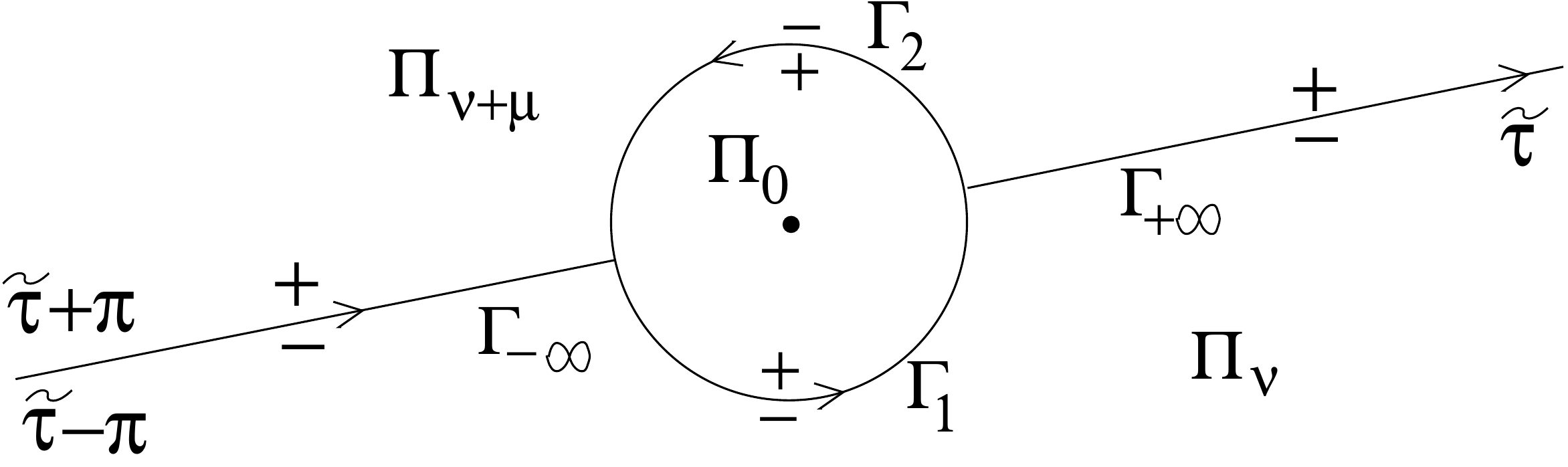}}
\caption{The contour $\Gamma_{-\infty}\cup \Gamma_1\cup \Gamma_2\cup \Gamma_{+\infty}$ of the Riemann-Hilbert problem, which divides the plane in regions $\Pi_\nu$, $\Pi_{\nu+\mu}$ and $\Pi_0$. The directional angles $\widetilde{\tau}$, $\widetilde{\tau}\pm \pi$ and the orientations are depicted.  } 
\label{RH-1}
\end{figure}

We formulate the following R-H, given the monodromy data. Consider the $z$-plane with the following  branch cut from $0$ to $\infty$: 
$$ 
\widetilde{\tau} -\pi <\arg z <\widetilde{\tau}+\pi.
$$
Consider a circle around $z=0$ of some radius $r$. The oriented  contour $\Gamma=\Gamma(\widetilde{\tau})$ of the R-H is the union of the following paths (see Figure \ref{RH-1}):
\beas
&&\Gamma_{-\infty}:
\quad
 \arg z=\widetilde{\tau}\pm \pi, 
 \quad
 ~|z|>r, 
 \quad~\hbox{ half-line coming from $\infty$ along the branch-cut}
\\
\\
&&\Gamma_{+\infty}:
\quad
  \arg z=\widetilde{\tau}, 
  \quad
  ~|z|>r,  
  \quad
  ~\hbox{ half-line going  to $\infty$  in direction $\widetilde{\tau}$},
\\
\\
&& \Gamma_1:
\quad
\widetilde{\tau} -\pi <\arg z \leq \widetilde{\tau}, 
\quad
~|z|=r,  
\quad
~\hbox{ half-circle in anti-clockwise sense},
\\
\\
&&
\Gamma_2:
\quad
\widetilde{\tau}  \leq \arg z  <\widetilde{\tau}+\pi,  
\quad
~|z|=r, 
\quad
~\hbox{ half-circle in anti-clockwise sense}.
\eeas 
Recalling that $\tau_\nu<\widetilde{\tau}<\tau_{\nu+\mu}$, we call: 
 
$\Pi_\nu$ the unbounded domain to the right of $\Gamma_{-\infty}\cup \Gamma_1 \cup \Gamma_{+\infty}$, 

$\Pi_0$ the ball inside the circle  $\Gamma_1\cup \Gamma_2$, 

 $\Pi_{\nu+\mu}$ the remaining unbounded region $\mathbb{C}\backslash \{\Pi_\nu\cup\Pi_0\cup \Gamma\}$.

\noindent
The R-H problem we need is as follows: 
\be
\label{19aprile2017-1}
\mathcal{G}_+(\zeta)=\mathcal{G}_-(\zeta) H(\zeta,u),
\quad
\quad
\zeta\in\Gamma,
\ee
where the jump $H(\zeta,u)$ is uniquely specified by assigning the monodromy data $\mathbb{S}_{\nu-\mu}$, $\mathbb{S}_\nu$, $B_1$,  $C^{(0)}_\nu$, $D^{(0)}$ and $L^{(0)}$ (i.e. $\mu_1,...,\mu_n$, $R^{(0)}$).   Since  $\Gamma_-$ lies along the branch-cut, we use  the symbol $\zeta_{\pm}$ if  $\arg \zeta=\widetilde{\tau} \pm \pi$.  Hence, $ H(\zeta,u)$ is 

\beas
 H(\zeta,u):= &&e^{Q(\zeta_-,u)}\mathbb{S}_{\nu-\mu}^{-1}e^{-Q(\zeta_-,u)} 
 \quad\quad
 \hbox{ along $\Gamma_{-\infty}$},
\\
&&e^{Q(\zeta,u)}\mathbb{S}_\nu e^{-Q(\zeta,u)}
\quad\quad
\hbox{ along $\Gamma_{+\infty}$},
\\
&&e^{Q(\zeta,u)} (C^{(0)}_\nu)^{-1} \zeta^{-L^{(0)}}\zeta^{-D^{(0}}
\quad
\quad
\hbox{ along $\Gamma_1$},
\\
&&e^{Q(\zeta,u)} \mathbb{S}_\nu^{-1}(C^{(0)}_\nu)^{-1} \zeta^{-L^{(0)}}\zeta^{-D^{(0}}
\quad
\quad
\hbox{ along $\Gamma_2$}.
\eeas
We require that the solution satisfies the conditions
\bea
\label{2marzo2017-2}
&&\mathcal{G}(z)\sim I+ \hbox{series in }z^{-1},
\quad
z\to\infty,
\quad
~z\in\Pi_\nu\cup \Pi_{\nu+\mu},
\\
\label{19aprile2017-2}
&&
\mathcal{G}(z) \hbox{ holomorphic in $\Pi_0$ and $\det(\mathcal{G}(0))\neq 0$}.
\eea

 By (\ref{17marzo2017-1})-(\ref{17marzo2017-2}),  our R-H  has  the following solution for $u\in \mathcal{V}$: \be
\label{28feb2017-4}
\mathcal{G}(z,u)=
\left\{
\begin{array}{c}
\mathcal{G}_0(z,u)  
\quad
~\hbox{ for $z\in \Pi_0$},
\\
\mathcal{G}_\nu(z,u) 
\quad
~ \hbox{ for $z\in \Pi_\nu$},
\\
 \mathcal{G}_{\nu+\mu}(z,u) 
 \quad
 ~ \hbox{ for $z\in \Pi_{\nu+\mu}$},
\end{array}
\right. 
\quad
~\hbox{holomorphic of $u\in \mathcal{V}$}.
\ee
  By the result of \cite{Miwa},  this solution can be analytically continued in $u$ as a meromorphic function on the universal covering of $\mathbb{C}^n\backslash \Delta_{\mathbb{C}_n}$.  Consider a loop around $\Delta$, as in (\ref{10maggio2016-1}), involving two coalescing coordinates $u_a,u_b$, starting from a point in $\mathcal{V}$.  We want to prove that the above continuation is single valued along this loop. As  in the proof of Theorem \ref{pincopallino}, we just need to consider the case when  $|u_a-u_b|$ is small and only $PR_{ab}$ and $PR_{ba}$ cross $l(\widetilde{\tau})$. Let 
$$ 
\varepsilon:=u_a-u_b.
$$
The lemma will be proved if we  prove that $\mathcal{G}$ in (\ref{28feb2017-4}) is holomorphic in a neighbourhood of $\varepsilon=0$, except at most for a finite number of  poles (the Malgrange's divisor).

\begin{figure}
\centerline{\includegraphics[width=0.6 \textwidth]{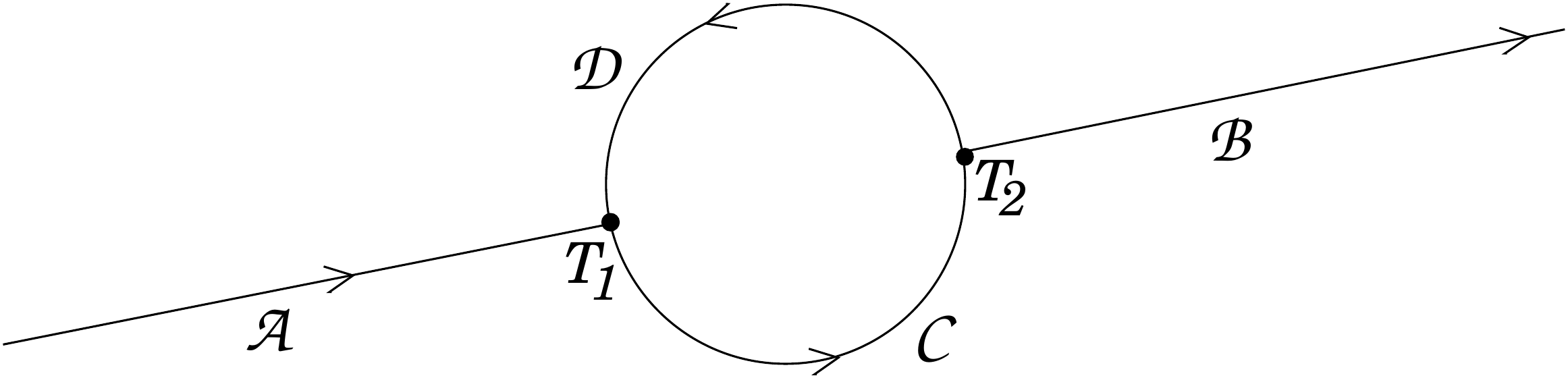}}
\caption{Jump matrices $\mathcal{A}$,  $\mathcal{B}$,  $\mathcal{C}$,  $\mathcal{D}$ along $\Gamma$, used in step 1.} 
\label{RH-2}
\end{figure}

In the following, we will drop  $u$ and only write the dependence on  $\varepsilon$. For example, we write   $H(\zeta,\varepsilon)$ instead of $H(\zeta,u)$. For our convenience, 
as in Figure \ref{RH-2} we call 
\beas
 H(\zeta,\varepsilon) &&=:\mathcal{A}(\zeta_-,\varepsilon)
 \quad
 \quad
 \hbox{ along $\Gamma_{-\infty}$},
\\
&&=:\mathcal{B}(\zeta,\varepsilon)
\quad
\quad
\hbox{ along $\Gamma_{+\infty}$},
\\
&&=:\mathcal{C}(\zeta,\varepsilon)
\quad
\quad
\hbox{ along $\Gamma_1$},
\\
&&=:\mathcal{D}(\zeta,\varepsilon)
\quad
\quad
\hbox{ along $\Gamma_2$}.
\eeas
 $\mathcal{A},...,\mathcal{D}$ are  {\it holomorphic} functions of $\varepsilon$. The following cyclic relations are easily verified:
\be
\label{17aprile2017-2}
\mathcal{A}(ze^{-2\pi i},\varepsilon)~\mathcal{D}(z,\varepsilon)~\mathcal{C}(ze^{-2\pi i},\varepsilon)^{-1}~=I
\quad
\quad
\mathcal{C}(z,\varepsilon)~\mathcal{D}(z,\varepsilon)^{-1}\mathcal{B}(z,\varepsilon)^{-1}~=I.
\ee
In particular, the following ``smoothness condition'' holds  at the points $T_1$ and $T_2$ of intersection of $\Gamma_{-\infty}$ and $\Gamma_{+\infty}$ with the circle $|z|=r$ respectively: 
$$
\mathcal{A}(\zeta_-,\varepsilon)~\mathcal{D}(\zeta_+,\varepsilon)~\mathcal{C}(\zeta_-,\varepsilon)^{-1}~=I 
\quad
~\hbox{ at $T_1$},
\quad
\quad
\mathcal{C}(\zeta,\varepsilon)~\mathcal{D}(\zeta,\varepsilon)^{-1}\mathcal{B}(\zeta,\varepsilon)^{-1}=I 
\quad
~ \hbox{ at $T_2$}.
$$

Indeed, 
\beas
&&\mathcal{A}(ze^{-2\pi i},\varepsilon)\mathcal{D}(z,\varepsilon)\mathcal{C}(ze^{-2\pi i},\varepsilon)^{-1}=
\\
&&
= e^{Q(ze^{-2\pi i})}\mathbb{S}_{\nu-\mu}^{-1}e^{-Q(ze^{-2\pi i})}\cdot e^{Q(z)}\mathbb{S}_\nu^{-1}(C_\nu^{(0)})^{-1}z^{-L^{(0)}}z^{-D^{(0)}}\cdot (ze^{-2\pi i})^{D^{(0)}}(ze^{-2\pi i})^{L^{(0)}} C_\nu^{(0)} e^{-Q(ze^{-2\pi i})}
\\
&&
=e^{-2\pi i B_1} e^{Q(z)}~\mathbb{S}_{\nu-\mu}^{-1}~e^{2\pi i B_1}~\mathbb{S}_\nu^{-1}~(C_\nu^{(0)})^{-1}z^{-L^{(0)}}z^{-D^{(0)}}\cdot z^{D^{(0)}}z^{L^{(0)}} e^{-2\pi i L^{(0)}} C_\nu^{(0)}~ e^{-Q(z)}~e^{2\pi i B_1}
\\
&&
=e^{-2\pi i B_1} e^{Q(z)}~\Bigl(\mathbb{S}_{\nu-\mu}^{-1}~e^{2\pi i B_1}~\mathbb{S}_\nu^{-1}~(C_\nu^{(0)})^{-1} e^{-2\pi i L^{(0)}} C_\nu^{(0)}\Bigr)~ e^{-Q(z)}~e^{2\pi i B_1}=I.
\eeas
In the last step, we have used (\ref{17aprile2017-1}). Moreover, 
$$
\mathcal{C}(\zeta,\varepsilon)\mathcal{D}(z,\varepsilon)^{-1}\mathcal{B}(z,\varepsilon)^{-1}=
e^{Q(z)}(C_\nu^{(0)})^{-1}z^{-L^{(0)}}z^{-D^{(0)}} \cdot z^{D^{(0)}}z^{L^{(0)}} C_\nu^{(0)} \mathbb{S}_\nu e^{-Q(z)}
\cdot e^{Q(z)}\mathbb{S}_\nu^{-1} e^{-Q(z)}=I.
$$
The last  result follows from simple cancellations.

In order to complete the proof, we need the theoretical  background, in particular the $L^p$ formulation of Riemann-Hilbert problems, found in  the  test-book \cite{Its}, the lecture notes \cite{ItsLect}  and the papers  \cite{Zhou} \cite{LpDZ} (see also     \cite{Deift}  \cite{DKMVZ} \cite{DKMVZ1} and \cite{Clancey} \cite{Pog}  \cite{Vekua}).  The proof is completed in the following steps, suggested to us by Marco Bertola. 

\vskip 0.2 cm 
$\bullet$ {\bf Step 1.} We contruct a naive solution $\mathfrak{S}(z,\varepsilon)$ to the R-H, which does not satisfy the asymptotic condition 
(\ref{2marzo2017-2}).  We start by defining $\mathfrak{S}(z,\varepsilon)=I$ in $\Pi_0$.  Then, keeping into account the jumps $\mathcal{C}$ and $\mathcal{B}$  along $\Gamma_1$ and $\Gamma_{+\infty}$ respecively (see Figure \ref{RH-2}), we have
\be
\label{17aprile2017-5}
\mathfrak{S}(z,\varepsilon)
=
\left\{
\begin{array}{cc}
I & \hbox{ for } z\in \Pi_0,
\\
\\
\mathcal{C}(z,\varepsilon)^{-1}&\hbox{ for } z\in \Pi_\nu,
\\
\\
\mathcal{C}(z,\varepsilon)^{-1}\mathcal{B}(z,\varepsilon) & \hbox{ for } z\in \Pi_{\nu+\mu},
\end{array}
\right.
\ee
On the other hand, starting with $\mathfrak{S}(z,\varepsilon)=I$ in $\Pi_0$ and keeping into account the jump $\mathcal{D}$  at $\Gamma_2$, we must have 
\be
\label{17aprile2017-3}
\mathfrak{S}(z,\varepsilon)= \mathcal{D}(z,\varepsilon)^{-1}
\quad
~\hbox{ for } z\in \Pi_{\nu+\mu}.
\ee
The second  relation in (\ref{17aprile2017-2}) ensures that (\ref{17aprile2017-3}) and the last expression in (\ref{17aprile2017-5}) coincide. 
Moreover, starting with $\mathfrak{S}(z,\varepsilon)=I$ in $\Pi_0$ and crossing $\Gamma_1$ and then $\Gamma_{-\infty}$ with jumps $\mathcal{C}$ and $\mathcal{A}$, we find a third representation of  $\mathfrak{S}(z,\varepsilon)$ for $z\in \Pi_{\nu+\mu}$, namely 
\be
\label{17aprile2017-4}
\mathfrak{S}(z,\varepsilon)=\mathcal{C}(ze^{-2\pi i},\varepsilon)^{-1} \mathcal{A}(ze^{-2\pi i},\varepsilon)
\quad
~\hbox{ for } z\in \Pi_{\nu+\mu}.
\ee
Now, the first relation in (\ref{17aprile2017-2}) ensures that (\ref{17aprile2017-3}) and (\ref{17aprile2017-4}) coincide.

\begin{figure}
\centerline{\includegraphics[width=0.6 \textwidth]{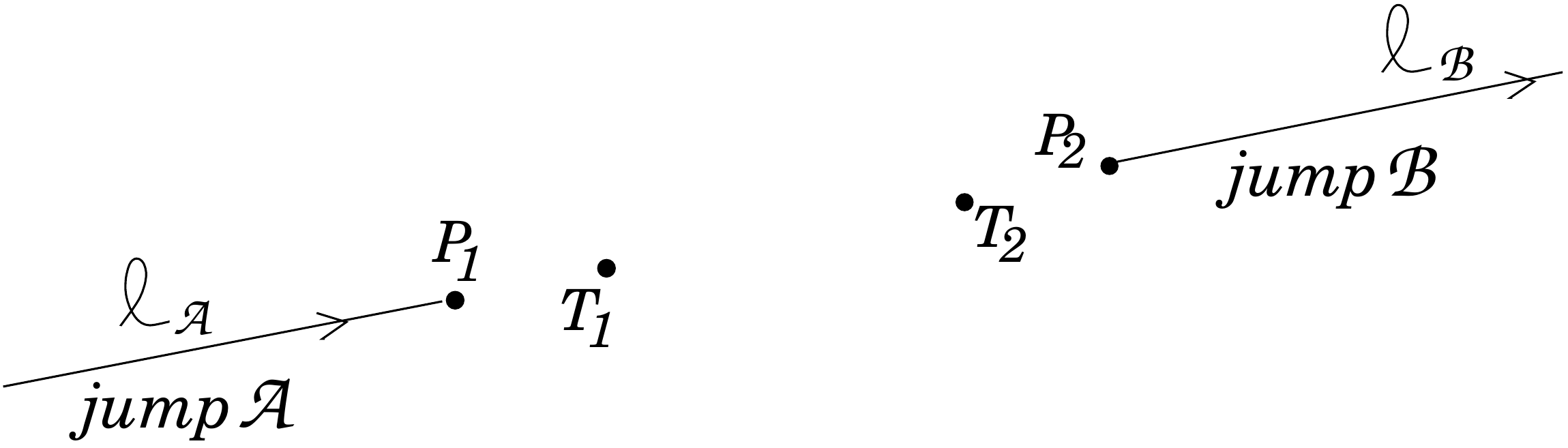}}
\caption{Step 2:  the auxiliary Riemann-Hilbert problem with contour $\ell_\mathcal{A}$ and $\ell_\mathcal{B}$.} 
\label{RH-3}
\end{figure}

\vskip 0.2 cm
$\bullet$ {\bf Step 2.}
We consider an auxiliary R-H  as in Figure \ref{RH-3}, whose boundary contour  is the union of a half line $\ell_{\mathcal{A}}$ contained in $\Gamma_{-\infty}$ from $\infty$ to a point  $P_1$ preceding $T_1$, and a half line $\ell_{\mathcal{B}}$ contained in $\Gamma_{+\infty}$ from a point $P_2$ following $T_2$ to $\infty$.  The jump along these half lines is  $H(\zeta,\varepsilon)$ (namely,  $\mathcal{A}(\zeta_-,\varepsilon)$ and $\mathcal{B}(\zeta,\varepsilon)$ on the two half lines respectively). The R-H is then 
$$
\Psi_+(\zeta)=\Psi_-(\zeta) H(\zeta,\varepsilon) 
\quad
\quad
\zeta\in \ell_{\mathcal{A}}\cup \ell_{\mathcal{B}},
$$
\be
\label{19aprile2017-4}
\Psi(z)\sim I+ \hbox{series in }z^{-1},
\quad
~z\to\infty,
\quad
~z\in\Pi_\nu\cup \Pi_{\nu+\mu}.
\ee
Keeping the above asymptotics  into account, the R-H  is rewritten as follows:  $$ 
\Psi(z)= I+\int_{\ell_{\mathcal{A}}\cup \ell_{\mathcal{B}}} \frac{\Psi_-(\zeta)( H(\zeta,\varepsilon)-I)}{\zeta -z }\frac{d\zeta}{2\pi i }.
$$
or, letting $\delta\Psi:=\Psi-I$ and $\delta H:=  H-I$, 
\be
\label{28feb2017-3}
\delta\Psi(z)= 
\int_{\ell_{\mathcal{A}}\cup \ell_{\mathcal{B}}}
\frac{\delta\Psi_-(\zeta)~\delta H(\zeta,\varepsilon)}{\zeta -z }\frac{d\zeta}{2\pi i }+\int_{\ell_{\mathcal{A}}\cup \ell_{\mathcal{B}}}
\frac{\delta H(\zeta,\varepsilon)}{\zeta -z }\frac{d\zeta}{2\pi i }
.
\ee
We solve the problem  by computing $\delta\Psi_-(\zeta)$, as the solution of the following integral equation (by taking the limit for $z\to z_{-}$ belonging to the ``$-$'' side of $\ell_{\mathcal{A}}\cup \ell_{\mathcal{B}}$): 
$$ 
\delta\Psi_{-}(z_{-})= 
\int_{\ell_{\mathcal{A}}\cup \ell_{\mathcal{B}}}
\frac{\delta\Psi_-(\zeta)~\delta H(\zeta,\varepsilon)}{\zeta -z_{-} }\frac{d\zeta}{2\pi i }+\int_{\ell_{\mathcal{A}}\cup \ell_{\mathcal{B}}}
\frac{\delta H(\zeta,\varepsilon)}{\zeta -z_{-} }\frac{d\zeta}{2\pi i }
$$
$$
=
C_{-}\Bigl[ \delta\Psi_-\delta H(\cdot,\varepsilon))\Bigr](z_{-}) + C_{-}\Bigl[\delta H(\cdot,\varepsilon)\Bigr](z_{-}).
$$ 
Here $C_{-}$ stands for the Cauchy boundary operator. We will  write $C_{-}\bigl[ \delta\Psi_-\delta H(\cdot,\varepsilon)\bigr]$ as $C_{-}\bigl[ \bullet~ \delta H(\cdot,\varepsilon)\bigr] \delta\Psi_{-}$, to represent the operator $C_{-}\bigl[ \bullet~ \delta H(\cdot,\varepsilon)\bigr] $ acting on $ \delta\Psi_{-}$.  We observe the following facts: 
 
 \begin{itemize}
\item[1.]
 If $u$ is in the cell containing $\mathcal{V}$, as $\zeta\to \infty$   along $\ell_{\mathcal{A}}$ and $\ell_{\mathcal{B}}$,  the off-diagonal  matrix entries of the jump are exponentially small. Indeed
\be
\label{28feb2017-1}
H_{ij}(z,\varepsilon)\equiv H_{ij}(\zeta,u)= s_{ij} \exp\Bigl\{(u_i-u_j)\zeta+((B_1)_{ii}-(B_1)_{jj})\ln \zeta\Bigr\}\longrightarrow \delta_{ij}.
\ee
This is due to the fact that $s_{ij}$ is either $(\mathcal{S}_\nu)_{ij}$ or $(\mathcal{S}_{\nu-\mu}^{-1})_{ij}$. Thus,  $\delta H_{ij}\in L^2(\ell_{\mathcal{A}} \cup \ell_{\mathcal{B}}, |d\zeta|)$, and $C_{-}\bigl[\delta H    \bigr]_{ij}\in L^2(\ell_{\mathcal{A}} \cup \ell_{\mathcal{B}}, |d\zeta|)$. Hence,  the problem is well posed in $L^2$, consisting in finding $\delta\Psi_-$ as the solution of 
\be
\label{18aprile2017-1}
\left( I- C_{-}\Bigl[ \bullet~\delta H(\cdot,\varepsilon)\Bigr]\right)\delta\Psi_-= C_{-}\Bigl[\delta H(\cdot,\varepsilon) \Bigr].
\ee

\item[2.] If $u$ is in the cell containing $\mathcal{V}$, by assumption both the operator   and the given term  in (\ref{18aprile2017-1}) depend holomorphically on $u$.  
 Along the loops $(u_i-u_j)\mapsto (u_i-u_j) e^{2\pi i}$, $1\leq i\neq j\leq n$, the property (\ref{28feb2017-1}) is lost,  because $u$ leaves the $\widetilde{\tau}$-cell containing $\mathcal{V}$, so that some Stokes rays cross the ray $R(\widetilde{\tau})$. On the other hand, if the vanishing condition (\ref{6gen2017-1-BBBIS}) holds, then  $s_{ab}=s_{ba}=0$.\footnote{No difficulty arises from the fact that $\mathbb{S}_{\nu-\mu}^{-1}$ appears. If for simplicity we take the labelling (\ref{29gen2016-1})-(\ref{29gen2016-3}), then $\mathbb{S}_{\nu-\mu}$ has diagonal blocks equal to $p_j\times p_j$ identity matrices. This structure  persists on taking the inverse.}   Thus, (\ref{28feb2017-1}) continues to hold along the loop $\varepsilon \mapsto
\varepsilon e^{2\pi i}$.   It follows that     $ I- C_{-}\bigl[ \bullet~ \delta  H(\cdot,\varepsilon)\bigr]$ is an  analytic operator in $\varepsilon$ and the term $C_{-}\bigl[\delta  H(\cdot,\varepsilon) \bigr]$ is also  analytic, for $\varepsilon$  belonging to a sufficiently small  closed ball $U$ centred at  $\varepsilon=0$.

\item[3.] If $P_1$ and $P_2$ are sufficiently far away from the origin, we can take $\| \delta H(\cdot ,\epsilon) \|_{\infty}=\sup_{\zeta \in \ell_{\mathcal{A}} \cup \ell_{\mathcal{B}}} | H(\zeta ,\epsilon)| $  so small  that the operator norm $\|\cdot \|$ in $L^2$   satisfies, for $\varepsilon \in U$,  
\be
\label{24aprile2017-3}
\left\|C_{-}\bigl[ \bullet~\delta H(\cdot,\varepsilon)\bigr]\right \|\leq \|C_{-}\| ~\|\delta H(\cdot ,\epsilon) \|_{\infty} <1. 
\ee
Here, $\|C_{-}\|$ is the operator norm of the Cauchy operator.\footnote{Here we use the simple estimate  $\|C_-(f\delta H)\|_{L^2}\leq \|C_-\| ~\|\delta H\|_{\infty} ~ \|f\|_{L^2}$, for any $f\in L^2$.}  
By (\ref{24aprile2017-3}),  the inverse exists: 
\be
\label{18aprile2017-2}
\left( I- C_{-}\Bigl[ \bullet~\delta H(\cdot,\varepsilon)\Bigr]\right)^{-1}=\sum_{k=1}^{+\infty}\left( C_{-}\bigl[ \bullet~\delta H(\cdot,\varepsilon)\bigr]\right)^k.
\ee
The series in the r.h.s. converges in operator norm and defines  an analytic operator in $\varepsilon\in U$. 
\end{itemize}

\noindent 
Using (\ref{18aprile2017-2}),  we find the unique $L^2$-solution of  (\ref{18aprile2017-1}) and then, substituting into (\ref{28feb2017-3}), we find the ordinary  solution $\Psi(z,\varepsilon)$ of the auxiliary problem, which is holomoprhic in $\varepsilon\in U$.

\begin{figure}
\centerline{\includegraphics[width=0.6 \textwidth]{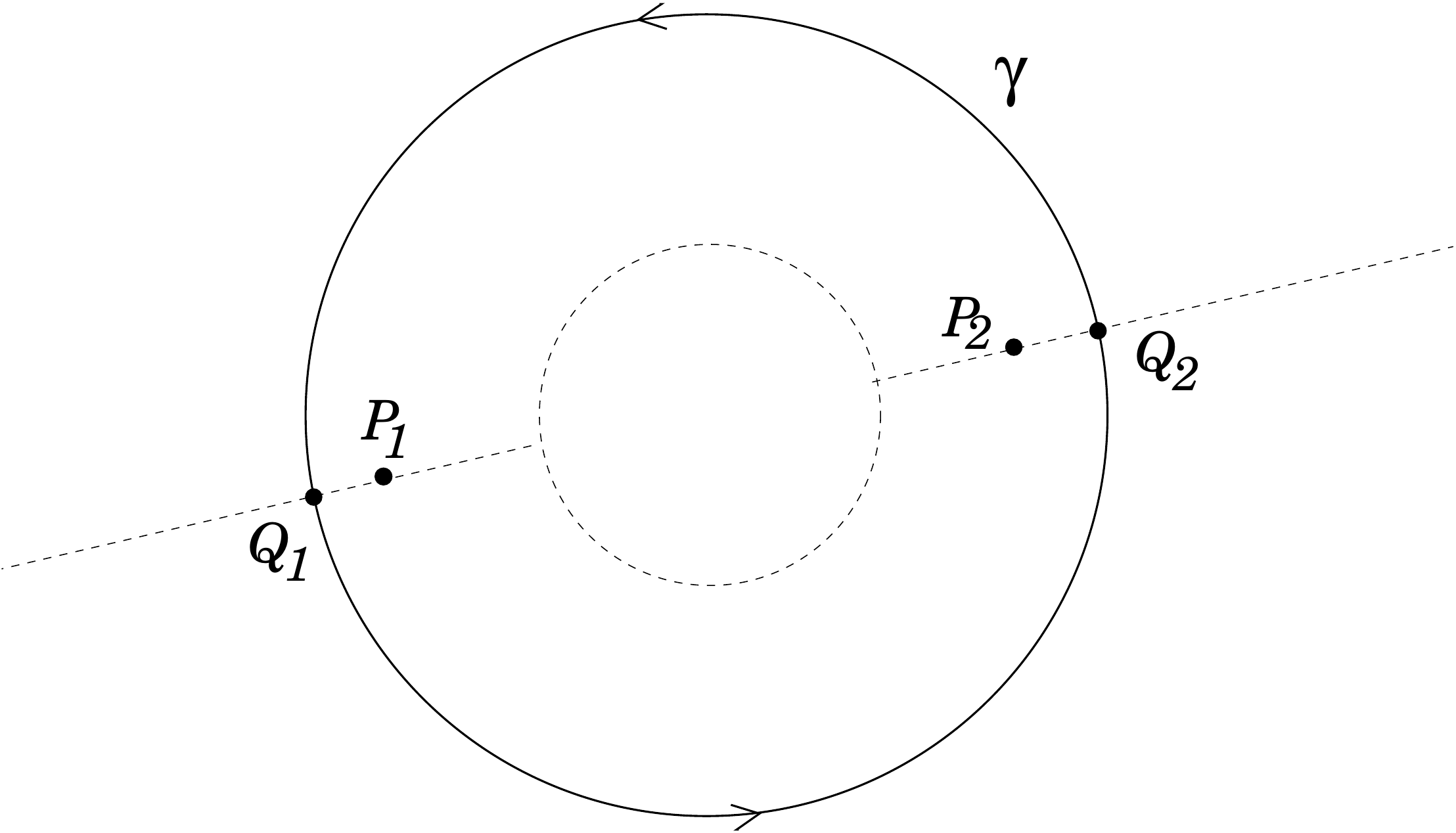}}
\caption{Step 3: the continuous Riemann-Hilbert problem on the circle $\gamma$, with jump $ \Psi(\zeta,\varepsilon)\mathfrak{S}(\zeta,\varepsilon)^{-1}$.} 
\label{RH-4}
\end{figure}

\vskip 0.2 cm 
$\bullet$ {\bf Step 3:}  We construct a R-H along a closed contour with a {\it continuous} jump. Consider a ``big'' counter-clockwise oriented circle $\gamma$ centered at the origin and intersecting $\Gamma_{-\infty}$ at  a point $Q_1$ preceding $P_1$,  $\Gamma_{+\infty}$ at a point $Q_2$ following  $P_2$.  See Figure \ref{RH-4}.  If $\mathcal{G}$ is the solution to the starting problem (\ref{19aprile2017-1}), (\ref{2marzo2017-2}), (\ref{19aprile2017-2}),  we construct a matrix-valued function $\Phi$ as follows: 
\bea
\label{19aprile2017-8}
\Phi:=&& \mathcal{G} \cdot \Psi(z,\varepsilon)^{-1},
\quad
~ \hbox{for $z$  outside $\gamma$}, 
\\
\label{19aprile2017-9}
&& 
\mathcal{G}\cdot \mathfrak{S}(z,\varepsilon)^{-1},
\quad
~ \hbox{for $z$  inside $\gamma$}.
\eea
By constriction,  $\Phi$ only has jumps along $\gamma$:
\be
\label{24aprile2017-1}
\Phi_{+}(\zeta)~=\Phi_-(\zeta) ~\widetilde{H}(\zeta,\varepsilon),
\quad
\quad
\widetilde{H}(\zeta,\varepsilon):= \Psi(\zeta,\varepsilon)\mathfrak{S}(\zeta,\varepsilon)^{-1}.
\ee
By construction, the jump matrix $\widetilde{H}(\zeta,\varepsilon)$ is continuous in $\zeta$ along $\gamma$, and is analytic in $\varepsilon\in U$. By  (\ref{19aprile2017-4}), then (\ref{2marzo2017-2}) is equivalent to
$$
\Phi(z)\sim I+ \hbox{series in }z^{-1},
\quad
~z\to\infty,
\quad
~z\in\Pi_\nu\cup \Pi_{\nu+\mu}.
$$
Therefore, the R-H for $\Phi$ is solved as in (\ref{18aprile2017-1}) and  (\ref{28feb2017-3}) by
\bea
\label{19aprile2017-5}
&&
\left( I- C_{-}\Bigl[ \bullet~\delta \widetilde{H}(\cdot,\varepsilon)\Bigr]\right)\delta\Phi_-= C_{-}\Bigl[\delta \widetilde{H}(\cdot,\varepsilon) \Bigr],
\\
\label{19aprile2017-6}
&&
\delta\Phi(z)= 
\int_\gamma
\frac{\delta\Phi_-(\zeta)~\delta \widetilde{H}(\zeta,\varepsilon)}{\zeta -z }\frac{d\zeta}{2\pi i }+\int_\gamma
\frac{\delta \widetilde{H}(\zeta,\varepsilon)}{\zeta -z }\frac{d\zeta}{2\pi i }.
\eea
Here $C_-$  is  Cauchy operator along $\gamma$. Since $\gamma$ is a closed contour and $\widetilde{H}(\zeta,\varepsilon)$ is continuous, the procedure and results of \cite{Zhou} \cite{LpDZ}  \cite{DKMVZ}  apply. The operator $C_{-}\Bigl[ \bullet~\delta \widetilde{H}(\cdot,\varepsilon)\Bigr]$ is Fredholm, $I-C_{-}[ \bullet~\delta \widetilde{H}(\cdot,\varepsilon)]$ has index $0$ and its kernel is $\{0\}$.  Therefore, the ``analytic Fredholm  alternative''  of \cite{Zhou}  holds. Namely,  either   $I-C_{-}[ \bullet~\delta \widetilde{H}(\cdot,\varepsilon)]$ can be inverted (and  (\ref{19aprile2017-5}) can be  solved) for every $\varepsilon \in U$, except for a finite number of isolated values, or is invertible for no   $\varepsilon$. In the first case, $ \left(I-C_{-}[ \bullet~\delta \widetilde{H}(\cdot,\varepsilon) ]\right)^{-1}$ is meromorphic, with poles at the isolated points in $U$. 

By (\ref{19aprile2017-8})-(\ref{19aprile2017-9}), solvability of the R-H  (\ref{24aprile2017-1}) is equivalent to the existence of the solution $\mathcal{G}(z,\varepsilon)\equiv \mathcal{G}(z,u)$ for the problem (\ref{19aprile2017-1}), (\ref{2marzo2017-2}), (\ref{19aprile2017-2}). By assumption (i.e. by the result of \cite{Miwa}) we know that locally in $u$ the solution $\mathcal{G}(z,u)$ exists. We therefore conclude that the ``Fredholm analytic alternative'' implies the existence of the solution $\Phi_-(\zeta,\varepsilon)$  of  (\ref{19aprile2017-5})   for every $\varepsilon \in U$, except for a finite number of poles, and that (\ref{19aprile2017-6}) gives an ordinary solution $\Phi(z,\varepsilon)$, meromoprhic as a function of $\varepsilon$ in  $U$. By (\ref{19aprile2017-8})-(\ref{19aprile2017-9}), the same conclusion holds for $\mathcal{G}(z,\varepsilon)\equiv \mathcal{G}(z,u)$. This proves the Lemma (as for  $\widehat{A}_1$, it suffices to note that $\widehat{A}_1(t)=z(Y^{-1}(z,t) dY(z,t)/dz-\Lambda(t))$). 
 $\Box$

\vskip 0.3 cm 
\noindent
{\bf Theorem  \ref{9gen2017-1}}  immediately follows from  Lemma \ref{8gen2017-9} and Lemma \ref{27feb2017-1}.
\vskip 0.5 cm
\hrule
\vskip 0.5 cm 
\centerline{\bf \Large PART V: Examples and  Applications}

\vskip 0.3 cm 

Our work is motivated both by the general deformation problems of linear systems with coalescing eigenvalues and by applications to Frobenius Manifolds and  Painlev\'e equations.  The applications  are sketched in the sections below, which  are a natural continuation of the Introduction, of which we keep the notations (for example, $\mathbb{S}_1$, $\mathbb{S}_2$ instead of $\mathbb{S}_{\nu+k\mu}$).

\section{Monodromy Data of Semisimple Frobenius Manifolds}
\label{16dicembre2016-5} 
  A  Frobenius manifold $M$ essentially is an analytic manifold with {a Frobenius algebra} structure {on} the tangent bundle and a deformed flat connection {(see \cite{Dub1} \cite{Dub2} for a precise definition)}. The manifold is called semisimple if the algebra is semisimple in an open dense subset, whose points are called semisimple points of $M$. In suitable coordinates $u=(u_1,...,u_n)$, called {\it canonical}, flatness  is translated into   $n+1$   {\it compatible} linear systems  of dimension  $n\times n$ 
\bea
\label{23gen2016-1}
&&
\frac{\partial Y}{\partial z}=\widehat{A}(z,u)Y,
\quad
\quad
\widehat{A}(z,u):=U+\frac{V(u)}{z}
 \\
\label{23gen2016-2}
&&
\frac{\partial Y}{\partial u_k}=\Omega_k(z,u)Y,
\quad
\quad
\Omega_k(z,u):=zE_k+V_k(u),
\quad
\quad
1\leq k \leq n.
\eea
Here  $E_k$ is the matrix with zero entries except for $(E_k)_{kk}=1$, 
 $
U=\hbox{diag}(u_1,...,u_n)$ and $V$ is skew-symmetric.  The system (\ref{23gen2016-1}) is of type (\ref{ourcase12}). If we write  $u=u(t)$ as in  (\ref{4gen2016-2}), then the following identification holds 
$$ 
U\equiv \Lambda(t)  ,
\quad
\quad
V(u(t))\equiv \widehat{A}_1(t).
$$
The matrices $V_k$ are defined by 
$$ 
V_k(u)=\frac{\partial \Psi(u)}{\partial u_k}~\Psi(u)^{-1}.
$$ 
The matrix $\Psi(u)$  gives the change of  basis between flat and canonical coordinates, according to the formulae in Exercise 3.2 of \cite{Dub2}.  It is crucial for our discussion that $\Psi(u)$ is always  holomorphic and invertible  at semisimple points,   also when $U$ has coalescing eigenvalues there. The proof is given in  \cite{OURPAPER}. Therefore,  the matrices $V_k(u)$ are holomorphic at semisimple points. 
  $\Psi(u)$ diagonalises $V(u)$, with constant   eigenvalues $\mu_1,...,\mu_n$ independent of the point of the manifold (see \cite{Dub1},\cite{Dub2}):
    $$ 
V(u)=\Psi(u)~\mu ~ \Psi(u)^{-1},
\quad
\quad
\mu:=\hbox{diag}(\mu_1,\mu_2,...,\mu_n). 
$$
Therefore, $V(u)$ is  holomorphically similar to  $\mu$ at semisimple points. 

The system (\ref{23gen2016-1}) admits a normal form at $z=0$ such that the corresponding   fundamental matrix,  denoted  
 \be
 \label{11novembre2016-3}
 Y_0(z,u)=\Bigl(\Psi(u)+\sum_{l=1}^\infty \Phi_l(u) z^l\Bigr)z^{\mu} z^{R},
 \ee
   has monodromy exponent  $R$  independent  of the point of the manifold.  $Y_0(z,u)$ is  holomorphic  of $u$ on the domain where  $V(u)$ is holomorphic.  In our notations, $R\equiv R^{(0)}$, and $Y_0\equiv Y^{(0)}$, as in (\ref{24gen2016-1}).

The system (\ref{23gen2016-1}), (\ref{23gen2016-2}) is the system  (\ref{ourcase12}), (\ref{3nov2015-3})  (let $t_k\mapsto u_k$ in (\ref{3nov2015-3})). The compatibility condition reads
\bea
\label{17gen2016-1} 
&&
[U,V_k]=[E_k,V],
\quad
~\Longrightarrow
\quad
~ (\delta_{ki}-\delta_{kj})V_{ij}=(u_i-u_j)(V_k)_{ij},
\quad
~1\leq i,j,k \leq n;
\\
&&
\label{23gen2016-3}
\frac{\partial V}{\partial u_k}=[V_k,V].
\eea 
Equations (\ref{23gen2016-3}) coincide with the isomonodromy deformation equations (\ref{23gen2016-4}) and $V_k(u)$ coincides with the matrix (\ref{3nov2015-2}).

Next, we establish the translation  between  our  Stokes and central connection matrices and those defined in \cite{Dub2}. 
Following \cite{Dub2}, Section 4,  we consider an oriented ray  $\ell_{+}(\phi):=\{z\in\mathcal{R}~|~\arg z=\phi\}$ and  (for $\epsilon>0$ small)  the following two sectors 
$$ 
\Pi_{\rm right}:=S(\phi-\pi-\epsilon,\phi+\epsilon),
\quad
\quad
\Pi_{\rm left}:=S(\phi-\epsilon,\phi+\pi+\epsilon).
$$
In \cite{OURPAPER}, we  introduce the open dense subset of points $p\in M$ such that the eigenvalues of $U$ at $p$ are pairwise distinct and  no Stokes rays associated with $U$ at $p$ coincide with $\ell_{+}(\phi)$, and  we call any connected component of this set an {\it $\ell$-chamber}. 
Let $\mathcal{V}$ be an open connected  domain such that $\overline{\mathcal{V}}$ is contained in an $\ell$-chamber. For suitable $\epsilon$, we can identify\footnote{In the notation used in the main body of the paper, 
$$ 
e^{-2\pi i}\Pi_{\rm left}= \mathcal{S}_{\nu}(\overline{\mathcal{V}}),
\quad
~~\Pi_{\rm right}= \mathcal{S}_{\nu+\mu}(\overline{\mathcal{V}}),
\quad
~~\Pi_{\rm left}= \mathcal{S}_{\nu+2\mu}(\overline{\mathcal{V}}),
\quad
~~\hbox{ for } \tau_\nu<\widetilde{\tau}<\tau_{\nu+1}.
$$
}
\be
\label{1giugno2017-1}
e^{-2\pi i}\Pi_{\rm left}= \mathcal{S}_1(\overline{\mathcal{V}}),
\quad
\quad
\Pi_{\rm right}= \mathcal{S}_2(\overline{\mathcal{V}}),
\quad
\quad
\Pi_{\rm left}= \mathcal{S}_3(\overline{\mathcal{V}}),
\ee
where $e^{-2\pi i}\Pi_{\rm left}:=\{z\in \mathcal{R}~| z=\zeta e^{-2\pi i}, ~\zeta\in \Pi_{\rm left}\}$, and $\mathcal{S}_r(\overline{\mathcal{V}})$ is defined in the Introduction. Let $Y_{\rm left}(z,u)$,  $Y_{\rm right}(z,u)$ be the unique fundamental matrix solutions having the canonical asymptotics 
$Y_F(z,u)=(I+O(1/z))e^{zU}$ in $\Pi_{\rm left}$  and $\Pi_{\rm right}$ respectively.  The Stokes matrices $S$ and $S_{-}$ of  \cite{Dub2} are defined by the relations,
\be
\label{23gen2016-6}
Y_{\rm left}(z,u)=Y_{\rm right}(z,u)S,~
\quad
\quad
Y_{\rm left}(e^{2\pi i}z,u)=Y_{\rm right}(z,u)S_{-},
\quad
\quad
z\in \mathcal{R}.
\ee
The symmetries of the system (\ref{23gen2016-1}) imply that $S_{-}=S^T$.
In our notations as in  (\ref{2maggio2016-4}),  the Stokes matrices are defined by 
\be
\label{23gen2016-7}
Y_3(z,u)=Y_2(z,u)\mathbb{S}_2,
\quad
\quad
Y_2(z,u)=Y_1(z,u) \mathbb{S}_1.
\ee
We identify 
\be
\label{13giugno2017-5}
Y_3(z,u)=Y_{\rm left}(z,u),
\quad
\quad
Y_2(z,u)=Y_{\rm right}(z,u)
\ee
Let $B_1$ denote the exponent of   formal monodromy\footnote{In general, a formal solution is $(I+\sum_{k=1}^\infty F_k(u)z^{-k})z^{B_1} e^{zU}$, but in case of Frobenius manifolds $B_1=0$.} at $z=\infty$, so that the  relation 
 $
Y_1(ze^{-2\pi i},u)=Y_3(z,u) e^{-2\pi i B_1}
$ holds.\footnote{In the notation of the main body of the paper, $Y_r\mapsto Y_{\nu+(r-1)\mu}$,  $r=1,2,3$, $\mathbb{S}_1\mapsto \mathbb{S}_{\nu}$, $\mathbb{S}_2 \mapsto \mathbb{S}_{\nu+\mu}$ and 
 $
Y_{\nu}(z_{(\nu)})=Y_{\nu+2\mu}(z_{(\nu+2\mu)})e^{-2\pi i L}$, where $z_{(\nu+(r-1)\mu)}\in \mathcal{S}_{\nu+(r-1)\mu}(\overline{\mathcal{V}})$ is seen as a point of $\mathcal{R}$ and not of $\mathbb{C}$.  
}
Since $V$ is skew symmetric and  $B_1=\hbox{ diag }(V)=0$,   the above relation  reduces to 
$$ 
Y_1(ze^{-2\pi i},u)=Y_{\rm left}(z,u).
$$
 Therefore (\ref{23gen2016-7}) coincides with (\ref{23gen2016-6}), with 
$$ 
S_{-}=\mathbb{S}_1^{-1},
\quad
~S=\mathbb{S}_2.
$$
The central connection matrix such that $Y_1=Y^{(0)}C^{(0)}$ was defined in (\ref{2maggio2016-5}) and in Definition \ref{13giugno2017-7}. In the theory of Frobenius manifolds, such as in  \cite{OURPAPER}, the central connection {matrix} is denoted  by $C$, defined by  
 $$Y_{\rm right}(z,u)=Y_0(z,u) C.
 $$
  Since  $Y_0=Y^{(0)}$,  $Y_{\rm right}=Y_2$,  $Y_2=Y_1 \mathbb{S}_1$, and $\mathbb{S}_1^{-1}=S^T$,  then 
  $$C^{(0)}=C \mathbb{S}_1^{-1}=CS^T
  .
  $$
  
Summarising,  monodromy data of a Frobenius manifold are $\mu$, $R$, $S$, $C$, versus the monodromy data $\mu_1,...,\mu_n$, $R^{(0)}$, $\mathbb{S}_1,\mathbb{S}_2$, $C^{(0)}$ of the present paper ($B_1=0$). 
\vskip 0.2 cm 
Coalescence points for $U$ in (\ref{23gen2016-1}) are singular points for the monodromy preserving deformation equations (\ref{17gen2016-1}) and  (\ref{23gen2016-3}). Their study is at the core of the analytic continuation of Frobenius structures. Our Theorem \ref{16dicembre2016-1} allows to extend the isomonodromic approach to Frobenius manifolds at coalescence points if the manifold is semisimple at these points.  Let $u^{(0)}=(u_1^{(0)},...,u_n^{(0)})$  denote  a  coalescence point. By a change $Y\mapsto PY$  in (\ref{23gen2016-1}),  given by a permutation matrix $P$,  there is no loss of generality in assuming that  
\beas
u_1^{(0)}=\cdots=u_{p_1}^{(0)}=:\lambda_1
\\
u_{p_1+1}^{(0)}=\cdots=u_{p_1+p_2}^{(0)}=:\lambda_2
\\
 \vdots&&
\\
u_{p_1+\cdots+p_{s-1}+1}^{(0)}=\cdots=u_{p_1+\cdots+p_{s-1}+p_s}^{(0)}=:\lambda_s,
\eeas
where  $p_1,...,p_s$ are  integers such that  $p_1+\cdots+p_s=n$, and $\lambda_j\neq \lambda_k$ for $j\neq k$. In order to have a correspondence with \cite{Dub2}, as in formula (\ref{1giugno2017-1}) and (\ref{13giugno2017-5}), we  take the ray $\ell_+(\phi)$ with 
\be
\label{24gen2016-7}
\phi=\widetilde{\tau}+\pi 
\quad
~\hbox{ mod } 2\pi,
\ee
 where $\widetilde{\tau}$ is the direction of an  admissible ray  for $U$ at the point $u^{(0)}$, as in Definition \ref{14feb2016-5}.  
Similarly to  (\ref{22marzo2016-10}), we consider positive numbers $\delta_k$ and $\epsilon_0$ as follows 
\be
\label{11novembre2016-1}
\delta_k:=\frac{1}{2}\min_{j\neq k} \left\{
 \left|
\lambda_k-\lambda_j +\rho e^{i\left(\frac{\pi}{2}-\phi\right)}\right|, ~\rho\in\mathbb{R}
\right\},
\quad
\quad
0<\epsilon_0<\min_{1\leq k\leq s} \delta_k.
\ee
Consider the neighbourhood (polydisc) of $u^{(0)}$ defined by
$$ 
\mathcal{U}_{\epsilon_0}(u^{(0)}):=
\Bigl\{ u\in\mathbb{C}^n~\Bigl|~ |u-u_0|\leq \epsilon_0\Bigr\}
$$
and denote by $\Delta$ the coalescence locus passing through $u^{(0)}$, namely
$$
\Delta:=\{u(p)\in \mathcal{U}_{\epsilon_0}(u^{(0)})~|~ u_i=u_j \hbox{ for some } i\neq j\}.
$$

 If  $u^{(0)}$ is a semisimple coalescence point, then the Frobenius Manifold $M$   is semisimple in $\mathcal{U}_{\epsilon_0}(u^{(0)})$ for sufficiently small $\epsilon_0$ (if necessary, we further restrict $\epsilon_0$ in (\ref{11novembre2016-1})). 
Given the above assumption of semisimplicity,  then $\Psi(u)$ is holomorphic at $\Delta$ and this implies that $V(u)$ is {holomorphically} similar to $\mu$. Equation (\ref{17gen2016-1}) for $k=i$ is $
V_{ij}=(u_i-u_j)(V_i)_{ij}
$, 
which implies that $V_{ij}(u)=0$ for $i\neq j$ and $u_i=u_j$. Therefore, recalling that $V(u^{(0)})$    corresponds to  $\widehat{A}_1(0)$, we conclude that  the vanishing condition (\ref{31luglio2016-1}) holds true and then our  Theorem \ref{16dicembre2016-1} applies. 
 We note that  $\hbox{diag}\left(V(u^{(0)})\right)=0$, then the diagonal entries of $\widehat{A}_1(0)$ do not differ by non-zero integers, so that  also  Corollary \ref{17dicembre2016-2} applies.  Then, the following holds: 

\vskip 0.2 cm 
\noindent
\begin{shaded}
{\bf Theorem \ref{16dicembre2016-1} as applied to Frobenius Manifolds:} [More details in \cite{OURPAPER}] {\it Let the Frobenius manifold $M$ be semisimple in a neighborhood $\mathcal{U}_{\epsilon_0}(u^{(0)})$ of a coalescence point $u^{(0)}$, where $\epsilon_0$ satisfies (\ref{11novembre2016-1}). Then the constant monodromy data $\mu$ $R$, $S$, $C$  of the manifold are well defined  in the whole  $\mathcal{U}_{\epsilon_1}(u^{(0)})$,  for any  $\epsilon_1<\epsilon_0$. In particular, they are well defined at $u^{(0)}$ and at all  points of $\Delta$.
 These data coincide with the data 
of the system (\ref{23gen2016-1}) computed at  fixed $u=u^{(0)}$, as explained above. 
 }
\end{shaded}

 We recall that  the monodromy data for the whole manifold can be computed  by an action of the braid group (see \cite{Dub1}, \cite{Dub2} and \cite{OURPAPER}) staring from  the data obtained in $\mathcal{U}_{\epsilon_1}(u^{(0)})$. Hence, our result allows to obtain the monodromy data for the whole  manifold  from  the data computed at a  coalescence point.  
 This relevant fact is important in the {following two} cases:
 
 \vskip 0.2 cm 
 
 a) The Frobenius structure (i.e. $V(u)$ in   (\ref{23gen2016-1}) ) is known everywhere, but the computation of monodromy data is  extremely difficult -- or impossible --   at generic semisimple points where  $U=\hbox{diag}(u_1,...,u_n)$ has distinct eigenvalues.  On the other hand,  the system  (\ref{23gen2016-1}) at a coalescence point simplifies, so that we may be able to explicitly  solve it  in terms of special functions and compute $S$ and $C$.   
  In \cite{OURPAPER} we give a detailed example of this kind for the Frobenius manifold associated {with} the Coxeter group $A_3$.

 \vskip 0.2 cm 

b) The Frobenius structure is explicitly known {\it   only at points where  $U$ has two or more non-distinct eigenvalues}.  
The quantum cohomology of Grassmannians  Gr$(k,n)$ are important examples of this case:  
 the explicit form of $V(u)$ is known only along the  small quantum cohomology, where two eigenvalues of $U$ may coincide, depending on $k$ and $n$. Indeed, coalescence of eigenvalues is the most frequent case \cite{Cotti1}.   $S$ and $C$ can be explicitly computed at the small quantum cohomology locus and Theorem \ref{16dicembre2016-1} allows  their extension to the whole manifold.  In \cite{OURPAPER} we do explicit computations for  Gr$(2,4)$. 
\vskip 0.2 cm


\section{Computation of Monodromy Data of Painlev\'e Transcendents.  Example of  the  Algebraic Solution associated {with}  $A_3$}
\label{31luglio2016-2}

Theorem \ref{16dicembre2016-1} provides an alternative to Jimbo's approach   for the computation  of the monodromy data associated {with}  
Painlev\'e 6 transcendents holomorphic at a critical point.  The  example below refers to the $A_3$-algebraic solution of 
 \cite{DM}.\footnote{The example,  reinterpreted in the framework of Frobenius manifolds,  
gives the analytic computation of  Stokes matrices for the $A_3$-Frobenius manifold.}

\vskip 0.2 cm 
Equations (\ref{17gen2016-1}) and skew symmetry of $V(u)$  imply  that   $\sum_i \partial_i V= \sum_i u_i
 \partial_i V =0$. Thus if $n=3$, 
  $$
  V(u_1,u_2,u_3)\equiv V(t),
  \quad
  \quad 
 t:={u_2-u_1\over u_3-u_1},
 $$
Write 
$$
V(t)=
\left( 
\begin{array}{ccc}
0 & \Omega_2 & -\Omega_3 \cr
                 - \Omega_2 & 0 & \Omega_1 \cr
                  \Omega_3 & -\Omega_1 & 0 \cr
\end{array}
\right).
$$
Suppose we want to study the coalescence $u_2-u_1\to 0$ in the system (\ref{23gen2016-1}), with $u_3-u_1\neq 0$. With the substitutions $Y(z)\mapsto e^{u_1 z}Y(z)$, and $z\to (u_3-u_1)z$,  (\ref{23gen2016-1}) becomes
  \be
\label{28feb2016-1}
\frac{dY}{d z}=\left[
\left(
\begin{array}{ccc}
0 & 0& 0
\\
0 & t & 0 
\\
0 & 0 & 1
\end{array}
\right)+\frac{V(t)}{z}\right]Y.
\ee 
 The coalescence $u_2-u_1\to 0$ corresponds to $t\to 0$. 

 In equations (\ref{23gen2016-3}), write $\partial V/\partial u_k=dV/dt \cdot \partial t/\partial u_k$, in order to obtain the following  equivalent equations
\be 
\label{equationV}
{d \Omega_1\over dt} = {1 \over t} ~ \Omega_2 \Omega_3, ,
\quad
\quad
{d \Omega_2\over dt} = {1 \over 1-t} ~ \Omega_1 \Omega_3, 
\quad
\quad
{d \Omega_3\over dt} = {1 \over t( t-1)} ~ \Omega_1 \Omega_2.
\ee
$V(t)$ can be expressed in terms of  transcendents $y(t)$  satisfying the following 
Painlev\'e 6 equation, called $PVI_\mu$ (see \cite{Dub2}, and also  \cite{davidG1} for an asymptotic study of (\ref{equationV})): 
{\small $$
{d^2y \over dt^2}={1\over 2}\left[ 
{1\over y}+{1\over y-1}+{1\over y-t}
\right]
           \left({dy\over dt}\right)^2
-\left[
{1\over t}+{1\over t-1}+{1\over y-t}
\right]{dy \over dt}
+{1\over 2}
{y(y-1)(y-t)\over t^2 (t-1)^2}
\left[
(2\mu-1)^2+{t(t-1)\over (y-t)^2}
\right]
,
$$}

\noindent
with parameter $\mu\in\mathbb{C}$. The eigenvalues of $V(t)$ are $\mu,0,-\mu$.  The following are the explicit formulae (see \cite{davidG0}).

\beas 
\Omega_1=i { \sqrt{y-1}\sqrt{y-t} \over \sqrt{t}} \left[
{A\over (y-1)(y-t)} + \mu\right],
&&
\Omega_2=i { \sqrt{y}\sqrt{y-t} \over \sqrt{1-t}} \left[
{A\over y(y-t)} + \mu\right],
\\
\\
\Omega_3=- { \sqrt{y}\sqrt{y-1} \over \sqrt{t}\sqrt{1-t}} \left[
{A\over y(y-1)} + \mu\right]
,
&&
 A:= {1\over 2} \left[ {dy\over dt} t(t-1)-y(y-1)\right].
\eeas
The branches (signs) in the square roots above are arbitrary. A change of the
 sign of one root  implies a change of two signs in
 $(\Omega_1,\Omega_2,\Omega_3)$, which still yields a solution of  (\ref{equationV}). 

The {``Painlev\'e transcendent"} corresponding to {the} $A_3$-Frobenius manifold is the following  algebraic solution of $PVI_\mu$, $\mu=-\frac{1}{4}$, obtained in \cite{DM} (there is a misprint  in $t(s)$ {in} \cite{DM}),
\be 
\label{para1}
y(s)={(1-s)^2~(1+3s)~(9s^2-5)^2
\over 
(1+s)~(243s^6+1539s^4-207s^2+25)}
,\quad\quad
t(s)=
{(1-s)^3~(1+3s)
\over (1+s)^3~(1-3s)
}~.
\ee
As it is shown in \cite{DM}, the Jimbo's monodromy data  of  the Jimbo-Miwa-Ueno   isomonodromic Fuchsian system associated {with} algebraic solutions of $PVI_\mu$ are  tr$(M_iM_j)=2-S_{ij}^2$,  $1\leq  i<j \leq 3$, where $S$ is the Stokes matrix (in upper triangular form)  of the corresponding Frobenius manifold. $S$ is well  known \cite{Dub1}, and $S+S^T$ is the Coxeter matrix of the reflection group $A_3$.    Moreover, Jimbo's isomonodromic method  \cite{Jimbo}, as applied in \cite{DM} (see also \cite{Kaneko}, \cite{D2} for holomorphic solutions) provides tr$(M_iM_j)$.  Here we apply Theorem \ref{16dicembre2016-1} and  obtain $S$ in  an alternative, and probably simpler, way. 

First, we take an holomorphic branch. It is obtained by letting $s\to -{1\over 3}$, which gives a  convergent Taylor expansion at $t=0$: 
{\small
\be 
\label{tayl}
y(t)= {1\over 2}t+{13\over 32}t^2+{13\over 64}t^3+{201\over 4096}t^4-{229\over 8192}t^5
-{101055\over 2097152}t^6-{167867\over 4194304}t^7
-{3235319\over 134217728}t^8 +O(t^9).
\ee
}

\noindent
Substitution of the parametric formulae (\ref{para1}) -- or  equivalently of (\ref{tayl}) --  into ({\ref{equationV}) yields  (two changes of signs   are allowed), 
\beas
&&
\Omega_1(t)=i\sqrt{2}\left({1\over 8}-{1\over 256}t
-{17\over 16384}t^2-{257\over 524288}t^3+O(t^4)\right),
\\
&&
\Omega_2(t)=-{1\over 32}t-{1\over 64}t^2-{173\over 16384}t^3+O(t^4),
\\
&&
\Omega_3(t)=i\sqrt{2}
\left({1\over 8}
+{1\over 256}t
+{47\over 16384}t^2+{1217\over 524288}t^3+O(t^4)\right).
\eeas
We observe that the following limits exist:
$$
\lim_{t\to 0 }\Omega_1(t)= {i\over 4\sqrt{2}}
,\quad
\quad
\lim_{t\to 0 }\Omega_2(t)= 0
,
\quad
\quad
\lim_{t\to 0 }\Omega_3(t)= {i\over 4\sqrt{2}}.
$$
Thus, the assumptions of Theorem \ref{16dicembre2016-1} hold, because   $\Omega_2(t)\to 0$ as $t\to 0$. Also the Corollary \ref{17dicembre2016-2} holds, because diag$(V)=(0,0,0)$.    Accordingly,  the Stokes matrices corresponding to (\ref{tayl}) for $|t|$ small  can be computed using  (\ref{28feb2016-1}) at $t=0$, namely:
\be
\label{28feb2016-3}
\frac{dY}{dz}=\left[
\left(
\begin{array}{ccc}
0 & 0& 0
\\
0 & 0 & 0 
\\
0 & 0 & 1
\end{array}
\right)
+
\frac{V(0)}{z}
\right]Y,
\quad
\quad
V(0)=
\left(
\begin{array}{ccc}
0 & 0& -i\sqrt{2}/8
\\
0 & 0 & i\sqrt{2}/8
\\
i\sqrt{2}/8 & -i\sqrt{2}/8 & 0
\end{array}
\right).
\ee
This system  is integrable, as follows. First, we do a gauge trasformation  $Y=\Psi  \widetilde{Y}$,  such that  
\footnote{
Each columns of $\Psi$ can be multiplied by a constant. We have chosen $\Psi$ such that 
 $\Psi^T\Psi=
 \left(
 \begin{array}{ccc}
0& 0 & 1
\\
0 & 1 &0
\\
1 & 0 & 0
\end{array}
\right)
$. This has  a meaning in the framework of Frobenius manifolds, but is of no importance for our computation. }
$$
\Psi^{-1}V(0)\Psi=\hbox{diag}(-1/4,~0,~1/4),
\quad
\quad
\Psi=\left(
\begin{array}{ccc}
i/2 & 1/\sqrt{2} & -i/2
\\
-i/2 & 1/\sqrt{2} & i/2 
\\
1/\sqrt{2} & 0 & 1/\sqrt{2}
\end{array}
\right)
.
$$
Hence (\ref{28feb2016-3}) becomes 
$$
\frac{d\widetilde{Y}}{dz}=\left[
\left(
\begin{array}{ccc}
1/2 & 0 & 1/2 
\\
0 & 0 & 0
\\
1/2 & 0 & 1/2 
\end{array}
\right)
+\frac{1}{z}
\left(
\begin{array}{ccc}
-1/4 & 0 & 0 
\\
0 & 0 & 0
\\
0 & 0 & 1/4 
\end{array}
\right)
\right] \widetilde{Y}.
$$
We consider a column $(y_1,y_2,y_3)^T$ of $\widetilde{Y}$, so that 
$$
y_1^\prime=\frac{1}{2}(y_1+y_3)-\frac{1}{4z}y_1,
\quad
\quad
y_2^\prime=0,
\quad
\quad
y_3^\prime=\frac{1}{2}(y_1+y_3)+\frac{1}{4z}y_3.
$$
By elimination of $y_3(z)$ and setting  $y_1(z)=\sqrt{z/2} ~e^{z/2} w\left( iz/2\right)$ we find that the system reduces to  the Bessel equation 
$$ 
\xi^2 \frac{d^2w}{d\xi^2}+\xi \frac{dw}{d\xi} +\left[\xi^2-\left(\frac{3}{4}\right)^2 \right]w=0.
$$
The last equation is integrated in terms of Hankel functions, so that we find general solutions $y_1(z)$ of the form  
$$
y_1(z,c_1,c_2,m,n)=\sqrt{\frac{z}{2}}~e^{\frac{z}{2}}\left( c_1 H_{\frac{3}{4}}^{(1)}\left(\frac{iz}{2}e^{im\pi}\right)+c_2 H_{\frac{3}{4}}^{(2)}\left(\frac{iz}{2}e^{in\pi}\right)\right),
\quad
\quad
c_1,~c_2\in\mathbb{C},
\quad
~m,~n\in\mathbb{Z}.
$$
The Stokes rays of the system (\ref{28feb2016-3}) are given by $\Re(iz)=0$,  namely $
\arg z=\frac{\pi}{2}+k\pi$, $ k\in\mathbb{Z}$.  
Consider three canonical sectors 
$$ 
\mathcal{S}_1=S\Bigl( -\frac{3\pi}{2},\frac{\pi}{2}\Bigr),
\quad~~~
\mathcal{S}_2=S\Bigl( -\frac{\pi}{2},\frac{3\pi}{2}\Bigr),
\quad~~~
\mathcal{S}_3=S\Bigl( \frac{\pi}{2},\frac{5\pi}{2}\Bigr).
$$
The asymptotic behaviour of  fundamental matrices $\widetilde{Y}_r(z)=\Psi^{-1}Y_r(z)$ corresponding to canonical asymptotics of $Y_r(z)$ for $z\to \infty$ in $\mathcal{S}_r$, $r=1,2,3$, is of the type
\be
\label{28feb2016-8} 
\widetilde{Y}(z)=\Psi^{-1}\left(I+\mathcal{O}\left(\frac{1}{z}\right)\right)
\left(
\begin{array}{ccc}
1 & 0 & 0 
\\ 0  & 1 & 0 
\\ 0 & 0 & e^z
\end{array}
\right)
= \left(I+\mathcal{O}\left(\frac{1}{z}\right)\right)\left(
\begin{array}{ccc}
-\frac{i}{2} & \frac{i}{2} & \frac{1}{\sqrt{2}}e^z
\\
\\
 \frac{1}{\sqrt{2}} &  \frac{1}{\sqrt{2}} & 0
 \\
 \\
 \frac{i}{2}  & -\frac{i}{2} & \frac{1}{\sqrt{2}}e^z
\end{array}
\right).
\ee
Let us  compute $\mathbb{S}_1$, such that $Y_2(z)=Y_1(z)\mathbb{S}_1$. 
The behaviour for  $z\to \infty$ of {Hankel} functions is 
\beas
&&
H_{\frac{3}{4}}^{(1)}(iz/2)=\frac{2e^{-i\frac{7\pi}{8}} }{\sqrt{\pi z}}~\left(1+\mathcal{O}\left(\frac{1}{z}\right)\right)~e^{-z/2},
\quad
\quad 
-\frac{3\pi}{2}<\arg z<\frac{3\pi}{2},
\\
&&
H_{\frac{3}{4}}^{(2)}(iz/2)=\frac{2e^{i\frac{3\pi}{8}}}{\sqrt{\pi z}}~\left(1+\mathcal{O}\left(\frac{1}{z}\right)\right)~e^{z/2}, 
\quad
\quad
-\frac{5\pi}{2}<\arg z<\frac{\pi}{2} , 
\eeas
 It impliies that
$$\tilde{Y}_1(z)
= \left(
\begin{array}{ccc}
y_1(z,c_1^{(-)},0,0,0) & y_1(z,c_1^{(+)},0,0,0) &y_1(z,0,c_2,0,0)
\\
\star & \star & \star
\\
\star & \star & \star
\end{array}
\right)
$$
 with
\be
\label{28feb2016-5}
y_1(z,c_1^{(-)},0,0,0) =c_1^{(-)} \sqrt{\frac{z}{2}} e^{z/2}~H_{\frac{3}{4}}^{(1)}\left(\frac{iz}{2}\right)=-\frac{i}{2} \left(1+\mathcal{O}\left(\frac{1}{z}\right)\right),
\quad
\quad
c_1^{(-)}:=-\frac{i}{2}\sqrt{ \frac{\pi}{2} }~e^{7\pi i/8},
\ee

\be
\label{28feb2016-6}
y_1(z,c_1^{(+)},0,0,0)=c_1^{(+)} \sqrt{\frac{z}{2}} e^{z/2}~H_{\frac{3}{4}}^{(1)}\left(\frac{iz}{2}\right)=\frac{i}{2} \left(1+\mathcal{O}\left(\frac{1}{z}\right)\right),
\quad
\quad
c_1^{(+)}:=\frac{i}{2}\sqrt{ \frac{\pi}{2} }~e^{7\pi i/8},
\ee

\be
\label{28feb2016-12}
y_1(z,0,c_2,0,0)=c_2\sqrt{\frac{z}{2}} e^{z/2}~H_{\frac{3}{4}}^{(2)}\left(\frac{iz}{2}\right)=\frac{e^z}{\sqrt{2}} \left(1+\mathcal{O}\left(\frac{1}{z}\right)\right),
\quad
\quad
c_2:=\frac{\sqrt{\pi}}{2}e^{-3\pi i/8}.
\ee
The asymptotics of $H_{\frac{3}{4}}^{(1)}(iz/2)$ extends up to $\arg z=3\pi/2$. Therefore, the first two matrix entries in the first row of  $Y_2(z)$ are the same of $Y_1(z)$, which implies 
$$ 
\mathbb{S}_1=
\left(
\begin{array}{ccc}
1& 0 & (\mathbb{S}_1)_{13} 
\\
0 & 1 & (\mathbb{S}_1)_{23}
\\
0 & 0 & 1
\end{array}
\right).
$$
To find the third entry,  we observe that $\mathcal{S}_2$ is obtained from  $\mathcal{S}_1$  by a rotation $z\mapsto z e^{-i\pi}$, and that   $H_{\frac{3}{4}}^{(1)}(ize^{-i\pi}/2)$ gives the correct asymptotics for $-\pi/2<\arg z<5\pi/2$. Therefore,
$$\tilde{Y}_2(z)
= \left(
\begin{array}{ccc}
y_1(z,c_1^{(-)},0,0,0) & y_1(z,c_1^{(+)},0,0,0) &y_1(z,\widehat{c}_1,0,-1,0)
\\
\star & \star & \star
\\
\star & \star & \star
\end{array}
\right),
$$
with,
\be
\label{28feb2016-11}
y_1(z,\widehat{c}_1,0,-1,0)=\widehat{c}_1\sqrt{\frac{z}{2}} e^{z/2}~H_{\frac{3}{4}}^{(1)}(ize^{-i\pi}/2)=\frac{e^z}{\sqrt{2}} \left(1+\mathcal{O}\left(\frac{1}{z}\right)\right),
\quad
\quad
\widehat{c}_1:=\frac{\sqrt{\pi}}{2}e^{3\pi i/8}.
\ee
Finally,  the cyclic relation (see \cite{Wasow}) 
$$ 
H_{\frac{3}{4}}^{(1)}\left(e^{-i\pi}~\frac{iz}{2}\right)=\sqrt{2}H_{\frac{3}{4}}^{(1)}\left(\frac{iz}{2}\right)+ e^{-{3\pi i/4}} H_{\frac{3}{4}}^{(2)}\left(\frac{iz}{2}\right),
$$
  together with (\ref{28feb2016-11}) and (\ref{28feb2016-12}),  implies that 
  \be
  \label{28feb2016-4}
  y_1(z,\widehat{c}_1,0,-1,0)=\sqrt{\frac{\pi}{2}}~e^{3\pi i /8}~\sqrt{\frac{z}{2}} e^{z/2}~H_{\frac{3}{4}}^{(1)}\left(\frac{iz}{2}\right)+   y_1(z,0,{c}_2,0,0).
  \ee
  On the other hand, from the definition of $\mathbb{S}_1$ we must have 
  \be
  \label{28feb2016-9}
   y_1(z,\widehat{c}_1,0,-1,0)=(\mathbb{S}_1)_{13}~y_1(z,c_1^{(-)},0,0,0)+ (\mathbb{S}_1)_{23}~ y_1(z,c_1^{(+)},0,0,0)+   y_1(z,0,{c}_2,0,0).
   \ee
   Clearly, (\ref{28feb2016-5}) (\ref{28feb2016-6}) and  (\ref{28feb2016-4}) are not enough to determine $(\mathbb{S}_1)_{13}$ and $(\mathbb{S}_1)_{23}$. 
  Thus, we analyse the second row of (\ref{28feb2016-8}), to which corresponds $y_2(z)$. Recall that  $y_2^\prime(z)=0$. Therefore,  we  choose $y_2(z)=1/\sqrt{2}$ for the first two entries, and $y_2(z)=0$ for  the third.  
   This gives, for the second row of $Y_2=Y_1\mathbb{S}_1$:
   $$
   \left[ \frac{1}{\sqrt{2}},~\frac{1}{\sqrt{2}},~0\right]=\left[ \frac{1}{\sqrt{2}},~\frac{1}{\sqrt{2}},~\Bigl((\mathbb{S}_1)_{13}+(\mathbb{S}_1)_{23}\Bigr) \frac{1}{\sqrt{2}}\right]
\quad
~\Longrightarrow 
\quad~ 
  (\mathbb{S}_1)_{23} = -  (\mathbb{S}_1)_{13} .
$$
Thus, (\ref{28feb2016-4}) and (\ref{28feb2016-9}) become
$$
(\mathbb{S}_1)_{13}~\Bigl(y_1(z,c_1^{(-)},0,0,0)-  y_1(z,c_1^{(+)},0,0,0)\Bigr)=\sqrt{\frac{\pi}{2}}~e^{3\pi i /8}~\sqrt{\frac{z}{2}} e^{z/2}~H_{\frac{3}{4}}^{(1)}\left(\frac{iz}{2}\right).
$$
Keeping into account (\ref{28feb2016-5}) and (\ref{28feb2016-6}) we find $(\mathbb{S}_1)_{13}=1$. 
Thus 
$$ 
\mathbb{S}_1=
\left(
\begin{array}{ccc}
1& 0 & ~1 
\\
0 & 1 & -1
\\
0 & 0 & ~1
\end{array}
\right)
.
$$
$\mathbb{S}_2$ can be computed in a similar way, by a further rotation. On the other hands, since $V$ is skew symmetric
$$
\mathbb{S}_2=\mathbb{S}_1^{-T}
=\left(
\begin{array}{ccc}
1& 0 & 0
\\
0 & 1 & 0
\\
-1 & 1 & 1
\end{array}
\right).
$$
The result is in accordance with Theorem \ref{16dicembre2016-1}, which predicts that the entry $(1,2)$ of $\mathbb{S}_1$ and the entry $(2,1)$ of $\mathbb{S}_2$ must be zero.  It is also  in accordance with the monodromy  data of $y(t)$ obtained  in  \cite{DM}.

\bre  
If we choose $V(0)$ with different signs, we obtain different signs in $\mathbb{S}_1$. For example, consider the  choice 
$$
\overline{V}(0)=\left(
\begin{array}{ccc}
0 & 0& -i\sqrt{2}/8
\\
0 & 0 &- i\sqrt{2}/8
\\
i\sqrt{2}/8 & i\sqrt{2}/8 & 0
\end{array}
\right)=JV(0)J,
\quad
\quad
J:=\hbox{diag}(1,-1,1).
$$
\ere
The same procedure as above yields
$$ 
\overline{\mathbb{S}}_1=
\left(
\begin{array}{ccc}
1& 0 & -1 
\\
0 & 1 & -1
\\
0 & 0 & ~1
\end{array}
\right)~\equiv  J\mathbb{S}_1J.
$$
This sign  freedom corresponds to the invariance of $U=\hbox{diag}(u_1,u_2,u_3)$, namely $JUJ\equiv U$.
The  result $\overline{\mathbb{S}}_1$ above  is in accordance with the known result of \cite{Dub1}  that the  Stokes matrix $S$   of the $A_3$ Frobenius manifold is (up to permutation, change of signs and action of the braid group)  the Stokes matrix $S$ such that $S+S^T$ is the Coxeter matrix  of the reflection group $A_3$.

\bre 
 Another Stokes matrix $S$ obtained by an action of the braid group from that computed above  exists with   entries $(S_{12},S_{13}, S_{23})=(1,1,1)$;  the corresponding  branch $y(t)$  at $t=0$ is obtained letting  $s\to 1$, yielding the    Puiseux series \cite{DM}
$$
y(t)=
{4^{2/3}\over 50} t^{2/3} +{1\over 2} t+{1941 \cdot 2^{2/3}\over 2500} t^{4/3} -{2^{1/3}\over 150} t^{5/3} +O(t^2),
$$
 to which corresponds the behaviour of $V(t)$,
  $$
  \Omega_1(t)=
  {2^{1/6}\over 12 \cdot t^{1/6}}
  -{5\cdot 2^{5/6}\over 96}~t^{1/6}+O(t^{1/2}),
  \quad
  \quad
  \Omega_2(t)= 
  {i\over 6}+{ 2^{1/3}~i\over 96}~t^{2/3}+O(t^{4/3})
  ,
$$
$$
\Omega_3(t)=  
  {2^{1/6}\over 12 \cdot t^{1/6}}
  +{5\cdot 2^{5/6}\over 96}~t^{1/6}+O(t^{1/2}).
  $$
  Thus, $V(t)$ has a branch point at $t=0$, no entry vanishes and both $\Omega_1(t)$ and  $\Omega_3(t)$ diverge, without contradiction with Theorem \ref{16dicembre2016-1}.
\ere

\vskip 0.5 cm 
\hrule

\section*{APPENDIX: Examples of Cell Decomposition}

 \bex
{\rm 
Let 
$$
\Lambda(t)=\hbox{diag}(u_1(t),u_2(t),u_3(t)):=\hbox{diag}(0,t,1).
$$
In this example, the coalescence locus  in a neighbourhood of $t=0$  is $\{0\}$, while the global coalescence locus in $\mathbb{C}$  is  $\{0,1\}$. 
 At $t=0$ we  have 
$$ 
\arg(u_1(0)-u_3(0))=\arg(0-1),
\quad
\quad
\arg(u_3(0)-u_1(0))=\arg(1-0).
$$ 
 We {\it choose} $
\widehat{\arg}(1)=0$, $\widehat{\arg}(-1)=\pi$. 
This implies that an admissible direction $\eta$ such that $\eta-2\pi<\widehat{\arg}(u_i(0)-u_j(0))<
\eta$ 
must satisfy 
$$ 
\eta-2\pi < 0 <\eta,
\quad
\quad
\eta-2\pi< \pi<\eta
\quad
\quad
\Longrightarrow
\quad
\quad
\pi<\eta<2\pi. 
$$ 
Therefore $\tau=3\pi/2-\eta$ satisfies 
$$ 
-\frac{\pi}{2}<\tau<\frac{\pi}{2}.
$$

\vskip 0.2 cm 
\noindent
-- At $t\neq 0$: $u_1(t)=u_1(0)$ and $u_3(t)=u_3(0)$, and 
$$ 
\arg(u_1(t)-u_2(t))=\arg(-t),
\quad
\quad
\arg(u_2(t)-u_1(t))=\arg(t),
$$
$$ 
\arg(u_3(t)-u_2(t))=\arg(1-t),
\quad
\quad
\arg(u_2(t)-u_3(t))=\arg(t-1).
$$
We {\it impose:}
$$
\eta-2\pi< \widehat{\arg}(-t)<\eta,
\quad
\quad
\eta-2\pi< \widehat{\arg}(t)<\eta,
$$
$$ 
\Downarrow
$$
$$
\quad
~\eta-2\pi<\widehat{\arg}(t)<\eta-\pi~\hbox{ out } ~\eta-\pi<\widehat{\arg}(t)<\eta. 
$$ 
The above gives the 2 cells of $\mathcal{U}_{\epsilon_0}(0)$ for $\epsilon_0<1$.  
$$ c(-):=\{t\in \mathcal{U}_{\epsilon_0}(0)~|~ \eta-2\pi<\arg(t)<\eta-\pi\},
\quad
\quad
 c(+):=\{t\in \mathcal{U}_{\epsilon_0}(0)~|~ \eta-\pi<\arg(t)<\eta\}.
$$

Since $u(t)$ is globally defined (and   $t=1$ is another coalescence point), one can globally divide the $t$-plane into cells. Accordingly, we also  impose the condition
$$
\eta-2\pi< \widehat{\arg}(1-t)<\eta,
\quad
\quad
\eta-2\pi< \widehat{\arg}(t-1)<\eta,
$$
$$ 
\Downarrow
$$
$$
\quad~\eta-2\pi<\widehat{\arg}(t-1)<\eta-\pi~\hbox{ out } ~\eta-\pi<\widehat{\arg}(t-1)<\eta. 
$$ 
Therefore,  the $t$ plane is globally  partitioned  into 3 cells by the above relation, as in  figure \ref{cellex1}.

\begin{figure}
\minipage{0.5\textwidth}
\centerline{\includegraphics[width=0.7\textwidth]{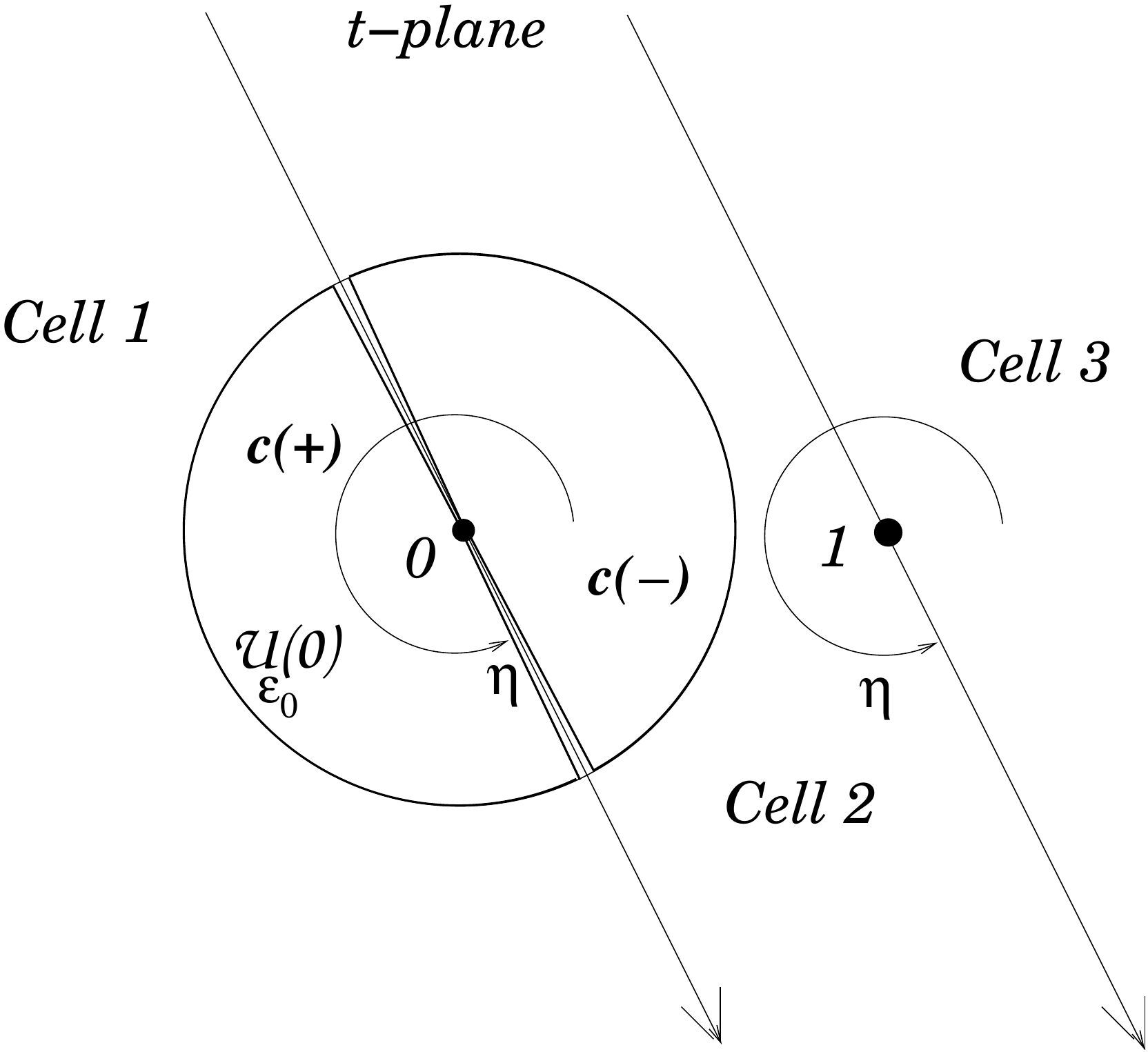}}
\caption{Cell partition (Cell 1, Cell 2, Cell 3) of the $t$-sheet  $\eta-2\pi<\arg(t)<\eta$ and $\eta-2\pi<\arg(t-1)<\eta$.  The neighbourhood $\mathcal{U}_{\epsilon_0}(0)$ (the disk) splits into two cells $c(+)$ and $c(-)$.}
\label{cellex1}
\endminipage\hfill
\minipage{0.45\textwidth}
\centerline{\includegraphics[width=0.65\textwidth]{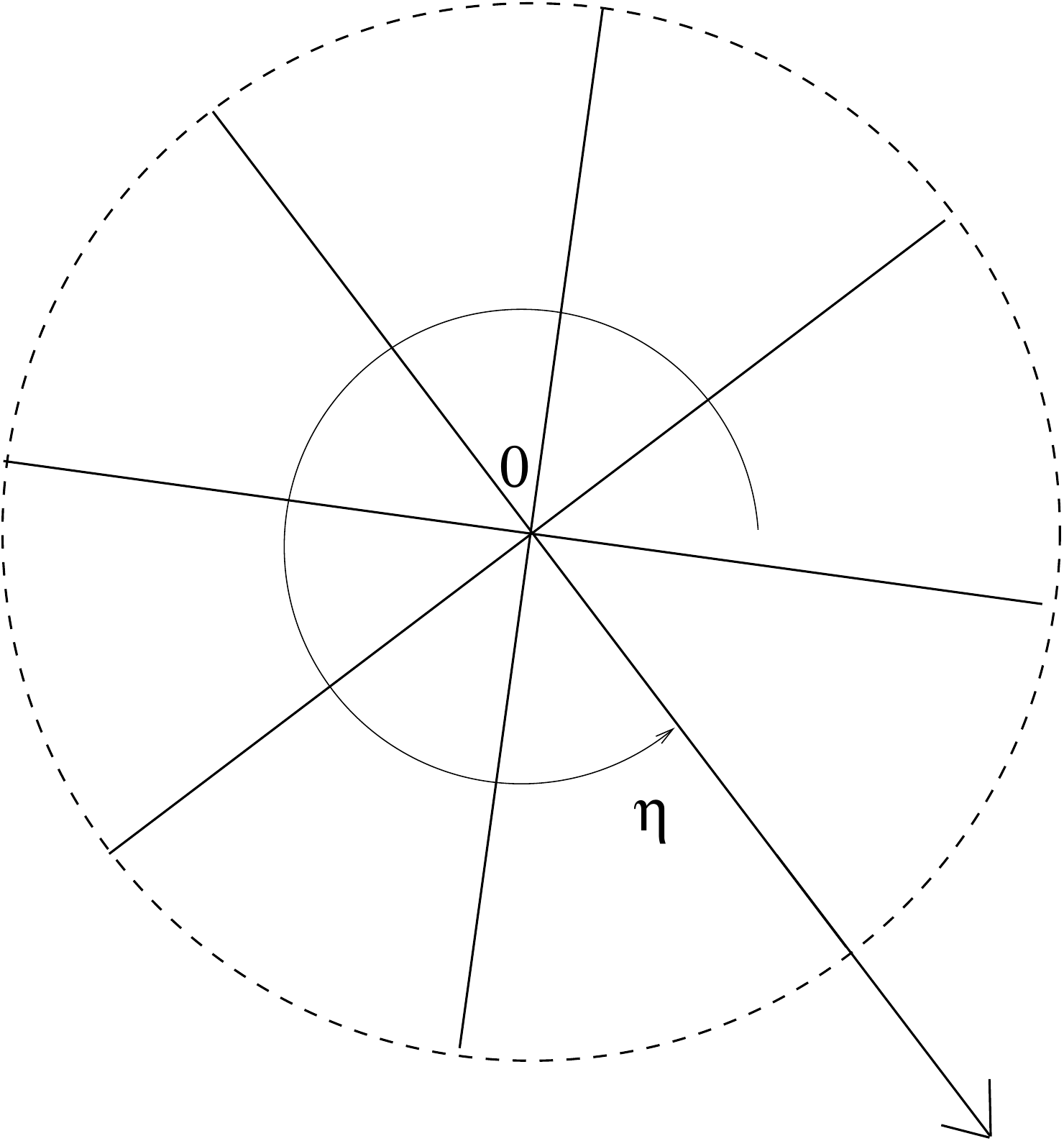}}
\caption{The cells of $\mathcal{U}_{\epsilon_0}(0)$ of Example \ref{18giugno-1}.}
\label{cellex3}
\endminipage
\end{figure}

\begin{figure}
\centerline{\includegraphics[width=0.7\textwidth]{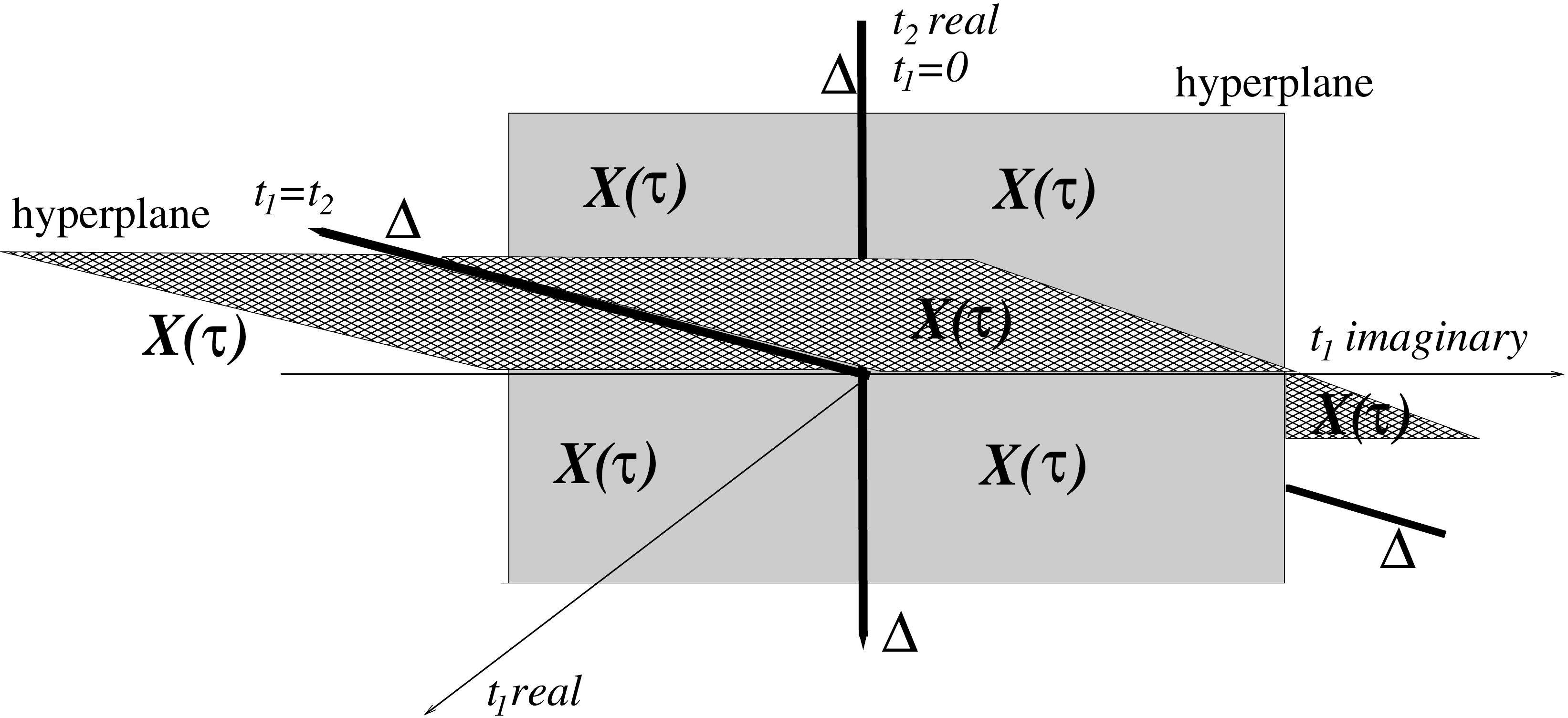}}
\caption{Example \ref{17feb2016-2}, with $\eta=3\pi/2$. The horizontal plane is $t_1\in\mathbb{C}$. The vertical axis is $t_2\in\mathbb{R}$. The  thick lines $t_1=t_2$ (real) and $t_1=0$ ($t_2$ real) are the projection of $\Delta_{\mathbb{C}^2}$. The  planes   (minus $\Delta$) are the projection of the crossing locus $X(\tau)$. The full planes (which include the thick lines)   are the projection of  $W(\tau)$. They disconnect  $ 
 \{t\in\mathbb{C}^2~|~t_2\in\mathbb{R}\}
 $. }
\label{cellex4}
\end{figure} 

}
\eex

\bex
\label{18giugno-1}
{\rm 
Let 
$$ 
\Lambda(t)=\hbox{ diag}\Bigl(u_1(t),u_2(t),u_3(t),u_4(t),u_5(t)\Bigr):=\hbox{ diag}\Bigl(0,~t,~te^{i\frac{\pi}{2}},~te^{i\pi},~te^{i\frac{3\pi}{2}}\Bigr).
$$
The coalescence locus is $t=0$.  
 The admissible direction $\eta$ can be chosen arbitrarily, because $\Lambda(0)=0$ has no Stokes rays. Once $\eta$ is fixed, we impose $\eta-2\pi <\widehat{\arg}(u_i(t)-u_j(t)) <\eta$. Thus, for $0\leq l,k\leq 3$:
$$ 
\eta-2\pi < \widehat{\arg}(t e^{i \frac{\pi}{2}k})<\eta,
\quad~
\eta-2\pi < \widehat{\arg}(-t e^{i \frac{\pi}{2}k})<\eta,
\quad~ 
\eta-2\pi < \widehat{\arg}\left(t(e^{i \frac{\pi}{2}l}- e^{i \frac{\pi}{2}k})\right)<\eta.
$$
The first two constraints imply
$$ 
\eta-2\pi-\frac{\pi}{2}k < \arg t <\eta-\pi-\frac{\pi}{2}k,
\quad
~\hbox{ or }
\quad~\eta-\pi-\frac{\pi}{2}k < \arg t <\eta-\frac{\pi}{2}k.
$$ 
By prosthaphaeresis formulas we have $ 
e^{i \frac{\pi}{2}l}- e^{i \frac{\pi}{2}k}=2i\sin\frac{\pi}{4}(l-k)~e^{i\frac{\pi}{4}(l+k)}
$. 
Therefore, the third constraint gives 
$$ 
\eta-2\pi-\frac{\pi}{4}(l+k) < \arg t <\eta-\pi-\frac{\pi}{4}(l+k),
\quad~
\hbox{ or }
\quad~\eta-\pi-\frac{\pi}{4}(l+k) < \arg t <\eta-\frac{\pi}{4}(l+k).
$$ 
It turns out that the cell-partition of $\mathcal{U}_{\epsilon_0}(0)$ is into 8 slices of angular width $\pi/4$, with angles determined by $\eta$. See figure \ref{cellex3}.

}
\eex


 \bex
\label{17feb2016-2}
 {\rm
 We consider $t=(t_1,t_2)\in\mathbb{C}^2$ and 
 $ 
 \Lambda(t)=\hbox{diag}(0,t_1,t_2)
  $.  
 The coalescence locus can be studied  globally on $\mathbb{C}^2$: 
 $$
 \Delta_{\mathbb{C}^2}=\{t\in\mathbb{C}^2~|~t_1=t_2\}\cup \{t\in\mathbb{C}^2~|~t_1=0\}\cup \{t\in\mathbb{C}^2~|~t_2=0\}.
 $$ 
 This is the union of complex lines (complex dimension $=1$) of complex co-dimension $=1$. In particular, $t=0$  is the point of maximal coalescence.  $\Lambda(0)=0$ has has no Stokes rays,  thus we choose $\eta$ freely.  The cell-partition for a chosen $\eta$ is given  (see previous examples) by:
 $$ 
 \eta-2\pi < \arg(t_i)<\eta-\pi, 
 \quad~\hbox{ or }
 \quad~\eta-\pi<\arg(t_i)<\eta,
 \quad\quad
 i=1,2,
 $$ 
 and
 $$ 
 \eta-2\pi < \arg(t_1-t_2)<\eta-\pi, 
 \quad~\hbox{ or }
 \quad~\eta-\pi<\arg(t_1-t_2)<\eta,
 \quad
 \quad
 i=1,2.
 $$
 In figure \ref{cellex4} we represent the projection of $\mathbb{C}^2$ onto the subspace
  $ 
 \{t\in\mathbb{C}^2~|~t_2\in\mathbb{R}\}
 $,  
for the choice  $\eta=3\pi/2$. 
 The two thick lines 
 $$ 
 t_1=t_2 \hbox{ real},
 \quad
 \quad
 t_1=0 \hbox{ with $t_2$ real},
 $$
 are the projection of $\Delta_{\mathbb{C}^2}$. The following  planes, without  the thick lines,       $$ 
\Bigl\{t~|~  \arg(t_1-t_2)=\frac{\pi}{2} \hbox{ or } \frac{3\pi}{2} \hbox{ mod } 2\pi\Bigr\} \cup  \Bigl\{t~|~  \arg(t_1)=\frac{\pi}{2} \hbox{ or } \frac{3\pi}{2} \hbox{ mod } 2\pi \Bigr\}
 $$ 
 are the projection  of the crossing locus $X(\tau)$. The planes, including the thick lines, are the projections of $W(\tau)$. 
 }
 \eex


   \end{document}